# SOFT NEUTROSOPHIC ALGEBRAIC STRUCTURES AND THEIR GENERALIZATION, VOL. 2

❖ **Mumtaz Ali**

E-mail: mumtazali770@yahoo.com

bloomy_boy2006@yahoo.com

❖ **Florentin Smarandache**

E-mail: smarand@unm.edu

fsmarandache@gmail.com

❖ **Muhammad Shabir**

E-mail: mshabirbhatti@yahoo.co.uk

# Contents

**Chapter No. 1**

**INTRODUCTION**

**Chapter No. 2**

**SOFT NEUTROSOPHIC GROUPOIDS AND THEIR GENERALIZATION**

**Chapter No. 3**

**SOFT NEUTROSOPHIC RINGS, SOFT NEUTROSOPHIC FIELDS AND THEIR GENERALIZATION**

**Chapter No. 4**

**SOFT NEUTROSOPHIC GROUP RINGS AND THEIR GENERALIZATION**

**Chapter No. 5**

**SOFT NEUTROSOPHIC SEMIGROUP RINGS AND THEIR GENERALIZATION**

**Chapter No. 6**

**SOFT MIXED NEUTROSOPHIC N-ALGEBRAIC STRUCTURES**

# DEDICATION by Mumtaz Ali

This book is dedicated to my dear Parents. Their kind support, praying and love always courage me to work hard all the time. I'm so much thankful to them with whole heartedly and dedicated this book to them. This is my humble way of paying homage to my parents.

# PREFACE

In this book we define some new notions of soft neutrosophic algebraic structures over  neutrosophic algebraic structures. We define some different soft neutrosophic algebraic structures but the main motivation is two-fold.  Firstly the classes of soft neutrosophic group ring and soft neutrosophic semigroup ring defined in this book is basically the generalization of two classes of rings: neutrosophic group rings and neutrosophic semigroup rings. These soft neutrosophic group rings and soft neutrosophic semigroup rings are defined over neutrosophic group rings and neutrosophic semigroup rings respectively. This is basically the collection of parameterized subneutrosophic group ring and subneutrosophic semigroup ring of a neutrosophic group ring and neutrosophic semigroup ring respectively. These soft neutrosophic algebraic structures are even bigger than the corresponding classical algebraic structures. It is interesting to see that it wider the area of research for the researcher of algebraic structures.

The organization of this book is as follows. This book has seven chapters. The first chapter is about introductory in nature, for it recalls some basic definitions and notions which are essential to make this book a self-contained one. Chapter two, introduces for the first time soft neutrosophic groupoid over neutrosophic groupoid and their generalization to catch up the whole theory in broader sense. In chapter three, we introduced soft neutrosophic rings which are basically defined over neutrosophic rings and we also generalized this study with some of its core properties. In this chapter, soft neutrosophic fields are also defined over neutrosophic fields with its generalization. Chapter four is about to introduce soft neutrosophic group rings over neutrosophic group rings and generalizes this theory with some of its fundamental properties. In chapter five, we kept soft neutrosophic semigroup rings over neutrosophic semigroup rings in the context of soft sets. The generalization of soft neutrosophic semigroup rings have also been made in the present chapter with some of their interesting properties. Chapter six introduces soft mixed neutrosophic N-algebraic structures over mixed neutrosophic N-algebraic structures with some of their basic properties. In the final chapter seven, we gave a number of suggested problems.

We are very thankful to dear Farzana Ahmad for her kind support, her precious timing, loving and caring.

**MUMTAZ ALI**

**FLORENTIN SMARANDACHE**

**MUHAMMAD SHABIR**

# Chapter No.1

# INTRODUCTION

In this introductory chapter, give certain basic concepts and notions. There are seven sections in this chapter. The first section contains the basic definitions and notions of neutrosophic groupoids, neutrosophic bigroupoids and neutrosophic N-groupoids. In section two, we mentioned neutrosophic rings, neutrosophic birings and neutrosophic N-rings. Section three is about the notions of neutrosophic fields and their generalization. The fourth section contains the basic theory and notions of neutrosophic group rings with its generalization. Similarly section five has covered the basic definitions of neutrosophic semigroup ring and their generalization. In the sixth section, we give the basic concepts about mixed neutrosophic N-algebraic structures. In the last section seven, some basic description and literature is discussed about soft sets and their related properties which will certainly help the readers to understand the notions of soft sets.

## 1.1 Neutrosophic Groupoids, Neutrosophic Bigroupoids, Neutrosophic N-groupoids and their Properties

Groupoids are the most generalized form of all the algebraic structures with a single binary operation. Neutrosophic groupoids are bigger in size than the classical groupoids with some kind of extra properties called indeterminacies. The definition of neutrosophic groupoid is as follows.

**Definition 1.1.1**. Let $G$ be a groupoid. The groupoid generated by $G$ and $I$ i.e. $G \cup I$ is denoted by $\langle G \cup I \rangle$ and defined to be a neutrosophic groupoid where $I$ is the indeterminacy element and termed as neutrosophic element with the property $I^2 = I$. For an integer $n, n + I$ and $nI$ are neutrosophic elements and $0.I = I$. Then inverse of $I$ is not defined and hence does not exist.

**Example 1.1.1.** Let $G = \{a, b \in \mathbb{Z}_3 : a * b = a + 2b \pmod{3}\}$ be a groupoid. Then $\langle G \cup I \rangle = \{0, 1, 2, I, 2I, 1 + I, 2 + I, 1 + 2I, 2 + 2I, *\}$ is a neutrosophic groupoid.

**Example 1.1.2.** Let $\left( \langle \mathbb{Z}^+ \cup I \rangle, * \right)$ be the set of positive integers with a binary operation $*$ where $a * b = 2a + 3b$ for all $a, b \in \langle \mathbb{Z}^+ \cup I \rangle$. Then $\left( \langle \mathbb{Z}^+ \cup I \rangle, * \right)$ is a neutrosophic groupoid.

**Definition 1.1.2.** Let $\langle G \cup I \rangle$ be a neutrosophic groupoid. A proper subset $P$ of $\langle G \cup I \rangle$ is said to be a neutrosophic subgroupoid if $P$ is a neutrosophic groupoid under the operations of $\langle G \cup I \rangle$.

A neutrosophic groupoid $\langle G \cup I \rangle$ is said to have a subgroupoid if $\langle G \cup I \rangle$ has a proper subset which is a groupoid under the operations of $\langle G \cup I \rangle$.

**Example 1.1.3.** Let $(\langle \mathbb{Z}_{10} \cup I \rangle, *)$ be a neutrosophic groupoid under the binary operation $*$ which is defined by $a * b = 3a + 2b \pmod{10}$ for all $a, b \in \langle \mathbb{Z}_{10} \cup I \rangle$. Then $P = \{0, 5, 5I, 5 + 5I\}$ is a neutrosophic subgroupoid of $(\langle \mathbb{Z}_{10} \cup I \rangle, *)$ and $Q = (\mathbb{Z}_{10}, *)$ is just a groupoid of $(\langle \mathbb{Z}_{10} \cup I \rangle, *)$.

**NOTE.** Let $\langle G \cup I \rangle$ be a neutrosophic groupoid. Suppose $P_1$ and $P_2$ be any two neutrosophic subgroupoids of $\langle G \cup I \rangle$. Then $P_1 \cup P_2$, the union of two neutrosophic subgroupoids in general need not be a neutrosophic subgroupoid.

**Definition 1.1.4.** A proper subset $P$ of $\langle G \cup I \rangle$ is said to be a neutrosophic strong subgroupoid if $P$ is a neutrosophic groupoid of $\langle G \cup I \rangle$ and all the elements of $P$ are neutrosophic elements.

**Definition 1.1.5.** Let $(\langle G \cup I \rangle, *)$ be a neutrosophic groupoid of finite order. A proper subset $P$ of $\langle G \cup I \rangle$ is said to be a Lagrange neutrosophic subloop, if $P$ is a neutrosophic subgroupoid under the operation $*$ and $o(P) / o\langle G \cup I \rangle$.

**Definition 1.1.6.** If every neutrosophic subgroupoid of $\langle G \cup I \rangle$ is Lagrange then we call $\langle G \cup I \rangle$ to be a Lagrange neutrosophic groupoid.

**Definition 1.1.7.** If $\langle G \cup I \rangle$ has no Lagrange neutrosophic subgroupoid then we call $\langle G \cup I \rangle$ to be a Lagrange free neutrosophic groupoid.

**Definition 1.1.8.** If $\langle G \cup I \rangle$ has at least one Lagrange neutrosophic subgroupoid then we call $\langle G \cup I \rangle$ to be a weakly Lagrange neutrosophic groupoid.

**Definition 1.1.9.** Let $\langle G \cup I \rangle$ be a neutrosophic groupoid under a binary operation $*$. $P$ be a proper subset of $\langle G \cup I \rangle$. $P$ is said to be a neutrosophic ideal of $\langle G \cup I \rangle$ if the following conditions are satisfied.

1. $P$ is a neutrosophic subgroupoid of $\langle G \cup I \rangle$.
2. For all $p \in P$ and for all $s \in \langle G \cup I \rangle$ we have $p * s$ and $s * p$ are in $P$.

We now proceed onto define neutrosophic bigroupoid.

**Neutrosophic Bigroupoid**

Here we give some definitions and notions of neutrosophic bigroupoid with basic properties and give some examples and to understand the theory of neutrosophic bigroupoid.

**Definition 1.1.10.** Let $(BN(\text{G}),*,\circ)$ be a non-empty set with two binary operations $*$ and $\circ$. $(BN(\text{G}),*,\circ)$ is said to be a neutrosophic bigroupoid if $BN(\text{G}) = P_1 \cup P_2$ where at least one of $(P_1,*)$ or $(P_2,\circ)$ is a neutrosophic groupoid and other is just a groupoid. $P_1$ and $P_2$ are proper subsets of $BN(\text{G})$.

If both $(P_1,*)$ and $(P_2,\circ)$ in the above definition are neutrosophic groupoids then we call $(BN(\text{G}),*,\circ)$ to be a strong neutrosophic bigroupoid. All strong neutrosophic bigroupoids are trivially neutrosophic bigroupoids.

**Example 1.1.4.** Let $\{B_N(\text{G}),*,\circ\}$ be a non-empty set with $B_N(\text{G}) = G_1 \cup G_2$, where $G_1 = \left\{ \left\langle Z_{10} \cup I \right\rangle \mid a * b = 2a + 3b (\text{mod}\,10); a,b \in \left\langle Z_{10} \cup I \right\rangle \right\}$ and

$G_2 = \left\{ \left\langle Z_4 \cup I \right\rangle \mid a \circ b = 2a + b (\text{mod}\,4); a,b \in \left\langle Z_4 \cup I \right\rangle \right\}$. Then $\{B_N(\text{G}),*,\circ\}$ is a neutrosophic bigroupoid.

**Example 1.1.5.** Let $BN(\text{G}) = \left\{ G_1 \cup G_2, *_1, *_2 \right\}$ where $G_1 = \left\{ \left\langle Z_4 \cup I \right\rangle \mid a \circ b = 2a + b (\text{mod}\,4); a,b \in \left\langle Z_4 \cup I \right\rangle \right\}$ is a neutrosophic groupoid and $G_2 = \left\{ \left\langle Z_{12} \cup I \right\rangle \mid a * b = 8a + 4b (\text{mod}\,12); a,b \in \left\langle Z_{12} \cup I \right\rangle \right\}$ is also another neutrosophic groupoid. Then clearly $BN(\text{G}) = \left\{ G_1 \cup G_2, *_1, *_2 \right\}$ is a neutrosophic groupoid.

**Definition 1.1.11.** Let $(BN(\text{G}) = P_1 \cup P_2; *,\circ)$ be a neutrosophic bigroupoid. A proper subset $(T,*,\circ)$ is said to be a neutrosophic subbigroupoid of $BN(\text{G})$ if

1. $T = T_1 \cup T_2$ where $T_1 = P_1 \cap T$ and $T_2 = P_2 \cap T$ and
2. At least one of $(T_1,\circ)$ or $(T_2,*)$ is a neutrosophic groupoid.

**Definition 1.1.12.** Let $(BN(\mathrm{G}) = \mathrm{P}_1 \cup P, *, \circ)$ be a neutrosophic strong bigroupoid. A proper subset $T$ of $BN(S)$ is called the strong neutrosophic subbigroupoid if $T = T_1 \cup T_2$ with $T_1 = P_1 \cap T$ and $T_2 = P_2 \cap T$ and if both $(T_1, *)$ and $(T_2, \circ)$ are neutrosophic subgroupoids of $(P_1, *)$ and $(P_2, \circ)$ respectively. We call $T = T_1 \cup T_2$ to be a neutrosophic strong subbigroupoid, if at least one of $(T_1, *)$ or $(T_2, \circ)$ is a groupoid then $T = T_1 \cup T_2$ is only a neutrosophic subgroupoid.

**Definition 1.1.13.** Let $(BN(\mathrm{G}) = \mathrm{P}_1 \cup P_2, *, \circ)$ be any neutrosophic bigroupoid. Let $J$ be a proper subset of $\mathrm{BN(G)}$ such that $J_1 = J \cap P_1$ and $J_2 = J \cap P_2$ are ideals of $P_1$ and $P_2$ respectively. Then $J$ is called the neutrosophic biideal of $BN(\mathrm{G})$.

**Definition 1.1.14.** Let $(BN(\mathrm{G}), *, \circ)$ be a strong neutrosophic bigroupoid where $BN(S) = \mathrm{P}_1 \cup P_2$ with $(P_1, *)$ and $(P_2, \circ)$ be any two neutrosophic groupoids. Let $J$ be a proper subset of $BN(\mathrm{G})$ where $I = I_1 \cup I_2$ with $I_1 = I \cap P_1$ and $I_2 = I \cap P_2$ are neutrosophic ideals of the neutrosophic groupoids $P_1$ and $P_2$ respectively. Then $I$ is called or defined as the strong neutrosophic biideal of $BN(\mathrm{G})$.

Union of any two neutrosophic biideals in general is not a neutrosophic biideal. This is true of neutrosophic strong biideals.

We now go further to generalize the concept of neutrosophic groupoids

**Neutrosophic $N$ -groupoid**

The basic definitions and notions of neutrosophic N-groupoid are presented here with illustrative examples for the interesting readers.

**Definition 1.1.15.** Let $\{N(G), *_1, ..., *_2\}$ be a non-empty set with $N$ -binary operations defined on it. We call $N(G)$ a neutrosophic $N$ -groupoid ( $N$ a positive integer) if the following conditions are satisfied.

1. $N(G) = G_1 \cup ... \cup G_N$ where each $G_i$ is a proper subset of $N(G)$ i.e. $G_i \subset G_j$ or $G_j \subset G_i$ if $i \neq j$ .

2. $(G_i, *_i)$ is either a neutrosophic groupoid or a groupoid for $i = 1, 2, 3, ..., N$ .

If all the $N$ -groupoids $(G_i, *_i)$ are neutrosophic groupoids (i.e. for $i = 1, 2, 3, ..., N$ ) then we call $N(G)$ to be a neutrosophic strong $N$ - groupoid.

**Example 1.1.6.** Let $N(G) = \{G_1 \cup G_2 \cup G_3, *_1, *_2, *_3\}$ be a neutrosophic 3-groupoid, where $G_1 = \{\langle \mathbb{Z}_{10} \cup I \rangle \mid a * b = 2a + 3b \pmod{10}; a, b \in \langle \mathbb{Z}_{10} \cup I \rangle\}$ is a neutrosophic groupoid, $G_2 = \{\mathbb{Z}_4 \mid a \circ b = 2a + b \pmod{4}; a, b \in \mathbb{Z}_4\}$ is just a groupoid and $G_3 = \{\langle \mathbb{Z}_{12} \cup I \rangle \mid a * b = 8a + 4b \pmod{12}; a, b \in \langle \mathbb{Z}_{12} \cup I \rangle\}$ .

**Example 1.1.7.** Let $N(G) = \{G_1 \cup G_2 \cup G_3, *_1, *_2, *_3\}$ be a neutrosophic 3-groupoid, where $G_1 = \{\langle \mathbb{Z}_{10} \cup I \rangle \mid a * b = 2a + 3b (\bmod 10); a, b \in \langle \mathbb{Z}_{10} \cup I \rangle\}$ is a neutrosophic groupoid, $G_2 = \{\langle \mathbb{Z}_4 \cup I \rangle \mid a \circ b = 2a + b (\bmod 4); a, b \in \langle \mathbb{Z}_4 \cup I \rangle\}$ is a neutrosophic groupoid and $G_3 = \{\langle \mathbb{Z}_{12} \cup I \rangle \mid a * b = 8a + 4b (\bmod 12); a, b \in \langle \mathbb{Z}_{12} \cup I \rangle\}$ is also another neutrosophic groupoid. In fact $N(G)$ is a neutrosophic strong 3-groupoid.

**Definition 1.1.16.** Let $N(G) = \{G_1 \cup G_2 \cup \ldots \cup G_N, *_1, *_2, \ldots, *_N\}$ be a neutrosophic $N$-groupoid. A proper subset $P = \{P_1 \cup P_2 \cup \ldots P_N, *_1, *_2, \ldots, *_N\}$ of $N(G)$ is said to be a neutrosophic $N$-subgroupoid if $P_i = P \cap G_i, i = 1, 2, \ldots, N$ are subgroupoids of $G_i$ in which at least some of the subgroupoids are neutrosophic subgroupoids.

**Definition 1.1.17.** Let $N(G) = \{G_1 \cup G_2 \cup \ldots \cup G_N, *_1, *_2, \ldots, *_N\}$ be a neutrosophic strong $N$-groupoid. A proper subset $T = \{T_1 \cup T_2 \cup \ldots \cup T_N, *_1, *_2, \ldots, *_N\}$ of $N(G)$ is said to be a neutrosophic strong sub $N$-groupoid if each $(T_i, *_i)$ is a neutrosophic subgroupoid of $(G_i, *_i)$ for $i = 1, 2, \ldots, N$ where $T_i = G_i \cap T$.
If only a few of the $(T_i, *_i)$ in $T$ are just subgroupoids of $(G_i, *_i)$. i.e. $(T_i, *_i)$ are not neutrosophic subgroupoids then we call $T$ to be a sub $N$-groupoid of $N(G)$.

**Definition 1.1.18.** Let $N(G) = \{G_1 \cup G_2 \cup \ldots \cup G_N, *_1, *_2, \ldots, *_N\}$ be a neutrosophic $N$-groupoid. A proper subset $P = \{P_1 \cup P_2 \cup \ldots \cup P_N, *_1, *_2, \ldots, *_N\}$ of $N(G)$ is said to be a neutrosophic $N$-subgroupoid, if the following conditions are true.

1. $P$ is a neutrosophic sub $N$-groupoid of. $N(G)$
2. Each $P_i = G \cap P_i, i = 1, 2, ..., N$ is an ideal of $G_i$.

Then $P$ is called or defined as the neutrosophic $N$-ideal of the neutrosophic $N$-groupoid $N(G)$.

**Definition 1.1.19.** Let $N(G) = \{G_1 \cup G_2 \cup .... G_N, *_1, *_2, ..., *_N\}$ be a neutrosophic strong $N$-groupoid. A proper subset $J = \{J_1 \cup J_2 \cup .... J_N, *_1, *_2, ..., *_N\}$ where $J_t = J \cap G_t$ for $t = 1, 2, ..., N$ is said to be a neutrosophic strong $N$-ideal of $N(G)$ if the following conditions are satisfied.

1. Each it is a neutrosophic subgroupoid of $G_t, t = 1, 2, ..., N$ i.e. It is a neutrosophic strong N-subgroupoid of $N(G)$.
2. Each it is a two sided ideal of $G_t$ for $t = 1, 2, ..., N$.

Similarly one can define neutrosophic strong $N$-left ideal or neutrosophic strong right ideal of $N(G)$.

**Definition 1.1.20.** A neutrosophic strong $N$-ideal is one which is both a neutrosophic strong $N$-left ideal and $N$-right ideal of $N(G)$.

## 1.2 Neutrosophic Rings, Neutrosophic Birings, Neutrosophic N-rings and their Properties

In this section, we give the fundamental definitions of neutrosophic rings and their generalization with some basic characterization.

We now proceed to define neutrosophic ring and give some examples.

**Definition 1.2.1.** Let $R$ be any ring. The neutrosophic ring denoted by $\langle R \cup I \rangle$ is also a ring generated by $R$ and $I$ under the operations of $R$. $I$ is called the neutrosophic element with the property $I^2 = I$. For an integer n , $n + I$ and $nI$ are neutrosophic elements and $0.I = 0$. $I^{-1}$, the inverse of $I$ is not defined and hence does not exist.

**Example 1.2.1.** Let $\mathbb{Z}$ be a ring of integers. Then $\langle Z \cup I \rangle$ is a neutrosophic ring of integers.

**Example 1.2.2.** Let $\langle \mathbb{R} \cup I \rangle$ and $\langle \mathbb{C} \cup I \rangle$ are neutrosophic rings of real numbers and complex numbers respectively.

**Definition 1.2.2.** Let $\langle R \cup I \rangle$ be a neutrosophic ring. A proper subset $P$ of $\langle R \cup I \rangle$ is said to be a neutrosophic subring if $P$ itself is a neutrosophic ring under the operations of $\langle R \cup I \rangle$.

**Definition 1.2.3.** Let $\langle R \cup I \rangle$ be any neutrosophic ring. A non-empty subset $P$ of $\langle R \cup I \rangle$ is defined to be a neutrosophic ideal of $\langle R \cup I \rangle$ if the following conditions are satisfied.

1. $P$ is a neutrosophic subring of $\langle R \cup I \rangle$.

2. For every $p \in P$ and $r \in \langle R \cup I \rangle$, $rp$ and $pr \in P$.

## Neutrosophic Biring

Here we give some basic definitions and notions about neutrosophic birings.

**Definition 1.2.4.** Let $(BN(\mathrm{R}), *, \circ)$ be a non-empty set with two binary operations $*$ and $\circ$. $(BN(\mathrm{R}), *, \circ)$ is said to be a neutrosophic biring if $BN(\mathrm{Rs}) = R_1 \cup R_2$ where at least one of $(\mathrm{R}_1, *, \circ)$ or $(\mathrm{R}_2, *, \circ)$ is a neutrosophic ring and other is just a ring. $R_1$ and $R_2$ are proper subsets of $BN(\mathrm{R})$.

To illustrate this situation, see the following example.

**Example 1.2.3.** Let $BN(\mathrm{R}) = (\mathrm{R}_1, *, \circ) \cup (\mathrm{R}_2, *, \circ)$ where $(\mathrm{R}_1, *, \circ) = (\langle \mathbb{Z} \cup I \rangle, +, \times)$ and $(\mathrm{R}_2, *, \circ) = (\mathbb{Q}, +, \times)$. Clearly $(\mathrm{R}_1, *, \circ)$ is a neutrosophic ring under addition and multiplication. $(\mathrm{R}_2, *, \circ)$ is just a ring. Thus $(BN(\mathrm{R}), *, \circ)$ is a neutrosophic biring.

**Definition 1.2.5.** Let $BN(\mathrm{R}) = (\mathrm{R}_1, *, \circ) \cup (\mathrm{R}_2, *, \circ)$ be a neutrosophic biring. Then $BN(\mathrm{R})$ is called a commutative neutrosophic biring if each $(\mathrm{R}_1, *, \circ)$ and $(\mathrm{R}_2, *, \circ)$ is a commutative neutrosophic ring.

We explain it in the following example.

**Example 1.2.4.** Let $BN(\mathrm{R}) = (\mathrm{R}_1, *, \circ) \cup (\mathrm{R}_2, *, \circ)$ where $(\mathrm{R}_1, *, \circ) = (\langle \mathbb{Z} \cup I \rangle, +, \times)$ and $(\mathrm{R}_2, *, \circ) = (\mathbb{Q}, +, \times)$. Clearly $(\mathrm{R}_1, *, \circ)$ is a commutative neutrosophic ring and $(\mathrm{R}_2, *, \circ)$ is also a commutative ring. Thus $(BN(\mathrm{R}), *, \circ)$ is a commutative neutrosophic biring.

**Definition 1.2.6.** Let $BN(\mathrm{R}) = (\mathrm{R}_1, *, \circ) \cup (\mathrm{R}_2, *, \circ)$ be a neutrosophic biring. Then $BN(\mathrm{R})$ is called a pseudo neutrosophic biring if each $(\mathrm{R}_1, *, \circ)$ and $(\mathrm{R}_2, *, \circ)$ is a pseudo neutrosophic ring.

**Example 1.2.5.** Let $BN(\mathrm{R}) = (\mathrm{R}_1, +, \times) \cup (\mathrm{R}_2, +, \circ)$ where $(\mathrm{R}_1, +, \times) = \{0, I, 2I, 3I\}$ is a pseudo neutrosophic ring under addition and multiplication modulo 4 and $(\mathrm{R}_2, +, \times) = \{0, \pm 1I, \pm 2I, \pm 3I, \ldots\}$ is another pseudo neutrosophic ring. Thus $(BN(\mathrm{R}), +, \times)$ is a pseudo neutrosophic biring.

**Definition 1.2.7.** Let $(BN(\mathrm{R}) = R_1 \cup R_2; *, \circ)$ be a neutrosophic biring. A proper subset $(T, *, \circ)$ is said to be a neutrosophic subbiring of $BN(\mathrm{R})$ if

1. $T = T_1 \cup T_2$ where $T_1 = R_1 \cap T$ and $T_2 = R_2 \cap T$,
2. At least one of $(T_1, \circ)$ or $(T_2, *)$ is a neutrosophic ring.

**Definition 1.2.8:** If both $(\mathrm{R}_1, *)$ and $(\mathrm{R}_2, \circ)$ in the above definition 1.2.4 are neutrosophic rings then we call $(BN(\mathrm{R}), *, \circ)$ to be a neutrosophic strong biring.

**Definition 1.2.9.** Let $(BN(\mathrm{R}) = R_1 \cup R_2; *, \circ)$ be a neutrosophic biring and let $(T, *, \circ)$ is a neutrosophic subbiring of $BN(\mathrm{R})$. Then $(T, *, \circ)$ is called a neutrosophic biideal of $BN(R)$ if

1. $T = T_1 \cup T_2$ where $T_1 = R_1 \cap T$ and $T_2 = R_2 \cap T$.
2. At least one of $(T_1, *, \circ)$ or $(T_2, *, \circ)$ is a neutrosophic ideal.

If both $(T_1,*,\circ)$ and $(T_2,*,\circ)$ in the above definition are neutrosophic ideals,

then we call $(T,*,\circ)$ to be a neutrosophic strong biideal of $BN(R)$.

**Definition 1.2.10.** Let $(BN(R) = R_1 \cup R_2; *,\circ)$ be a neutrosophic biring and let $(T,*,\circ)$ is a neutrosophic subbiring of $BN(R)$. Then $(T,*,\circ)$ is called a pseudo neutrosophic biideal of $BN(R)$ if

1. $T = T_1 \cup T_2$ where $T_1 = R_1 \cap T$ and $T_2 = R_2 \cap T$.
2. $(T_1,*,\circ)$ and $(T_2,*,\circ)$ are pseudo neutrosophic ideals.

We now generalize the concept of neutrosophic rings.

**Neutrosophic $N$-ring**

The definitions and notions of neutrosophic N-rings are presented here with the illustrative examples.

**Definition 1.2.11.** Let $\{N(R),*_1,...,*_2,\circ_1,\circ_2,...,\circ_N\}$ be a non-empty set with $N$-binary operations defined on it. We call $N(R)$ a neutrosophic $N$-ring ($N$ a positive integer) if the following conditions are satisfied.

1. $N(R) = R_1 \cup R_2 \cup ... \cup R_N$ where each $R_i$ is a proper subset of $N(R)$ i.e. $R_i \not\subset R_j$ or $R_j \not\subset R_i$ if $i \neq j$.
2. $(R_i,*_i,\circ_i)$ is either a neutrosophic ring or a ring for $i = 1,2,3,...,N$.

This situation can be explained in the following example.

**Example 1.2.6.** Let $N(R) = (R_1,*,\circ) \cup (R_2,*,\circ) \cup (R_3,*,\circ)$ where $(R_1,*,\circ) = (\langle \mathbb{Z} \cup I \rangle,+,\times)$, $(R_2,*,\circ) = (\mathbb{Q},+,\times)$ and $(R_3,*,\circ) = (Z_{12},+,\times)$. Thus

$(N(\mathrm{R}),*,\circ)$ is a neutrosophic $N$-ring.

**Definition 1.2.12.** Let $\mathrm{N(R)} = \{\mathrm{R}_1 \cup R_2 \cup .... \cup \mathrm{R}_N, *_1, *_2, ..., *_N, \circ_1, \circ_2, ..., \circ_N\}$ be a neutrosophic N-ring. Then $N(\mathrm{R})$ is called a pseudo neutrosophic N-ring if each $(\mathrm{R}_i, *_i)$ is a pseudo neutrosophic ring where $i = 1, 2, ..., \mathrm{N}$.

**Example 1.2.7.** Let $N(\mathrm{R}) = (\mathrm{R}_1, +, \times) \cup (\mathrm{R}_2, +, \times) \cup (\mathrm{R}_3, +, \times)$ where $(\mathrm{R}_1, +, \times) = \{0, I, 2I, 3I\}$ is a pseudo neutrosophic ring under addition and multiplication modulo 4, $(\mathrm{R}_2, +, \times) = \{0, \pm1I, \pm2I, \pm3I, ...\}$ is a pseudo neutrosophic ring and $(\mathrm{R}_3, +, \times) = \{0, \pm2I, \pm4I, \pm6I...\}$. Thus $(N(\mathrm{R}), +, \times)$ is a pseudo neutrosophic 3-ring.

**Definition 1.2.13**. If all the $N$-rings $(\mathrm{R}_i, *_i)$ in definition 1.2.11 are neutrosophic rings (i.e. for $i = 1, 2, 3, ..., N$) then we call $\mathrm{N(R)}$ to be a neutrosophic strong $N$-ring.

**Example 1.2.8.** Let $N(\mathrm{R}) = (\mathrm{R}_1, *, \circ) \cup (\mathrm{R}_2, *, \circ) \cup (\mathrm{R}_3, *, \circ)$ where $(\mathrm{R}_1, *, \circ) = (\langle \mathbb{Z} \cup I \rangle, +, \times)$, $(\mathrm{R}_2, *, \circ) = (\langle \mathbb{Q} \cup I \rangle, +, \times)$ and $(\mathrm{R}_3, *, \circ) = (\langle \mathbb{Z}_{12} \cup I \rangle, +, \times)$. Thus $(N(\mathrm{R}), *, \circ)$ is a strong neutrosophic $N$-ring.

**Definition 1.2.14.** Let $\mathrm{N(R)} = \{\mathrm{R}_1 \cup R_2 \cup .... \cup \mathrm{R}_N, *_1, *_2, ..., *_N, \circ_1, \circ_2, ..., \circ_N\}$ be a neutrosophic $N$-ring. A proper subset $P = \{\mathrm{P}_1 \cup P_2 \cup .... \mathrm{P}_N, *_1, *_2, ..., *_N\}$ of $\mathrm{N(R)}$ is said to be a neutrosophic $N$-subring if $P_i = P \cap R_i, i = 1, 2, ..., N$ are subrings of $R_i$ in which at least some of the subrings are neutrosophic subrings.

**Definition 1.2.15.** Let $\mathrm{N(R)} = \{\mathrm{R}_1 \cup R_2 \cup .... \cup \mathrm{R}_N, *_1, *_2, ..., *_N, \circ_1, \circ_2, ..., \circ_N\}$ be a neutrosophic $N$-ring. A proper subset $T = \{\mathrm{T}_1 \cup T_2 \cup .... \cup \mathrm{T}_N, *_1, *_2, ..., *_N, \circ_1, \circ_2, ..., \circ_N\}$ of $N(R)$ is said to be a neutrosophic strong sub $N$-ring if each $(T_i, *_i)$ is a neutrosophic subring of $(R_i, *_i, \circ_i)$ for $i = 1, 2, ..., N$ where $T_i = R_i \cap T$.

**Definition 1.2.16.** Let $N(R) = \{R_1 \cup R_2 \cup .... \cup R_N, *_1, *_2, ..., *_N, \circ_1, \circ_2, ..., \circ_N\}$ be a neutrosophic $N$-ring. A proper subset $P = \{P_1 \cup P_2 \cup .... \cup P_N, *_1, *_2, ..., *_N, \circ_1, \circ_2, ..., \circ_N\}$ where $P_t = P \cap R_t$ for $t = 1, 2, ..., N$ is said to be a neutrosophic $N$-ideal of $N(R)$ if the following conditions are satisfied.

1. Each it is a neutrosophic subring of $R_t, t = 1, 2, ..., N$.
2. Each it is a two sided ideal of $R_t$ for $t = 1, 2, ..., N$.

If $(P_i, *_i, \circ_i)$ in the above definition are neutrosophic ideals, then we call $(P_i, *_i, \circ_i)$ to be a neutrosophic strong N-ideal of $N(R)$.

**Definition 1.2.17.** Let $N(R) = \{R_1 \cup R_2 \cup .... \cup R_N, *_1, *_2, ..., *_N, \circ_1, \circ_2, ..., \circ_N\}$ be a neutrosophic $N$-ring. A proper subset $P = \{P_1 \cup P_2 \cup .... \cup P_N, *_1, *_2, ..., *_N, \circ_1, \circ_2, ..., \circ_N\}$ where $P_t = P \cap R_t$ for $t = 1, 2, ..., N$ is said to be a pseudo neutrosophic $N$-ideal of $N(R)$ if the following conditions are satisfied.

1. Each it is a neutrosophic subring of $R_t, t = 1, 2, ..., N$.
2. Each $(P_i, *_i, \circ_i)$ is a pseudo neutrosophic ideal.

In the next section, we give the basic and fundamental definitions and notions about neutrosophic fields and their generalization.

## 1.3   Neutrosophic Fields, Neutrosophic Bifields, Neutrosophic N-fields and their Properties

In this section, we give the definitions of neutrosophic fields, neutrosophic bifields and neutrosophic N-fields.

**Neutrosophic Field**

**Definition 1.3.1.** Let $K$ be a field . We call the field generated by $K \cup I$ to be the neutrosophic field for it involves the indeterminacy factor in it. We define $I^2 = I$ , $I + I = 2I$ i.e. $I + I +, ..., + I = nI$ , and if $k \in K$ then $k.I = kI, 0I = 0$ . We denote the neutrosophic field by $K(I)$ which is generated by $K \cup I$ that is $K(I) = \langle K \cup I \rangle$ . $\langle K \cup I \rangle$ denotes the field generated by $K$ and $I$ .

**Example 1.3.1.** Let $\mathbb{C}$ be a field of complex numbers. Then $\mathbb{C}(I) = \langle \mathbb{C} \cup I \rangle$ is a neutrosophic field of complex numbers.

**Definition 1.3.2.** Let $K(I)$ be a neutrosophic field, $P \subset K(I)$ is a neutrosophic subfield of $K(I)$ if $P$ itself is a neutrosophic field under the operation of $K(I)$ .

**Neutrosophic Bi-field**

**Definition 1.3.3.** Let $(BN(\text{F}), *, \circ)$ be a non-empty set with two binary operations $*$ and $\circ$ . $(BN(\text{F}), *, \circ)$ is said to be a neutrosophic bifield if

$BN(\text{F}) = F_1 \cup F_2$ where at least one of $(\text{F}_1, *, \circ)$ or $(\text{F}_2, *, \circ)$ is a neutrosophic field and other is just a field. $F_1$ and $F_2$ are proper subsets of $BN(\text{F})$.

If in the above definition both $(\text{F}_1, *, \circ)$ and $(\text{F}_2, *, \circ)$ are neutrosophic fields, then we call $(BN(\text{F}), *, \circ)$ to be a neutrosophic strong bifield.

**Example 1.3.2.** Let $BN(\text{F}) = (\text{F}_1, *, \circ) \cup (\text{F}_2, *, \circ)$ where $(\text{F}_1, *, \circ) = (\langle \mathbb{C} \cup I \rangle, +, \times)$ and $(\text{F}_2, *, \circ) = (\mathbb{Q}, +, \times)$. Clearly $(\text{F}_1, *, \circ)$ is a neutrosophic field and $(\text{F}_2, *, \circ)$ is just a field. Thus $(BN(\text{F}), *, \circ)$ is a neutrosophic bifield.

**Example 1.3.3.** Let $BN(\text{F}) = (\text{F}_1, *, \circ) \cup (\text{F}_2, *, \circ)$ where $(\text{F}_1, *, \circ) = (\langle \mathbb{C} \cup I \rangle, +, \times)$ and $(F_2, *, \circ) = (\langle \mathbb{Q} \cup I \rangle, +, \times)$. Clearly $(\text{F}_1, *, \circ)$ is a neutrosophic field and $(\text{F}_2, *, \circ)$ is also a neutrosophic field. Thus $(BN(\text{F}), *, \circ)$ is a neutrosophic strong bifield.

**Definition 1.3.4.** Let $BN(\text{F}) = (\text{F}_1 \cup F_2, *, \circ)$ be a neutrosophic bifield. A proper subset $(T, *, \circ)$ is said to be a neutrosophic subbifield of $BN(\text{F})$ if

1. $T = T_1 \cup T_2$ where $T_1 = F_1 \cap T$ and $T_2 = F_2 \cap T$ and
2. At least one of $(T_1, \circ)$ or $(T_2, *)$ is a neutrosophic field and the other is just a field.

**Example 1.3.4.** Let $BN(\text{F}) = (\text{F}_1, *, \circ) \cup (\text{F}_2, *, \circ)$ where $(\text{F}_1, *, \circ) = (\langle \mathbb{R} \cup I \rangle, +, \times)$ and $(\text{F}_2, *, \circ) = (\mathbb{C}, +, \times)$. Let $P = P_1 \cup P_2$ be a proper subset of $BN(\text{F})$, where $P_1 = (\mathbb{Q}, +, \times)$ and $P_2 = (\mathbb{R}, +, \times)$. Clearly $(P, +, \times)$ is a neutrosophic subbifield of $BN(\text{F})$.

**Neutrosophic N-field**

**Definition 1.3.5.** Let $\{N(F), *_1, ..., *_2, \circ_1, \circ_2, ..., \circ_N\}$ be a non-empty set with $N$ -binary operations defined on it. We call $N(F)$ a neutrosophic $N$ -field ( $N$ a positive integer) if the following conditions are satisfied.

1. $N(F) = F_1 \cup F_2 \cup ... \cup F_N$ where each $F_i$ is a proper subset of $N(F)$ i.e. $R_i \not\subset R_j$ or $R_j \not\subset R_i$ if $i \neq j$ .

2. $(F_i, *_i, \circ_i)$ is either a neutrosophic field or just a field for $i = 1, 2, 3, ..., N$ .

If in the above definition each $(F_i, *_i, \circ_i)$ is a neutrosophic field, then we call $N(F)$ to be a neutrosophic strong N-field.

**Example 1.3.5.** Let $N(F) = \{F_1 \cup F_2 \cup F_3, *_1, *_2, *_3, \circ_1, \circ_2, \circ_3\}$ where $(F_1, *_1, \circ_1) = (\langle \mathbb{R} \cup I \rangle, +, \times)$ , $(F_2, *_2, \circ_2) = (\langle \mathbb{C} \cup I \rangle, +, \times)$ and $(F_3, *_3, \circ_3) = (\langle \mathbb{Q} \cup I \rangle, +, \times)$ . Clearly $N(F)$ is a neutrosophic 3-field.

**Definition 1.3.6.** Let $N(F) = \{F_1 \cup F_2 \cup .... \cup F_N, *_1, *_2, ..., *_N, \circ_1, \circ_2, ..., \circ_N\}$ be a neutrosophic $N$ -field. A proper subset $T = \{T_1 \cup T_2 \cup .... \cup T_N, *_1, *_2, ..., *_N, \circ_1, \circ_2, ..., \circ_N\}$ of $N(F)$ is said to be a neutrosophic $N$ -subfield if each $(T_i, *_i)$ is a neutrosophic subfield of $(F_i, *_i, \circ_i)$ for $i = 1, 2, ..., N$ where $T_i = F_i \cap T$ .

Next we give some basic material of neutrosophic group ring and their related for the readers to understand the theory of neutrosophic group ring.

## 1.4  Neutrosophic Group Ring, Neutrosophic Bigroup Biring and Neutrosophic N-group N-ring

**Neutrosophic Group Ring**

In this section, we define the neutrosophic group ring and displayed some of their basic and fundamental properties and notions.

**Definition 1.4.1.** Let $\langle G \cup I \rangle$ be any neutrosophic group and $R$ be any ring with $1$ which is commutative or field. We define the neutrosophic group ring $R\langle G \cup I \rangle$ of the neutrosophic group $\langle G \cup I \rangle$ over the ring $R$ as follows:

1. $R\langle G \cup I \rangle$ consists of all finite formal sum of the form $\alpha = \sum_{i=1}^{n} r_i g_i$, $n < \infty$, $r_i \in R$ and $g_i \in \langle G \cup I \rangle$ $\left( \alpha \in R\langle G \cup I \rangle \right)$.

2. Two elements $\alpha = \sum_{i=1}^{n} r_i g_i$ and $\beta = \sum_{i=1}^{m} s_i g_i$ in $R\langle G \cup I \rangle$ are equal if and only if $r_i = s_i$ and $n = m$.

3. Let $\alpha = \sum_{i=1}^{n} r_i g_i, \beta = \sum_{i=1}^{m} s_i g_i \in R\langle G \cup I \rangle$; $\alpha + \beta = \sum_{i=1}^{n} (\alpha_i + \beta_i) g_i \in R\langle G \cup I \rangle$, as $\alpha_i, \beta_i \in R$, so $\alpha_i + \beta_i \in R$ and $g_i \in \langle G \cup I \rangle$.

4. $0 = \sum_{i=1}^{n} 0 g_i$ serve as the zero of $R\langle G \cup I \rangle$.

5. Let $\alpha = \sum_{i=1}^{n} r_i g_i \in R\langle G \cup I \rangle$ then $-\alpha = \sum_{i=1}^{n} (-\alpha_i) g_i$ is such that
$$\begin{aligned} \alpha + (-\alpha) &= 0 \\ &= \sum_{i=1}^{n} (\alpha_i + (-\alpha_i)) g_i \\ &= \sum 0 g_i \end{aligned}$$
Thus we see that $R\langle G \cup I \rangle$ is an abelian group under $+$.

6. The product of two elements $\alpha, \beta$ in $R\langle G \cup I \rangle$ is follows:

Let $\alpha = \sum_{i=1}^{n} \alpha_i g_i$ and $\beta = \sum_{j=1}^{m} \beta_j h_j$. Then $\alpha.\beta = \sum_{\substack{1 \leq i \leq n \\ 1 \leq j \leq m}} \alpha_i.\beta_j g_i h_j$

$$= \sum_k y_k t_k$$

Where $y_k = \sum \alpha_i \beta_j$ with $g_i h_j = t_k$, $t_k \in \langle G \cup I \rangle$ and $y_k \in R$.

Clearly $\alpha.\beta \in R\langle G \cup I \rangle$.

7. Let $\alpha = \sum_{i=1}^{n} \alpha_i g_i$ and $\beta = \sum_{j=1}^{m} \beta_j h_j$ and $\gamma = \sum_{k=1}^{p} \delta_k l_k$.

Then clearly $\alpha(\beta + \gamma) = \alpha\beta + \alpha\gamma$ and $(\beta + \gamma)\alpha = \beta\alpha + \gamma\alpha$ for all $\alpha, \beta, \gamma \in R\langle G \cup I \rangle$.

Hence $R\langle G \cup I \rangle$ is a ring under the binary operations $+$ and $\cdot$.

**Example 1.4.1.** Let $\mathbb{Q}\langle G \cup I \rangle$ be a neutrosophic group ring, where $\mathbb{Q} =$ field of rationals and $\langle G \cup I \rangle = \{1, g, g^2, g^3, g^4, g^5, I, gI, ..., g^5 I : g^6 = 1, I^2 = I\}$.

**Example 1.4.2.** Let $\mathbb{C}\langle G \cup I \rangle$ be a neutrosophic group ring, where $\mathbb{C} =$ field of complex numbers and

$$\langle G \cup I \rangle = \{1, g, g^2, g^3, g^4, g^5, I, gI, ..., g^5 I : g^6 = 1, I^2 = I\}.$$

**Definition 1.4.2.** Let $R\langle G \cup I \rangle$ be a neutrosophic group ring and let $P$ be a proper subset of $R\langle G \cup I \rangle$. Then $P$ is called a subneutrosophic group ring of $R\langle G \cup I \rangle$ if $P = R\langle H \cup I \rangle$ or $S\langle G \cup I \rangle$ or $T\langle H \cup I \rangle$. In

$P = R\langle H \cup I \rangle$, $R$ is a ring and $\langle H \cup I \rangle$ is a proper neutrosophic subgroup of $\langle G \cup I \rangle$ or in $S\langle G \cup I \rangle$, $S$ is a proper subring with 1 of $R$ and $\langle G \cup I \rangle$ is a neutrosophic group and if $P = T\langle H \cup I \rangle$, $T$ is a subring of $R$ with unity and $\langle H \cup I \rangle$ is a proper neutrosophic subgroup of $\langle G \cup I \rangle$.

**Definition 1.4.3.** Let $R\langle G \cup I \rangle$ be a neutrosophic group ring. A proper subset $P$ of $R\langle G \cup I \rangle$ is called a neutrosophic subring if $P = \langle S \cup I \rangle$ where $S$ is a subring of $RG$ or $R$.

**Definition 1.4.4.** Let $R\langle G \cup I \rangle$ be a neutrosophic group ring. A proper subset $T$ of $R\langle G \cup I \rangle$ which is a pseudo neutrosophic subring. Then we call $T$ to be a pseudo neutrosophic subring.

**Definition 1.4.5.** Let $R\langle G \cup I \rangle$ be a neutrosophic group ring. A proper subset $P$ of $R\langle G \cup I \rangle$ is called a subgroup ring if $P = SH$ where $S$ is a subring of $R$ and $H$ is a subgroup of $G$. $SH$ is the group ring of the subgroup $H$ over the subring $S$.

**Definition 1.4.6.** Let $R\langle G \cup I \rangle$ be a neutrosophic group ring. A proper subset $P$ of $R\langle G \cup I \rangle$ is called a subring but $P$ should not have the group ring structure is defined to be a subring of $R\langle G \cup I \rangle$.

**Definition 1.4.7.** Let $R\langle G \cup I \rangle$ be a neutrosophic group ring. A proper subset $P$ of $R\langle G \cup I \rangle$ is called a neutrosophic ideal of $R\langle G \cup I \rangle$,

1. if $P$ is a neutrosophic subring or subneutrosophic group ring of $R\langle G \cup I\rangle$.

2. For all $p \in P$ and $\alpha \in R\langle G \cup I\rangle$, $\alpha p$ and $p\alpha \in P$.

**Definition 1.4.8.** Let $R\langle G \cup I\rangle$ be a neutrosophic group ring. A proper subset $P$ of $R\langle G \cup I\rangle$ is called a neutrosophic ideal of $R\langle G \cup I\rangle$,

1. if $P$ is a pseudo neutrosophic subring or pseudo subneutrosophic group ring of $R\langle G \cup I\rangle$.

2. For all $p \in P$ and $\alpha \in R\langle G \cup I\rangle$, $\alpha p$ and $p\alpha \in P$.

Further we give some basic concepts about the neutrosophic bigroup biring for the interested readers.

**Neutrosophic Bigroup Biring**

Here the definitions and notions about neutrosophic bigroup biring is introduced to develop some idea on it. We now proceed to define neutrosophic bigroup biring as follows.

**Definition 1.4.9.** Let $R_B\langle G \cup I\rangle = \left\{ R_1\langle G_1 \cup I\rangle \cup R_2\langle G_2 \cup I\rangle, *_1, *_2 \right\}$ be a non-empty set with two binary operations on $R_B\langle G \cup I\rangle$. Then $R_B\langle G \cup I\rangle$ is called a neutrosophic bigroup biring if

1. $R_B \langle G \cup I \rangle = R_1 \langle G_1 \cup I \rangle \cup R_2 \langle G_2 \cup I \rangle$ , where $R_1 \langle G_1 \cup I \rangle$ and $R_2 \langle G_2 \cup I \rangle$ are proper subsets of $R_B \langle G \cup I \rangle$ .

2. $\left( R_1 \langle G_1 \cup I \rangle, *_1, *_2 \right)$ be a neutrosophic group ring and

3. $\left( R_2 \langle G_2 \cup I \rangle, *_1, *_2 \right)$ be a neutrosophic group ring.

**Example 1.4.3.** Let $R_B \langle G \cup I \rangle = \mathbb{Q} \langle G_1 \cup I \rangle \cup \mathbb{R} \langle G_2 \cup I \rangle$ be a neutrosophic bigroup biring, where $R_B = \mathbb{Q} \cup \mathbb{R}$ and $\langle G_1 \cup I \rangle = \{1, g, g^2, g^3, g^4, g^5, I, gI, ..., g^5 I : g^6 = 1, I^2 = I\}$ and $\langle G_2 \cup I \rangle = \{1, g, g^2, g^3, I, gI, g^2 I, g^3 I : g^4 = 1, I^2 = I\}$ .

**Example 1.4.4.** Let $R_B \langle G \cup I \rangle = \mathbb{R} \langle G_1 \cup I \rangle \cup \mathbb{C} \langle G_2 \cup I \rangle$ be a neutrosophic bigroup biring, where $\langle G_1 \cup I \rangle = \{1, g, g^2, g^3, g^4, g^5, I, gI, ..., g^5 I : g^6 = 1, I^2 = I\}$ and $\langle G_2 \cup I \rangle = \{1, g, g^2, g^3, I, gI, g^2 I, g^3 I : g^4 = 1, I^2 = I\}$ .

**Definition 1.4.10.** Let $R_B \langle G \cup I \rangle = \{ R_1 \langle G_1 \cup I \rangle \cup R_2 \langle G_2 \cup I \rangle, *_1, *_2 \}$ be a neutrosophic bigroup biring and let $P = \{ R_1 \langle H_1 \cup I \rangle \cup R_2 \langle H_2 \cup I \rangle \}$ . Then $P$ is called subneutrosophic bigroup biring if $\left( R_1 \langle H_1 \cup I \rangle, *_1, *_2 \right)$ is a subneutrosophic group ring of $\left( R_1 \langle G_1 \cup I \rangle, *_1, *_2 \right)$ and $\left( R_2 \langle H_2 \cup I \rangle, *_1, *_2 \right)$ is also a subneutrosophic group ring of $\left( R_2 \langle G_2 \cup I \rangle, *_1, *_2 \right)$ .

**Definition 1.4.11.** Let $R_B \langle G \cup I \rangle = \{ R_1 \langle G_1 \cup I \rangle \cup R_2 \langle G_2 \cup I \rangle, *_1, *_2 \}$ be a neutrosophic bigroup biring and let $P = \{ \langle S_1 \cup I \rangle \cup \langle S_2 \cup I \rangle \}$ . Then $P$ is called neutrosophic subbiring if $\left( \langle S_1 \cup I \rangle, *_1, *_2 \right)$ is a neutrosophic subring of $\left( R_1 \langle G_1 \cup I \rangle, *_1, *_2 \right)$ and $\left( \langle S_2 \cup I \rangle, *_1, *_2 \right)$ is also a neutrosophic subring of $\left( R_2 \langle G_2 \cup I \rangle, *_1, *_2 \right)$.

**Definition 1.4.12.** Let $R_B \langle G \cup I \rangle = \left\{ R_1 \langle G_1 \cup I \rangle \cup R_2 \langle G_2 \cup I \rangle, *_1, *_2 \right\}$ be a neutrosophic bigroup biring and let $P = \left\{ R_1 \langle H_1 \cup I \rangle \cup R_2 \langle H_2 \cup I \rangle \right\}$. Then $P$ is called pseudo neutrosophic subbiring if $\left( R_1 \langle H_1 \cup I \rangle, *_1, *_2 \right)$ is a pseudo neutrosophic subring of $\left( R_1 \langle G_1 \cup I \rangle, *_1, *_2 \right)$ and $\left( R_2 \langle H_2 \cup I \rangle, *_1, *_2 \right)$ is also a pseudo neutrosophic subring of $\left( R_2 \langle G_2 \cup I \rangle, *_1, *_2 \right)$.

**Definition 1.4.13.** Let $R_B \langle G \cup I \rangle = \left\{ R_1 \langle G_1 \cup I \rangle \cup R_2 \langle G_2 \cup I \rangle, *_1, *_2 \right\}$ be a neutrosophic bigroup biring and let $P = \left\{ R_1 \langle H_1 \cup I \rangle \cup R_2 \langle H_2 \cup I \rangle \right\}$. Then $P$ is called neutrosophic subbigroup ring if $\left( R_1 H_1, *_1, *_2 \right)$ is a neutrosophic subgroup ring of $\left( R_1 \langle G_1 \cup I \rangle, *_1, *_2 \right)$ and $\left( R_2 H_2, *_1, *_2 \right)$ is also a neutrosophic subgroup ring of $\left( R_2 \langle G_2 \cup I \rangle, *_1, *_2 \right)$.

**Definition 1.4.14.** Let $R_B \langle G \cup I \rangle = \left\{ R_1 \langle G_1 \cup I \rangle \cup R_2 \langle G_2 \cup I \rangle, *_1, *_2 \right\}$ be a neutrosophic bigroup biring and let $P = \left\{ P_1 \cup P_2 \right\}$. Then $P$ is called subbiring if $\left( P_1, *_1, *_2 \right)$ is a subring of $\left( R_1 \langle G_1 \cup I \rangle, *_1, *_2 \right)$ and $\left( P_2, *_1, *_2 \right)$ is also a subring of $\left( R_2 \langle G_2 \cup I \rangle, *_1, *_2 \right)$.

**Definition 1.4.15.** Let $R_B \langle G \cup I \rangle = \left\{ R_1 \langle G_1 \cup I \rangle \cup R_2 \langle G_2 \cup I \rangle, *_1, *_2 \right\}$ be a neutrosophic bigroup biring and let $J = \left\{ J_1 \cup J_2 \right\}$. Then $J$ is called neutrosophic biideal if $(J_1, *_1, *_2)$ is a neutrosophic ideal of $(R_1 \langle G_1 \cup I \rangle, *_1, *_2)$ and $(J_2, *_1, *_2)$ is also a neutrosophic ideal of $(R_2 \langle G_2 \cup I \rangle, *_1, *_2)$.

**Definition 1.4.16.** Let $R_B \langle G \cup I \rangle = \{ R_1 \langle G_1 \cup I \rangle \cup R_2 \langle G_2 \cup I \rangle, *_1, *_2 \}$ be a neutrosophic bigroup biring and let $J = \{ J_1 \cup J_2 \}$. Then $J$ is called pseudo neutrosophic biideal if $(J_1, *_1, *_2)$ is a pseudo neutrosophic ideal of $(R_1 \langle G_1 \cup I \rangle, *_1, *_2)$ and $(J_2, *_1, *_2)$ is also a pseudo neutrosophic ideal of $(R_2 \langle G_2 \cup I \rangle, *_1, *_2)$.

In the next, we generalize this concept and give basic definition of neutrosophic N-group N-ring.

## Neutrosophic N-group N-ring

We now give the definition of neutrosophic N-group N-ring and give some of their related properties.

**Definition 1.4.17.** Let $N(R \langle G \cup I \rangle) = \{ R_1 \langle G_1 \cup I \rangle \cup R_2 \langle G_2 \cup I \rangle \cup ... \cup R_n \langle G_n \cup I \rangle, *_1, *_2, ..., *_n \}$ be a non-empty set with $n$ binary operations on $N(R \langle G \cup I \rangle)$. Then $N(R \langle G \cup I \rangle)$ is called a neutrosophic N-group N-ring if

1. $N\left(R\langle G\cup I\rangle\right)=R_1\langle G_1\cup I\rangle\cup R_2\langle G_2\cup I\rangle\cup...\cup R_n\langle G_n\cup I\rangle$, where $R_i\langle G_i\cup I\rangle$ is a proper subset of $N\left(R\langle G\cup I\rangle\right)$ for all $i$.

2. $\left(R_i\langle G_i\cup I\rangle,*_i,*_i\right)$ be a neutrosophic group ring for all $i$.

**Example 1.4.5.** Let $N\left(R\langle G\cup I\rangle\right)=\mathbb{Q}\langle G_1\cup I\rangle\cup\mathbb{R}\langle G_2\cup I\rangle\cup\mathbb{Z}\langle G_3\cup I\rangle$ be a neutrosophic 3-group 3-ring, where $R=\mathbb{Q}\cup\mathbb{R}\cup\mathbb{Z}$ and
$\langle G_1\cup I\rangle=\{1,g,g^2,g^3,g^4,g^5,I,gI,...,g^5I:g^6=1,I^2=I\}$,
$\langle G_2\cup I\rangle=\{1,g,g^2,g^3,I,gI,g^2I,g^3I:g^4=1,I^2=I\}$ and
$\langle G_3\cup I\rangle=\{1,g,g^2,g^3,g^4,...g^7,I,gI,g^2I,...,g^7I:g^8=1,I^2=I\}$.

**Example 1.4.6.** Let $N\left(R\langle G\cup I\rangle\right)=\mathbb{R}\langle G_1\cup I\rangle\cup\mathbb{C}\langle G_2\cup I\rangle\cup\mathbb{Z}\langle G_3\cup I\rangle$ be a neutrosophic 3-group 3-ring, where
$\langle G_1\cup I\rangle=\{1,g,g^2,g^3,g^4,g^5,I,gI,...,g^5I:g^6=1,I^2=I\}$,
$\langle G_2\cup I\rangle=\{1,g,g^2,g^3,I,gI,g^2I,g^3I:g^4=1,I^2=I\}$ and
$\langle G_3\cup I\rangle=\{1,g,g^2,g^3,g^4,g^5,g^6,g^7,I,gI,...,g^7I:g^8=1,I^2=I\}$.

**Definition 1.4.18.** Let
$N\left(R\langle G\cup I\rangle\right)=\{R_1\langle G_1\cup I\rangle\cup R_2\langle G_2\cup I\rangle\cup...\cup R_n\langle G_n\cup I\rangle,*_1,*_2,...,*_n\}$ be a neutrosophic N-group N-ring and let
$P=\{R_1\langle H_1\cup I\rangle\cup R_2\langle H_2\cup I\rangle\cup...\cup R_n\langle H_n\cup I\rangle,*_1,*_2,...,*_n\}$ be a proper subset of

$N\left(R\langle G\cup I\rangle\right)$. Then $P$ is called subneutrosophic N-group N-ring if $\left(R_i\langle H_i\cup I\rangle, *_i, *_i\right)$ is a subneutrosophic group ring of $N\left(R\langle G\cup I\rangle\right)$ for all $i$.

**Definition 1.4.19.** Let
$N\left(R\langle G\cup I\rangle\right) = \left\{R_1\langle G_1\cup I\rangle \cup R_2\langle G_2\cup I\rangle \cup...\cup R_n\langle G_n\cup I\rangle, *_1, *_2,...,*_n\right\}$ be a neutrosophic N-group N-ring and let $P = \left\{\langle S_1\cup I\rangle \cup \langle S_2\cup I\rangle \cup...\cup \langle S_n\cup I\rangle\right\}$. Then $P$ is called neutrosophic sub N-ring if $\left(\langle S_i\cup I\rangle, *_i, *_i\right)$ is a neutrosophic sub N-ring of $\left(R_i\langle G_i\cup I\rangle, *_i, *_i\right)$ for all $i$.

**Definition 1.4.20.** Let
$N\left(R\langle G\cup I\rangle\right) = \left\{R_1\langle G_1\cup I\rangle \cup R_2\langle G_2\cup I\rangle \cup...\cup R_n\langle G_n\cup I\rangle, *_1, *_2,...,*_n\right\}$ be a neutrosophic N-group N-ring and let
$P = \left\{R_1\langle H_1\cup I\rangle \cup R_2\langle H_2\cup I\rangle \cup...\cup R_n\langle H_n\cup I\rangle\right\}$. Then $P$ is called pseudo neutrosophic subbiring if $\left(R_i\langle H_i\cup I\rangle, *_i, *_i\right)$ is a pseudo neutrosophic sub N-ring of $\left(R_i\langle G_i\cup I\rangle, *_i, *_i\right)$ for all $i$

**Definition 1.4.21.** Let
$N\left(R\langle G\cup I\rangle\right) = \left\{R_1\langle G_1\cup I\rangle \cup R_2\langle G_2\cup I\rangle \cup...\cup R_n\langle G_n\cup I\rangle, *_1, *_2,...,*_n\right\}$ be a neutrosophic N-group N-ring and let $P = \left\{R_1 H_1 \cup R_2 H_2 \cup...\cup R_n H_n\right\}$. Then $P$ is called neutrosophic sub N-group ring if $\left(R_i H_i, *_i, *_i\right)$ is a neutrosophic subgroup ring of $\left(R_i\langle G_i\cup I\rangle, *_i, *_i\right)$ for all $i$.

**Definition 1.4.22.** Let
$N\left(R\langle G\cup I\rangle\right) = \left\{R_1\langle G_1\cup I\rangle \cup R_2\langle G_2\cup I\rangle \cup...\cup R_n\langle G_n\cup I\rangle, *_1, *_2,...,*_n\right\}$ be a

neutrosophic N-group N-ring and let $P = \{P_1 \cup P_2 \cup ... \cup P_n\}$. Then $P$ is called sub N-ring if $(P_i, *_i, *_i)$ is a subring of $(R_i \langle G_i \cup I \rangle, *_i, *_i)$ for all $i$.

**Definition 1.4.23.** Let
$N(R \langle G \cup I \rangle) = \{R_1 \langle G_1 \cup I \rangle \cup R_2 \langle G_2 \cup I \rangle \cup ... \cup R_n \langle G_n \cup I \rangle, *_1, *_2, ..., *_n\}$ be a neutrosophic N-group N-ring and let $J = \{J_1 \cup J_2 \cup ... \cup J_n\}$. Then $J$ is called neutrosophic N-ideal if $(J_i, *_i, *_i)$ is a neutrosophic ideal of $(R_i \langle G_i \cup I \rangle, *_i, *_i)$ for all $i$.

**Definition 1.4.24.** Let
$N(R \langle G \cup I \rangle) = \{R_1 \langle G_1 \cup I \rangle \cup R_2 \langle G_2 \cup I \rangle \cup ... \cup R_n \langle G_n \cup I \rangle, *_1, *_2, ..., *_n\}$ be a neutrosophic N-group N-ring and let $J = \{J_1 \cup J_2 \cup ... J_n\}$. Then $J$ is called pseudo neutrosophic N-ideal if $(J_i, *_i, *_i)$ is a pseudo neutrosophic ideal of $(R_i \langle G_i \cup I \rangle, *_i, *_i)$ for all $i$.

In the next section, we give some important definitions and notions about the theory of neutrosophic semigroup ring and their generalization.

## 1.5 Neutrosophic Semigroup Ring, Neutrosophic Bisemigroup, Neutrosophic N-semigroup and their Properties

**Neutrosophic Semigroup Ring**

In this section, we give the definitions and notions of neutrosophic semigroup rings and their generalization.

We now proceed to define neutrosophic semigroup ring.

**Definition 1.5.1.** Let $\langle S \cup I \rangle$ be any neutrosophic semigroup. $R$ be any ring with 1 which is commutative or field. We define the neutrosophic semigroup ring $R\langle S \cup I \rangle$ of the neutrosophic semigroup $\langle S \cup I \rangle$ over the ring $R$ as follows:

1. $R\langle S \cup I \rangle$ consists of all finite formal sum of the form $\alpha = \sum_{i=1}^{n} r_i g_i$, $n < \infty$, $r_i \in R$ and $g_i \in \langle S \cup I \rangle \left( \alpha \in R\langle S \cup I \rangle \right)$.

2. Two elements $\alpha = \sum_{i=1}^{n} r_i g_i$ and $\beta = \sum_{i=1}^{m} s_i g_i$ in $R\langle S \cup I \rangle$ are equal if and only if $r_i = s_i$ and $n = m$.

3. Let $\alpha = \sum_{i=1}^{n} r_i g_i, \beta = \sum_{i=1}^{m} s_i g_i \in R\langle S \cup I \rangle$. Then $\alpha + \beta = \sum_{i=1}^{n} (\alpha_i + \beta_i) g_i \in R\langle S \cup I \rangle$, as $\alpha_i, \beta_i \in R$, so $\alpha_i + \beta_i \in R$ and $g_i \in \langle S \cup I \rangle$.

4. $0 = \sum_{i=1}^{n} 0 g_i$ serve as the zero of $R\langle S \cup I \rangle$.

5. Let $\alpha = \sum_{i=1}^{n} r_i g_i \in R\langle S \cup I \rangle$ then $-\alpha = \sum_{i=1}^{n} (-\alpha_i) g_i$ is such that

$$\alpha + (-\alpha) = 0$$

$$= \sum_{i=1}^{n} \left( \alpha_i + (-\alpha_i) \right) g_i$$

$$= \sum 0 g_i$$

Thus we see that $R \langle S \cup I \rangle$ is an abelian group under $+$.

6. The product of two elements $\alpha, \beta$ in $R \langle S \cup I \rangle$ is follows:

Let $\alpha = \sum_{i=1}^{n} \alpha_i g_i$ and $\beta = \sum_{j=1}^{m} \beta_j h_j$. Then $\alpha.\beta = \sum_{\substack{1 \leq i \leq n \\ 1 \leq j \leq m}}^{n} \alpha_i.\beta_j g_i h_j$

$$= \sum_{k} y_k t_k$$

Where $y_k = \sum \alpha_i \beta_j$ with $g_i h_j = t_k$, $t_k \in \langle S \cup I \rangle$ and $y_k \in R$.

Clearly $\alpha.\beta \in R \langle S \cup I \rangle$.

7. Let $\alpha = \sum_{i=1}^{n} \alpha_i g_i$, $\beta = \sum_{j=1}^{m} \beta_j h_j$ and $\gamma = \sum_{k=1}^{p} \delta_k l_k$.

Then clearly $\alpha(\beta + \gamma) = \alpha\beta + \alpha\gamma$ and $(\beta + \gamma)\alpha = \beta\alpha + \gamma\alpha$ for all $\alpha, \beta, \gamma \in R \langle S \cup I \rangle$.

Hence $R \langle S \cup I \rangle$ is a ring under the binary operations $+$ and $\cdot$.

**Example 1.5.1.** Let $\mathbb{Q} \langle Z^+ \cup \{0\} \cup \{I\} \rangle$ be a neutrosophic semigroup ring, where $\mathbb{Q} =$ field of rationals and $\langle S \cup I \rangle = \langle Z^+ \cup \{0\} \cup \{I\} \rangle$ be a neutrosophic semigroup under $+$. Let $\langle H_1 \cup I \rangle = \langle 2Z^+ \cup I \rangle$ and $\langle H_2 \cup I \rangle = \langle 3Z^+ \cup I \rangle$.

**Example 1.5.2.** Let $\mathbb{C}\langle Z^+ \cup \{0\} \cup \{I\}\rangle$ be a neutrosophic semigroup ring, where $\mathbb{C}$ = field of complex numbers and $\langle Z^+ \cup \{0\} \cup \{I\}\rangle$ is a neutrosophic semigroup ring under $+$.

**Definition 1.5.2.** Let $R\langle S \cup I\rangle$ be a neutrosophic semigroup ring and let $P$ be a proper subset of $R\langle S \cup I\rangle$. Then $P$ is called a subneutrosophic semigroup ring of $R\langle S \cup I\rangle$ if $P = R\langle H \cup I\rangle$ or $Q\langle S \cup I\rangle$ or $T\langle H \cup I\rangle$. In $P = R\langle H \cup I\rangle$, $R$ is a ring and $\langle H \cup I\rangle$ is a proper neutrosophic subsemigroup of $\langle S \cup I\rangle$ or in $Q\langle S \cup I\rangle$, $Q$ is a proper subring with $1$ of $R$ and $\langle S \cup I\rangle$ is a neutrosophic semigroup and if $P = T\langle H \cup I\rangle$, $T$ is a subring of $R$ with unity and $\langle H \cup I\rangle$ is a proper neutrosophic subsemigroup of $\langle S \cup I\rangle$.

**Definition 1.5.3.** Let $R\langle S \cup I\rangle$ be a neutrosophic semigroup ring. A proper subset $P$ of $R\langle S \cup I\rangle$ is called a neutrosophic subring if $P = \langle S_1 \cup I\rangle$ where $S_1$ is a subring of $RS$ or $R$.

**Definition 1.5.4.** Let $R\langle S \cup I\rangle$ be a neutrosophic semigroup ring. A proper subset $T$ of $R\langle S \cup I\rangle$ which is a pseudo neutrosophic subring. Then we call $T$ to be a pseudo neutrosophic subring.

**Definition 1.5.5.** Let $R\langle S \cup I\rangle$ be a neutrosophic group ring. A proper subset $P$ of $R\langle G \cup I\rangle$ is called a subgroup ring if $P = SH$ where $S$ is a subring of $R$ and $H$ is a subgroup of $G$.

$SH$ is the group ring of the subgroup $H$ over the subring $S$.

**Definition 1.5.6.** Let $R\langle G \cup I \rangle$ be a neutrosophic semigroup ring. A proper subset $P$ of $R\langle S \cup I \rangle$ is called a subring but $P$ should not have the semigroup ring structure and is defined to be a subring of $R\langle S \cup I \rangle$.

**Definition 1.5.7.** Let $R\langle S \cup I \rangle$ be a neutrosophic semigroup ring. A proper subset $P$ of $R\langle S \cup I \rangle$ is called a neutrosophic ideal of $R\langle S \cup I \rangle$,

1. if $P$ is a neutrosophic subring or subneutrosophic semigroup ring of $R\langle S \cup I \rangle$.
2. For all $p \in P$ and $\alpha \in R\langle S \cup I \rangle$, $\alpha p$ and $p\alpha \in P$.

**Definition 1.5.8.** Let $R\langle S \cup I \rangle$ be a neutrosophic semigroup ring. A proper subset $P$ of $R\langle S \cup I \rangle$ is called a neutrosophic ideal of $R\langle S \cup I \rangle$,

1. if $P$ is a pseudo neutrosophic subring or pseudo subneutrosophic semigroup ring of $R\langle S \cup I \rangle$.
2. For all $p \in P$ and $\alpha \in R\langle S \cup I \rangle$, $\alpha p$ and $p\alpha \in P$.

**Neutrosophic Bisemigroup Biring**

Here we present some dfinitions and notions about neutrosophic bisemigroup biring.

**Definition 1.5.9.** Let $R_B \langle S \cup I \rangle = \{R_1 \langle S_1 \cup I \rangle \cup R_2 \langle S_2 \cup I \rangle, *_1, *_2\}$ be a non-empty set with two binary operations on $R_B \langle S \cup I \rangle$. Then $R_B \langle S \cup I \rangle$ is called a neutrosophic bisemigroup biring if

1. $R_B \langle S \cup I \rangle = R_1 \langle S_1 \cup I \rangle \cup R_2 \langle S_2 \cup I \rangle$, where $R_1 \langle S_1 \cup I \rangle$ and $R_2 \langle S_2 \cup I \rangle$ are proper subsets of $R_B \langle S \cup I \rangle$.
2. $(R_1 \langle S_1 \cup I \rangle, *_1, *_2)$ is a neutrosophic semigroup ring and
3. $(R_2 \langle S_2 \cup I \rangle, *_1, *_2)$ is a neutrosophic semigroup ring.

The following examples illustrates this situation.

**Example 1.5.3.** Let $R_B \langle S \cup I \rangle = \mathbb{Q} \langle S_1 \cup I \rangle \cup \mathbb{R} \langle S_2 \cup I \rangle$ be a neutrosophic bisemigroup biring, where $R_B = \mathbb{Q} \cup \mathbb{R}$ and $\langle S_1 \cup I \rangle = \{1, g, g^2, g^3, g^4, g^5, I, gI, ..., g^5 I : g^6 = 1, I^2 = I\}$ is a neutrosophic semigroup and $\langle S_2 \cup I \rangle = \langle Z^+ \setminus \{0\} \cup I \rangle$ be a neutrosophic semigroup under +.

**Example 1.5.4.** Let $R_B \langle S \cup I \rangle = \mathbb{R} \langle S_1 \cup I \rangle \cup \mathbb{C} \langle S_2 \cup I \rangle$ be a neutrosophic bisemigroup biring, where $\langle S_1 \cup I \rangle = \{1, g, g^2, g^3, g^4, g^5, I, gI, ..., g^5 I : g^6 = 1, I^2 = I\}$ and $\langle S_2 \cup I \rangle = \{1, g, g^2, g^3, I, gI, g^2 I, g^3 I : g^4 = 1, I^2 = I\}$ are neutrosophic semigroups.

**Definition 1.5.10.** Let $R_B \langle S \cup I \rangle = \{R_1 \langle S_1 \cup I \rangle \cup R_2 \langle S_2 \cup I \rangle, *_1, *_2\}$ be a neutrosophic bisemigroup biring and let $P = \{R_1 \langle H_1 \cup I \rangle \cup R_2 \langle H_2 \cup I \rangle\}$ be a proper subset of $R_B \langle S \cup I \rangle$. Then $P$ is called subneutrosophic bisemigroup biring if $\left(R_1\langle H_1 \cup I\rangle, *_1, *_2\right)$ is a subneutrosophic semigroup ring of $\left(R_1\langle S_1 \cup I\rangle, *_1, *_2\right)$ and $\left(R_2\langle H_2 \cup I\rangle, *_1, *_2\right)$ is also a subneutrosophic semigroup ring of $\left(R_2\langle S_2 \cup I\rangle, *_1, *_2\right)$.

**Definition 1.5.11.** Let $R_B\langle S \cup I\rangle = \left\{R_1\langle S_1 \cup I\rangle \cup R_2\langle S_2 \cup I\rangle, *_1, *_2\right\}$ be a neutrosophic bigroup biring and let $P = \left\{\langle S_1 \cup I\rangle \cup \langle S_2 \cup I\rangle\right\}$ be a proper subset of $R_B\langle S \cup I\rangle$. Then $P$ is called neutrosophic subbiring if $\left(\langle S_1 \cup I\rangle, *_1, *_2\right)$ is a neutrosophic subring of $\left(R_1\langle S_1 \cup I\rangle, *_1, *_2\right)$ and $\left(\langle S_2 \cup I\rangle, *_1, *_2\right)$ is also a neutrosophic subring of $\left(R_2\langle S_2 \cup I\rangle, *_1, *_2\right)$.

**Definition 1.5.12.** Let $R_B\langle S \cup I\rangle = \left\{R_1\langle S_1 \cup I\rangle \cup R_2\langle S_2 \cup I\rangle, *_1, *_2\right\}$ be a neutrosophic bisemigroup biring and let $P = \left\{R_1\langle H_1 \cup I\rangle \cup R_2\langle H_2 \cup I\rangle\right\}$ be a proper subset of $R_B\langle S \cup I\rangle$. Then $P$ is called pseudo neutrosophic subbiring if $\left(R_1\langle H_1 \cup I\rangle, *_1, *_2\right)$ is a pseudo neutrosophic subring of $\left(R_1\langle S_1 \cup I\rangle, *_1, *_2\right)$ and $\left(R_2\langle H_2 \cup I\rangle, *_1, *_2\right)$ is also a pseudo neutrosophic subring of $\left(R_2\langle S_2 \cup I\rangle, *_1, *_2\right)$.

**Definition 1.5.12.** Let $R_B\langle S \cup I\rangle = \left\{R_1\langle S_1 \cup I\rangle \cup R_2\langle S_2 \cup I\rangle, *_1, *_2\right\}$ be a neutrosophic bigroup biring and let $P = \left\{R_1\langle H_1 \cup I\rangle \cup R_2\langle H_2 \cup I\rangle\right\}$ be a proper subset of $R_B\langle S \cup I\rangle$. Then $P$ is called neutrosophic subbigroup ring if $\left(R_1 H_1, *_1, *_2\right)$ is a neutrosophic subgroup ring of $\left(R_1\langle S_1 \cup I\rangle, *_1, *_2\right)$ and $\left(R_2 H_2, *_1, *_2\right)$ is also a neutrosophic subgroup ring of $\left(R_2\langle S_2 \cup I\rangle, *_1, *_2\right)$.

**Definition 1.5.13.** Let $R_B \langle S \cup I \rangle = \{R_1 \langle S_1 \cup I \rangle \cup R_2 \langle S_2 \cup I \rangle, *_1, *_2 \}$ be a neutrosophic bisemigroup biring and let $P = \{P_1 \cup P_2 \}$ be a proper subset of $R_B \langle S \cup I \rangle$. Then $P$ is called subbiring if $(P_1, *_1, *_2)$ is a subring of $(R_1 \langle S_1 \cup I \rangle, *_1, *_2)$ and $(P_2, *_1, *_2)$ is also a subring of $(R_2 \langle S_2 \cup I \rangle, *_1, *_2)$.

**Definition 1.5.14.** Let $R_B \langle S \cup I \rangle = \{R_1 \langle S_1 \cup I \rangle \cup R_2 \langle S_2 \cup I \rangle, *_1, *_2 \}$ be a neutrosophic bisemigroup biring and let $J = \{J_1 \cup J_2 \}$ be a proper subset of $R_B \langle S \cup I \rangle$. Then $J$ is called neutrosophic biideal if $(J_1, *_1, *_2)$ is a neutrosophic ideal of $(R_1 \langle S_1 \cup I \rangle, *_1, *_2)$ and $(J_2, *_1, *_2)$ is also a neutrosophic ideal of $(R_2 \langle S_2 \cup I \rangle, *_1, *_2)$.

**Definition 1.5.15.** Let $R_B \langle S \cup I \rangle = \{R_1 \langle S_1 \cup I \rangle \cup R_2 \langle S_2 \cup I \rangle, *_1, *_2 \}$ be a neutrosophic bisemigroup biring and let $J = \{J_1 \cup J_2 \}$. Then $J$ is called pseudo neutrosophic biideal if $(J_1, *_1, *_2)$ is a pseudo neutrosophic ideal of $(R_1 \langle S_1 \cup I \rangle, *_1, *_2)$ and $(J_2, *_1, *_2)$ is also a pseudo neutrosophic ideal of $(R_2 \langle S_2 \cup I \rangle, *_1, *_2)$.

Now we generalize the concept of neutrosophic semigroup ring and give the definition of neutrosophic N-semigroup N-ring.

**Neutrosophic N-semigroup N-ring**

Here the definitions and other important notions about neutrosophic N-semigroup N-ring are introduced.

**Definition 1.5.16.** Let
$N(R\langle S \cup I\rangle) = \{R_1\langle S_1 \cup I\rangle \cup R_2\langle S_2 \cup I\rangle \cup ... \cup R_n\langle S_n \cup I\rangle, *_1, *_2, ..., *_n\}$ be a non-empty set with $n$ binary operations on $N(R\langle S \cup I\rangle)$. Then $N(R\langle S \cup I\rangle)$ is called a neutrosophic N-semigroup N-ring if

1. $N(R\langle S \cup I\rangle) = R_1\langle S_1 \cup I\rangle \cup R_2\langle S_2 \cup I\rangle \cup ... \cup R_n\langle S_n \cup I\rangle$, where $R_i\langle S_i \cup I\rangle$ is a proper subset of $N(R\langle S \cup I\rangle)$ for all $i$.

2. $(R_i\langle S_i \cup I\rangle, *_i, *_i)$ be a neutrosophic group ring for all $i$.

**Example 1.5.5.** Let $N(R\langle S \cup I\rangle) = \mathbb{Q}\langle S_1 \cup I\rangle \cup \mathbb{R}\langle S_2 \cup I\rangle \cup \mathbb{Z}\langle S_3 \cup I\rangle$ be a neutrosophic 3-semigroup 3-ring, where $R = \mathbb{Q} \cup \mathbb{R} \cup \mathbb{Z}$ and

$\langle S_1 \cup I\rangle = \{1, g, g^2, g^3, g^4, g^5, I, gI, ..., g^5 I : g^6 = 1, I^2 = I\}$,

$\langle S_2 \cup I\rangle = \{1, g, g^2, g^3, I, gI, g^2 I, g^3 I : g^4 = 1, I^2 = I\}$ and $\langle S_3 \cup I\rangle = \langle Z^+ \setminus \{0\} \cup I\rangle$ are neutrosophic semigroups.

**Example 1.5.6.** Let $N(R\langle S \cup I\rangle) = \mathbb{R}\langle S_1 \cup I\rangle \cup \mathbb{C}\langle S_2 \cup I\rangle \cup \mathbb{Z}\langle S_3 \cup I\rangle$ be a neutrosophic 3-semigroup 3-ring, where
$\langle S_1 \cup I\rangle = \{1, g, g^2, g^3, g^4, g^5, I, gI, ..., g^5 I : g^6 = 1, I^2 = I\}$,

$\langle S_2 \cup I\rangle = \{1, g, g^2, g^3, I, gI, g^2 I, g^3 I : g^4 = 1, I^2 = I\}$ and $\langle G_3 \cup I\rangle = \langle Z^+ \setminus \{0\} \cup I\rangle$ are neutrosophic semigroups.

**Definition 1.5.17.** Let

$N\left(R\langle S\cup I\rangle\right)=\left\{R_1\langle S_1\cup I\rangle\cup R_2\langle S_2\cup I\rangle\cup...\cup R_n\langle S_n\cup I\rangle,*_1,*_2,...,*_n\right\}$ be a neutrosophic N-semigroup N-ring and let

$P=\left\{R_1\langle H_1\cup I\rangle\cup R_2\langle H_2\cup I\rangle\cup...\cup R_n\langle H_n\cup I\rangle,*_1,*_2,...,*_n\right\}$ be a proper subset of $N\left(R\langle S\cup I\rangle\right)$. Then $P$ is called subneutrosophic N-semigroup N-ring if $\left(R_i\langle H_i\cup I\rangle,*_i,*_i\right)$ is a subneutrosophic semigroup ring of $N\left(R\langle S\cup I\rangle\right)$ for all $i$.

**Definition 1.5.18.** Let

$N\left(R\langle S\cup I\rangle\right)=\left\{R_1\langle S_1\cup I\rangle\cup R_2\langle S_2\cup I\rangle\cup...\cup R_n\langle S_n\cup I\rangle,*_1,*_2,...,*_n\right\}$ be a neutrosophic N-semigroup N-ring and let $P=\left\{\langle P_1\cup I\rangle\cup\langle P_2\cup I\rangle\cup...\cup\langle P_n\cup I\rangle\right\}$. Then $P$ is called neutrosophic sub N-ring if $\left(\langle P_i\cup I\rangle,*_i,*_i\right)$ is a neutrosophic sub N-ring of $\left(R_i\langle S_i\cup I\rangle,*_i,*_i\right)$ for all $i$.

**Definition 1.5.19.** Let

$N\left(R\langle S\cup I\rangle\right)=\left\{R_1\langle S_1\cup I\rangle\cup R_2\langle S_2\cup I\rangle\cup...\cup R_n\langle S_n\cup I\rangle,*_1,*_2,...,*_n\right\}$ be a neutrosophic N-semigroup N-ring and let $P=\left\{R_1\langle H_1\cup I\rangle\cup R_2\langle H_2\cup I\rangle\cup...\cup R_n\langle H_n\cup I\rangle\right\}$. Then $P$ is called pseudo neutrosophic subbiring if $\left(R_i\langle H_i\cup I\rangle,*_i,*_i\right)$ is a pseudo neutrosophic sub N-ring of $\left(R_i\langle S_i\cup I\rangle,*_i,*_i\right)$ for all $i$

**Definition 1.5.20.** Let

$N\left(R\langle S\cup I\rangle\right)=\left\{R_1\langle S_1\cup I\rangle\cup R_2\langle S_2\cup I\rangle\cup...\cup R_n\langle S_n\cup I\rangle,*_1,*_2,...,*_n\right\}$ be a neutrosophic N-semigroup N-ring and let $P=\left\{R_1H_1\cup R_2H_2\cup...\cup R_nH_n\right\}$. Then $P$ is called neutrosophic sub N-semigroup ring if $(R_i H_i, *_i, *_i)$ is a neutrosophic subgroup ring of $(R_i \langle S_i \cup I \rangle, *_i, *_i)$ for all $i$.

**Definition 1.5.21.** Let
$N\left(R\langle S \cup I\rangle\right) = \{R_1\langle S_1 \cup I\rangle \cup R_2\langle S_2 \cup I\rangle \cup ... \cup R_n\langle S_n \cup I\rangle, *_1, *_2, ..., *_n\}$ be a neutrosophic N-semigroup N-ring and let $P = \{P_1 \cup P_2 \cup ... \cup P_n\}$. Then $P$ is called sub N-ring if $(P_i, *_i, *_i)$ is a subring of $(R_i\langle S_i \cup I\rangle, *_i, *_i)$ for all $i$.

**Definition 1.5.22.** Let
$N\left(R\langle S \cup I\rangle\right) = \{R_1\langle S_1 \cup I\rangle \cup R_2\langle S_2 \cup I\rangle \cup ... \cup R_n\langle S_n \cup I\rangle, *_1, *_2, ..., *_n\}$ be a neutrosophic N-semigroup N-ring and let $J = \{J_1 \cup J_2 \cup ... \cup J_n\}$. Then $J$ is called neutrosophic N-ideal if $(J_i, *_i, *_i)$ is a neutrosophic ideal of $(R_i\langle S_i \cup I\rangle, *_i, *_i)$ for all $i$.

**Definition 1.5.23.** Let
$N\left(R\langle S \cup I\rangle\right) = \{R_1\langle S_1 \cup I\rangle \cup R_2\langle S_2 \cup I\rangle \cup ... \cup R_n\langle S_n \cup I\rangle, *_1, *_2, ..., *_n\}$ be a neutrosophic N-semigroup N-ring and let $J = \{J_1 \cup J_2 \cup ... J_n\}$. Then $J$ is called pseudo neutrosophic N-ideal if $(J_i, *_i, *_i)$ is a pseudo neutrosophic ideal of $(R_i\langle S_i \cup I\rangle, *_i, *_i)$ for all $i$.

In the next section of this chapter, we give some theory about mixed neutrosophic N-algebraic structures to develop some idea of the readers.

## 1.6 Mixed Neutrosophic N-algebraic structures

In this section, the authors define mixed neutrosophic N-algebraic structures and give examples for better illustration.

We now give the definition of mixed neutrosophic N-algebraic structure.

**Definition 1.6.1.** Let $\{\langle M \cup I \rangle = M_1 \cup M_2 \cup ... \cup M_N, *_1, *_2, ..., *_N\}$ such that $N \geq 5$. Then $\langle M \cup I \rangle$ is called a mixed neutrosophic $N$-algebraic structure if

1. $\langle M \cup I \rangle = M_1 \cup M_2 \cup ... \cup M_N$, where each $M_i$ is a proper subset of $\langle M \cup I \rangle$.

2. Some of $(M_i, *_i)$ are neutrosophic groups.

**49**

3. Some of $(M_j, *_j)$ are neutrosophic loops.

4. Some of $(M_k, *_k)$ are neutrosophic groupoids.

5. Some of $(M_r, *_r)$ are neutrosophic semigroups.

6. Rest of $(M_t, *_t)$ can be loops or groups or semigroups or groupoids. ( 'or' not used in the mutually exclusive sense).

The following example illustrates this fact.

**Example 1.6.1.** Let $\{\langle M \cup I \rangle = M_1 \cup M_2 \cup M_3 \cup M_4 \cup M_5, *_1, *_2, *_3, *_4, *_5\}$ be a mixed neutosophic $5$-structure, where $M_1 = \langle \mathbb{Z}_3 \cup I \rangle$, a neutrosophic group under multiplication modolu $3$, $M_2 = \langle \mathbb{Z}_6 \cup I \rangle$, a neutrosophic semigroup under multiplication modolu $6$, $M_3 = \{0,1,2,3,1I,2I,3I,$, a neutrosophic groupoid under multiplication modolu $4\}$, $M_4 = S_3$, and $M_5 = \{Z_{10},$ a semigroup under multiplication modolu $10\}$.

**Definition 1.6.2.** Let $\{\langle D \cup I \rangle = D_1 \cup D_2 \cup ... \cup D_N, *_1, *_2, ..., *_N\}$. Then $\langle D \cup I \rangle$ is called a mixed dual neutrosophic $N$-algebraic structure if

1. $\langle D \cup I \rangle = D_1 \cup D_2 \cup ... \cup D_N$, where each $D_i$ is a proper subset of $\langle D \cup I \rangle$.
2. Some of $(D_i, *_i)$ are groups.
3. Some of $(D_j, *_j)$ are loops.
4. Some of $(D_k, *_k)$ are groupoids.
5. Some of $(D_r, *_r)$ are semigroups.

6. Rest of $(D_t, *_t)$ can be neutrosophic loops or neutrosophic groups or neutrosophic semigroups or neutrosophic groupoids. ('or' not used in the mutually exclusive sense).

**Example 1.6.2.** Let $\{\langle D \cup I \rangle = D_1 \cup D_2 \cup D_3 \cup D_4 \cup D_5, *_1, *_2, *_3, *_4, *_5\}$ be a mixed dual neutosophic $5$-structure, where $D_1 = L_7(4)$, $D_2 = S_4$, $D_3 = \{Z_{10},$ a semigroup under multiplication modulo $10\}$, $D_4 = \{0,1,2,3,$ a groupoid under multiplication modulo $4\}, D_5 = \langle L_7(4) \cup I \rangle$.

**Definition 1.6.3** Let $\{\langle W \cup I \rangle = W_1 \cup W_2 \cup ... \cup W_N, *_1, *_2, ..., *_N\}$. Then $\langle W \cup I \rangle$ is called a weak mixed neutrosophic $N$-algebraic structure if

1. $\langle W \cup I \rangle = W_1 \cup W_2 \cup \ldots \cup W_N$, where each $W_i$ is a proper subset of $\langle W \cup I \rangle$.

2. Some of $(W_i, *_i)$ are neutrosophic groups or neutrosophic loops.

3. Some of $(W_j, *_j)$ are neutrosophic groupoids or neutrosophic semigroups.

4. Rest of $(W_t, *_t)$ are loops or groups or semigroups or groupoids. i.e in the collection $\{W_i, *_i\}$, all the 4 algebraic neutrosophic structures may not be present.

At most 3 algebraic neutrosophic structures are present and at least 2 algebraic neutrosophic structures are present. Rest being non-neutrosophic algebraic structures.

**Definition 1.6.4.** Let $\{\langle V \cup I \rangle = V_1 \cup V_2 \cup \ldots \cup V_N, *_1, *_2, \ldots, *_N\}$. Then $\langle V \cup I \rangle$ is called a weak mixed dual neutrosophic $N$-algebraic structure if

1. $\langle V \cup I \rangle = V_1 \cup V_2 \cup \ldots \cup V_N$, where each $V_i$ is a proper subset of $\langle V \cup I \rangle$.

**51**

2. Some of $(V_i, *_i)$ are groups or loops.

3. Some of $(V_j, *_j)$ are groupoids or semigroups.

4. Rest of $(V_t, *_t)$ are neutrosophic loops or neutrosophic groups or neutrosophic semigroups or neutrosophic groupoids.

**Definition 1.6.5.** Let $\{\langle M \cup I \rangle = M_1 \cup M_2 \cup \ldots \cup M_N, *_1, *_2, \ldots, *_N\}$ be a neutrosophic $N$-algebraic structure. A proper subset $\{\langle P \cup I \rangle = P_1 \cup P_2 \cup \ldots \cup P_N, *_1, *_2, \ldots, *_N\}$ is called a mixed neutrosophic sub $N$-structure if $\langle P \cup I \rangle$ is a mixed neutrosophic $N$-structure under the operation of $\langle M \cup I \rangle$.

**Definition 1.6.6.** Let $\{\langle W \cup I \rangle = W_1 \cup W_2 \cup ... \cup W_N, *_1, *_2, ..., *_N\}$ be a neutrosophic $N$-algebraic structure. We call a finite non-empty subset $P$ of $\langle W \cup I \rangle$ to be a weak mixed deficit neutrosophic sub $N$-algebraic structure if $\{P = P_1 \cup P_2 \cup ... \cup P_t, *_1, *_2, ..., *_t\}$, $1 < t < N$ with $P_i = P \cap L_k$, $1 \le i \le t$ and $1 \le k \le N$ and some $P_i's$ are neutrosophic groups or neutrosophic loops, some $P_j's$ are neutrosophic groupoids or neutrosophic semigroups and the rest of $P_k's$ are groups or loops or groupoids or semigroups.

**Definition 1.6.7.** Let $\{\langle M \cup I \rangle = M_1 \cup M_2 \cup ... \cup M_N, *_1, *_2, ..., *_N\}$ be a neutrosophic $N$-algebraic structure of finite order. A proper mixed neutrosophic sub $N$-structure $P$ of $\langle M \cup I \rangle$ is called Lagrange mixed neutrosophic sub $N$-structure if $\dfrac{o(P)}{o(\langle M \cup I \rangle)}$.

**Definition 1.6.8.** If every mixed neutrosophic sub $N$-structure of $\langle M \cup I \rangle$ is a Lagrange mixed neutrosophic sub $N$-structures. Then $\langle M \cup I \rangle$ is said to be a Lagrange mixed neutrosophic $N$-structure.

**Definition 1.6.9.** If some mixed neutrosophic sub $N$-structure of $\langle M \cup I \rangle$ are Lagrange mixed neutrosophic sub $N$-structures. Then $\langle M \cup I \rangle$ is said to be a weak Lagrange mixed neutrosophic $N$-structure.

**Definition 1.6.10.** If every mixed neutrosophic sub $N$-structure of $\langle M \cup I \rangle$ is not a Lagrange mixed neutrosophic sub $N$-structures. Then $\langle M \cup I \rangle$ is said to be a Lagrange free mixed neutrosophic $N$-structure.

In this final section of this chapter, we give some basic and fundamental theory of soft sets.

## 1.7   Soft Sets

Throughout this section $U$ refers to an initial universe, $E$ is a set of parameters, $P(U)$ is the power set of $U$, and $A, B \subset E$. Molodtsov defined the soft set in the following manner:

**Definition 1.7.1.** A pair $(F, A)$ is called a soft set over $U$ where $F$ is a mapping given by $F : A \to P(U)$.

In other words, a soft set over $U$ is a parameterized family of subsets of the universe $U$. For $a \in A$, $F(\mathrm{a})$ may be considered as the set of $a$-elements of the soft set $(F, A)$, or as the set of $a$-approximate elements of the soft set.

This situation can be explained in the following example.

**Example 1.7.1.** Suppose that $U$ is the set of shops. $E$ is the set of parameters and each parameter is a word or sentence. Let

$$E = \begin{Bmatrix} \text{high rent,normal rent,} \\ \text{in good condition,in bad condition} \end{Bmatrix}.$$

Let us consider a soft set $(F, A)$ which describes the attractiveness of shops that Mr. $Z$ is taking on rent. Suppose that there are five houses in the universe $U = \{s_1, s_2, s_3, s_4, s_5\}$ under consideration, and that $A = \{a_1, a_2, a_3\}$ be the set of parameters where

$a_1$ stands for the parameter 'high rent,

$a_2$ stands for the parameter 'normal rent,

$a_3$ stands for the parameter 'in good condition.

Suppose that

$$F(a_1) = \{s_1, s_4\} \ ,$$

$$F(a_2) = \{s_2, s_5\} \ ,$$

$$F(a_3) = \{s_3\}.$$

**54**

The soft set $(F, A)$ is an approximated family $\{F(a_i), i = 1, 2, 3\}$ of subsets of the set $U$ which gives us a collection of approximate description of an object. Then $(F, A)$ is a soft set as a collection of approximations over $U$, where

$$F(a_1) = high \ rent = \{s_1, s_2\},$$

$$F(a_2) = normal \ rent = \{s_2, s_5\},$$

$$F(a_3) = in \ good \ condition = \{s_3\}.$$

**Definition 1.7.2.** For two soft sets $(F, A)$ and $(H, B)$ over $U$, $(F, A)$ is called a soft subset of $(H, B)$ if

1. $A \subseteq B$ and
2. $F(a) \subseteq H(a)$, for all $a \in A$.

This relationship is denoted by $(F, A) \subset (H, B)$. Similarly $(F, A)$ is called a soft superset of $(H, B)$ if $(H, B)$ is a soft subset of $(F, A)$ which is denoted by $(F, A) \supset (H, B)$.

**Definition 1.7.3.** Two soft sets $(F, A)$ and $(H, B)$ over $U$ are called soft equal if $(F, A)$ is a soft subset of $(H, B)$ and $(H, B)$ is a soft subset of $(F, A)$.

**55**

**Definition 1.7.4.** Let $(F, A)$ and $(K, B)$ be two soft sets over a common universe $U$ such that $A \cap B \neq \phi$. Then their restricted intersection is denoted by $(F, A) \cap_R (K, B) = (H, C)$ where $(H, C)$ is defined as $H(c) = F(c) \cap K(c)$ for all $c \in C = A \cap B$.

**Definition 1.7.5.** The extended intersection of two soft sets $(F, A)$ and $(K, B)$ over a common universe $U$ is the soft set $(H, C)$, where $C = A \cup B$,

and for all $c \in C$ , $H(c)$ is defined as

$$H(c) = \begin{cases} F(c) & \text{if } c \in A - B, \\ G(c) & \text{if } c \in B - A, \\ F(c) \cap G(c) & \text{if } c \in A \cap B. \end{cases}$$

We write $(F,A) \cap_{\varepsilon} (K,B) = (H,C)$.

**Definition 1.7.6.** The restricted union of two soft sets $(F,A)$ and $(K,B)$ over a common universe $U$ is the soft set $(H,C)$, where $C = A \cup B$ , and for all $c \in C$ , $H(c)$ is defined as $H(c) = F(c) \cup G(c)$ for all $c \in C$ .
We write it as $(F,A) \cup_R (K,B) = (H,C)$.

**Definition 1.7.7.** The extended union of two soft sets $(F,A)$ and $(K,B)$ over a common universe $U$ is the soft set $(H,C)$, where $C = A \cup B$ , and for all $c \in C$ , $H(c)$ is defined as

$$H(c) = \begin{cases} F(c) & \text{if } c \in A - B, \\ G(c) & \text{if } c \in B - A, \\ F(c) \cup G(c) & \text{if } c \in A \cap B. \end{cases}$$

We write $(F,A) \cup_{\varepsilon} (K,B) = (H,C)$.

# Chapter No.2

# Soft Neutrosophic Groupoid AND Their Generaliaztion

In this chapter, we introduce soft neutrosophic groupoid over a neutrosophic groupoid with its generalization. We also give sufficient amount of illustrative examples and establish some of their basic and fundamental properties and characteristics.

## 2.1    Soft Neutrosophic Groupoid

The definitions and notions of soft neutrosophic groupoid over a neutrosophic groupoid is introduce here with examples and basic properties.

**Definition 2.1.1.** Let $\{\langle G \cup I \rangle, *\}$ be a neutrosophic groupoid and $(F, A)$ be a soft set over $\{\langle G \cup I \rangle, *\}$. Then $(F, A)$ is called soft neutrosophic groupoid if and only if $F(a)$ is neutrosophic subgroupoid of $\{\langle G \cup I \rangle, *\}$ for all $a \in A$.

This situation can be illustrated in the following examples.

**Example 2.1.1.**  Let $\langle Z_{10} \cup I \rangle = \begin{Bmatrix} 0,\ 1,\ 2,\ 3,\ \dots,\ 9,\ I,\ 2I,\ \dots,\ 9I, \\ 1+\ I,\ 2\ +\ I,\ \dots,\ 9\ +\ 9I \end{Bmatrix}$ be a

neutrosophic groupoid where $*$ is defined on $\langle Z_{10} \cup I \rangle$ by $a*b = 3a + 2b (\mathrm{mod}\,10)$ for all $a,b \in \langle Z10 \cup I \rangle$. Let $A = \{a_1, a_2\}$ be a set of parameters.

Then $(F, A)$ is a soft neutrosophic groupoid over $\{\langle Z10 \cup I \rangle, *\}$, where

$$F(a_1) = \{0, 5, 5I, 5 + 5I\},$$

$$F(a_2) = (Z_{10}, *).$$

**Example 2.1.2.** Let $\langle Z_4 \cup I \rangle = \left\{ \begin{array}{l} 0,1,2,3,I,2I,3I,1+I,1+2I,1+3I \\ 2+I,2+2I,2+3I,3+I,3+2I,3+3I \end{array} \right\}$ be a neutrosophic groupoid with respect to the operation $*$ where $*$ is defined as $a*b = 2a + b (\mathrm{mod}\,4)$ for all $a,b \in \langle Z_4 \cup I \rangle$. Let $A = \{a_1, a_2, a_3\}$ be a set of parameters.

Then $(F, A)$ is a soft neutrosophic groupoid over $\langle Z_4 \cup I \rangle$, where

$$F(a_1) = \{0, 2, 2I, 2 + 2I\},$$

$$F(a_2) = \{0, 2, 2 + 2I\},$$

$$F(a_3) = \{0, 2 + 2I\}.$$

**Theorem 2.1.1.** A soft neutrosophic groupoid over $\{\langle G \cup I \rangle, *\}$ always contain a soft groupoid over $(G, *)$.

**Proof.** The proof is left as an exercise for the readers.

**Theorem 2.1.2.** Let $(F,A)$ and $(H,A)$ be two soft neutrosophic groupoids over $\{\langle G \cup I \rangle, *\}$. Then their intersection $(F,A) \cap (H,A)$ is again a soft neutrosophic groupoid over $\{\langle G \cup I \rangle, *\}$.

**Proof.** The proof is straightforward.

**Theorem 2.1.3.** Let $(F,A)$ and $(H,B)$ be two soft neutrosophic groupoids over $\{\langle G \cup I \rangle, *\}$. If $A \cap B = \phi$, then $(F,A) \cup (H,B)$ is a soft neutrosophic groupoid over $\{\langle G \cup I \rangle, *\}$.

**Remark 2.1.1.** The extended union of two soft neutrosophic groupoids $(F,A)$ and $(K,B)$ over a neutrosophic groupoid $\{\langle G \cup I \rangle, *\}$ is not a soft neutrosophic groupoid over $\{\langle G \cup I \rangle, *\}$.

**Proposition 2.1.1.** The extended intersection of two soft neutrosophic groupoids over a neutrosophic groupoid $\{\langle G \cup I \rangle, *\}$ is a soft neutrosophic groupoid over $\{\langle G \cup I \rangle, *\}$.

**Remark 2.1.2.** The restricted union of two soft neutrosophic groupoids $(F,A)$ and $(K,B)$ over $\{\langle G \cup I \rangle, *\}$ is not a soft neutrosophic groupoid over $\{\langle G \cup I \rangle, *\}$.

**59**

**Proposition 2.1.2.** The restricted intersection of two soft neutrosophic groupoids over $\{\langle G \cup I \rangle, *\}$ is a soft neutrosophic groupoid over $\{\langle G \cup I \rangle, *\}$.

**Proposition 2.1.3.** The *AND* operation of two soft neutrosophic groupoids over $\{\langle G \cup I \rangle, *\}$ is a soft neutrosophic groupoid over $\{\langle G \cup I \rangle, *\}$.

**Remark 2.1.3.** The *OR* operation of two soft neutosophic groupoids over $\{\langle G \cup I \rangle, *\}$ is not a soft nuetrosophic groupoid over $\{\langle G \cup I \rangle, *\}$.

**Definition 21.2.** Let $(F, A)$ be a soft neutrosophic groupoid over $\{\langle G \cup I \rangle, *\}$. Then $(F, A)$ is called an absolute-soft neutrosophic groupoid over $\{\langle G \cup I \rangle, *\}$ if $F(a) = \{\langle G \cup I \rangle, *\}$, for all $a \in A$.

**Theorem 2.1.4.** Every absolute-soft neutrosophic groupoid over $\{\langle G \cup I \rangle, *\}$ always contain absolute soft groupoid over $\{G, *\}$.

**Definition 2.1.3.** Let $(F, A)$ and $(H, B)$ be two soft neutrosophic groupoids over $\{\langle G \cup I \rangle, *\}$. Then $(H, B)$ is a soft neutrosophic subgroupoid of $(F, A)$, if

1. $B \subseteq A$.
2. $H(a)$ is neutrosophic subgroupoid of $F(a)$, for all $a \in B$.

**Example 2.1.3.** Let $\langle Z_4 \cup I \rangle = \left\{ \begin{array}{l} 0,1,2,3,I,2I,3I,1+I,1+2I,1+3I \\ 2+I,2+2I,2+3I,3+I,3+2I,3+3I \end{array} \right\}$ be a neutrosophic groupoid with respect to the operation $*$ where $*$ is defined as $a*b = 2a+b \pmod 4$ for all $a,b \in \langle Z_4 \cup I \rangle$. Let $A = \{a_1, a_2, a_3\}$ be a set of parameters. Then $(F, A)$ is a soft neutrosophic groupoid over $\langle Z_4 \cup I \rangle$, where

$$F(a_1) = \{0, 2, 2I, 2+2I\},$$

$$F(a_2) = \{0, 2, 2+2I\},$$

$$F(a_3) = \{0, 2+2I\}.$$

Let $B = \{a_1, a_2\} \subseteq A$. Then $(H, B)$ is a soft neutrosophic subgroupoid of $(F, A)$, where

$$H(a_1) = \{0, 2+2I\},$$

$$H(a_2) = \{0, 2+2I\}.$$

**Definition 3.1.4.** Let $\langle\langle G \cup I\rangle, *\rangle$ be a neutrosophic groupoid and $(F, A)$ be a soft neutrosophic groupoid over $\langle\langle G \cup I\rangle, *\rangle$. Then $(F, A)$ is called soft Lagrange neutrosophic groupoid if and only if $F(a)$ is a Lagrange neutrosophic subgroupoid of $\langle\langle G \cup I\rangle, *\rangle$ for all $a \in A$.

This situation can be explained in the next example.

**Example 2.1.4.** Let $\langle Z_4 \cup I\rangle = \left\{ \begin{array}{l} 0,1,2,3,I,2I,3I,1+I,1+2I,1+3I \\ 2+I,2+2I,2+3I,3+I,3+2I,3+3I \end{array} \right\}$ be a neutrosophic groupoid of order 16 with respect to the operation $*$ where $*$ is defined as $a*b = 2a + b \pmod 4$ for all $a, b \in \langle Z_4 \cup I\rangle$. Let $A = \{a_1, a_2\}$ be a set of parameters. Then $(F, A)$ is a soft Lagrange neutrosophic groupoid over $\langle Z_4 \cup I\rangle$, where

$$F(a_1) = \{0, 2, 2I, 2+2I\},$$

$$F(a_2) = \{0, 2+2I\}.$$

**Theorem 2.1.5.** Every soft Lagrange neutrosophic groupoid over $\{\langle G \cup I \rangle, *\}$ is a soft neutrosophic groupoid over $\{\langle G \cup I \rangle, *\}$ but the converse is not true.
We can easily show the converse by the help of example.

**Theorem 2.1.6.** If $\{\langle G \cup I \rangle, *\}$ is a Lagrange neutrosophic groupoid, then $(F, A)$ over $\{\langle G \cup I \rangle, *\}$ is a soft Lagrange neutrosophic groupoid but the converse is not true.

**Remark 2.1.4.** Let $(F, A)$ and $(K, C)$ be two soft Lagrange neutrosophic groupoids over $\{\langle G \cup I \rangle, *\}$. Then

1. Their extended intersection $(F, A) \cap_E (K, C)$ may not be a soft Lagrange neutrosophic groupoid over $\{\langle G \cup I \rangle, *\}$.

2. Their restricted intersection $(F, A) \cap_R (K, C)$ may not be a soft Lagrange neutrosophic groupoid over $\{\langle G \cup I \rangle, *\}$.

3. Their *AND* operation $(F, A) \wedge (K, C)$ may not be a soft Lagrange neutrosophic groipoid over $\{\langle G \cup I \rangle, *\}$.

4. Their extended union $(F, A) \cup_E (K, C)$ may not be a soft Lagrange neutrosophic groupoid over $\{\langle G \cup I \rangle, *\}$.

5. Their restricted union $(F, A) \cup_R (K, C)$ may not be a soft Lagrange neutrosophic groupoid over $\{\langle G \cup I \rangle, *\}$.

6. Their *OR* operation $(F, A) \vee (K, C)$ may not be a soft Lagrange

neutrosophic groupoid over $\{\langle G \cup I \rangle, *\}$.

One can easily verify $(1), (2), (3), (4), (5)$ and $(6)$ by the help of examples.

**Definition 2.1.5.** Let $\{\langle G \cup I \rangle, *\}$ be a neutrosophic groipoid and $(F, A)$ be a soft neutrosophic groupoid over $\{\langle G \cup I \rangle, *\}$. Then $(F, A)$ is called soft weak Lagrange neutrosophic groupoid if at least one $F(a)$ is not a Lagrange neutrosophic subgroupoid of $\{\langle G \cup I \rangle, *\}$ for some $a \in A$.

**Example 2.1.5.** Let $\langle Z_4 \cup I \rangle = \left\{ \begin{array}{l} 0,1,2,3,I,2I,3I,1+I,1+2I,1+3I \\ 2+I,2+2I,2+3I,3+I,3+2I,3+3I \end{array} \right\}$ be a neutrosophic groupoid of order 16 with respect to the operation $*$ where $*$ is defined as $a * b = 2a + b \pmod{4}$ for all $a, b \in \langle Z_4 \cup I \rangle$. Let $A = \{a_1, a_2, a_3\}$ be a set of parameters. Then $(F, A)$ is a soft weak Lagrange neutrosophic groupoid over $\langle Z_4 \cup I \rangle$, where

$$F(a_1) = \{0, 2, 2I, 2 + 2I\},$$

$$F(a_2) = \{0, 2, 2 + 2I\},$$

$$F(a_3) = \{0, 2 + 2I\}.$$

**Theorem 2.1.7.** Every soft weak Lagrange neutrosophic groupoid over $\{\langle G \cup I \rangle, *\}$ is a soft neutrosophic groupoid over $\{\langle G \cup I \rangle, *\}$ but the converse is not true.

**63**

**Theorem 2.1.8.** If $\{\langle G \cup I \rangle, *\}$ is weak Lagrange neutrosophic groupoid, then $(F, A)$ over $\{\langle G \cup I \rangle, *\}$ is also soft weak Lagrange neutrosophic groupoid but the converse is not true.

**Remark 2.1.5.** Let $(F,A)$ and $(K,C)$ be two soft weak Lagrange neutrosophic groupoids over $\{\langle G \cup I \rangle, *\}$. Then

1. Their extended intersection $(F,A) \cap_E (K,C)$ is not a soft weak Lagrange neutrosophic groupoid over $\{\langle G \cup I \rangle, *\}$.

2. Their restricted intersection $(F,A) \cap_R (K,C)$ is not a soft weak Lagrange neutrosophic groupoid over $\{\langle G \cup I \rangle, *\}$.

3. Their *AND* operation $(F,A) \wedge (K,C)$ is not a soft weak Lagrange neutrosophic groupoid over $\{\langle G \cup I \rangle, *\}$.

4. Their extended union $(F,A) \cup_E (K,C)$ is not a soft weak Lagrnage neutrosophic groupoid over $\{\langle G \cup I \rangle, *\}$.

5. Their restricted union $(F,A) \cup_R (K,C)$ is not a soft weak Lagrange neutrosophic groupoid over $\{\langle G \cup I \rangle, *\}$.

6. Their *OR* operation $(F,A) \vee (K,C)$ is not a soft weak Lagrange neutrosophic groupoid over $\{\langle G \cup I \rangle, *\}$.

One can easily verify $(1),(2),(3),(4),(5)$ and $(6)$ by the help of examples.

**Definition 2.1.6.** Let $\{\langle G \cup I \rangle, *\}$ be a neutrosophic groupoid and $(F,A)$ be a soft neutrosophic groupoid over $\{\langle G \cup I \rangle, *\}$. Then $(F,A)$ is called soft Lagrange free neutrosophic groupoid if $F(a)$ is not a lagrange neutrosophic subgroupoid of $\{\langle G \cup I \rangle, *\}$ for all $a \in A$.

**Example 2.1.6.** Let $\langle Z_4 \cup I \rangle = \left\{ \begin{array}{l} 0,1,2,3,I,2I,3I,1+I,1+2I,1+3I \\ 2+I,2+2I,2+3I,3+I,3+2I,3+3I \end{array} \right\}$ be a neutrosophic groupoid of order 16 with respect to the operation $*$ where

$*$ is defined as $a*b = 2a+b \pmod 4$ for all $a,b \in \langle Z_4 \cup I \rangle$. Let $A = \{a_1, a_2, a_3\}$ be a set of parameters. Then $(F,A)$ is a soft Lagrange free neutrosophic groupoid over $\langle Z_4 \cup I \rangle$, where

$$F(a_1) = \{0, 2I, 2+2I\},$$

$$F(a_2) = \{0, 2, 2+2I\}.$$

**Theorem 2.1.9.** Every soft Lagrange free neutrosophic groupoid over $\{\langle G \cup I \rangle, *\}$ is trivially a soft neutrosophic groupoid over $\{\langle G \cup I \rangle, *\}$ but the converse is not true.

**Theorem 2.1.10.** If $\{\langle G \cup I \rangle, *\}$ is a Lagrange free neutrosophic groupoid, then $(F,A)$ over $\{\langle G \cup I \rangle, *\}$ is also a soft Lagrange free neutrosophic groupoid but the converse is not true.

**Remark 2.1.6.** Let $(F,A)$ and $(K,C)$ be two soft Lagrange free neutrosophic groupoids over $\{\langle G \cup I \rangle, *\}$. Then

1. Their extended intersection $(F,A) \cap_E (K,C)$ is not a soft Lagrange free neutrosophic groupoid over $\{\langle G \cup I \rangle, *\}$.
2. Their restricted intersection $(F,A) \cap_R (K,C)$ is not a soft Lagrange free neutrosophic groupoid over $\{\langle G \cup I \rangle, *\}$.
3. Their *AND* operation $(F,A) \wedge (K,C)$ is not a soft Lagrange free

**65**

neutrosophic groupoid over $\{\langle G \cup I \rangle, *\}$.

4. Their extended union $(F,A) \cup_E (K,C)$ is not a soft Lagrnage free

neutrosophic groupoid over $\{\langle G \cup I \rangle, *\}$.

5. Their restricted union $(F,A) \cup_R (K,C)$ is not a soft Lagrange free neutrosophic groupoid over $\{\langle G \cup I \rangle, *\}$.

6. Their *OR* operation $(F,A) \vee (K,C)$ is not a soft Lagrange free neutrosophic groupoid over $\{\langle G \cup I \rangle, *\}$.

One can easily verify (1),(2),(3),(4),(5) and (6) by the help of examples.

**Definition 2.1.7.** $(F,A)$ is called soft neutrosophic ideal over $\{\langle G \cup I \rangle, *\}$ if $F(a)$ is a neutrosophic ideal of $\{\langle G \cup I \rangle, *\}$, for all $a \in A$.

**Theorem 2.1.11.** Every soft neutrosophic ideal $(F,A)$ over $\{\langle G \cup I \rangle, *\}$ is trivially a soft neutrosophic subgroupid but the converse may not be true.

**Proposition 2.1.4.** Let $(F,A)$ and $(K,B)$ be two soft neutrosophic ideals over $\{\langle G \cup I \rangle, *\}$. Then

1. Their extended intersection $(F,A) \cap_E (K,B)$ is soft neutrosophic ideal over $\{\langle G \cup I \rangle, *\}$.

2. Their restricted intersection $(F,A) \cap_R (K,B)$ is soft neutrosophic ideal over $\{\langle G \cup I \rangle, *\}$.

3. Their *AND* operation $(F,A) \wedge (K,B)$ is soft neutrosophic ideal over $\{\langle G \cup I \rangle, *\}$.

**66**

**Remark 2.1.7.** Let $(F,A)$ and $(K,B)$ be two soft neutrosophic ideal over $\{\langle G \cup I \rangle, *\}$. Then

1. Their extended union $(F,A) \cup_E (K,B)$ is not soft neutrosophic ideal over $\{\langle G \cup I \rangle, *\}$.

2. Their restricted union $(F,A) \cup_R (K,B)$ is not soft neutrosophic ideal over $\{\langle G \cup I \rangle, *\}$.

3. Their *OR* operation $(F,A) \vee (K,B)$ is not soft neutrosophic ideal over $\{\langle G \cup I \rangle, *\}$.

One can easily proved (1),(2), and (3) by the help of examples.

**Theorem 2.1.12.** Let $(F,A)$ be a soft neutrosophic ideal over $\{\langle G \cup I \rangle, *\}$ and $\{(H_i, B_i) : i \in J\}$ is a non-empty family of soft neutrosophic ideals of $(F,A)$. Then

1. $\underset{i \in J}{\cap}(H_i, B_i)$ is a soft neutrosophic ideal of $(F,A)$.

2. $\underset{i \in J}{\wedge}(H_i, B_i)$ is a soft neutrosophic ideal of $\underset{i \in J}{\wedge}(F,A)$.

## Soft Neutrosophic Strong Groupoid

The soft neutrosophic strong groupoid over a neutrosophic groupoid is introduced here which is of pure neutrosophic character.

**Definition 2.1.8.** Let $\{\langle G \cup I \rangle, *\}$ be a neutrosophic groupoid and $(F,A)$ be a soft set over $\{\langle G \cup I \rangle, *\}$. Then $(F,A)$ is called soft neutrosophic strong groupoid if and only if $F(a)$ is a neutrosophic strong subgroupoid of $\{\langle G \cup I \rangle, *\}$ for all $a \in A$.

**Example 2.1.7.** Let $\langle Z_4 \cup I \rangle = \begin{Bmatrix} 0,1,2,3,I,2I,3I,1+I,1+2I,1+3I \\ 2+I,2+2I,2+3I,3+I,3+2I,3+3I \end{Bmatrix}$ be a neutrosophic groupoid with respect to the operation $*$ where $*$ is defined as $a * b = 2a + b (\bmod 4)$ for all $a, b \in \langle Z_4 \cup I \rangle$. Let $A = \{a_1, a_2, a_3\}$ be a set of parameters. Then $(F,A)$ is a soft neutrosophic strong groupoid over $\langle Z_4 \cup I \rangle$, where

$$F(a_1) = \{0, 2I, 2+2I\},$$

$$F(a_2) = \{0, 2+2\mathrm{I}\}.$$

**Proposition 2.1.5.** Let $(F,A)$ and $(K,\mathrm{C})$ be two soft neutrosophic strong groupoids over $\{\langle G \cup I \rangle, *\}$. Then

1. Their extended intersection $(F,A) \cap_E (K,\mathrm{C})$ is a soft neutrosophic strong groupoid over $\{\langle G \cup I \rangle, *\}$.
2. Their restricted intersection $(F,A) \cap_R (K,\mathrm{C})$ is a soft neutrosophic strong groupoid over $\{\langle G \cup I \rangle, *\}$.
3. Their $AND$ operation $(F,A) \wedge (K,\mathrm{C})$ is a soft neutrosophic strong groupoid over $\{\langle G \cup I \rangle, *\}$.

**Remark 2.1.8.** Let $(F,A)$ and $(K,\mathrm{C})$ be two soft neutrosophic strong

groupoids over $\{\langle G \cup I \rangle, *\}$. Then

1. Their extended union $(F,A) \cup_E (K,C)$ is a soft neutrosophic strong groupoid over $\{\langle G \cup I \rangle, *\}$.
2. Their restricted union $(F,A) \cup_R (K,C)$ is a soft neutrosophic strong groupoid over $\{\langle G \cup I \rangle, *\}$.
3. Their *OR* operation $(F,A) \vee (K,C)$ is a soft neutrosophic strong groupoid over $\{\langle G \cup I \rangle, *\}$.

**Definition 2.1.9.** Let $(F,A)$ and $(H,C)$ be two soft neutrosophic strong groupoids over $\{\langle G \cup I \rangle, *\}$. Then $(H,C)$ is called soft neutrosophic strong sublgroupoid of $(F,A)$, if

1. $C \subseteq A$.
2. $H(a)$ is a neutrosophic strong subgroupoid of $F(a)$ for all $a \in A$.

**Definition 2.1.10.** Let $\{\langle G \cup I \rangle, *\}$ be a neutrosophic strong groupoid and $(F,A)$ be a soft neutrosophic groupoid over $\{\langle G \cup I \rangle, *\}$. Then $(F,A)$ is called soft Lagrange neutrosophic strong groupoid if and only if $F(a)$ is a Lagrange neutrosophic strong subgroupoid of $\{\langle G \cup I \rangle, *\}$ for all $a \in A$.

**Theorem 2.1.13.** Every soft Lagrange neutrosophic strong groupoid over $\{\langle G \cup I \rangle, *\}$ is a soft neutrosophic groupoid over $\{\langle G \cup I \rangle, *\}$ but the converse is not true.

**Theorem 2.1.14.** If $\{\langle G \cup I \rangle, *\}$ is a Lagrange neutrosophic strong

groupoid, then $(F,A)$ over $\{\langle G \cup I \rangle, *\}$ is a soft Lagrange neutrosophic groupoid but the converse is not true.

**Remark 2.1.9.** Let $(F,A)$ and $(K,C)$ be two soft Lagrange neutrosophic strong groupoids over $\{\langle G \cup I \rangle, *\}$. Then

1. Their extended intersection $(F,A) \cap_E (K,C)$ may not be a soft Lagrange neutrosophic strong groupoid over $\{\langle G \cup I \rangle, *\}$.
2. Their restricted intersection $(F,A) \cap_R (K,C)$ may not be a soft Lagrange strong neutrosophic groupoid over $\{\langle G \cup I \rangle, *\}$.
3. Their *AND* operation $(F,A) \wedge (K,C)$ may not be a soft Lagrange neutrosophic strong groupoid over $\{\langle G \cup I \rangle, *\}$.
4. Their extended union $(F,A) \cup_E (K,C)$ may not be a soft Lagrange neutrosophic strong groupoid over $\{\langle G \cup I \rangle, *\}$.
5. Their restricted union $(F,A) \cup_R (K,C)$ may not be a soft Lagrange neutrosophic strong groupoid over $\{\langle G \cup I \rangle, *\}$.
6. Their *OR* operation $(F,A) \vee (K,C)$ may not be a soft Lagrange neutrosophic strong groupoid over $\{\langle G \cup I \rangle, *\}$.

One can easily verify $(1),(2),(3),(4),(5)$ and $(6)$ by the help of examples.

**Definition 2.1.11.** Let $\{\langle G \cup I \rangle, *\}$ be a neutrosophic strong groupoid and $(F,A)$ be a soft neutrosophic groupoid over $\{\langle G \cup I \rangle, *\}$. Then $(F,A)$ is called soft weak Lagrange neutrosophic strong groupoid if at least one $F(a)$ is not a Lagrange neutrosophic strong subgroupoid of $\{\langle G \cup I \rangle, *\}$ for some $a \in A$.

**Theorem 2.1.15.** Every soft weak Lagrange neutrosophic strong groupoid over $\{\langle G \cup I \rangle, *\}$ is a soft neutrosophic groupoid over $\{\langle G \cup I \rangle, *\}$ but the converse is not true.

**Theorem 2.1.16.** If $\{\langle G \cup I \rangle, *\}$ is weak Lagrange neutrosophic strong groupoid, then $(F, A)$ over $\{\langle G \cup I \rangle, *\}$ is also soft weak Lagrange neutrosophic strong groupoid but the converse is not true.

**Remark 2.1.10.** Let $(F, A)$ and $(K, C)$ be two soft weak Lagrange neutrosophic strong groupoids over $\{\langle G \cup I \rangle, *\}$. Then

1. Their extended intersection $(F, A) \cap_E (K, C)$ is not a soft weak Lagrange neutrosophic strong groupoid over $\{\langle G \cup I \rangle, *\}$.
2. Their restricted intersection $(F, A) \cap_R (K, C)$ is not a soft weak Lagrange neutrosophic strong groupoid over $\{\langle G \cup I \rangle, *\}$.
3. Their $AND$ operation $(F, A) \wedge (K, C)$ is not a soft weak Lagrange neutrosophic strong groupoid over $\{\langle G \cup I \rangle, *\}$.
4. Their extended union $(F, A) \cup_E (K, C)$ is not a soft weak Lagrnage neutrosophic strong groupoid over $\{\langle G \cup I \rangle, *\}$.
5. Their restricted union $(F, A) \cup_R (K, C)$ is not a soft weak Lagrange neutrosophic strong groupoid over $\{\langle G \cup I \rangle, *\}$.
6. Their $OR$ operation $(F, A) \vee (K, C)$ is not a soft weak Lagrange neutrosophic strong groupoid over $\{\langle G \cup I \rangle, *\}$.

One can easily verify (1),(2),(3),(4),(5) and (6) by the help of examples.

**Definition 2.1.12.** Let $\langle G \cup I \rangle$ be a neutrosophic strong groupoid and

$(F, A)$ be a soft neutrosophic groupoid over $\langle G \cup I \rangle$. Then $(F, A)$ is called soft Lagrange free neutrosophic strong groupoid if $F(a)$ is not a Lagrange neutrosophic strong subgroupoid of $\langle G \cup I \rangle$ for all $a \in A$.

**Theorem 2.1.17.** Every soft Lagrange free neutrosophic strong groupoid over $\langle L \cup I \rangle$ is a soft neutrosophic groupoid over $\{\langle G \cup I \rangle, *\}$ but the converse is not true.

**Theorem 2.1.18.** If $\{\langle G \cup I \rangle, *\}$ is a Lagrange free neutrosophic strong groupoid, then $(F, A)$ over $\{\langle G \cup I \rangle, *\}$ is also a soft Lagrange free neutrosophic strong groupoid but the converse is not true.

**Remark 2.1.11.** Let $(F, A)$ and $(K, C)$ be two soft Lagrange free neutrosophic strong groupoids over $\langle L \cup I \rangle$. Then

1. Their extended intersection $(F, A) \cap_E (K, C)$ is not a soft Lagrange free neutrosophic strong groupoid over $\{\langle G \cup I \rangle, *\}$.

2. Their restricted intersection $(F, A) \cap_R (K, C)$ is not a soft Lagrange free neutrosophic strong groupoid over $\{\langle G \cup I \rangle, *\}$.

3. Their *AND* operation $(F, A) \wedge (K, C)$ is not a soft Lagrange free neutrosophic strong groupoid over $\{\langle G \cup I \rangle, *\}$.

4. Their extended union $(F, A) \cup_E (K, C)$ is not a soft Lagrange free neutrosophic strong groupoid over $\{\langle G \cup I \rangle, *\}$.

5. Their restricted union $(F, A) \cup_R (K, C)$ is not a soft Lagrange free neutrosophic groupoid over $\{\langle G \cup I \rangle, *\}$.

**72**

6. Their *OR* operation $(F, A) \vee (K, C)$ is not a soft Lagrange free

neutrosophic strong groupoid over $\{\langle G \cup I \rangle, *\}$.

One can easily verify (1),(2),(3),(4),(5) and (6) by the help of examples.

**Definition 2.1.13.** $(F, A)$ is called soft neutrosophic strong ideal over $\{\langle G \cup I \rangle, *\}$ if $F(a)$ is a neutrosophic strong ideal of $\{\langle G \cup I \rangle, *\}$, for all $a \in A$.

**Theorem 2.1.19.** Every soft neutrosophic strong ideal $(F, A)$ over $\{\langle G \cup I \rangle, *\}$ is trivially a soft neutrosophic strong groupoid.

**Theorem 2.1.20.** Every soft neutrosophic strong ideal $(F, A)$ over $\{\langle G \cup I \rangle, *\}$ is trivially a soft neutrosophic ideal.

**Proposition 2.1.6.** Let $(F, A)$ and $(K, B)$ be two soft neutrosophic strong ideals over $\{\langle G \cup I \rangle, *\}$. Then

1. Their extended intersection $(F, A) \cap_E (K, B)$ is soft neutrosophic strong ideal over $\{\langle G \cup I \rangle, *\}$.
2. Their restricted intersection $(F, A) \cap_R (K, B)$ is soft neutrosophic strong ideal over $\{\langle G \cup I \rangle, *\}$.
3. Their *AND* operation $(F, A) \wedge (K, B)$ is soft neutrosophic strong ideal over $\{\langle G \cup I \rangle, *\}$.

**Remark 2.1.12.** Let $(F, A)$ and $(K, B)$ be two soft neutrosophic strong

ideal over $\{\langle G \cup I \rangle, *\}$. Then

1. Their extended union $(F, A) \cup_E (K, B)$ is not soft neutrosophic strong ideal over $\{\langle G \cup I \rangle, *\}$.

2. Their restricted union $(F, A) \cup_R (K, B)$ is not soft neutrosophic strong ideal over $\{\langle G \cup I \rangle, *\}$.

3. Their *OR* operation $(F, A) \vee (K, B)$ is not soft neutrosophic strong ideal over $\{\langle G \cup I \rangle, *\}$.

One can easily proved $(1), (2),$ and $(3)$ by the help of examples.

**Theorem 2.1.21.** Let $(F, A)$ be a soft neutrosophic strong ideal over $\{\langle G \cup I \rangle, *\}$ and $\{(H_i, B_i) : i \in J\}$ is a non-empty family of soft neutrosophic strong ideals of $(F, A)$. Then

1. $\underset{i \in J}{\cap}(H_i, B_i)$ is a soft neutrosophic strong ideal of $(F, A)$.

2. $\underset{i \in J}{\wedge}(H_i, B_i)$ is a soft neutrosophic strong ideal of $\underset{i \in J}{\wedge}(F, A)$.

In the proceeding section, we define soft neutrosophic bigroupoid and soft neutrosophic strong bigroupoid over a neutrosophic bigroupoid.

## 2.2 Soft Neutrosophic Bigroupoid

In this section soft neutrosophic bigroupoid over a neutrosophic bigroupoid is presented here. We also introduced soft neutrosophic strong bigroupoid over a neutrosophic bigroupoid in the current section.

Some basic properties are also given with sufficient amount of illustrative examples.

We now proceed on to define soft neutrosophic bigroupoid over a neutrosophic bigroupoid.

**Definition 2.2.1.** Let $\{B_N(G), *, \circ\}$ be a neutrosophic bigroupoid and $(F, A)$ be a soft set over $\{B_N(G), *, \circ\}$. Then $(F, A)$ is called soft neutrosophic bigroupoid if and only if $F(a)$ is neutrosophic sub bigroupoid of $\{B_N(G), *, \circ\}$ for all $a \in A$.

The following examples illustrates this fact.

**Example 2.2.1.** Let $\{B_N(G), *, \circ\}$ be a neutrosophic groupoid with $B_N(G) = G_1 \cup G_2$, where

$G_1 = \left\{ \langle Z_{10} \cup I \rangle \mid a*b = 2a + 3b(\mathrm{mod}\, 10); a, b \in \langle Z_{10} \cup I \rangle \right\}$ and

$G_2 = \left\{ \langle Z_4 \cup I \rangle \mid a \circ b = 2a + b(\mathrm{mod}\, 4); a, b \in \langle Z_4 \cup I \rangle \right\}$.

Let $A = \{a_1, a_2\}$ be a set of parameters. Then $(F, A)$ is a soft neutrosophic bigroupoid over $\{B_N(G), *, \circ\}$, where

$$F(a_1) = \{0, 5, 5I, 5+5I\} \cup \{0, 2, 2I, 2+2I\},$$

$$F(a_2) = (Z_{10}, *) \cup \{0, 2+2I\}.$$

**Theorem 2.2.1.** Let $(F, A)$ and $(H, A)$ be two soft neutrosophic bigroupoids over $\{B_N(G),*,\circ\}$. Then their intersection $(F,A) \cap (H,A)$ is again a soft neutrosophic groupoid over $\{B_N(G),*,\circ\}$.

**Proof.** The proof is straightforward.

**Theorem 2.2.2.** Let $(F,A)$ and $(H,B)$ be two soft neutrosophic groupoids over $\{\langle G \cup I \rangle,*\}$. If $A \cap B = \phi$, then $(F,A) \cup (H,B)$ is a soft neutrosophic groupoid over $\{\langle G \cup I \rangle,*\}$.

**Proposition 2.2.1.** Let $(F,A)$ and $(K,C)$ be two soft neutrosophic bigroupoids over $\{B_N(G),*,\circ\}$. Then

1. Their extended intersection $(F,A) \cap_E (K,C)$ is a soft neutrosophic bigroupoid over $\{B_N(G),*,\circ\}$.
2. Their restricted intersection $(F,A) \cap_R (K,C)$ is a soft neutrosophic bigroupoid over $\{B_N(G),*,\circ\}$.
3. Their *AND* operation $(F,A) \wedge (K,C)$ is a soft neutrosophic bigroupoid over $\{B_N(G),*,\circ\}$.

**Remark 2.2.1.** Let $(F,A)$ and $(K,C)$ be two soft neutrosophic biloops over $\{B_N(G),*,\circ\}$. Then

1. Their extended union $(F,A) \cup_E (K,C)$ is not a soft neutrosophic bigroupoid over $\{B_N(G),*,\circ\}$.
2. Their restricted union $(F,A) \cup_R (K,C)$ is not a soft neutrosophic bigroupoid over $\{B_N(G),*,\circ\}$.

**76**

3. Their *OR* operation $(F,A) \vee (K,C)$ is not a soft neutrosophic bigroupoid over $\{B_N(G), *, \circ\}$.

One can easily verify $(1), (2),$ and $(3)$ by the help of examples.

**Definition 2.2.2.** Let $(F,A)$ be a soft neutrosophic bigroupoid over $\{B_N(G), *, \circ\}$. Then $(F,A)$ is called an absolute soft neutrosophic bigroupoid over $\{B_N(G), *, \circ\}$ if $F(a) = \{B_N(G), *, \circ\}$ for all $a \in A$.

**Definition 2.2.3.** Let $(F,A)$ and $(H,C)$ be two soft neutrosophic bigroupoids over $\{B_N(G), *, \circ\}$. Then $(H,C)$ is called soft neutrosophic sub bigroupoid of $(F,A)$, if

1. $C \subseteq A$.
2. $H(a)$ is a neutrosophic sub bigroupoid of $F(a)$ for all $a \in A$.

**Example 2.2.2.** Let $\{B_N(G), *, \circ\}$ be a neutrosophic groupoid with $B_N(G) = G_1 \cup G_2$, where

$G_1 = \left\{ \left\langle Z_{10} \cup I \right\rangle \mid a * b = 2a + 3b(\text{mod}10); a, b \in \left\langle Z_{10} \cup I \right\rangle \right\}$ and

$G_2 = \left\{ \left\langle Z_4 \cup I \right\rangle \mid a \circ b = 2a + b(\text{mod}4); a, b \in \left\langle Z_4 \cup I \right\rangle \right\}$.

Let $A = \{a_1, a_2\}$ be a set of parameters. Let $(F,A)$ is a soft neutrosophic bigroupoid over $\{B_N(G), *, \circ\}$, where

$$F(a_1) = \{0, 5, 5I, 5 + 5I\} \cup \{0, 2, 2I, 2 + 2I\},$$

$$F(a_2) = (Z_{10}, *) \cup \{0, 2 + 2I\}.$$

Let $B = \{a_1\} \subseteq A$. Then $(H, B)$ is a soft neutrosophic sub bigroupoid of $(F, A)$, where

$$H(a_1) = \{0, 5\} \cup \{0, 2 + 2I\}.$$

**Definition 2.2.4.** Let $\{B_N(G), *, \circ\}$ be a neutrosophic strong bigroupoid and $(F, A)$ be a soft neutrosophic bigroupoid over $\{B_N(G), *, \circ\}$. Then $(F, A)$ is called soft Lagrange neutrosophic bigroupoid if and only if $F(a)$ is a Lagrange neutrosophic sub bigroupoid of $\{B_N(G), *, \circ\}$ for all $a \in A$.

**Theorem 2.2.3.** Every soft Lagrange neutrosophic bigroupoid over $\{B_N(G), *, \circ\}$ is a soft neutrosophic bigroupoid over $\{B_N(G), *, \circ\}$ but the converse is not true.

One can easily see the converse by the help of examples.

**Theorem 2.2.4.** If $\{B_N(G), *, \circ\}$ is a Lagrange neutrosophic bigroupoid, then $(F, A)$ over $\{B_N(G), *, \circ\}$ is a soft Lagrange neutrosophic bigroupoid but the converse is not true.

 **Remark 2.2.2.** Let $(F, A)$ and $(K, C)$ be two soft Lagrange neutrosophic bigroupoids over $\{B_N(G), *, \circ\}$. Then

1. Their extended intersection $(F, A) \cap_E (K, C)$ may not be a soft Lagrange neutrosophic bigroupoid over $\{B_N(G), *, \circ\}$.
2. Their restricted intersection $(F, A) \cap_R (K, C)$ may not be a soft Lagrange neutrosophic bigroupoid over $\{B_N(G), *, \circ\}$.

3. Their $AND$ operation $(F,A) \wedge (K,C)$ may not be a soft Lagrange neutrosophic bigroupoid over $\{B_N(G), *, \circ\}$.

4. Their extended union $(F,A) \cup_E (K,C)$ may not be a soft Lagrange neutrosophic bigroupoid over $\{B_N(G), *, \circ\}$.

5. Their restricted union $(F,A) \cup_R (K,C)$ may not be a soft Lagrange neutrosophic bigroupoid over $\{B_N(G), *, \circ\}$.

6. Their $OR$ operation $(F,A) \vee (K,C)$ may not be a soft Lagrange neutrosophic bigroupoid over $\{B_N(G), *, \circ\}$.

One can easily verify $(1), (2), (3), (4), (5)$ and $(6)$ by the help of examples.

**Definition 2.2.5.** Let $\{B_N(G), *, \circ\}$ be a neutrosophic bigroupoid and $(F,A)$ be a soft neutrosophic bigroupoid over $\{B_N(G), *, \circ\}$. Then $(F,A)$ is called soft weak Lagrange neutrosophic bigroupoid if at least one $F(a)$ is not a Lagrange neutrosophic sub bigroupoid of $\{B_N(G), *, \circ\}$ for some $a \in A$.

**Theorem 2.2.5.** Every soft weak Lagrange neutrosophic bigroupoid over $\{B_N(G), *, \circ\}$ is a soft neutrosophic groupoid over $\{B_N(G), *, \circ\}$ but the converse is not true.

**Theorem 2.2.6.** If $\{B_N(G), *, \circ\}$ is weak Lagrange neutrosophic bigroupoid, then $(F,A)$ over $\{B_N(G), *, \circ\}$ is also soft weak Lagrange neutrosophic bigroupoid but the converse is not true.

**Remark 2.2.3.** Let $(F,A)$ and $(K,C)$ be two soft weak Lagrange neutrosophic bigroupoids over $\{B_N(G), *, \circ\}$. Then

1. Their extended intersection $(F,A) \cap_E (K,C)$ is not a soft weak

**79**

Lagrange neutrosophic bigroupoid over $\{B_N(G), *, \circ\}$.

2. Their restricted intersection $(F, A) \cap_R (K, C)$ is not a soft weak Lagrange neutrosophic bigroupoid over $\{B_N(G), *, \circ\}$.

3. Their *AND* operation $(F, A) \wedge (K, C)$ is not a soft weak Lagrange neutrosophic bigroupoid over $\{B_N(G), *, \circ\}$.

4. Their extended union $(F, A) \cup_E (K, C)$ is not a soft weak Lagrnage neutrosophic bigroupoid over $\{B_N(G), *, \circ\}$.

5. Their restricted union $(F, A) \cup_R (K, C)$ is not a soft weak Lagrange neutrosophic bigroupoid over $\{B_N(G), *, \circ\}$.

6. Their *OR* operation $(F, A) \vee (K, C)$ is not a soft weak Lagrange neutrosophic bigroupoid over $\{B_N(G), *, \circ\}$.

One can easily verify (1),(2),(3),(4),(5) and (6) by the help of examples.

**Definition 2.2.6.** Let $\{B_N(G), *, \circ\}$ be a neutrosophic bigroupoid and $(F, A)$ be a soft neutrosophic groupoid over $\{B_N(G), *, \circ\}$. Then $(F, A)$ is called soft Lagrange free neutrosophic bigroupoid if $F(a)$ is not a Lagrange neutrosophic sub bigroupoid of $\{B_N(G), *, \circ\}$ for all $a \in A$.

**Theorem 2.2.7.** Every soft Lagrange free neutrosophic bigroupoid over $\{B_N(G), *, \circ\}$ is a soft neutrosophic bigroupoid over $\{B_N(G), *, \circ\}$ but the converse is not true.

**Theorem 2.2.8.** If $\{B_N(G), *, \circ\}$ is a Lagrange free neutrosophic bigroupoid, then $(F, A)$ over $\{B_N(G), *, \circ\}$ is also a soft Lagrange free neutrosophic bigroupoid but the converse is not true.

**Remark 2.2.4.** Let $(F,A)$ and $(K,C)$ be two soft Lagrange free neutrosophic bigroupoids over $\{B_N(G), *, \circ\}$. Then

1. Their extended intersection $(F,A) \cap_E (K,C)$ is not a soft Lagrange free neutrosophic bigroupoid over $\{B_N(G), *, \circ\}$.
2. Their restricted intersection $(F,A) \cap_R (K,C)$ is not a soft Lagrange free neutrosophic bigroupoid over $\{B_N(G), *, \circ\}$.
3. Their *AND* operation $(F,A) \wedge (K,C)$ is not a soft Lagrange free neutrosophic bigroupoid over $\{B_N(G), *, \circ\}$.
4. Their extended union $(F,A) \cup_E (K,C)$ is not a soft Lagrange free neutrosophic bigroupoid over $\{B_N(G), *, \circ\}$.
5. Their restricted union $(F,A) \cup_R (K,C)$ is not a soft Lagrange free neutrosophic bigroupoid over $\{B_N(G), *, \circ\}$.
6. Their *OR* operation $(F,A) \vee (K,C)$ is not a soft Lagrange free neutrosophic bigroupoid over $\{B_N(G), *, \circ\}$.

One can easily verify $(1), (2), (3), (4), (5)$ and $(6)$ by the help of examples.

**Definition 2.2.7.** $(F,A)$ is called soft neutrosophic biideal over $\{B_N(G), *, \circ\}$ if $F(a)$ is a neutrosophic biideal of $\{B_N(G), *, \circ\}$, for all $a \in A$.

We now give some characterization of soft neutrosophic biideals.

**Theorem 2.2.9.** Every soft neutrosophic biideal $(F,A)$ over $\{B_N(G), *, \circ\}$ is a soft neutrosophic bigroupoid.

**Proposition 2.2.2.** Let $(F,A)$ and $(K,B)$ be two soft neutrosophic biideals over $\{B_N(G),*,\circ\}$. Then

1. Their extended intersection $(F,A) \cap_E (K,B)$ is soft neutrosophic biideal over $\{B_N(G),*,\circ\}$.
2. Their restricted intersection $(F,A) \cap_R (K,B)$ is soft neutrosophic biideal over $\{B_N(G),*,\circ\}$.
3. Their *AND* operation $(F,A) \wedge (K,B)$ is soft neutrosophic biideal over $\{B_N(G),*,\circ\}$.

**Remark 2.2.5.** Let $(F,A)$ and $(K,B)$ be two soft neutrosophic biideals over $\{B_N(G),*,\circ\}$. Then

1. Their extended union $(F,A) \cup_E (K,B)$ is not soft neutrosophic biideals over $\{B_N(G),*,\circ\}$.
2. Their restricted union $(F,A) \cup_R (K,B)$ is not soft neutrosophic biidleals over $\{B_N(G),*,\circ\}$.
3. Their *OR* operation $(F,A) \vee (K,B)$ is not soft neutrosophic biideals over $\{B_N(G),*,\circ\}$.

One can easily proved $(1),(2),$ and $(3)$ by the help of examples

**Theorem 2.2.10.** Let $(F,A)$ be a soft neutrosophic biideal over $\{B_N(G),*,\circ\}$ and $\{(H_i,B_i):i \in J\}$ is a non-empty family of soft neutrosophic biideals of $(F,A)$. Then

1. $\underset{i \in J}{\cap}(H_i,B_i)$ is a soft neutrosophic biideal of $(F,A)$.

2. $\underset{i \in J}{\wedge}(H_i, B_i)$ is a soft neutrosophic biideal of $\underset{i \in J}{\wedge}(F, A)$.

We now define soft neutrosophic strong bigroupoid over a neutrosophic bigroupoid.

## Soft Neutrosophic Strong Bigroupoid

Here the definition and other important notions of soft neutrosophic strong bigroupoid over a neutrosophic bigroupoid is introduced. Soft neutrosophic strong bigroupoid over a neutrosophic groupoid is a parameterized family of neutrosophic strong subbigroupoid which is of purely neutrosophic character.

**Definition 2.2.8.** Let $\{B_N(G), *, \circ\}$ be a neutrosophic bigroupoid and $(F, A)$ be a soft set over $\{B_N(G), *, \circ\}$. Then $(F, A)$ is called soft neutrosophic strong bigroupoid if and only if $F(a)$ is neutrosophic strong subbigroupoid of $\{B_N(G), *, \circ\}$ for all $a \in A$.

**Example 2.2.3.** Let $\{B_N(G), *, \circ\}$ be a neutrosophic groupoid with $B_N(G) = G_1 \cup G_2$, where $G_1 = \left\{ \langle Z_{10} \cup I \rangle \mid a * b = 2a + 3b \pmod{10}; a, b \in \langle Z_{10} \cup I \rangle \right\}$ and $G_2 = \left\{ \langle Z_4 \cup I \rangle \mid a \circ b = 2a + b \pmod{4}; a, b \in \langle Z_4 \cup I \rangle \right\}$.

Let $A = \{a_1, a_2\}$ be a set of parameters. Then $(F, A)$ is a soft neutrosophic strong bigroupoid over $\{B_N(G), *, \circ\}$, where

$$F(a_1) = \{0, 5 + 5I\} \cup \{0, 2 + 2I\},$$

$$F(a_2) = \{0, 5I\} \cup \{0, 2 + 2I\} \ .$$

**Theorem 2.2.11.** Let $(F,A)$ and $(H,A)$ be two soft neutrosophic strong bigroupoids over $\{B_N(G), *, \circ\}$. Then their intersection $(F,A) \cap (H,A)$ is again a soft neutrosophic strong bigroupoid over $\{B_N(G), *, \circ\}$.

The proof is left as an exercise for the interested readers.

**Theorem 2.2.12.** Let $(F,A)$ and $(H,B)$ be two soft neutrosophic strong bigroupoids over $\{B_N(G), *, \circ\}$. If $A \cap B = \phi$, then $(F,A) \cup (H,B)$ is a soft neutrosophic strong bigroupoid over $\{B_N(G), *, \circ\}$.

**Proposition 2.2.3.** Let $(F,A)$ and $(K,C)$ be two soft neutrosophic strong bigroupoids over $\{B_N(G), *, \circ\}$. Then

1.  Their extended intersection $(F,A) \cap_E (K,C)$ is a soft neutrosophic strong bigroupoid over $\{B_N(G), *, \circ\}$.
2.  Their restricted intersection $(F,A) \cap_R (K,C)$ is a soft neutrosophic strong bigroupoid over $\{B_N(G), *, \circ\}$.
3.  Their *AND* operation $(F,A) \wedge (K,C)$ is a soft neutrosophic strong bigroupoid over $\{B_N(G), *, \circ\}$.

**Remark 2.2.6.** Let $(F,A)$ and $(K,C)$ be two soft neutrosophic strong bigroupoids over $\{B_N(G), *, \circ\}$. Then

1.  Their extended union $(F,A) \cup_E (K,C)$ is not a soft neutrosophic

strong bigroupoid over $\{B_N(G),*,\circ\}$.

2. Their restricted union $(F,A) \cup_R (K,C)$ is not a soft neutrosophic strong bigroupoid over $\{B_N(G),*,\circ\}$.

3. Their *OR* operation $(F,A) \vee (K,C)$ is not a soft neutrosophic strong bigroupoid over $\{B_N(G),*,\circ\}$.

One can easily verify $(1), (2),$ and $(3)$ by the help of examples.

**Definition 2.2.9.** Let $(F,A)$ and $(H,C)$ be two soft neutrosophic strong bigroupoids over $\{B_N(G),*,\circ\}$. Then $(H,C)$ is called soft neutrosophic strong sub bigroupoid of $(F,A)$, if

1. $C \subseteq A$.
2. $H(a)$ is a neutrosophic strong sub bigroupoid of $F(a)$ for all $a \in A$.

**Definition 2.2.10.** Let $\{B_N(G),*,\circ\}$ be a neutrosophic strong bigroupoid and $(F,A)$ be a soft neutrosophic strong bigroupoid over $\{B_N(G),*,\circ\}$. Then $(F,A)$ is called soft Lagrange neutrosophic strong bigroupoid if and only if $F(a)$ is a Lagrange neutrosophic strong sub bigroupoid of $\{B_N(G),*,\circ\}$ for all $a \in A$.

**Theorem 2.2.13.** Every soft Lagrange neutrosophic strong bigroupoid over $\{B_N(G),*,\circ\}$ is a soft neutrosophic strong bigroupoid over $\{B_N(G),*,\circ\}$ but the converse is not true.

One can easily see the converse by the help of examples.

**85**

**Theorem 2.2.14.** If $\{B_N(G), *, \circ\}$ is a Lagrange neutrosophic strong bigroupoid, then $(F, A)$ over $\{B_N(G), *, \circ\}$ is a soft Lagrange neutrosophic strong bigroupoid but the converse is not true.

**Remark 2.2.7.** Let $(F, A)$ and $(K, C)$ be two soft Lagrange neutrosophic strong bigroupoids over $\{B_N(G), *, \circ\}$. Then

1. Their extended intersection $(F, A) \cap_E (K, C)$ may not be a soft Lagrange neutrosophic strong bigroupoid over $\{B_N(G), *, \circ\}$.

2. Their restricted intersection $(F, A) \cap_R (K, C)$ may not be a soft Lagrange neutrosophic strong bigroupoid over $\{B_N(G), *, \circ\}$.

3. Their *AND* operation $(F, A) \wedge (K, C)$ may not be a soft Lagrange neutrosophic strong bigroupoid over $\{B_N(G), *, \circ\}$. Their extended union $(F, A) \cup_E (K, C)$ may not be a soft Lagrange neutrosophic strong bigroupoid over $\{B_N(G), *, \circ\}$.

4. Their restricted union $(F, A) \cup_R (K, C)$ may not be a soft Lagrange neutrosophic strong bigroupoid over $\{B_N(G), *, \circ\}$.

5. Their extended union $(F, A) \cap_E (K, C)$ may not be a soft Lagrange neutrosophic strong bigroupoid over $\{B_N(G), *, \circ\}$.

6. Their *OR* operation $(F, A) \vee (K, C)$ may not be a soft Lagrange neutrosophic strong bigroupoid over $\{B_N(G), *, \circ\}$.

One can easily verify (1),(2),(3),(4),(5) and (6) by the help of examples.

**Definition 2.2.11.** Let $\{B_N(G), *, \circ\}$ be a neutrosophic strong bigroupoid and $(F, A)$ be a soft neutrosophic strong bigroupoid over $\{B_N(G), *, \circ\}$. Then $(F, A)$ is called soft weak Lagrange neutrosophic strong bigroupoid if at least one $F(a)$ is not a Lagrange neutrosophic strong sub bigroupoid of $\{B_N(G), *, \circ\}$ for some $a \in A$.

**Theorem 2.2.15.** Every soft weak Lagrange neutrosophic strong bigroupoid over $\{B_N(G), *, \circ\}$ is a soft neutrosophic strong bigroupoid over $\{B_N(G), *, \circ\}$ but the converse is not true.

**Theorem 2.2.16.** If $\{B_N(G), *, \circ\}$ is weak Lagrange neutrosophic strong bigroupoid, then $(F, A)$ over $\{B_N(G), *, \circ\}$ is also soft weak Lagrange neutrosophic strong bigroupoid but the converse is not true.

**Remark 2.2.8.** Let $(F, A)$ and $(K, C)$ be two soft weak Lagrange neutrosophic strong bigroupoids over $\{B_N(G), *, \circ\}$. Then

1. Their extended intersection $(F, A) \cap_E (K, C)$ is not a soft weak Lagrange neutrosophic strong bigroupoid over $\{B_N(G), *, \circ\}$.
2. Their restricted intersection $(F, A) \cap_R (K, C)$ is not a soft weak Lagrange neutrosophic strong bigroupoid over $\{B_N(G), *, \circ\}$.
3. Their *AND* operation $(F, A) \wedge (K, C)$ is not a soft weak Lagrange neutrosophic strong bigroupoid over $\{B_N(G), *, \circ\}$.
4. Their extended union $(F, A) \cup_E (K, C)$ is not a soft weak Lagrnage neutrosophic strong bigroupoid over $\{B_N(G), *, \circ\}$.
5. Their restricted union $(F, A) \cup_R (K, C)$ is not a soft weak Lagrange neutrosophic strong bigroupoid over $\{B_N(G), *, \circ\}$.
6. Their *OR* operation $(F, A) \vee (K, C)$ is not a soft weak Lagrange neutrosophic strong bigroupoid over $\{B_N(G), *, \circ\}$.

One can easily verify $(1), (2), (3), (4), (5)$ and $(6)$ by the help of examples.

**Definition 2.2.12.** Let $\{B_N(G), *, \circ\}$ be a neutrosophic strong bigroupoid and $(F, A)$ be a soft neutrosophic strong bigroupoid over $\{B_N(G), *, \circ\}$. Then $(F, A)$ is called soft Lagrange free neutrosophic strong bigroupoid if $F(a)$ is not a Lagrange neutrosophic strong sub bigroupoid of $\{B_N(G), *, \circ\}$ for all $a \in A$.

**Theorem 2.2.17.** Every soft Lagrange free neutrosophic strong bigroupoid over $\{B_N(G), *, \circ\}$ is a soft neutrosophic strong bigroupoid over $\{B_N(G), *, \circ\}$ but the converse is not true.

**Theorem 2.2.18.** If $\{B_N(G), *, \circ\}$ is a Lagrange free neutrosophic strong bigroupoid, then $(F, A)$ over $\{B_N(G), *, \circ\}$ is also a soft Lagrange free neutrosophic strong bigroupoid but the converse is not true.

**Remark 2.2.9.** Let $(F, A)$ and $(K, C)$ be two soft Lagrange free neutrosophic strong bigroupoids over $\{B_N(G), *, \circ\}$. Then

1. Their extended intersection $(F, A) \cap_E (K, C)$ is not a soft Lagrange free neutrosophic strong bigroupoid over $\{B_N(G), *, \circ\}$.
2. Their restricted intersection $(F, A) \cap_R (K, C)$ is not a soft Lagrange free neutrosophic strong bigroupoid over $\{B_N(G), *, \circ\}$.
3. Their *AND* operation $(F, A) \wedge (K, C)$ is not a soft Lagrange free neutrosophic strong bigroupoid over $\{B_N(G), *, \circ\}$.
4. Their extended union $(F, A) \cup_E (K, C)$ is not a soft Lagrange free neutrosophic strong bigroupoid over $\{B_N(G), *, \circ\}$.
5. Their restricted union $(F, A) \cup_R (K, C)$ is not a soft Lagrange free

neutrosophic strong bigroupoid over $\{B_N(G), *, \circ\}$.

6. Their *OR* operation $(F,A) \vee (K,C)$ is not a soft Lagrange free neutrosophic strong bigroupoid over $\{B_N(G), *, \circ\}$.

One can easily verify (1),(2),(3),(4),(5) and (6) by the help of examples.

**Definition 2.2.13.** $(F,A)$ is called soft neutrosophic strong biideal over $\{B_N(G), *, \circ\}$ if $F(a)$ is a neutrosophic strong biideal of $\{B_N(G), *, \circ\}$, for all $a \in A$.

**Theorem 2.2.19.** Every soft neutrosophic strong biideal $(F,A)$ over $\{B_N(G), *, \circ\}$ is a soft neutrosophic strong bigroupoid.

**Proposition 2.2.4.** Let $(F,A)$ and $(K,B)$ be two soft neutrosophic strong biideals over $\{B_N(G), *, \circ\}$. Then

1. Their extended intersection $(F,A) \cap_E (K,B)$ is soft neutrosophic strong biideal over $\{B_N(G), *, \circ\}$.
2. Their restricted intersection $(F,A) \cap_R (K,B)$ is soft neutrosophic strong biideal over $\{B_N(G), *, \circ\}$.
3. Their *AND* operation $(F,A) \wedge (K,B)$ is soft neutrosophic strong biideal over $\{B_N(G), *, \circ\}$.

**Remark 2.2.10.** Let $(F,A)$ and $(K,B)$ be two soft neutrosophic strong biideals over $\{B_N(G), *, \circ\}$. Then

1. Their extended union $(F,A) \cup_E (K,B)$ is not soft neutrosophic strong biideals over $\{B_N(G), *, \circ\}$.

**89**

2. Their restricted union $(F,A) \cup_R (K,B)$ is not soft neutrosophic strong biidleals over $\{B_N(G), *, \circ\}$.

3. Their $OR$ operation $(F,A) \vee (K,B)$ is not soft neutrosophic strong biideals over $\{B_N(G), *, \circ\}$.

One can easily proved $(1),(2),$ and $(3)$ by the help of examples

**Theorem 2.2.20.** Let $(F,A)$ be a soft neutrosophic strong biideal over $\{B_N(G), *, \circ\}$ and $\{(H_i, B_i): i \in J\}$ is a non-empty family of soft neutrosophic strong biideals of $(F,A)$. Then

1. $\underset{i \in J}{\cap}(H_i, B_i)$ is a soft neutrosophic strong biideal of $(F,A)$.

2. $\underset{i \in J}{\wedge}(H_i, B_i)$ is a soft neutrosophic strong biideal of $\underset{i \in J}{\wedge}(F,A)$.

In this last section of chapter two, we generalize the concept of soft neutrosophic groupoid.

## 2.3  Soft Neutrosophic N-groupoid

We now extend the definition of soft neutrosophic groupoid over a neutrosophic groupoid to soft neutrosophic N-groupoid over neutrosophic N-groupoid. Further more, we also define the strong part of soft neutrosophic N-groupoid at the end of this section. The basic and fundamental properties are established with sufficient amount of examples.

We proceed to define soft neutrosophic N-groupoid over a neutrosophic N-groupoid.

**Definition 2.3.1.** Let $N(G) = \{G_1 \cup G_2 \cup \ldots \cup G_N, *_1, *_2, \ldots, *_N\}$ be a neutrosophic N-groupoid and $(F, A)$ be a soft set over $N(G) = \{G_1 \cup G_2 \cup \ldots \cup G_N, *_1, *_2, \ldots, *_N\}$. Then $(F, A)$ is called soft neutrosophic N-groupoid if and only if $F(a)$ is neutrosophic sub N-groupoid of $N(G) = \{G_1 \cup G_2 \cup \ldots \cup G_N, *_1, *_2, \ldots, *_N\}$ for all $a \in A$.

Lets take a look to the following example for illustration.

**Example 2.3.1.** Let $N(G) = \{G_1 \cup G_2 \cup G_3, *_1, *_2, *_3\}$ be a neutrosophic 3-groupoid, where $G_1 = \{\langle Z_{10} \cup I \rangle \mid a * b = 2a + 3b (\mathrm{mod}\, 10); a, b \in \langle Z_{10} \cup I \rangle\}$, $G_2 = \{\langle Z_4 \cup I \rangle \mid a \circ b = 2a + b (\mathrm{mod}\, 4); a, b \in \langle Z_4 \cup I \rangle\}$ and $G_3 = \{\langle Z_{12} \cup I \rangle \mid a * b = 8a + 4b (\mathrm{mod}\, 12); a, b \in \langle Z_{12} \cup I \rangle\}$.

Let $A = \{a_1, a_2\}$ be a set of parameters. Then $(F, A)$ is a soft neutrosophic N-groupoid over $N(G) = \{G_1 \cup G_2 \cup G_3, *_1, *_2, *_3\}$, where

$$F(a_1) = \{0, 5, 5I, 5 + 5I\} \cup \{0, 2, 2I, 2 + 2I\} \cup \{0, 2\},$$

$$F(a_2) = (Z_{10}, *) \cup \{0, 2 + 2I\} \cup \{0, 2I\}.$$

**Theorem 2.3.1.** Let $(F, A)$ and $(H, A)$ be two soft neutrosophic N-groupoids over $N(G)$. Then their intersection $(F, A) \cap (H, A)$ is again a soft neutrosophic N-groupoid over $N(G)$.

The proof is straightforward so left as an exercise for the interested readers.

**Theorem 2.3.2.** Let $(F,A)$ and $(H,B)$ be two soft neutrosophic N-groupoids over $N(G)$. If $A \cap B = \phi$, then $(F,A) \cup (H,B)$ is a soft neutrosophic N-groupoid over $N(G)$.

**Proposition 2.3.1.** Let $(F,A)$ and $(K,C)$ be two soft neutrosophic N-groupoids over $N(G)$. Then

1. Their extended intersection $(F,A) \cap_E (K,C)$ is a soft neutrosophic N-groupoid over $N(G)$.
2. Their restricted intersection $(F,A) \cap_R (K,C)$ is a soft neutrosophic N-groupoid over $N(G)$.
3. Their *AND* operation $(F,A) \wedge (K,C)$ is a soft neutrosophic N-groupoid over $N(G)$.

**Remark 2.3.1.** Let $(F,A)$ and $(K,C)$ be two soft neutrosophic N-groupoids over $N(G)$. Then

1. Their extended union $(F,A) \cup_E (K,C)$ is not a soft neutrosophic N-groupoid over $N(G)$.
2. Their restricted union $(F,A) \cup_R (K,C)$ is not a soft neutrosophic N-groupoid over $N(G)$.
3. Their *OR* operation $(F,A) \vee (K,C)$ is not a soft neutrosophic N-groupoid over $N(G)$.

One can easily verify $(1),(2),$ and $(3)$ by the help of examples.

**92**

**Definition 2.3.2.** Let $(F,A)$ be a soft neutrosophic N-groupoid over $N(G)$. Then $(F,A)$ is called an absolute soft neutrosophic N-groupoid over $N(G)$ if $F(a) = N(G)$ for all $a \in A$.

**Definition 2.3.3.** Let $(F,A)$ and $(H,C)$ be two soft neutrosophic N-groupoids over $N(G)$. Then $(H,C)$ is called soft neutrosophic sub N-groupoid of $(F,A)$, if

1. $C \subseteq A$.
2. $H(a)$ is a neutrosophic sub bigroupoid of $F(a)$ for all $a \in A$.

**Example 2.3.2.** Let $N(G) = \{G_1 \cup G_2 \cup G_3, *_1, *_2, *_3\}$ be a neutrosophic 3-groupoid, where $G_1 = \{\langle Z_{10} \cup I \rangle \mid a*b = 2a + 3b \pmod{10}; a,b \in \langle Z_{10} \cup I \rangle\}$, $G_2 = \{\langle Z_4 \cup I \rangle \mid a \circ b = 2a + b \pmod{4}; a,b \in \langle Z_4 \cup I \rangle\}$ and $G_3 = \{\langle Z_{12} \cup I \rangle \mid a*b = 8a + 4b \pmod{12}; a,b \in \langle Z_{12} \cup I \rangle\}$.

Let $A = \{a_1, a_2\}$ be a set of parameters. Then $(F,A)$ is a soft neutrosophic N-groupoid over $N(G) = \{G_1 \cup G_2 \cup G_3, *_1, *_2, *_3\}$, where

$$F(a_1) = \{0,5,5I,5+5I\} \cup \{0,2,2I,2+2I\} \cup \{0,2\},$$

$$F(a_2) = (Z_{10},*) \cup \{0,2+2I\} \cup \{0,2I\}.$$

Let $B = \{a_1\} \subseteq A$. Then $(H,B)$ is a soft neutrosophic sub N-groupoid of $(F,A)$, where

$$H(a_1) = \{0,5\} \cup \{0,2+2I\} \cup \{0,2\}.$$

**Definition 2.3.4.** Let $N(G)$ be a neutrosophic N-groupoid and $(F,A)$ be a soft neutrosophic N-groupoid over $N(G)$. Then $(F,A)$ is called soft Lagrange neutrosophic N-groupoid if and only if $F(a)$ is a Lagrange neutrosophic sub N-groupoid of $N(G)$ for all $a \in A$.

The interested readers can construct a lot of examples easily for better understanding.

**Theorem 2.3.3.** Every soft Lagrange neutrosophic N-groupoid over $N(G)$ is a soft neutrosophic N-groupoid over $N(G)$ but the converse may not be true.

One can easily see the converse by the help of examples.

**Theorem 2.3.4.** If $N(G)$ is a Lagrange neutrosophic N-groupoid, then $(F,A)$ over $N(G)$ is a soft Lagrange neutrosophic N-groupoid but the converse is not true.

**Remark 2.3.2.** Let $(F,A)$ and $(K,C)$ be two soft Lagrange neutrosophic N-groupoids over $N(G)$. Then

1. Their extended intersection $(F,A) \cap_E (K,C)$ may not be a soft Lagrange neutrosophic N-groupoid over $N(G)$.
2. Their restricted intersection $(F,A) \cap_R (K,C)$ may not be a soft Lagrange neutrosophic N-groupoid over $N(G)$.
3. Their $AND$ operation $(F,A) \wedge (K,C)$ may not be a soft Lagrange neutrosophic N-groupoid over $N(G)$.
4. Their extended union $(F,A) \cup_E (K,C)$ may not be a soft Lagrange neutrosophic N-groupoid over $N(G)$.

**94**

5. Their restricted union $(F,A) \cup_R (K,C)$ may not be a soft Lagrange neutrosophic N-groupoid over $N(G)$.
6. Their *OR* operation $(F,A) \vee (K,C)$ may not be a soft Lagrange neutrosophic N-groupoid over $N(G)$.

One can easily verify $(1),(2),(3),(4),(5)$ and $(6)$ by the help of examples.

**Definition 2.3.5.** Let $N(G)$ be a neutrosophic N-groupoid and $(F,A)$ be a soft neutrosophic N-groupoid over $N(G)$. Then $(F,A)$ is called soft weak Lagrange neutrosophic N-groupoid if at least one $F(a)$ is not a Lagrange neutrosophic sub N-groupoid of $N(G)$ for some $a \in A$.

**Theorem 2.3.5.** Every soft weak Lagrange neutrosophic N-groupoid over $N(G)$ is a soft neutrosophic N-groupoid over $N(G)$ but the converse is not true.

**Theorem 2.3.6.** If $N(G)$ is weak Lagrange neutrosophic N-groupoid, then $(F,A)$ over $N(G)$ is also a soft weak Lagrange neutrosophic bigroupoid but the converse is not true.

**Remark 2.3.3.** Let $(F,A)$ and $(K,C)$ be two soft weak Lagrange neutrosophic N-groupoids over $N(G)$. Then

1. Their extended intersection $(F,A) \cap_E (K,C)$ may not be a soft weak Lagrange neutrosophic N-groupoid over $N(G)$.
2. Their restricted intersection $(F,A) \cap_R (K,C)$ may not be a soft weak Lagrange neutrosophic N-groupoid over $N(G)$.
3. Their *AND* operation $(F,A) \wedge (K,C)$ may not be a soft weak

Lagrange neutrosophic N-groupoid over $N(G)$.

4. Their extended union $(F, A) \cup_E (K, C)$ may not be a soft weak Lagrnage neutrosophic N-groupoid over $N(G)$.

5. Their restricted union $(F, A) \cup_R (K, C)$ may not be a soft weak Lagrange neutrosophic N-groupoid over $N(G)$.

6. Their *OR* operation $(F, A) \vee (K, C)$ may not be a soft weak Lagrange neutrosophic N-groupoid over $N(G)$.

One can easily verify $(1), (2), (3), (4), (5)$ and $(6)$ by the help of examples.

**Definition 2.3.6.** Let $N(G)$ be a neutrosophic N-groupoid and $(F, A)$ be a soft neutrosophic N-groupoid over $N(G)$. Then $(F, A)$ is called soft Lagrange free neutrosophic N-groupoid if $F(a)$ is not a Lagrange neutrosophic sub N-groupoid of $N(G)$ for all $a \in A$.

**Theorem 2.3.7.** Every soft Lagrange free neutrosophic N-groupoid over $N(G)$ is a soft neutrosophic N-groupoid over $N(G)$ but the converse is not true.

**Theorem 2.3.8.** If $N(G)$ is a Lagrange free neutrosophic N-groupoid, then $(F, A)$ over $N(G)$ is also a soft Lagrange free neutrosophic N-groupoid but the converse is not true.

**Remark 2.3.4.** Let $(F, A)$ and $(K, C)$ be two soft Lagrange free neutrosophic N-groupoids over $N(G)$. Then

1. Their extended intersection $(F, A) \cap_E (K, C)$ is not a soft Lagrange free neutrosophic N-groupoid over $N(G)$.

2. Their restricted intersection $(F,A) \cap_R (K,C)$ is not a soft Lagrange free neutrosophic N-groupoid over $N(G)$.

3. Their *AND* operation $(F,A) \wedge (K,C)$ is not a soft Lagrange free neutrosophic N-groupoid over $N(G)$.

4. Their extended union $(F,A) \cup_E (K,C)$ is not a soft Lagrange free neutrosophic N-groupoid over $N(G)$.

5. Their restricted union $(F,A) \cup_R (K,C)$ is not a soft Lagrange free neutrosophic N-groupoid over $N(G)$.

6. Their *OR* operation $(F,A) \vee (K,C)$ is not a soft Lagrange free neutrosophic N-groupoid over $N(G)$.

One can easily verify $(1),(2),(3),(4),(5)$ and $(6)$ by the help of examples.

**Definition 2.3.7.** $(F,A)$ is called soft neutrosophic N-ideal over $N(G)$ if and only if $F(a)$ is a neutrosophic N-ideal of $N(G)$, for all $a \in A$.

The interested readers can construct a lot of examples easily for better understanding.

**Theorem 2.3.8.** Every soft neutrosophic N-ideal $(F,A)$ over $N(G)$ is a soft neutrosophic N-groupoid.

**Proposition 2.3.2.** Let $(F,A)$ and $(K,B)$ be two soft neutrosophic N-ideals over $N(G)$. Then

1. Their extended intersection $(F,A) \cap_E (K,B)$ is soft neutrosophic N-ideal over $N(G)$.

2. Their restricted intersection $(F,A) \cap_R (K,B)$ is soft neutrosophic N-

**97**

ideal over $N(G)$.

3. Their *AND* operation $(F,A) \wedge (K,B)$ is soft neutrosophic N-ideal over $N(G)$.

**Remark 2.3.5.** Let $(F,A)$ and $(K,B)$ be two soft neutrosophic N-ideals over $N(G)$. Then

1. Their extended union $(F,A) \cup_E (K,B)$ is not a soft neutrosophic N-ideal over $N(G)$.

2. Their restricted union $(F,A) \cup_R (K,B)$ is not a soft neutrosophic N-idleal over $N(G)$.

3. Their *OR* operation $(F,A) \vee (K,B)$ is not a soft neutrosophic N-ideal over $N(G)$.

One can easily proved (1),(2), and (3) by the help of examples

**Theorem 2.3.9.** Let $(F,A)$ be a soft neutrosophic N-ideal over $N(G)$ and $\{(H_i, B_i) : i \in J\}$ be a non-empty family of soft neutrosophic N-ideals of $(F,A)$. Then

1. $\underset{i \in J}{\cap}(H_i, B_i)$ is a soft neutrosophic N-ideal of $(F,A)$.

2. $\underset{i \in J}{\wedge}(H_i, B_i)$ is a soft neutrosophic N-ideal of $\underset{i \in J}{\wedge}(F,A)$.

## Soft Neutrosophic Strong N-groupoid

Here we give the definition of soft neutrosophic strong N-groupoid over a neutrosophic N-groupoid. Some important facts can also be studied over here.

**Definition 2.3.8.** Let $N(G) = \{G_1 \cup G_2 \cup ... \cup G_N, *_1, *_2, ..., *_N\}$ be a neutrosophic N-groupoid and $(F, A)$ be a soft set over $N(G) = \{G_1 \cup G_2 \cup ... \cup G_N, *_1, *_2, ..., *_N\}$. Then $(F, A)$ is called soft neutrosophic strong N-groupoid if and only if $F(a)$ is neutrosophic strong sub N-groupoid of $N(G) = \{G_1 \cup G_2 \cup ... \cup G_N, *_1, *_2, ..., *_N\}$ for all $a \in A$.

**Example 2.3.3.** Let $N(G) = \{G_1 \cup G_2 \cup G_3, *_1, *_2, *_3\}$ be a neutrosophic 3-groupoid, where $G_1 = \{\langle Z_{10} \cup I \rangle \mid a*b = 2a + 3b \pmod{10}; a, b \in \langle Z_{10} \cup I \rangle\}$, $G_2 = \{\langle Z_4 \cup I \rangle \mid a \circ b = 2a + b \pmod 4; a, b \in \langle Z_4 \cup I \rangle\}$ and $G_3 = \{\langle Z_{12} \cup I \rangle \mid a*b = 8a + 4b \pmod{12}; a, b \in \langle Z_{12} \cup I \rangle\}$. Let $A = \{a_1, a_2\}$ be a set of parameters.

Then $(F, A)$ is a soft neutrosophic N-groupoid over $N(G) = \{G_1 \cup G_2 \cup G_3, *_1, *_2, *_3\}$, where

$$F(a_1) = \{0, 5I\} \cup \{0, 2I\} \cup \{0, 2I\},$$

$$F(a_2) = \{0, 5 + 5I\} \cup \{0, 2 + 2I\} \cup \{0, 2 + 2I\}.$$

**Theorem 2.3.10.** Let $(F, A)$ and $(H, A)$ be two soft neutrosophic strong N-groupoids over $N(G)$. Then their intersection $(F, A) \cap (H, A)$ is again a soft neutrosophic strong N-groupoid over $N(G)$.

The proof is straightforward so left as an exercise for the readers.

**Theorem 2.3.11.** Let $(F,A)$ and $(H,B)$ be two soft neutrosophic strong N-groupoids over $N(G)$. If $A \cap B = \phi$, then $(F,A) \cup (H,B)$ is a soft neutrosophic strong N-groupoid over $N(G)$.

**Theorem 2.3.12.** If $N(G)$ is a neutrosophic strong N-groupoid, then $(F,A)$ over $N(G)$ is also a soft neutrosophic strong N-groupoid.

**Proposition 2.3.3.** Let $(F,A)$ and $(K,C)$ be two soft neutrosophic strong N-groupoids over $N(G)$. Then

1. Their extended intersection $(F,A) \cap_E (K,C)$ is a soft neutrosophic strong N-groupoid over $N(G)$.
2. Their restricted intersection $(F,A) \cap_R (K,C)$ is a soft neutrosophic strong N-groupoid over $N(G)$.
3. Their $AND$ operation $(F,A) \wedge (K,C)$ is a soft neutrosophic strong N-groupoid over $N(G)$.

**Remark 2.3.6.** Let $(F,A)$ and $(K,C)$ be two soft neutrosophic strong N-groupoids over $N(G)$. Then

1. Their extended union $(F,A) \cup_E (K,C)$ is not a soft neutrosophic strong N-groupoid over $N(G)$.
2. Their restricted union $(F,A) \cup_R (K,C)$ is not a soft neutrosophic strong N-groupoid over $N(G)$.

**100**

3. Their *OR* operation $(F,A) \vee (K,C)$ is not a soft neutrosophic strong N-groupoid over $N(G)$.

One can easily verify $(1),(2),$ and $(3)$ by the help of examples.

We now give the definition of soft neutrosophic strong sub N-groupoid of a soft neutrosophic strong N-groupoid over a neutrosophic N-groupoid.

**Definition 43.** Let $(F,A)$ and $(H,C)$ be two soft neutrosophic strong N-groupoids over $N(G)$. Then $(H,C)$ is called soft neutrosophic strong sub N-groupoid of $(F,A)$, if

1. $C \subseteq A$.
2. $H(a)$ is a neutrosophic sub bigroupoid of $F(a)$ for all $a \in A$.

**Definition 2.3.9.** Let $N(G)$ be a neutrosophic strong N-groupoid and $(F,A)$ be a soft neutrosophic strong N-groupoid over $N(G)$. Then $(F,A)$ is called soft Lagrange neutrosophic strong N-groupoid if and only if $F(a)$ is a Lagrange neutrosophic sub N-groupoid of $N(G)$ for all $a \in A.$

The interested readers can construct a lot of examples easily for better understanding.

**Theorem 2.3.13.** Every soft Lagrange neutrosophic strong N-groupoid over $N(G)$ is a soft neutrosophic N-groupoid over $N(G)$ but the converse may not be true.

One can easily see the converse by the help of examples.

**Theorem 2.3.14.** Every soft Lagrange neutrosophic strong N-groupoid over $N(G)$ is a soft Lagrange neutrosophic N-groupoid over $N(G)$ but the converse may not be true.

**Theorem 2.3.15.** If $N(G)$ is a Lagrange neutrosophic strong N-groupoid, then $(F,A)$ over $N(G)$ is a soft Lagrange neutrosophic strong N-groupoid but the converse is not true.

**Remark 2.3.7.** Let $(F,A)$ and $(K,C)$ be two soft Lagrange neutrosophic strong N-groupoids over $N(G).$ Then

1. Their extended intersection $(F,A) \cap_E (K,C)$ may not be a soft Lagrange neutrosophic strong N-groupoid over $N(G)$.
2. Their restricted intersection $(F,A) \cap_R (K,C)$ may not be a soft Lagrange neutrosophic strong N-groupoid over $N(G)$.
3. Their $AND$ operation $(F,A) \wedge (K,C)$ may not be a soft Lagrange neutrosophic strong N-groupoid over $N(G)$.
4. Their extended union $(F,A) \cup_E (K,C)$ may not be a soft Lagrange neutrosophic strong N-groupoid over $N(G)$.
5. Their restricted union $(F,A) \cup_R (K,C)$ may not be a soft Lagrange neutrosophic strong N-groupoid over $N(G)$.
6. Their $OR$ operation $(F,A) \vee (K,C)$ may not be a soft Lagrange neutrosophic strong N-groupoid over $N(G)$.

One can easily verify $(1),(2),(3),(4),(5)$ and $(6)$ by the help of examples.

**Definition 2.3.10.** Let $N(G)$ be a neutrosophic strong N-groupoid and $(F,A)$ be a soft neutrosophic strong N-groupoid over $N(G)$. Then $(F,A)$ is called soft weak Lagrange neutrosophic strong N-groupoid if at least one $F(a)$ is not a Lagrange neutrosophic sub N-groupoid of $N(G)$ for some $a \in A$.

**Theorem 2.3.16.** Every soft weak Lagrange neutrosophic strong N-groupoid over $N(G)$ is a soft neutrosophic N-groupoid over $N(G)$ but the converse is not true.

**Theorem 2.3.17.** Every soft weak Lagrange neutrosophic strong N-groupoid over $N(G)$ is a soft weak Lagrange neutrosophic N-groupoid over $N(G)$ but the converse is not true.

**Remark 2.3.8.** Let $(F,A)$ and $(K,C)$ be two soft weak Lagrange neutrosophic strong N-groupoids over $N(G)$. Then

1. Their extended intersection $(F,A) \cap_E (K,C)$ may not be a soft weak Lagrange neutrosophic strong N-groupoid over $N(G)$.
2. Their restricted intersection $(F,A) \cap_R (K,C)$ may not be a soft weak Lagrange neutrosophic strong N-groupoid over $N(G)$.
3. Their $AND$ operation $(F,A) \wedge (K,C)$ may not be a soft weak Lagrange neutrosophic strong N-groupoid over $N(G)$.
4. Their extended union $(F,A) \cup_E (K,C)$ may not be a soft weak Lagrnage neutrosophic strong N-groupoid over $N(G)$.
5. Their restricted union $(F,A) \cup_R (K,C)$ may not be a soft weak Lagrange neutrosophic strong N-groupoid over $N(G)$.
6. Their $OR$ operation $(F,A) \vee (K,C)$ may not be a soft weak Lagrange

**103**

neutrosophic strong N-groupoid over $N(G)$.

One can easily verify (1),(2),(3),(4),(5) and (6) by the help of examples.

**Definition 2.3.11.** Let $N(G)$ be a neutrosophic strong N-groupoid and $(F,A)$ be a soft neutrosophic strong N-groupoid over $N(G)$. Then $(F,A)$ is called soft Lagrange free neutrosophic strong N-groupoid if $F(a)$ is not a Lagrange neutrosophic sub N-groupoid of $N(G)$ for all $a \in A$.

**Theorem 2.3.18.** Every soft Lagrange free neutrosophic strong N-groupoid over $N(G)$ is a soft neutrosophic N-groupoid over $N(G)$ but the converse is not true.

**Theorem 2.3.19.** Every soft Lagrange free neutrosophic strong N-groupoid over $N(G)$ is a soft Lagrange neutrosophic N-groupoid over $N(G)$ but the converse is not true.

**Theorem 2.3.20.** If $N(G)$ is a Lagrange free neutrosophic strong N-groupoid, then $(F,A)$ over $N(G)$ is also a soft Lagrange free neutrosophic strong N-groupoid but the converse is not true.

**Remark 2.3.9.** Let $(F,A)$ and $(K,C)$ be two soft Lagrange free neutrosophic N-groupoids over $N(G)$. Then

1. Their extended intersection $(F,A) \cap_E (K,C)$ is not a soft Lagrange free neutrosophic strong N-groupoid over $N(G)$.
2. Their restricted intersection $(F,A) \cap_R (K,C)$ is not a soft Lagrange free neutrosophic strong N-groupoid over $N(G)$.

**104**

3. Their *AND* operation $(F,A) \wedge (K,C)$ is not a soft Lagrange free neutrosophic strong N-groupoid over $N(G)$.

4. Their extended union $(F,A) \cup_E (K,C)$ is not a soft Lagrange free neutrosophic strong N-groupoid over $N(G)$.

5. Their restricted union $(F,A) \cup_R (K,C)$ is not a soft Lagrange free neutrosophic strong N-groupoid over $N(G)$.

6. Their *OR* operation $(F,A) \vee (K,C)$ is not a soft Lagrange free neutrosophic stong N-groupoid over $N(G)$.

One can easily verify (1),(2),(3),(4),(5) and (6) by the help of examples.

**Definition 2.3.12.** $(F,A)$ is called soft neutrosophic strong N-ideal over $N(G)$ if and only if $F(a)$ is a neutrosophic strong N-ideal of $N(G)$, for all $a \in A$.

The interested readers can construct a lot of examples easily for better understanding.

**Theorem 2.3.21.** Every soft neutrosophic strong N-ideal $(F,A)$ over $N(G)$ is a soft neutrosophic N-groupoid.

**Theorem 2.3.22.** Every soft neutrosophic strong N-ideal $(F,A)$ over $N(G)$ is a soft neutrosophic N-ideal but the converse is not true.

**Proposition 2.3.4.** Let $(F,A)$ and $(K,B)$ be two soft neutrosophic strong N-ideals over $N(G)$. Then

1. Their extended intersection $(F,A) \cap_E (K,B)$ is soft neutrosophic

strong N-ideal over $N(G)$.

2. Their restricted intersection $(F,A) \cap_R (K,B)$ is soft neutrosophic strong N-ideal over $N(G)$.

3. Their *AND* operation $(F,A) \wedge (K,B)$ is soft neutrosophicstrong N-ideal over $N(G)$.

**Remark 2.3.10.** Let $(F,A)$ and $(K,B)$ be two soft neutrosophic strong N-ideals over $N(G)$. Then

1. Their extended union $(F,A) \cup_E (K,B)$ is not a soft neutrosophic strong N-ideal over $N(G)$.

2. Their restricted union $(F,A) \cup_R (K,B)$ is not a soft neutrosophic strong N-idleal over $N(G)$.

3. Their *OR* operation $(F,A) \vee (K,B)$ is not a soft neutrosophic strong N-ideal over $N(G)$.

One can easily proved (1),(2), and (3) by the help of examples

**Theorem 2.3.23.** Let $(F,A)$ be a soft neutrosophic strong N-ideal over $N(G)$ and $\{(H_i, B_i) : i \in J\}$ be a non-empty family of soft neutrosophic strong N-ideals of $(F,A)$. Then

1. $\underset{i \in J}{\cap} (H_i, B_i)$ is a soft neutrosophic strong N-ideal of $(F,A)$.

2. $\underset{i \in J}{\wedge} (H_i, B_i)$ is a soft neutrosophic strong N-ideal of $\underset{i \in J}{\wedge} (F,A)$.

# Chapter No. 3

# Soft Nuetrosophic Rings, Soft Neutrosophic Fields and Their Generalization

In this chapter we take the important neutrosophic algebraic structurs neutrosophic rings and neutrosophic fields to define soft neutrosophic rings and soft neutrosophic fields respectively. We also extend this theory and give the general concept of soft neutrosophic N-rings and soft neutrosophic N-fields over neutrosophic N-rings and neutrosophic N-fields respectively. Some core properties are also studied here and many illustrative examples are given.

We proceed to define soft neutrosophic ring over a neutrosophic ring.

## 3.1    Soft Nuetrosophic Ring

In this section, the authors introduced soft neutrosophic ring over a neutrosophic ring for the first time. This means that soft neutrosophic ring is an approximated collection of neutrosophic subrings of the neutrosophic ring. We also coined here some basic and fundamental properties of soft neutrosophic ring over a neutrosophic ring with many cleared examples.

**Definition 3.1.1.** Let $\langle R \cup I \rangle$ be a neutrosophic ring and $(F, A)$ be a soft set over $\langle R \cup I \rangle$. Then $(F, A)$ is called soft neutrosophic ring if and only if $F(a)$ is a neutrosophic subring of $\langle R \cup I \rangle$ for all $a \in A$.

See the following examples which illustrate this situation.

**Example 3.1.1.** Let $\langle Z \cup I \rangle$ be a neutrosophic ring of integers and let $(F, A)$ be a soft set over $\langle Z \cup I \rangle$. Let $A = \{a_1, a_2, a_3, a_4\}$ be a set of parameters. Then clearly $(F, A)$ is a soft neutrosophic ring over $\langle Z \cup I \rangle$, where

$$F(a_1) = \langle 2Z \cup I \rangle, F(a_2) = \langle 3Z \cup I \rangle$$

$$F(a_3) = \langle 5Z \cup I \rangle, F(a_4) = \langle 6Z \cup I \rangle.$$

**Example.3.1.2.** Let $\langle \mathbb{C} \cup I \rangle$ be a neutrosophic ring of complex numbers and let $(F, A)$ be a soft set over $\langle \mathbb{C} \cup I \rangle$. Let $A = \{a_1, a_2, a_3, a_4\}$ be a set of parameters. Then clearly $(F, A)$ is a soft neutrosophic ring over $\langle \mathbb{C} \cup I \rangle$, where

$$F(a_1) = \langle \mathbb{R} \cup I \rangle, F(a_2) = \langle \mathbb{Q} \cup I \rangle$$

$$F(a_3) = \langle \mathbb{Z} \cup I \rangle, F(a_4) = \langle 2\mathbb{Z} \cup I \rangle.$$

**108**

We now give some characterization of soft neutrosophic ring over a neutrosophic ring.

**Theorem 3.1.1.** Let $(F, A)$ and $(H, A)$ be two soft neutrosophic rings over $\langle R \cup I \rangle$. Then their intersection $(F, A) \cap (H, A)$ is again a soft neutrosophic ring over $\langle R \cup I \rangle$.

The proof is straightforward, so left as an exercise for the interested readers.

**Theorem 3.1.2.** Let $(F, A)$ and $(H, B)$ be two soft neutrosophic rings over $\langle R \cup I \rangle$. If $A \cap B = \phi$, then $(F, A) \cup (H, B)$ is a soft neutrosophic ring over $\langle R \cup I \rangle$.

The readers can easily prove this theorem.

**Remark 3.1.1.** The extended union of two soft neutrosophic rings $(F, A)$ and $(K, B)$ over $\langle R \cup I \rangle$ is not a soft neutrosophic ring over $\langle R \cup I \rangle$.

We check this by the help of following example.

**Example 3.1.3.** Let $\langle Z \cup I \rangle$ be a neutrosophic ring of integers. Let $(F, A)$ and $(K, B)$ be two soft neutrosophic rings over $\langle Z \cup I \rangle$, where

$$F(a_1) = \langle 2Z \cup I \rangle, F(a_2) = \langle 3Z \cup I \rangle, F(a_3) = \langle 4Z \cup I \rangle,$$

And

$$K(a_1) = \left\langle 5Z \cup I \right\rangle, K(a_3) = \left\langle 7Z \cup I \right\rangle.$$

Their extended union $(F, A) \cup_E (K, B) = (H, C)$, where

$$H(a_1) = \left\langle 2Z \cup I \right\rangle \cup \left\langle 5Z \cup I \right\rangle,$$

$$H(a_2) = \left\langle 3Z \cup I \right\rangle,$$

$$H(a_3) = \left\langle 5Z \cup I \right\rangle \cup \left\langle 7Z \cup I \right\rangle.$$

Thus clearly $H(a_1) = \left\langle 2Z \cup I \right\rangle \cup \left\langle 5Z \cup I \right\rangle$, $H(a_3) = \left\langle 5Z \cup I \right\rangle \cup \left\langle 7Z \cup I \right\rangle$ is not a neutrosophic subrings of $\left\langle Z \cup I \right\rangle$.

**Remark 3.1.2.** The restricted union of two soft neutrosophic rings $F, A$ and $K, B$ over $\left\langle R \cup I \right\rangle$ is not a soft neutrosophic ring over $\left\langle R \cup I \right\rangle$.

This remark can easily be seen in the above Example 3.1.3.

**Remark 3.1.3.** The $OR$ operation of two soft neutrosophic rings over $\left\langle R \cup I \right\rangle$ may not be a soft neutrosophic ring over $\left\langle R \cup I \right\rangle$.

One can easily check this remark in the Example 3.1.3.

**Proposition 3.1.1.** The extended intersection of two soft neutrosophic rings over $\left\langle R \cup I \right\rangle$ is soft neutrosophic ring over $\left\langle R \cup I \right\rangle$.

The proof is straightforward.

**Proposition 3.1.2.** The restricted intersection of two soft neutrosophic rings over $\langle R \cup I \rangle$ is a soft neutrosophic ring over $\langle R \cup I \rangle$.

**Proposition 3.1.3.** The *AND* operation of two soft neutrosophic rings over $\langle R \cup I \rangle$ is a soft neutrosophic ring over $\langle R \cup I \rangle$.

**Definition 3.1.2.** Let $F, A$ be a soft set over a neutrosophic ring $\langle R \cup I \rangle$. Then $(F, A)$ is called an absolute soft neutrosophic ring if $F(a) = \langle R \cup I \rangle$ for all $a \in A$.

**Definition 3.1.3.** Let $(F, A)$ be a soft set over a neutrosophic ring $\langle R \cup I \rangle$. Then $(F, A)$ is called soft neutrosophic ideal over $\langle R \cup I \rangle$ if and only if $F(a)$ is a neutrosophic ideal of $\langle R \cup I \rangle$ for all $a \in A$.

**Example 3.1.4.** Let $\langle Z_{12} \cup I \rangle$ be a neutrosophic ring. Let $A = \{a_1, a_2\}$ be a set of parameters and $(F, A)$ be a soft set over $\langle Z_{12} \cup I \rangle$. Then clearly $(F, A)$ is a soft neutrosophic ideal over $\langle R \cup I \rangle$, where

$$F(a_1) = \{0, 6, 2I, 4I, 6I, 8I, 10I, 6 + 2I, ..., 6 + 10I\},$$

$$F(a_2) = \{0, 6, 6I, 6 + 6I\}.$$

**Theorem 3.1.3.** Every soft neutrosophic ideal $(F, A)$ over a neutrosophic ring $\langle R \cup I \rangle$ is a soft neutrosophic ring.

**Proof.** Let $(F, A)$ be a soft neutrosophic ideal over a neutrosophic ring $\langle R \cup I \rangle$. Then by definition $F(a)$ is a neutrosophic ideal for all $a \in A$. Since we know that every neutrosophic ideal is a neutrosophic subring. It follows that $F(a)$ is a neutrosophic subring of $\langle R \cup I \rangle$.

Thus by definition of soft neutrosophic ring, this implies that $(F, A)$ is a soft neutrosophic ring.

The converse of the above theorem is not true.

To check the converse, we see the following example.

**Example 3.1.5.** Let $\langle Z_{10} \cup I \rangle$ be a neutrosophic ring. Let $A = \{a_1, a_2\}$ be a set of parameters and $(F, A)$ be a soft neutrosophic ring over $\langle Z_{10} \cup I \rangle$, where

$$F(a_1) = \{0, 2, 4, 6, 8, 2I, 4I, 6I, 8I\},$$

$$\mathrm{F}(a_2) = \{0, 2I, 4I, 6I, 8I\}.$$

Then clearly $(F, A)$ is not a soft neutrosophic ideal over $\langle Z_{10} \cup I \rangle$.

**Proposition 3.1.4.** Let $(F, A)$ and $(K, B)$ be two soft neutosophic ideals over a neutrosophic ring $\langle R \cup I \rangle$. Then

1. Their extended intersection $(F, A) \cap_E (K, B)$ is again a soft neutrosophic ideal over $\langle R \cup I \rangle$.
2. Their restricted intersection $(F, A) \cap_R (K, B)$ is again a soft neutrosophic ideal over $\langle R \cup I \rangle$.

**112**

3. Their *AND* operation $(F,A) \vee (K,B)$ is again a soft neutrosophic ideal over $\langle R \cup I \rangle$.

**Remark 3.1.4.** Let $(F,A)$ and $(K,B)$ be two soft neutosophic ideals over a neutrosophic ring $\langle R \cup I \rangle$. Then

1. Their extended union $(F,A) \cup_E (K,B)$ is not a soft neutrosophic ideal over $\langle R \cup I \rangle$.
2. Their restricted union $(F,A) \cup_R (K,B)$ is not a soft neutrosophic ideal over $\langle R \cup I \rangle$.
3. Their *OR* operation $(F,A) \vee (K,B)$ is not a soft neutrosophic ideal over $\langle R \cup I \rangle$.

One can easily see these remarks with the help of examples.

We now introduce the substructures of soft neutrosophic ring over a neutrosophic ring.

**Definition 3.1.4.** Let $(F,A)$ and $(K,B)$ be two soft neutrosophic rings over $\langle R \cup I \rangle$. Then $(K,B)$ is called soft neutrosophic subring of $(F,A)$, if

1. $B \subseteq A$, and
2. $K(a)$ is a neutrosophic subring of $F(a)$ for all $a \in A$.

**Example 3.1.6.** Let $\langle C \cup I \rangle$ be the neutrosophic ring of complex numbers. Let $A = \{a_1, a_2, a_3\}$ be a set of parameters. Then $(F,A)$ is a soft neutrosophic ring over $\langle C \cup I \rangle$, where

**113**

$$F(a_1) = \langle Z \cup I \rangle, F(a_2) = \langle Q \cup I \rangle,$$

$$F(a_3) = \langle R \cup I \rangle.$$

Where $\langle Z \cup I \rangle$, $\langle Q \cup I \rangle$ and $\langle R \cup I \rangle$ are neutrosophic rings of integers, rational numbers, and real numbers respectively.

Let $B = \{a_2, a_3\}$ be a set of parmeters . Let $(K, B)$ be the neutrosophic subring of $(F, A)$ over $\langle C \cup I \rangle$, where

$$K(a_2) = \langle Z \cup I \rangle, K(a_3) = \langle Q \cup I \rangle.$$

**Theorem 3.1.4.** Every soft ring $(H, B)$ over a ring $R$ is a soft neutrosophic subring of a soft neutrosophic ring $(F, A)$ over the corresponding neutrosophic ring $\langle R \cup I \rangle$ if $B \subseteq A$.

This is straightforward and the proof is left for the interested readers as an exercise.

**Definition 3.1.5.** Let $(F, A)$ and $(K, B)$ be two soft neutrosophic rings over $\langle R \cup I \rangle$. Then $(K, B)$ is called soft neutrosophic ideal of $(F, A)$, if

1. $B \subseteq A$, and
2. $K(a)$ is a neutrosophic ideal of $F(a)$ for all $a \in A$.

**Example 3.1.7.** Let $\langle Z_{12} \cup I \rangle$ be a neutrosophic ring. Let $A = \{a_1, a_2\}$ be a set of parameters and $(F, A)$ be a soft set over $\langle Z_{12} \cup I \rangle$. Then clearly $(F, A)$ is a soft neutrosophic ring over $\langle Z_{12} \cup I \rangle$, where

**114**

$$F(a_1) = \{0, 6, 2I, 4I, 6I, 8I, 10I, 6 + 2I, ..., 6 + 10I\}\,,$$

$$F(a_2) = \{0, 2, 4, 6, 8, 2I, 4I, 6I, 8I\}\,.$$

Let $B = \{a_1, a_2\}$ be a set of parameters. Then clearly $(H, B)$ is a soft neutrosophic ideal of $(F, A)$ over $\langle Z_{12} \cup I \rangle$, where

$$H(a_1) = \{0, 6, 6 + 6I\}\,,$$

$$H(a_2) = \{0, 2, 4, 6, 8\}\,.$$

**Proposition 3.1.5.** All soft neutrosophic ideals are trivially soft neutrosophic subrings.

Now we proceed furthere to the generalization of soft neutrosophic ring.

## 3.2   Soft Neutrosophic Birings

In this section, soft neutrosophic biring over a neutrosophic biring is introduced. This is the parameterized collection of neutrosophic subrings of the neutrosophic ring. We also give some core properties of it with sufficient amount of examples.

**Definition 3.2.1.** Let $(BN(\mathbb{R}), *, \circ)$ be a neutrosophic biring and $(F, A)$ be a soft set over $(BN(\mathbb{R}), *, \circ)$. Then $(F, A)$ is called soft neutrosophic biring if and only if $F(a)$ is a neutrosophic subbiring of $(BN(\mathbb{R}), *, \circ)$ for all $a \in A$.

The following example illustrates this facts.

**Example 3.2.1.** Let $BN(\mathbb{R}) = (\mathbb{R}_1, *, \circ) \cup (\mathbb{R}_2, *, \circ)$ be a neutrosophic biring, where $(\mathbb{R}_1, *, \circ) = (\langle \mathbb{Z} \cup I \rangle, +, \times)$ and $(\mathbb{R}_2, *, \circ) = (\mathbb{Q}, +, \times)$.

Let $A = \{a_1, a_2, a_3, a_4\}$ be a set of parameters. Then clearly $(F, A)$ is a soft neutrosophic biring over $BN(R)$, where

$$F(a_1) = \langle 2\mathbb{Z} \cup I \rangle \cup \mathbb{R}, F(a_2) = \langle 3\mathbb{Z} \cup I \rangle \cup \mathbb{Q},$$

$$F(a_3) = \langle 5\mathbb{Z} \cup I \rangle \cup \mathbb{Z}, F(a_4) = \langle 6\mathbb{Z} \cup I \rangle \cup 2\mathbb{Z}.$$

**Example 3.2.2.** Let $BN(\mathbb{R}) = (\mathbb{R}_1, *, \circ) \cup (\mathbb{R}_2, *, \circ)$ be a neutrosophic biring, where $(\mathbb{R}_1, *, \circ) = (\langle \mathbb{C} \cup I \rangle, +, \times)$ and $(\mathbb{R}_2, *, \circ) = (\mathbb{R}, +, \times)$. Let $A = \{a_1, a_2, a_3, a_4\}$ be a set of parameters. Then clearly $(F, A)$ is a soft neutrosophic biring over $BN(R)$, where

$$F(a_1) = \langle \mathbb{R} \cup I \rangle \cup \mathbb{Q}, F(a_2) = \langle \mathbb{Q} \cup I \rangle \cup \mathbb{Z}$$

$$F(a_3) = \langle \mathbb{Z} \cup I \rangle \cup 2\mathbb{Z}, F(a_4) = \mathbb{Z} \cup 3\mathbb{Z}.$$

**Theorem 3.2.1 .** Let $(F, A)$ and $(H, A)$ be two soft neutrosophic birings over $BN(R)$. Then their intersection $(F, A) \cap (H, A)$ is again a soft neutrosophic biring over $BN(R)$ .

The proof is straightforward.

**Theorem 3.2.2.** Let $F, A$ and $H, B$ be two soft neutrosophic birings over $BN(R)$. If $A \cap B = \phi$, then $F, A \cup H, B$ is a soft neutrosophic biring over $BN(R)$.

**Remark 3.2.1.** The extended union of two soft neutrosophic birings $F, A$ and $K, B$ over $BN(R)$ is not a soft neutrosophic biring over $BN(R)$.

We can easily check this remark by the help of examples.

**Remark 3.2.2.** The restricted union of two soft neutrosophic rings $F, A$ and $K, B$ over $\langle R \cup I \rangle$ is not a soft neutrosophic ring over $\langle R \cup I \rangle$.

We can easily check this remark by the help of examples.

**Remark 3.2.3.** The $OR$ operation of two soft neutrosophic rings over $\langle R \cup I \rangle$ may not be a soft neutrosophic ring over $\langle R \cup I \rangle$.

One can easily check this remark with the help of example.

**Proposition 3.2.1.** The extended intersection of two soft neutrosophic birings over $BN(R)$ is soft neutrosophic biring over $BN(R)$.

The proof is straightforward, so left as an exercise for the readers.

**Proposition 3.2.2.** The restricted intersection of two soft neutrosophic birings over $BN(R)$ is a soft neutrosophic biring over $BN(R)$.

**Proposition 3.2.3.** The $AND$ operation of two soft neutrosophic birings over $BN(R)$ is soft neutrosophic biring over $BN(R)$.

**Proof.** The proof is easy.

**Definition 3.2.2.** Let $F, A$ be a soft set over a neutrosophic biring over $BN(R)$. Then $(F, A)$ is called an absolute soft neutrosophic biring if $F(a) = BN(R)$ for all $a \in A$.

**Definition 3.2.3.** Let $(F, A)$ be a soft set over a neutrosophic ring $BN(R)$. Then $(F, A)$ is called soft neutrosophic biideal over $BN(R)$ if and only if $F(a)$ is a neutrosophic biideal of $BN(R)$ for all $a \in A$.

**Example 3.2.3.** Let $BN(\mathrm{R}) = (\mathrm{R}_1, *, \circ) \cup (\mathrm{R}_2, *, \circ)$ be a neutrosophic biring, where $(\mathrm{R}_1, *, \circ) = (\langle \mathbb{Z} \cup I \rangle, +, \times)$ and $(\mathrm{R}_2, *, \circ) = (\mathbb{C}, +, \times)$. Let $A = \{a_1, a_2\}$ be a set of parameters. Then clearly $(F, A)$ is a soft neutrosophic biideal over $BN(R)$, where

$$F(a_1) = \langle 2\mathbb{Z} \cup I \rangle \cup \mathbb{R}, F(a_2) = \langle 4\mathbb{Z} \cup I \rangle \cup \mathbb{Q}.$$

**Theorem 3.2.3.** Every soft neutrosophic biideal $(F, A)$ over a neutrosophic biring $BN(R)$ is trivially a soft neutrosophic biring but the converse may not be true.

The converse is left as an exercise for the readers.

**Proposition 3.2.4.** Let $(F,A)$ and $(K,B)$ be two soft neutosophic biideals over a neutrosophic biring $BN(R)$. Then

1. Their extended intersection $(F,A) \cap_E (K,B)$ is again a soft neutrosophic biideal over $BN(R)$.
2. Their restricted intersection $(F,A) \cap_R (K,B)$ is again a soft neutrosophic biideal over $BN(R)$.
3. Their $AND$ operation $(F,A) \vee (K,B)$ is again a soft neutrosophic biideal over $BN(R)$.

**Remark 3.2.4.** Let $(F,A)$ and $(K,B)$ be two soft neutosophic biideals over a neutrosophic biring $BN(R)$. Then

1. Their extended union $(F,A) \cup_E (K,B)$ is not a soft neutrosophic biideal over $BN(R)$.
2. Their restricted union $(F,A) \cup_R (K,B)$ is not a soft neutrosophic biideal over $BN(R)$.
3. Their $OR$ operation $(F,A) \vee (K,B)$ is not a soft neutrosophic biideal over $BN(R)$.

One can easily check these remarks by the help of examples.

**Definition 3.2.4.** Let $(F,A)$ and $(K,B)$ be two soft neutrosophic birings over $BN(R)$. Then $(K,B)$ is called soft neutrosophic subbiring of $(F,A)$, if

1. $B \subseteq A$, and
2. $K(a)$ is a neutrosophic subbiring of $F(a)$ for all $a \in A$.

**119**

**Theorem 3.2.4.** Every soft biring over a biring is a soft neutrosophic subbiring of a soft neutrosophic biring over the corresponding neutrosophic biring if $B \subseteq A$.

The proof is straightforward, so left as an exercise for the readers.

**Definition 3.2.5.** Let $(F, A)$ and $(K, B)$ be two soft neutrosophic birings over $BN(R)$. Then $(K, B)$ is called a soft neutrosophic biideal of $(F, A)$, if

1. $B \subseteq A$, and
2. $K(a)$ is a neutrosophic biideal of $F(a)$ for all $a \in A$.

**Proposition 3.2.5.** All soft neutrosophic biideals are trivially soft neutrosophic subbirings.

**Proof.** It is left as an exercise for the readers.

In the next section of this chapter, we finally extend soft neutrosophic ring to soft neutrosophic N-ring.

## 3.3   Soft Neutrosophic N-Ring

In this section of the current chapter, we introduce soft neutrosophic N-ring over a neutrosophic N-ring. This is an approximated collection of neutrosophic sub N-rings of the neutrosophic ring. We also give some basic properties and illustrate it with different examples.

We now procced to define soft neutrosophic N-ring over a neutrosophic N-ring.

**Definition 3.3.1.** Let $(N(\mathbb{R}), *_1, *_2, ..., *_N)$ be a neutrosophic N-ring and $(F, A)$ be a soft set over $N(R)$. Then $(F, A)$ is called soft neutrosophic N-ring if and only if $F(a)$ is a neutrosophic sub N-ring of $N(R)$ for all $a \in A$.

We illustrate it in the following example.

**Example 3.3.1.** Let $N(\mathbb{R}) = (\mathbb{R}_1, *, \circ) \cup (\mathbb{R}_2, *, \circ) \cup (\mathbb{R}_3, *, \circ)$ be a neutrosophic 3-ring, where $(\mathbb{R}_1, *, \circ) = (\langle \mathbb{Z} \cup I \rangle, +, \times)$, $(\mathbb{R}_2, *, \circ) = (\mathbb{C}, +, \times)$ and $(\mathbb{R}_3, *, \circ) = (\mathbb{R}, +, \times)$. Let $A = \{a_1, a_2, a_3, a_4\}$ be a set of parameters. Then clearly $(F, A)$ is a soft neutrosophic N-ring over $N(R)$, where

$$F(a_1) = \langle 2\mathbb{Z} \cup I \rangle \cup \mathbb{R} \cup \mathbb{Q}, F(a_2) = \langle 3\mathbb{Z} \cup I \rangle \cup \mathbb{Q} \cup \mathbb{Z},$$

$$F(a_3) = \langle 5\mathbb{Z} \cup I \rangle \cup \mathbb{Z} \cup 2\mathbb{Z}, F(a_4) = \langle 6\mathbb{Z} \cup I \rangle \cup 2\mathbb{Z} \cup \mathbb{R}.$$

The interested readers can easily construct more examples to understand this situation.

**Theorem 3.3.1.** Let $F, A$ and $(H, A)$ be two soft neutrosophic N-rings over $N(R)$. Then their intersection $F, A \cap H, A$ is again a soft neutrosophic N-ring over $N(R)$.

The proof is straightforward, so left as an exercise for the readers.

**Theorem 3.3.2.** Let $F, A$ and $H, B$ be two soft neutrosophic N-rings over $N(R)$. If $A \cap B = \phi$, then $F, A \cup H, B$ is a soft neutrosophic

**121**

N-ring over $N(R)$.

**Remark 3.3.1.** The extended union of two soft neutrosophic N-rings $F, A$ and $K, B$ over $BN(R)$ is not a soft neutrosophic ring over $N(R)$.

We can check this by the help of example.

**Remark 3.3.2.** The restricted union of two soft neutrosophic N-rings $F, A$ and $K, B$ over $N(R)$ is not a soft neutrosophic N-ring over $N(R)$.

One can easily verify this remark by the help of example.

**Remark 3.3.3.** The $OR$ operation of two soft neutrosophic N-rings over $N(R)$ may not be a soft neutrosophic N-ring over $N(R)$.

One can easily check these remarks with the help of Examples.

**Proposition 3.3.1.** The extended intersection of two soft neutrosophic N-rings over $N(R)$ is soft neutrosophic Nring over $N(R)$.

**Proof.** The proof is left as an exercise.

**Proposition 3.3.2.** The restricted intersection of two soft neutrosophic N-rings over $N(R)$ is soft neutrosophic N-ring over $N(R)$.

**Proof.** It is obvious.

**Proposition 3.3.3.** The *AND* operation of two soft neutrosophic N-rings over $N(R)$ is soft neutrosophic N-ring over $N(R)$.

**Definition 3.3.2.** Let $F, A$ be a soft set over a neutrosophic N-ring $N(R)$. Then $(F, A)$ is called an absolute soft neutrosophic N-ring if $F(a) = N(R)$ for all $a \in A$.

**Definition 3.3.3.** Let $(F, A)$ be a soft set over a neutrosophic N-ring $N(R)$. Then $(F, A)$ is called soft neutrosophic N-ideal over $N(R)$ if and only if $F(a)$ is a neutrosophic N-ideal of $N(R)$ for all $a \in A$.

**Example 3.3.1.** Let $N(R) = (R_1, *, \circ) \cup (R_2, *, \circ) \cup (R_3, *, \circ)$ be a neutrosophic 3-ring, where $(R_1, *, \circ) = (\langle \mathbb{Z} \cup I \rangle, +, \times)$, $(R_2, *, \circ) = (\mathbb{C}, +, \times)$ and $(R_3, *, \circ) = (\mathbb{R}, +, \times)$. Let $A = \{a_1, a_2\}$ be a set of parameters. Then clearly $(F, A)$ is a soft neutrosophic N-ideal over $N(R)$, where

$$F(a_1) = \langle 2\mathbb{Z} \cup I \rangle \cup 3\mathbb{Z} \cup \mathbb{Q}, F(a_2) = \langle 3\mathbb{Z} \cup I \rangle \cup \mathbb{Q} \cup 5\mathbb{Z}.$$

**Theorem 3.3.2.** Every soft neutrosophic N-ideal $(F, A)$ over a neutrosophic N-ring $N(R)$ is trivially a soft neutrosophic N-ring but the converse may not be true.

One can see the converse by the help of example.

**Proposition 3.3.4.** Let $(F, A)$ and $(K, B)$ be two soft neutosophic N-ideals over a neutrosophic N-ring $N(R)$. Then

1. Their extended intersection $(F,A) \cap_E (K,B)$ is again a soft neutrosophic N-ideal over $N(R)$.
2. Their restricted intersection $(F,A) \cap_R (K,B)$ is again a soft neutrosophic N-ideal over $N(R)$.
3. Their *AND* operation $(F,A) \vee (K,B)$ is again a soft neutrosophic N-ideal over $N(R)$.

**Remark 3.3.4:** Let $(F,A)$ and $(K,B)$ be two soft neutosophic N-ideals over a neutrosophic N-ring $N(R)$. Then

1. Their extended union $(F,A) \cup_E (K,B)$ is not a soft neutrosophic N-ideal over $N(R)$.
2. Their restricted union $(F,A) \cup_R (K,B)$ is not a soft neutrosophic N-ideal over $N(R)$.
3. Their *OR* operation $(F,A) \vee (K,B)$ is not a soft neutrosophic N-ideal over $N(R)$.

One can easily see these by the help of examples.

**Definition 3.3.4.** Let $(F,A)$ and $(K,B)$ be two soft neutrosophic N-rings over $N(R)$. Then $(K,B)$ is called soft neutrosophic sub N-ring of $(F,A)$, if

1. $B \subseteq A$, and
2. $K(a)$ is a neutrosophic sub N-ring of $F(a)$ for all $a \in A$.

**Example 3.32.** Let $N(R) = (R_1, *, \circ) \cup (R_2, *, \circ) \cup (R_3, *, \circ)$ be a neutrosophic 3-ring, where $(R_1, *, \circ) = (\langle \mathbb{Z} \cup I \rangle, +, \times)$, $(R_2, *, \circ) = (\mathbb{C}, +, \times)$ and $(R_3, *, \circ) = (\mathbb{R}, +, \times)$. Let $A = \{a_1, a_2, a_3, a_4\}$ be a set of parameters. Then clearly $(F, A)$ is a soft neutrosophic N-ring over $N(R)$, where

$$F(a_1) = \langle 2\mathbb{Z} \cup I \rangle \cup \mathbb{R} \cup \mathbb{Q}, F(a_2) = \langle 3\mathbb{Z} \cup I \rangle \cup \mathbb{Q} \cup \mathbb{Z},$$

$$F(a_3) = \langle 5\mathbb{Z} \cup I \rangle \cup \mathbb{Z} \cup 2\mathbb{Z}, F(a_4) = \langle 6\mathbb{Z} \cup I \rangle \cup 2\mathbb{Z} \cup \mathbb{R}.$$

Let $B = \{a_1, a_3\} \subset A$ and let $(H, B)$ be another soft neutrosophic N-ring over $N(R)$, where

$$H(a_1) = \langle 4\mathbb{Z} \cup I \rangle \cup \mathbb{Q} \cup \mathbb{Z},$$

$$H(a_3) = \langle 6\mathbb{Z} \cup I \rangle \cup \mathbb{Z} \cup 2\mathbb{Z}.$$

Then clearly $(H, B)$ is a soft neutrosophic sub N-ring of $(F, A)$ over $N(R)$.

**Theorem 3.3.3.** Every soft N-ring over a N-ring is a soft neutrosophic sub N-ring of a soft neutrosophic N-ring over the corresponding neutrosophic N-ring if $B \subseteq A$.

The proof is straightforward, so it is left as an exercise for the readers.

**Definition 3.3.5.** Let $(F, A)$ and $(K, B)$ be two soft neutrosophic N-rings over $N(R)$. Then $(K, B)$ is called a soft neutrosophic N-ideal of $(F, A)$, if

1. $B \subseteq A$, and
2. $K(a)$ is a neutrosophic N-ideal of $F(a)$ for all $a \in A$.

**Proposition 3.3.5.** All soft neutrosophic N-ideals are trivially soft neutrosophic sub N-rings.

**Proof.** It is straightforward.

In this section, we studies the important aspect neutrosophic fields and introduced soft neutrosophic fields and their generalization with its fundamental properties and examples.

### 3.4   Soft Neutrosophic Field

In this last section, we introduced soft neutrosophic fields over neutrosophic fields and also construct its generalization to soft neutrosophic N-fields over neutrosophic N-fields respectively. We also studied some fundamental properites of soft neutrosophic fields over neutrosophic fields with its generalization.

**Defintion 3.4.1.** Let $K(I) = \langle K \cup I \rangle$ be a neutrosophic field and let $(F, A)$ be a soft set over $K(I)$. Then $(F, A)$ is said to be soft neutrosophic field if and only if $F(a)$ is a neutrosophic subfield of $K(I)$ for all $a \in A$.

**Example 3.4.1.** Let $\langle C \cup I \rangle$ be a neutrosophic field of complex numbers. Let $A = \{a_1, a_2\}$ be a set of parameters and let $(F, A)$ be a soft set of $\langle C \cup I \rangle$. Then $(\text{F}, \text{A})$ is a soft neutrosophic field over $\langle C \cup I \rangle$, where

$$F(a_1) = \langle R \cup I \rangle, F(a_2) = \langle Q \cup I \rangle.$$

Where $\langle R \cup I \rangle$ and $\langle Q \cup I \rangle$ are the neutosophic fields of real numbers and rational numbers respectively.

**Proposition 3.4.1.** Every soft neutrosophic field is trivially a soft neutrosophic ring.

**Proof.** The proof is trivial.

**Remark 3.4.1.** The converse of the above proposition is not true.

To see the converse, lets take a look to the following example.

**Example 3.4.2.** Let $\langle Z \cup I \rangle$ be a neutrosophic ring of integers. Let $A = \{a_1, a_2, a_3, a_4\}$ be a set of parameters and let $(F, A)$ be a soft set over $\langle Z \cup I \rangle$. Then $(F, A)$ is a soft neutrosophic ring over $\langle Z \cup I \rangle$, where

$$F(a_1) = \langle 2Z \cup I \rangle, F(a_2) = \langle 3Z \cup I \rangle$$

$$F(a_3) = \langle 5Z \cup I \rangle, F(a_4) = \langle 6Z \cup I \rangle.$$

Clearly $(F, A)$ is not a soft neutrosophic field over $\langle Z \cup I \rangle$.

**Definition 3.4.2.** Let $(F, A)$ be a soft neutrosophic field over a neutrosophic field $\langle K \cup I \rangle$. Then $(F, A)$ is called an absolute soft neutrosophic field if $F(a) = \langle K \cup I \rangle$, for all $a \in A$.

**Soft Neutrosophic Bifield**

Here we give the deifiniton of soft neutrosophic bifield over a neutrosophic bifield and give some of their properties with illustrative examples.

**Defintion 3.4.3.** Let $BN(K)$ be a neutrosophic bifield and let $(F, A)$ be a soft set over $BN(K)$. Then $(F, A)$ is said to be soft neutrosophic bifield if and only if $F(a)$ is a neutrosophic subbifield of $BN(K)$ for all $a \in A$.

**Example 3.4.3.** Let $BN(K) = \langle \mathbb{C} \cup I \rangle \cup \mathbb{R}$ be a neutrosophic bifield of neutrosophic complex numbers and real numbers. Let $A = \{a_1, a_2\}$ be a set of parameters and let $(F, A)$ be a soft set of $BN(K)$. Then $(F, A)$ is a soft neutrosophic bifield over $BN(K)$, where

$$F(a_1) = \langle \mathbb{R} \cup I \rangle \cup \mathbb{Q}, F(a_2) = \langle \mathbb{Q} \cup I \rangle \cup \mathbb{Q}.$$

Where $\langle \mathbb{R} \cup I \rangle$ and $\langle \mathbb{Q} \cup I \rangle$ are the neutosophic fields of real numbers and rational numbers.

**Proposition 3.4.2.** Every soft neutrosophic bifield is trivially a soft neutrosophic biring.

**Proof.** The proof is trivial.

**Definition 3.4.4.** Let $(F, A)$ be a soft neutrosophic bifield over a neutrosophic bifield $BN(K)$. Then $(F, A)$ is called an absolute soft neutrosophic bifield if $F(a) = BN(K)$, for all $a \in A$.

Next we generalizes the concept of soft neutrosophic fields over neutrosophic fields to introduce soft neutrosophic N-fields over neutrosophic N-fields.

**Soft Neutrosophic N-field**

Here soft neutrosophic N-fields are introduced over neutrosophic N-fields. This is an approximated collection of neutrosophic sub N-fields of neutrosophic N-field.

**Defintion 3.4.5.** Let $N(K)$ be a neutrosophic N-field and let $(F, A)$ be a soft set over $N(K)$. Then $(F, A)$ is said to be a soft neutrosophic N-field if and only if $F(a)$ is a neutrosophic sub N-field of $N(K)$ for all $a \in A$.

**Proposition 3.4.3.** Every soft neutrosophic N-field is trivially a soft neutrosophic N-ring.

The proof is trivial, so left as an exercise for the readers.

**Definition 3.4.6.** Let $(F, A)$ be a soft neutrosophic N-field over a neutrosophic N-field $N(K)$. Then $(F, A)$ is called an absolute soft neutrosophic N-field if $F(a) = N(K)$, for all $a \in A$.

# Chapter No. 4

## Soft Neutrosophic Group Ring and Their Generalization

In this chapter, we introduce soft neutrosophic group ring over a neutrosophic group ring. We also extend it to its generalization and defined soft neutrosophic N-group N-ring over neutrosophic N-group N-ring. We give sufficient amount of illustrative examples and establish some of their basic and fundamental properties and characteristics.

### 4.1    Soft Neutrosophic Group Ring

In this section, the authors for the first time introduced soft neutrosophic group ring over a neutrosophic group ring which is basically a parameterized collection of subneutrosophic group ring of the neutrosophic group ring. We also give some core properties and give many examples to illustrate it.

We now define soft neutrosophic group ring over a neutrosophic group ring as follows.

**Definition 4.1.1.** Let $R\langle G \cup I\rangle$ be a neutrosophic group ring. Let $A$ be a set of parameters and $(F,A)$ be a soft set over $R\langle G \cup I\rangle$. Then $(F,A)$ is called soft neutosophic group ring if and only if $F(a)$ is subneutrosophic group ring of $R\langle G \cup I\rangle$ for all $a \in A$.

This situation can be explained in the following example.

**Example 4.1.1.** Let $\mathbb{Q}\langle G \cup I\rangle$ be a neutrosophic group ring, where $\mathbb{Q} =$ field of rationals and $\langle G \cup I\rangle = \{1, g, g^2, g^3, g^4, g^5, I, gI, ..., g^5 I : g^6 = 1, I^2 = I\}$. Let $\langle H_1 \cup I\rangle = \{1, g^3 : g^6 = 1\}$, $\langle H_2 \cup I\rangle = \{1, g^3, I, g^3 I : g^6 = 1, I^2 = I\}$, $\langle H_3 \cup I\rangle = \{1, g^2, g^4 : g^6 = 1, I^2 = I\}$ and $\langle H_4 \cup I\rangle = \{1, g^2, g^4, I, g^2 I, g^4 I : g^6 = 1, I^2 = I\}$.

Let $A = \{a_1, a_2, a_3, a_4\}$ be a set of parameters. Then $(F,A)$ is a soft neutrosophic group ring over $R\langle G \cup I\rangle$, where

$$F(a_1) = \mathbb{Q}\langle H_1 \cup I\rangle,$$

$$F(a_2) = \mathbb{Q}\langle H_2 \cup I\rangle,$$

$$F(a_3) = \mathbb{Q}\langle H_3 \cup I\rangle,$$

$$F(a_4) = \mathbb{Q}\langle H_4 \cup I\rangle.$$

**Theorem 4.1.1 .** Let $F,A$ and $(H,A)$ be two soft neutrosophic group rings over $R\langle G \cup I\rangle$. Then their intersection $F,A \cap H,A$ is again a soft neutrosophic group ring over $R\langle G \cup I\rangle$.

The proof is straightforward, so left as an exercise for the readers.

**Theorem 4.1.2.** Let $F, A$ and $H, B$ be two soft neutrosophic group rings over $R\langle G \cup I\rangle$. If $A \cap B = \phi$, then $F, A \cup H, B$ is a soft neutrosophic group ring over $R\langle G \cup I\rangle$.

**Theorem 4.1.3:** If $R\langle G \cup I\rangle$ is a commutative neutrosophic group ring. Then the soft neutrosophic group ring $(F, A)$ is also commutative soft neutrosophic group ring.

**Proof:** It is trivial.

**Proposition 4.1.1.** Let $(F, A)$ and $(K, B)$ be two soft neutosophic group rings over a neutrosophic grpup ring $R\langle G \cup I\rangle$. Then

1. Their extended intersection $(F, A) \cap_E (K, B)$ is again a soft neutrosophic group ring over $R\langle G \cup I\rangle$.
2. Their restricted intersection $(F, A) \cap_R (K, B)$ is again a soft neutrosophic group ring over $R\langle G \cup I\rangle$.
3. Their $AND$ operation $(F, A) \vee (K, B)$ is again a soft neutrosophic group ring over $R\langle G \cup I\rangle$.

**Proof.** This is straightforward.

**Remark 4.1.1.** Let $(F, A)$ and $(K, B)$ be two soft neutosophic group rings over a neutrosophic group ring $R\langle G \cup I\rangle$. Then

1. Their restricted union $(F, A) \cup_R (K, B)$ is not a soft neutrosophic group ring over $R\langle G \cup I\rangle$.

**132**

2. Their extended union $(F, A) \cup_E (K, B)$ is not a soft neutrosophic group ring over $R\langle G \cup I \rangle$.

3. Their *OR* operation $(F, A) \vee (K, B)$ is not a soft neutrosophic group ring over $R\langle G \cup I \rangle$.

To establish these remarks, see the following example.

**Example 4.1.2.** Let $\mathbb{Q}\langle G \cup I \rangle$ be a neutrosophic group ring, where $\mathbb{Q} =$ field of rationals and $\langle G \cup I \rangle = \{1, g, g^2, g^3, g^4, g^5, I, gI, ..., g^5 I : g^6 = 1, I^2 = I\}$. Let $A = \{a_1, a_2, a_3, a_4\}$ be a set of parameters. Then $(F, A)$ is a soft neutrosophic group ring over $R\langle G \cup I \rangle$, where

$$F(a_1) = \mathbb{Q}\langle H_1 \cup I \rangle,$$

$$F(a_2) = \mathbb{Q}\langle H_2 \cup I \rangle,$$

$$F(a_3) = \mathbb{Q}\langle H_3 \cup I \rangle,$$

$$F(a_4) = \mathbb{Q}\langle H_4 \cup I \rangle.$$

Let $B = \{a_1\}$ and let $(H, B)$ be another soft neutrosophic group ring over $\mathbb{Q}\langle G \cup I \rangle$, where

$$H(a_1) = \mathbb{Q}\langle H_4 \cup I \rangle.$$

Then $C = A \cap B = \{a_1\}$. Their restricted union is $(F, A) \cup_R (H, B) = (K, C)$, where $K(a_1) = F(a_1) \cup H(a_1) = \mathbb{Q}\langle H_1 \cup I \rangle \cup \mathbb{Q}\langle H_4 \cup I \rangle$. Clearly $K(a_1)$ is not a subneutrosophic group ring of $\mathbb{Q}\langle G \cup I \rangle$. Hence $(K,C)$ is not a soft neutrosophic group ring over $\mathbb{Q}\langle G \cup I \rangle$.

Similarly we can show $2$ and $3$ by the help of examples.

**Definition 4.1.2.** Let $R\langle G \cup I \rangle$ be a neutrosophic group ring. Let $A$ be a set of parameters and $(F,A)$ be a soft set over $R\langle G \cup I \rangle$. Then $(F,A)$ is called soft neutosophic group subring if and only if $F(a)$ is subneutrosophic group subring of $R\langle G \cup I \rangle$ for all $a \in A$.

**Example 4.1.3.** Let $\mathbb{C}\langle G \cup I \rangle$ be a neutrosophic group ring, where $\mathbb{C} =$ field of complex numbers and
$\langle G \cup I \rangle = \{1, g, g^2, g^3, g^4, g^5, I, gI, ..., g^5 I : g^6 = 1, I^2 = I\}$.

Let $A = \{a_1, a_2, a_3\}$ be a set of parameters. Then $(F,A)$ is a soft neutrosophic group subring over $\mathbb{C}\langle G \cup I \rangle$, where

$$F(a_1) = \mathbb{R}\langle G \cup I \rangle,$$

$$F(a_2) = \mathbb{Q}\langle G \cup I \rangle,$$

$$F(a_3) = \mathbb{Z}\langle G \cup I \rangle.$$

**Theorem 4.1.4 .** Let $F,A$ and $(H,A)$ be two soft neutrosophic group subrings over $R\langle G \cup I \rangle$. Then their intersection $F,A \cap H,A$ is again a soft neutrosophic group subrings over $R\langle G \cup I \rangle$.

**Proof.** The proof is straightforward.

**Theorem 4.1.5.** Let $F, A$ and $H, B$ be two soft neutrosophic group subrings over $R\langle G \cup I\rangle$. If $A \cap B = \phi$, then $F, A \cup H, B$ is a soft neutrosophic group subring over $R\langle G \cup I\rangle$.

This is straightforward, so left as an exercise for the interested readers.

**Remark 4.1.2.** Let $(F, A)$ and $(K, B)$ be two soft neutosophic group subrings over $R\langle G \cup I\rangle$. Then

1. Their extended union $(F, A) \cup_E (K, B)$ is not a soft neutrosophic group subring over $R\langle G \cup I\rangle$.
2. Their restricted union $(F, A) \cup_R (K, B)$ is not a soft neutrosophic group subring over $R\langle G \cup I\rangle$.
3. Their $OR$ operation $(F, A) \vee (K, B)$ is not a soft neutrosophic group subring over $R\langle G \cup I\rangle$.

One can easily show these remaks by the help of examples.

**Proposition 4.1.2.** Let $(F, A)$ and $(K, B)$ be two soft neutosophic group subrings over $R\langle G \cup I\rangle$. Then

1. Their restricted intersection $(F, A) \cap_R (K, B)$ is a soft neutrosophic group suring over $R\langle G \cup I\rangle$.
2. Their extended intersection $(F, A) \cap_E (K, B)$ is a soft neutrosophic group subring over $R\langle G \cup I\rangle$.
3. Their $AND$ operation $(F, A) \vee (K, B)$ is a soft neutrosophic group subring over $R\langle G \cup I\rangle$.

**135**

These are straightforward and the readers can easily attempt to prove it.

**Definition 4.1.3.** Let $R\langle G \cup I \rangle$ be a neutrosophic group ring. Let $A$ be a set of parameters and $(F, A)$ be a soft set over $R\langle G \cup I \rangle$. Then $(F, A)$ is called soft mixed neutosophic group ring if for some $a \in A$, $F(a)$ is subneutrosophic group subring and for the remaining $a \in A$, $F(a)$ is a subneutrosophic group ring of $R\langle G \cup I \rangle$.

**Example 4.1.4.** Let $\mathbb{C}\langle G \cup I \rangle$ be a neutrosophic group ring, where $\mathbb{C} = $ field of complex numbers and
$\langle G \cup I \rangle = \{1, g, g^2, g^3, g^4, g^5, I, gI, ..., g^5 I : g^6 = 1, I^2 = I\}$. Let $A = \{a_1, a_2, a_3, a_4, a_5, a_6, a_7\}$ be a set of parameters. Then $(F, A)$ is a soft mixed neutrosophic group ring over $\mathbb{C}\langle G \cup I \rangle$, where

$$F(a_1) = \mathbb{R}\langle G \cup I \rangle,$$

$$F(a_2) = \mathbb{Q}\langle G \cup I \rangle,$$

$$F(a_3) = \mathbb{Z}\langle G \cup I \rangle,$$

$$F(a_4) = \mathbb{C}\langle H_1 \cup I \rangle,$$

$$F(a_5) = \mathbb{C}\langle H_2 \cup I \rangle,$$

$$F(a_6) = \mathbb{C}\langle H_3 \cup I \rangle,$$

$$F(a_7) = \mathbb{C}\langle H_4 \cup I \rangle.$$

Where $\langle H_1 \cup I \rangle = \{1, g^3 : g^6 = 1\}$, $\langle H_2 \cup I \rangle = \{1, g^3, I, g^3 I : g^6 = 1, I^2 = I\}$,

$\langle H_3 \cup I \rangle = \{1, g^2, g^4 : g^6 = 1, I^2 = I\}$ and $\langle H_4 \cup I \rangle = \{1, g^2, g^4, I, g^2 I, g^4 I : g^6 = 1, I^2 = I\}$.

**Theorem 4.1.6.** Let $F, A$ and $(H, A)$ be two soft mixed neutrosophic group rings over $R\langle G \cup I \rangle$. Then their intersection $F, A \cap H, A$ is again a soft mixed neutrosophic group ring over $R\langle G \cup I \rangle$.

**Theorem 4.1.7.** Let $F, A$ and $H, B$ be two soft mixed neutrosophic group rings over $R\langle G \cup I \rangle$. If $A \cap B = \phi$, then $F, A \cup H, B$ is a soft mixed neutrosophic group ring over $R\langle G \cup I \rangle$.

**Proof.** This is straightforward.

**Remark 4.1.3.** Let $R\langle G \cup I \rangle$ be a neutrosophic group ring. Then $R\langle G \cup I \rangle$ can have soft neutrosophic group ring, soft neutrosophic group subring and soft mixed neutrosophic group ring over $R\langle G \cup I \rangle$.

**Proposition 4.1.3.** Let $(F, A)$ and $(K, B)$ be two soft mixed neutosophic group rings over $R\langle G \cup I \rangle$. Then

1. Their extended intersection $(F, A) \cap_E (K, B)$ is a soft mixed neutrosophic group ring over $R\langle G \cup I \rangle$.
2. Their restricted intersection $(F, A) \cap_R (K, B)$ is a soft mixed neutrosophic group ring over $R\langle G \cup I \rangle$.
3. Their $AND$ operation $(F, A) \vee (K, B)$ is a soft mixed neutrosophic group ring over $R\langle G \cup I \rangle$.

**Remark 4.1.4.** Let $(F,A)$ and $(K,B)$ be two soft mixed neutosophic group rings over $R\langle G \cup I\rangle$. Then

1. Their restricted union $(F,A) \cup_R (K,B)$ is not a soft mixed neutrosophic group ring over $R\langle G \cup I\rangle$.
2. Their extended union $(F,A) \cup_E (K,B)$ is not a soft mixed neutrosophic group ring over $R\langle G \cup I\rangle$.
3. Their $OR$ operation $(F,A) \vee (K,B)$ is a soft mixed neutrosophic group ring over $R\langle G \cup I\rangle$.

One can easily show these remarks by the help of examples.

**Definition 4.1.4.** Let $R\langle G \cup I\rangle$ be a neutrosophic group ring. Let $A$ be a set of parameters and $(F,A)$ be a soft set over $R\langle G \cup I\rangle$. Then $(F,A)$ is called soft neutosophic subring if and only if $F(a)$ is neutrosophic subring of $R\langle G \cup I\rangle$ for all $a \in A$.

This situation can be explained in the following example.

**Example 4.1.5.** Let $\mathbb{C}\langle G \cup I\rangle$ be a neutrosophic group ring, where $\mathbb{C} = $ field of complex numbers and
$\langle G \cup I\rangle = \{1, g, g^2, g^3, g^4, g^5, I, gI, ..., g^5I : g^6 = 1, I^2 = I\}$.

Let $A = \{a_1, a_2, a_3\}$ be a set of parameters. Then $(F,A)$ is a soft neutrosophic subring over $\mathbb{C}\langle G \cup I\rangle$, where

$$F(a_1) = \langle \mathbb{R} \cup I\rangle,$$

$$F(a_2) = \langle \mathbb{Q} \cup I \rangle,$$

$$F(a_3) = \langle \mathbb{Z} \cup I \rangle.$$

**Theorem 4.1.8.** Let $F, A$ and $(H,A)$ be two soft neutrosophic subrings over $R \langle G \cup I \rangle$. Then their intersection $F, A \cap H, A$ is again a soft neutrosophic subring over $R \langle G \cup I \rangle$.

The proof is straightforward, so it is left as an exercise for the readers.

**Theorem 4.1.9.** Let $F, A$ and $H, B$ be two soft neutrosophic subrings over $R \langle G \cup I \rangle$. If $A \cap B = \phi$, then $F, A \cup H, B$ is also a soft neutrosophic subring over $R \langle G \cup I \rangle$.

**Proposition 4.1.4.** Let $(F, A)$ and $(K, B)$ be two soft neutosophic subrings over $R \langle G \cup I \rangle$. Then

1. Their extended intersection $(F, A) \cap_E (K, B)$ is a soft neutrosophic subring over $R \langle G \cup I \rangle$.
2. Their restricted intersection $(F, A) \cap_R (K, B)$ is a soft neutrosophic subring over $R \langle G \cup I \rangle$.
3. Their $AND$ operation $(F, A) \vee (K, B)$ is a soft neutrosophic subring over $R \langle G \cup I \rangle$.

These are straightforward, so left as an exercise for the interested readers.

**139**

**Remark 4.1.5.** Let $(F,A)$ and $(K,B)$ be two soft neutosophic subrings over $R\langle G \cup I\rangle$. Then

1. Their restricted union $(F,A) \cup_R (K,B)$ is not a soft neutrosophic subring over $R\langle G \cup I\rangle$.
2. Their extended union $(F,A) \cup_E (K,B)$ is not a soft neutrosophic subring over $R\langle G \cup I\rangle$.
3. Their $OR$ operation $(F,A) \vee (K,B)$ is a soft neutrosophic subring over $R\langle G \cup I\rangle$.

One can easily show these remarks by the help of examples.

**Definition 4.1.5.** Let $R\langle G \cup I\rangle$ be a neutrosophic group ring. Let $A$ be a set of parameters and $(F,A)$ be a soft set over $R\langle G \cup I\rangle$. Then $(F,A)$ is called soft pseudo neutosophic subring if and only if $F(a)$ is a pseudo neutrosophic subring of $R\langle G \cup I\rangle$ for all $a \in A$.

**Example 4.1.6.** Let $\mathbb{Z}_6\langle G \cup I\rangle$ be a neutrosophic group ring of the neutrosophic group $\langle G \cup I\rangle$ over $\mathbb{Z}_6$. Let $A = \{a_1, a_2\}$ be a set of parameters. Then $(F,A)$ is a soft pseudo neutrosophic subring over $\mathbb{Z}_6\langle G \cup I\rangle$, where

$$F(a_1) = \{0, 3I\}, F(a_2) = \{0, 2I, 4I\}.$$

**Theorem 4.1.10.** Let $F,A$ and $(H,A)$ be two soft pseudo neutrosophic subrings over $R\langle G \cup I\rangle$. Then their intersection $F,A \cap H,A$ is again a soft pseudo neutrosophic subring over $R\langle G \cup I\rangle$.

**Proof.** The proof is straightforward.

**Theorem 4.1.11.** Let $F, A$ and $H, B$ be two soft pseudo neutrosophic subrings over $R\langle G \cup I \rangle$. If $A \cap B = \phi$, then $F, A \cup H, B$ is a soft pseudo neutrosophic subring over $R\langle G \cup I \rangle$.

**Proposition 4.1.5.** Let $(F, A)$ and $(K, B)$ be two soft pseudo neutosophic subrings over $R\langle G \cup I \rangle$. Then

1. Their extended intersection $(F, A) \cap_E (K, B)$ is a soft pseudo neutrosophic subring over $R\langle G \cup I \rangle$.
2. Their restricted intersection $(F, A) \cap_R (K, B)$ is a soft pseudo neutrosophic subring over $R\langle G \cup I \rangle$.
3. Their $AND$ operation $(F, A) \vee (K, B)$ is a soft pseudo neutrosophic subring over $R\langle G \cup I \rangle$.

**Proof.** These are straightforward.

**Remark 4.1.6.** Let $(F, A)$ and $(K, B)$ be two soft pseudo neutosophic subrings over $R\langle G \cup I \rangle$. Then

1. Their restricted union $(F, A) \cup_R (K, B)$ is not a soft pseudo neutrosophic subring over $R\langle G \cup I \rangle$.
2. Their extended union $(F, A) \cup_E (K, B)$ is not a soft pseudo neutrosophic subring over $R\langle G \cup I \rangle$.
3. Their $OR$ operation $(F, A) \vee (K, B)$ is a soft pseudo neutrosophic subring over $R\langle G \cup I \rangle$.

**141**

One can easily show these remarks by the help of examples.

**Definition 4.1.6.** Let $R\langle G \cup I\rangle$ be a neutrosophic group ring. Let $A$ be a set of parameters and $(F, A)$ be a soft set over $R\langle G \cup I\rangle$. Then $(F, A)$ is called soft neutosophic subgroup ring if and only if $F(a)$ is neutrosophic subgroup ring of $R\langle G \cup I\rangle$ for all $a \in A$ .

**Example 4.1.7.** Let $\mathbb{C}\langle G \cup I\rangle$ be a neutrosophic group ring, where $\mathbb{C} =$ field of complex numbers and $\langle G \cup I\rangle = \{1, g, g^2, g^3, g^4, g^5, I, gI, ..., g^5 I : g^6 = 1, I^2 = I\}$. Let $H_1 = \{1, g^2, g^4 : g^6 = 1\}$ and $H_2 = \{1, g^3 : g^6 = 1\}$ be subgroups of the neutrosophic group $\langle G \cup I\rangle$.

Let $A = \{a_1, a_2\}$ be a set of parameters. Then $(F, A)$ is a soft neutrosophic subgroup ring over $\mathbb{C}\langle G \cup I\rangle$, where

$$F(a_1) = \mathbb{R}H_1, \ F(a_2) = \mathbb{Q}H_2 .$$

**Theorem 4.1.12.** Let $F, A$ and $(H, A)$ be two soft neutrosophic subgroup rings over $R\langle G \cup I\rangle$. Then their intersection $F, A \cap H, A$ is again a soft neutrosophic subgroup ring over $R\langle G \cup I\rangle$ .

**Proof.** The proof is straightforward.

**Theorem 4.1.13.** Let $F, A$ and $H, B$ be two soft neutrosophic subgroup rings over $R\langle G \cup I\rangle$. If $A \cap B = \phi$, then $F, A \cup H, B$ is a soft neutrosophic subgroup ring over $R\langle G \cup I\rangle$.

**Proposition 4.1.6.** Let $(F, A)$ and $(K, B)$ be two soft neutosophic subgroup rings over $R\langle G \cup I \rangle$. Then

1. Their extended intersection $(F, A) \cap_E (K, B)$ is a soft neutrosophic subgroup ring over $R\langle G \cup I \rangle$.
2. Their restricted intersection $(F, A) \cap_R (K, B)$ is a soft neutrosophic subgroup ring over $R\langle G \cup I \rangle$.
3. Their $AND$ operation $(F, A) \vee (K, B)$ is a soft neutrosophic subgroup ring over $R\langle G \cup I \rangle$.

This is straightforward, so left as an exercise for the readers.

**Remark 4.1.7.** Let $(F, A)$ and $(K, B)$ be two soft neutosophic subgroup rings over $R\langle G \cup I \rangle$. Then

1. Their restricted union $(F, A) \cup_R (K, B)$ is not a soft neutrosophic subgroup ring over $R\langle G \cup I \rangle$.
2. Their extended union $(F, A) \cup_E (K, B)$ is not a soft neutrosophic subgroup ring over $R\langle G \cup I \rangle$.
3. Their $OR$ operation $(F, A) \vee (K, B)$ is a soft neutrosophic subgroup ring over $R\langle G \cup I \rangle$.

One can easily show these remarks by the help of examples.

**Definition 4.1.7.** Let $R\langle G \cup I \rangle$ be a neutrosophic group ring. Let $A$ be a set of parameters and $(F,A)$ be a soft set over $R\langle G \cup I \rangle$. Then $(F,A)$ is called a soft subring over $R\langle G \cup I \rangle$ if and only if $F(a)$ is a subring of $R\langle G \cup I \rangle$ for all $a \in A$.

**Example 4.1.8.** Let $\mathbb{Z}_2\langle G \cup I \rangle$ be a neutrosophic group ring, where $\langle G \cup I \rangle = \{1, g, g^2, g^3, I, gI, g^2I, g^3I : g^4 = 1, I^2 = I\}$. Let $A = \{a_1, a_2\}$ be a set of parameters. Then $(F,A)$ is a soft subring over $\mathbb{Z}_2\langle G \cup I \rangle$, where

$$F(a_1) = \{0, 1 + g^2\},$$

$$F(a_2) = \{0, 1 + g, g + g^3, 1 + g^3\}.$$

**Theorem 4.1.14.** Let $F,A$ and $(H,A)$ be two soft subrings over $R\langle G \cup I \rangle$. Then their intersection $F,A \cap H,A$ is again a soft subring over $R\langle G \cup I \rangle$.

**Proof.** The proof is straightforward.

**Theorem 4.1.15.** Let $F,A$ and $H,B$ be two soft subrings over $R\langle G \cup I \rangle$. If $A \cap B = \phi$, then $F,A \cup H,B$ is a soft subring over $R\langle G \cup I \rangle$.

**Proof.** This is straightforward.

**Proposition 4.1.7.** Let $(F,A)$ and $(K,B)$ be two soft subrings over $R\langle G \cup I \rangle$. Then

**144**

1. Their extended intersection $(F,A) \cap_E (K,B)$ is a soft subring over $R\langle G \cup I \rangle$.

2. Their restricted intersection $(F,A) \cap_R (K,B)$ is a soft subring over $R\langle G \cup I \rangle$.

3. Their *AND* operation $(F,A) \vee (K,B)$ is a soft subring over $R\langle G \cup I \rangle$.

**Proof.** These are straightforward.

**Remark 4.1.8.** Let $(F,A)$ and $(K,B)$ be two soft subrings over $R\langle G \cup I \rangle$. Then

1. Their restricted union $(F,A) \cup_R (K,B)$ is not a soft subring over $R\langle G \cup I \rangle$.

2. Their extended union $(F,A) \cup_E (K,B)$ is not a soft subring over $R\langle G \cup I \rangle$.

3. Their *OR* operation $(F,A) \vee (K,B)$ is a soft subring over $R\langle G \cup I \rangle$.

One can easily show these remarks by the help of examples.

**Definition 4.1.8.** Let $R\langle G \cup I \rangle$ be a neutrosophic group ring. Then $(F,A)$ is called soft absolute neutosophic group subring if $F(a) = R\langle G \cup I \rangle$, for all $a \in A$.

**Definition 4.1.9.** Let $(F,A)$ and $(K,B)$ be two soft neutrosophic group rings over $R\langle G \cup I \rangle$. Then $(K,B)$ is called soft subneutrosophic group ring of $(F,A)$, if

1. $B \subseteq A$, and
2. $K(a)$ is a subneutrosophic group ring of $F(a)$ for all $a \in B$.

**Example 4.1.9.** Let $\mathbb{C}\langle G \cup I \rangle$ be a neutrosophic group ring, where $\mathbb{C} =$ field of complex numbers and

$$\langle G \cup I \rangle = \{1, g, g^2, g^3, g^4, g^5, I, gI, \dots, g^5 I : g^6 = 1, I^2 = I\}.$$

Let $A = \{a_1, a_2, a_3, a_4, a_5, a_6, a_7\}$ be a set of parameters. Then $(F, A)$ is a soft neutrosophic group ring over $\mathbb{C}\langle G \cup I \rangle$, where

$$F(a_1) = \mathbb{R}\langle G \cup I \rangle,$$

$$F(a_2) = \mathbb{Q}\langle G \cup I \rangle,$$

$$F(a_3) = \mathbb{Z}\langle G \cup I \rangle,$$

$$F(a_4) = \mathbb{C}\langle H_1 \cup I \rangle,$$

$$F(a_5) = \mathbb{C}\langle H_2 \cup I \rangle,$$

$$F(a_6) = \mathbb{C}\langle H_3 \cup I \rangle,$$

$$F(a_7) = \mathbb{C}\langle H_4 \cup I \rangle.$$

Where $\langle H_1 \cup I \rangle = \{1, g^3 : g^6 = 1\}$, $\langle H_2 \cup I \rangle = \{1, g^3, I, g^3 I : g^6 = 1, I^2 = I\}$, $\langle H_3 \cup I \rangle = \{1, g^2, g^4 : g^6 = 1, I^2 = I\}$ and $\langle H_4 \cup I \rangle = \{1, g^2, g^4, I, g^2 I, g^4 I : g^6 = 1, I^2 = I\}$.

Let $B = \{a_1, a_2, a_5, a_7\}$ be a set of parameters. Then $(H, B)$ is a soft subneutrosophic group ring of $(F, A)$, where

$$\text{H}(a_1) = \mathbb{Q}\langle G \cup I \rangle,$$

$$\text{H}(a_2) = \mathbb{Z}\langle G \cup I \rangle,$$

$$\text{H}(a_5) = \mathbb{C}\langle H_2 \cup I \rangle,$$

$$\text{H}(a_7) = \mathbb{C}\langle H_4 \cup I \rangle.$$

**Definition 4.1.10.** Let $R\langle G \cup I \rangle$ be a neutrosophic group ring. Let $A$ be a set of parameters and $(F,A)$ be a soft set over $R\langle G \cup I \rangle$. Then $(F,A)$ is called soft neutosophic ideal if and only if $F(a)$ is a neutrosophic ideal of $R\langle G \cup I \rangle$ for all $a \in A$.

Similarly one can define soft neutrosophic left ideal and soft neutrosophic right ideal over $R\langle G \cup I \rangle$.

**Example 4.1.10.** Let $R\langle G \cup I \rangle$ be a neutrosophic group ring, where $R = \mathbb{Z}$. Let $A = \{a_1, a_2, a_3\}$ be a set of parameters. Then $(F,A)$ is a soft neutrosophic ideal over $R\langle G \cup I \rangle$, where

$$F(a_1) = 2\mathbb{Z}\langle G \cup I \rangle,$$

$$F(a_2) = 4\mathbb{Z}\langle G \cup I \rangle,$$

$$F(a_3) = 6\mathbb{Z}\langle G \cup I \rangle.$$

**Theorem 4.1.16.** All soft neutrosophic ideals are trivially soft neutrosophic group rings but the converse is not true in general.

We can easily establish the converse by the help of example.

**Proposition 4.1.8.** Let $(F, A)$ and $(K, B)$ be two soft neutosophic ideals over $R \langle G \cup I \rangle$. Then

1. Their extended intersection $(F, A) \cap_E (K, B)$ is a soft neutrosophic ideal over $R \langle G \cup I \rangle$.
2. Their restricted intersection $(F, A) \cap_R (K, B)$ is a soft neutrosophic ideal over $R \langle G \cup I \rangle$.
3. Their *AND* operation $(F, A) \vee (K, B)$ is a soft neutrosophic ideal over $R \langle G \cup I \rangle$.

**Proof.** This is straightforward.

**Remark 4.1.9.** Let $(F, A)$ and $(K, B)$ be two soft neutosophic ideal over $R \langle G \cup I \rangle$. Then

1. Their restricted union $(F, A) \cup_R (K, B)$ is not a soft neutrosophic ideal over $R \langle G \cup I \rangle$.
2. Their extended union $(F, A) \cup_E (K, B)$ is not a soft neutrosophic ideal over $R \langle G \cup I \rangle$.
3. Their *OR* operation $(F, A) \vee (K, B)$ is a soft neutrosophic ideal over $R \langle G \cup I \rangle$.

One can easily show these remarks by the help of examples.

**Definition 4.1.11.** Let $R \langle G \cup I \rangle$ be a neutrosophic group ring. Let $A$ be a set of parameters and $(F, A)$ be a soft set over $R \langle G \cup I \rangle$. Then $(F, A)$ is called soft pseudo neutosophic ideal if and only if $F(a)$ is a pseudo neutrosophic ideal of $R\langle G \cup I \rangle$ for all $a \in A$.

**Example 4.1.11.** Let $R\langle G \cup I \rangle$ be a neutrosophic group ring, where $R = Z_2 = \{0,1\}$ and $\langle G \cup I \rangle = \{1, g, g^2, g^3, I, gI, g^2I, g^3I : g^4 = 1, I^2 = I\}$. Let $A = \{a_1, a_2, a_3\}$ be a set of parameters. Then $(F, A)$ is a soft pseudo neutrosophic ideal over $R\langle G \cup I \rangle$, where

$$F(a_1) = \left\langle 0, (1 + g + g^2 + g^3) \right\rangle,$$

$$F(a_2) = \left\langle I + gI + g^2I + g^3I \right\rangle,$$

$$F(a_3) = \left\langle 1 + g + g^2 + g^3 + I + gI + g^2I + g^3I \right\rangle.$$

**Theorem 4.1.17.** All soft pseudo neutrosophic ideals are trivially soft neutrosophic group rings but the converse is not true in general.

We can easily establish the converse by the help of example.

**Theorem 4.1.18.** All soft pseudo neutrosophic ideals are trivially soft neutrosophic ideals but the converse is not true in general.

We can easily establish the converse by the help of example.

**Proposition 4.1.9.** Let $(F, A)$ and $(K, B)$ be two soft pseudo neutosophic ideals over $R\langle G \cup I \rangle$. Then

1. Their extended intersection $(F, A) \cap_E (K, B)$ is a soft pseudo neutrosophic ideals over $R\langle G \cup I \rangle$.
2. Their restricted intersection $(F, A) \cap_R (K, B)$ is a soft pseudo

**149**

neutrosophic ideals over $R\langle G \cup I \rangle$.

3. Their *AND* operation $(F,A) \vee (K,B)$ is a soft pseudo neutrosophic ideals over $R\langle G \cup I \rangle$.

**Proof.** This is straightforward.

**Remark 4.1.10.** Let $(F,A)$ and $(K,B)$ be two soft pseudo neutosophic ideals over $R\langle G \cup I \rangle$. Then

1. Their restricted union $(F,A) \cup_R (K,B)$ is not a soft pseudo neutrosophic ideals over $R\langle G \cup I \rangle$.

2. Their extended union $(F,A) \cup_E (K,B)$ is not a soft pseudo neutrosophic ideals over $R\langle G \cup I \rangle$.

3. Their *OR* operation $(F,A) \vee (K,B)$ is a soft pseudo neutrosophic ideals over $R\langle G \cup I \rangle$.

One can easily show these remarks by the help of examples.

**Definition 4.1.12.** Let $(F,A)$ and $(K,B)$ be two soft neutrosophic group rings over $R\langle G \cup I \rangle$. Then $(K,B)$ is called soft neutrosophic ideal of $(F,A)$, if

1. $B \subseteq A$, and
2. $K(a)$ is a neutrosophic ideal of $F(a)$ for all $a \in B$.

**Example 4.1.12.** Let $R\langle G \cup I\rangle$ be a neutrosophic group ring, where $R = \mathbb{Z}$. Let $A = \{a_1, a_2, a_3\}$ be a set of parameters. Then $(F, A)$ is a soft neutrosophic group ring over $R\langle G \cup I\rangle$, where

$$F(a_1) = 2\mathbb{Z}\langle G \cup I\rangle,$$

$$F(a_2) = 4\mathbb{Z}\langle G \cup I\rangle,$$

$$F(a_3) = 6\mathbb{Z}\langle G \cup I\rangle.$$

Let $B = \{a_1, a_2\}$ be another set of parameters. Then clearly $(H, B)$ is a soft neutrosophic ideal of $t(F, A)$, where

$$H(a_1) = 8\mathbb{Z}\langle G \cup I\rangle,$$

$$H(a_2) = 12\mathbb{Z}\langle G \cup I\rangle.$$

**Theorem 4.1.19.** Every soft neutrosophic ideal of the soft neutrosophic group ring over a neutrosophic group ring is trivially a soft subneutrosophic group ring.

In the next section, we define soft neutrosophic bigroup biring over a neutrosophic bigroup biring and give some fundamental results and illustrated examples.

## 4.2   Soft Neutrosophic Bigroup Biring

In this section, the authors introduced soft neutrosophic bigroup biring over a neutrosophic bigroup biring. Basically this is the extension of soft neutrosophic group ring over neutrosophic group ring. We also studied some basic and fundamental properties with sufficient amount of examples.

**Definition 4.2.1.** Let $R_B \langle G \cup I \rangle = R_1 \langle G_1 \cup I \rangle \cup R_2 \langle G_2 \cup I \rangle$ be a neutrosophic bigroup biring over a biring $R_B = R_1 \cup R_2$. Let $A$ be a set of parameters and let $(F, A)$ be a soft set over $R_B \langle G \cup I \rangle$. Then $(F, A)$ is called a soft neutrosophic bigroup biring over $R_B \langle G \cup I \rangle$ if and only if $F(a)$ is a subneutrosophic bigroup biring of $R_B \langle G \cup I \rangle$ for all $a \in A$.

This fact can be explained in the following example.

**Example 4.2.1.** Let $R_B \langle G \cup I \rangle = \mathbb{Q} \langle G_1 \cup I \rangle \cup \mathbb{R} \langle G_2 \cup I \rangle$ be a neutrosophic bigroup biring, where $R_B = \mathbb{Q} \cup \mathbb{R}$ and

$\langle G_1 \cup I \rangle = \{1, g, g^2, g^3, g^4, g^5, I, gI, \ldots, g^5 I : g^6 = 1, I^2 = I\}$,

$\langle G_2 \cup I \rangle = \{1, g, g^2, g^3, I, gI, g^2 I, g^3 I : g^4 = 1, I^2 = I\}$.

Let $\langle H_1^{'} \cup I \rangle \cup \langle H_2^{'} \cup I \rangle = \{1, g^3 : g^6 = 1\} \cup \{1, g^2 : g^4 = 1\}$,

$\langle H_1^{''} \cup I \rangle \cup \langle H_2^{''} \cup I \rangle = \{1, g^3, I, g^3 I : g^6 = 1, I^2 = I\} \cup \{1, g^2, I, g^2 I\}$,

$$\left\langle H^{\sim}_1 \cup I\right\rangle \cup \left\langle H^{\sim}_2 \cup I\right\rangle = \{1, g^2, g^4 : g^6 = 1, I^2 = I\} \cup \{1, I, g^2 I : g^4 = 1, I^2 = I\} \text{ and}$$

$$\left\langle H^{\sim}_1 \cup I\right\rangle \cup \left\langle H^{\sim}_2 \cup I\right\rangle = \{1, g^2, g^4, I, g^2 I, g^4 I : g^6 = 1, I^2 = I\} \cup \{1, g^2 : g^4 = 1\}.$$

Let $A = \{a_1, a_2, a_3, a_4\}$ be a set of parameters. Then $(F, A)$ is a soft neutrosophic bigroup biring over $R_B \left\langle G \cup I\right\rangle$, where

$$F(a_1) = \mathbb{Q}\left\langle H^{'}_1 \cup I\right\rangle \cup \mathbb{R}\left\langle H^{''}_2 \cup I\right\rangle,$$

$$F(a_2) = \mathbb{Q}\left\langle H^{''}_1 \cup I\right\rangle \cup \mathbb{R}\left\langle H^{'}_2 \cup I\right\rangle,$$

$$F(a_3) = \mathbb{Q}\left\langle H^{\sim}_1 \cup I\right\rangle \cup \mathbb{R}\left\langle H^{''}_2 \cup I\right\rangle,$$

$$F(a_4) = \mathbb{Q}\left\langle H^{\sim}_1 \cup I\right\rangle \cup \mathbb{R}\left\langle H^{\sim}_2 \cup I\right\rangle.$$

We now study some properties and characterization of soft neutrosophic bigroup birings.

**Theorem 4.2.1.** Let $F, A$ and $(H, A)$ be two soft neutrosophic bigroup birings over $R_B \left\langle G \cup I\right\rangle$. Then their intersection $F, A \cap H, A$ is again a soft neutrosophic bigroup biring over $R_B \left\langle G \cup I\right\rangle$.

**Theorem 4.2.2.** Let $F, A$ and $H, B$ be two soft neutrosophic bigroup birings over $R_B \left\langle G \cup I\right\rangle$. If $A \cap B = \phi$, then $F, A \cup H, B$ is a soft neutrosophic bigroup biring over $R_B \left\langle G \cup I\right\rangle$.

This is straightforward, so left for the readers as an exercise.

**Proposition 4.2.1.** Let $(F, A)$ and $(K, B)$ be two soft neutosophic bigroup

birings over a neutrosophic bigroup biring $R_B \langle G \cup I \rangle$. Then

1. Their extended intersection $(F,A) \cap_E (K,B)$ is again a soft neutrosophic bigroup biring over $R_B \langle G \cup I \rangle$.
2. Their restricted intersection $(F,A) \cap_R (K,B)$ is again a soft neutrosophic bigroup biring over $R_B \langle G \cup I \rangle$.
3. Their $AND$ operation $(F,A) \vee (K,B)$ is again a soft neutrosophic bigroup biring over $R_B \langle G \cup I \rangle$.

These are straightforward.

**Remark 4.2.1.** Let $(F,A)$ and $(K,B)$ be two soft neutosophic bigroup birings over $R_B \langle G \cup I \rangle$. Then

1. Their restricted union $(F,A) \cup_R (K,B)$ is not a soft neutrosophic bigroup biring over $R_B \langle G \cup I \rangle$.
2. Their extended union $(F,A) \cup_E (K,B)$ is not a soft neutrosophic bigroup biring over $R_B \langle G \cup I \rangle$.
3. Their $OR$ operation $(F,A) \vee (K,B)$ is not a soft neutrosophic bigroup biring over $R_B \langle G \cup I \rangle$.

One can establish these remarks by the help of following examples.

**Definition 4.2.2.** Let $R_B \langle G \cup I \rangle = R_1 \langle G_1 \cup I \rangle \cup R_2 \langle G_2 \cup I \rangle$ be a neutrosophic bigroup biring. Let $A$ be a set of parameters and $(F,A)$ be a soft set over $R_B \langle G \cup I \rangle$. Then $(F, A)$ is called soft neutosophic bigroup subbiring over $R_B \langle G \cup I \rangle$ if and only if $F(a)$ is subneutrosophic bigroup subbiring of $R_B \langle G \cup I \rangle$ for all $a \in A$.

**Example 4.2.2.** Let $R_B \langle G \cup I \rangle = \mathbb{R} \langle G_1 \cup I \rangle \cup \mathbb{C} \langle G_2 \cup I \rangle$ be a neutrosophic bigroup biring, where $\langle G_1 \cup I \rangle = \{1, g, g^2, g^3, g^4, g^5, I, gI, ..., g^5 I : g^6 = 1, I^2 = I\}$ and $\langle G_2 \cup I \rangle = \{1, g, g^2, g^3, I, gI, g^2 I, g^3 I : g^4 = 1, I^2 = I\}$. Let $A = \{a_1, a_2\}$ be a set of parameters. Then $(F, A)$ is a soft neutrosophic bigroup subbiring over $R_B \langle G \cup I \rangle$, where

$$F(a_1) = \mathbb{Q} \langle G_1 \cup I \rangle \cup \mathbb{R} \langle G_2 \cup I \rangle,$$

$$F(a_2) = \mathbb{Z} \langle G_1 \cup I \rangle \cup \mathbb{Q} \langle G_2 \cup I \rangle.$$

**Theorem 4.2.3.** Let $F, A$ and $(H, A)$ be two soft neutrosophic bigroup subbirings over $R_B \langle G \cup I \rangle$. Then their intersection $F, A \cap H, A$ is again a soft neutrosophic bigroup subbirings over $R_B \langle G \cup I \rangle$.

**Proof.** The proof is straightforward.

**Theorem 4.2.4.** Let $F, A$ and $H, B$ be two soft neutrosophic bigroup subbirings over $R_B \langle G \cup I \rangle$. If $A \cap B = \phi$, then $F, A \cup H, B$ is a soft neutrosophic bigroup subbiring over $R_B \langle G \cup I \rangle$.

**Proof.** This is straightforward.

**Remark 4.2.2.** Let $(F,A)$ and $(K,B)$ be two soft neutosophic bigroup subbirings over $R_B \langle G \cup I \rangle$. Then

1. Their extended union $(F,A) \cup_E (K,B)$ is not a soft neutrosophic bigroup subbiring over $R_B \langle G \cup I \rangle$.
2. Their restricted union $(F,A) \cup_R (K,B)$ is not a soft neutrosophic bigroup subbiring over $R_B \langle G \cup I \rangle$.
3. Their $OR$ operation $(F,A) \vee (K,B)$ is not a soft neutrosophic bigroup subbiring over $R_B \langle G \cup I \rangle$.

One can easily show these remaks by the help of examples.

**Proposition 4.2.2.** Let $(F,A)$ and $(K,B)$ be two soft neutrosophic bigroup subbirings over $R_B \langle G \cup I \rangle$. Then

1. Their restricted intersection $(F,A) \cap_R (K,B)$ is a soft neutrosophic bigroup subbring over $R_b \langle G \cup I \rangle$.
2. Their extended intersection $(F,A) \cap_E (K,B)$ is a soft neutrosophic bigroup subbiring over $R_B \langle G \cup I \rangle$.
3. Their $AND$ operation $(F,A) \vee (K,B)$ is a soft neutrosophic bigroup subbring over $R_B \langle G \cup I \rangle$.

**Proof.** These are straightforward.

**Definition 4.2.3.** Let $R_B \langle G \cup I \rangle = R_1 \langle G_1 \cup I \rangle \cup R_2 \langle G_2 \cup I \rangle$ be a neutrosophic bigroup biring. Let $A$ be a set of parameters and $(F,A)$ be a soft set over $R_B \langle G \cup I \rangle$. Then $(F,A)$ is called soft mixed neutosophic bigroup biring over $R_B \langle G \cup I \rangle$ if for some $a \in A$, $F(a)$ is subneutrosophic bigroup subbiring and for the remaining $a \in A$, $F(a)$ is a subneutrosophic bigroup biring of $R_B \langle G \cup I \rangle$.

**Example 4.2.3.** Let $R_b \langle G \cup I \rangle = \mathbb{Z} \langle G_1 \cup I \rangle \cup \mathbb{C} \langle G_2 \cup I \rangle$ be a neutrosophic bigroup biring, where $\langle G_1 \cup I \rangle = \{1, g, g^2, g^3, g^4, g^5, I, gI, ..., g^5 I : g^6 = 1, I^2 = I\}$ and $\langle G_2 \cup I \rangle = \{1, g, g^2, g^3, I, gI, g^2 I, g^{53} I : g^4 = 1, I^2 = I\}$. Let $A = \{a_1, a_2, a_3, a_4\}$ be a set of parameters. Then $(F, A)$ is a soft mixed neutrosophic bigroup biring over $R_B \langle G \cup I \rangle$, where

$$F(a_1) = 2\mathbb{Z} \langle G_1 \cup I \rangle \cup \mathbb{R} \langle G_2 \cup I \rangle,$$

$$F(a_2) = 4\mathbb{Z} \langle G_1 \cup I \rangle \cup \mathbb{Q} \langle G_2 \cup I \rangle,$$

$$F(a_3) = \mathbb{Z} \langle H_1^{'} \cup I \rangle \mathbb{C} \langle H_2^{'} \cup I \rangle,$$

$$F(a_4) = \mathbb{Z} \langle H_1^{''} \cup I \rangle \mathbb{C} \langle H_2^{''} \cup I \rangle.$$

Where $\langle H_1^{'} \cup I \rangle = \{1, g^2 : g^4 = 1\}$, $\langle H_2^{'} \cup I \rangle = \{1, g^3, I, g^3 I : g^6 = 1, I^2 = I\}$, $\langle H_1^{''} \cup I \rangle = \{1, g^2, I, g^2 I : g^4 = 1, I^2 = I\}$ and $\langle H_2^{''} \cup I \rangle = \{1, g^2, g^4, I, g^2 I, g^4 I : g^6 = 1, I^2 = I\}$.

**Theorem 4.2.5.** Let $F, A$ and $(H, A)$ be two soft mixed neutrosophic bigroup birings over $R_B \langle G \cup I \rangle$. Then their intersection $F, A \cap H, A$ is again a soft mixed neutrosophic bigroup biring over $R_B \langle G \cup I \rangle$.

**Proof.** The proof is straightforward.

**Theorem 4.2.6.** Let $F, A$ and $H, B$ be two soft mixed neutrosophic bigroup birings over $R_B \langle G \cup I \rangle$. If $A \cap B = \phi$, then $F, A \cup H, B$ is a soft mixed neutrosophic bigroup biring over $R_B \langle G \cup I \rangle$.

**Proof.** This is straightforward.

**Remark 4.2.3.** Let $R_B \langle G \cup I \rangle$ be a neutrosophic bigroup biring. Then $R_B \langle G \cup I \rangle$ can have soft neutrosophic bigroup biring, soft neutrosophic bigroup subbiring and soft mixed neutrosophic bigroup biring over $R_B \langle G \cup I \rangle$.

**Proof:** It is obvious.

**Proposition 4.2.3.** Let $(F, A)$ and $(K, B)$ be two soft mixed neutosophic bigroup birings over $R_B \langle G \cup I \rangle$. Then

1. Their extended intersection $(F, A) \cap_E (K, B)$ is a soft mixed neutrosophic bigroup biring over $R_B \langle G \cup I \rangle$.
2. Their restricted intersection $(F, A) \cap_R (K, B)$ is a soft mixed neutrosophic bigroup biring over $R_B \langle G \cup I \rangle$.
3. Their $AND$ operation $(F, A) \vee (K, B)$ is a soft mixed neutrosophic bigroup biring over $R_B \langle G \cup I \rangle$.

**Proof.** This is straightforward.

**Remark 4.2.4.** Let $(F,A)$ and $(K,B)$ be two soft mixed neutosophic bigroup birings over $R_B \langle G \cup I \rangle$. Then

1. Their restricted union $(F,A) \cup_R (K,B)$ is not a soft mixed neutrosophic bigroup biring over $R_B \langle G \cup I \rangle$.
2. Their extended union $(F,A) \cup_E (K,B)$ is not a soft mixed neutrosophic bigroup biring over $R_B \langle G \cup I \rangle$.
3. Their $OR$ operation $(F,A) \vee (K,B)$ is a soft mixed neutrosophic bigroup biring over $R_B \langle G \cup I \rangle$.

One can easily show these remarks by the help of examples.

**Definition 4.2.3.** Let $R_{B_s} \langle G \cup I \rangle = R_1 \langle G_1 \cup I \rangle \cup R_2 \langle G_2 \cup I \rangle$ be a neutrosophic bigroup biring. Let $A$ be a set of parameters and $(F,A)$ be a soft set over $R_B \langle G \cup I \rangle$. Then $(F,A)$ is called soft neutosophic subbiring if and only if $F(a)$ is neutrosophic subbiring of $R_B \langle G \cup I \rangle$ for all $a \in A$ .

**Example 4.2.3.** Let $R_B \langle G \cup I \rangle = \mathbb{Z} \langle G_1 \cup I \rangle \cup \mathbb{C} \langle G_1 \cup I \rangle$ be a neutrosophic bigroup biring. Let $A = \{a_1, a_2, a_3\}$ be a set of parameters. Then $(F,A)$ is a soft neutrosophic subbiring over $R_B \langle G \cup I \rangle$, where

$$F(a_1) = \langle 2\mathbb{Z} \cup I \rangle \cup \langle \mathbb{R} \cup I \rangle,$$

$$F(a_2) = \langle 4\mathbb{Z} \cup I \rangle \cup \langle \mathbb{Q} \cup I \rangle,$$

$$F(a_3) = \langle 6\mathbb{Z} \cup I \rangle \cup \langle \mathbb{Z} \cup I \rangle.$$

**Theorem 4.2.7.** Let $F, A$ and $(H, A)$ be two soft neutrosophic subbirings over $R_B \langle G \cup I \rangle$. Then their intersection $F, A \cap H, A$ is again a soft neutrosophic subbiring over $R_B \langle G \cup I \rangle$.

**Proof.** The proof is straightforward.

**Theorem 4.2.8.** Let $F, A$ and $H, B$ be two soft neutrosophic subbirings over $R_B \langle G \cup I \rangle$. If $A \cap B = \phi$, then $F, A \cup H, B$ is also a soft neutrosophic subbiring over $R_B \langle G \cup I \rangle$.

**Proof.** This is straightforward.

**Proposition 4.2.4.** Let $(F, A)$ and $(K, B)$ be two soft neutosophic subbirings over $R_B \langle G \cup I \rangle$. Then

1. Their extended intersection $(F, A) \cap_E (K, B)$ is a soft neutrosophic subbiring over $R_B \langle G \cup I \rangle$.
2. Their restricted intersection $(F, A) \cap_R (K, B)$ is a soft neutrosophic subbiring over $R_B \langle G \cup I \rangle$.
3. Their $AND$ operation $(F, A) \vee (K, B)$ is a soft neutrosophic subbiring over $R_B \langle G \cup I \rangle$.

**Proof.** This is straightforward.

**Remark 4.2.5.** Let $(F, A)$ and $(K, B)$ be two soft neutosophic subbirings over $R_B \langle G \cup I \rangle$. Then

1. Their restricted union $(F, A) \cup_R (K, B)$ is not a soft neutrosophic subbiring over $R_B \langle G \cup I \rangle$.

2. Their extended union $(F, A) \cup_E (K, B)$ is not a soft neutrosophic subbiring over $R_B \langle G \cup I \rangle$.

3. Their $OR$ operation $(F, A) \vee (K, B)$ is a soft neutrosophic subbiring over $R_B \langle G \cup I \rangle$.

One can easily show these remarks by the help of examples.

**Definition 4.2.4.** Let $R_B \langle G \cup I \rangle = R_1 \langle G_1 \cup I \rangle \cup R_2 \langle G_2 \cup I \rangle$ be a neutrosophic bigroup biring. Let $A$ be a set of parameters and $(F, A)$ be a soft set over $R_B \langle G \cup I \rangle$. Then $(F, A)$ is called soft pseudo neutosophic subbiring if and only if $F(a)$ is a pseudo neutrosophic subbiring of $R_B \langle G \cup I \rangle$ for all $a \in A$.

**Theorem 4.2.9.** Let $F, A$ and $(H, A)$ be two soft pseudo neutrosophic subbirings over $R_B \langle G \cup I \rangle$. Then their intersection $F, A \cap H, A$ is again a soft pseudo neutrosophic subbiring over $R_B \langle G \cup I \rangle$.

**Proof.** The proof is straightforward.

**Theorem 4.2.10.** Let $F, A$ and $H, B$ be two soft pseudo neutrosophic subbirings over $R_B \langle G \cup I \rangle$. If $A \cap B = \phi$, then $F, A \cup H, B$ is a soft pseudo neutrosophic subbiring over $R_B \langle G \cup I \rangle$.

**Proposition 4.2.5.** Let $(F,A)$ and $(K,B)$ be two soft pseudo neutosophic subbirings over $R_B \langle G \cup I \rangle$. Then

1. Their extended intersection $(F,A) \cap_E (K,B)$ is a soft pseudo neutrosophic subbiring over $R_B \langle G \cup I \rangle$.
2. Their restricted intersection $(F,A) \cap_R (K,B)$ is a soft pseudo neutrosophic subbiring over $R_b \langle G \cup I \rangle$.
3. Their $AND$ operation $(F,A) \vee (K,B)$ is a soft pseudo neutrosophic subbiring over $R_B \langle G \cup I \rangle$.

**Proof.** This is straightforward.

**Remark 4.2.6.** Let $(F,A)$ and $(K,B)$ be two soft pseudo neutosophic subbirings over $R_B \langle G \cup I \rangle$. Then

1. Their restricted union $(F,A) \cup_R (K,B)$ is not a soft pseudo neutrosophic subbiring over $R_B \langle G \cup I \rangle$.
2. Their extended union $(F,A) \cup_E (K,B)$ is not a soft pseudo neutrosophic subbbiring over $R_B \langle G \cup I \rangle$.
3. Their $OR$ operation $(F,A) \vee (K,B)$ is a soft pseudo neutrosophic subbiring over $R_B \langle G \cup I \rangle$.

One can easily show these remarks by the help of examples.

**Definition 4.2.5.** Let $R_B \langle G \cup I \rangle = R_1 \langle G_1 \cup I \rangle \cup R_2 \langle G_2 \cup I \rangle$ be a neutrosophic bigroup biring. Let $A$ be a set of parameters and $(F,A)$ be a soft set over $R_B \langle G \cup I \rangle$. Then $(F, A)$ is called soft neutosophic subbigroup biring if and only if $F(a)$ is neutrosophic subbigroup biring of $R_B \langle G \cup I \rangle$ for all $a \in A$ .

**Example 4.2.4.** Let $R_B \langle G \cup I \rangle = \mathbb{Z} \langle G_1 \cup I \rangle \cup \mathbb{C} \langle G_2 \cup I \rangle$ be a neutrosophic bigroup biring, where $\langle G_1 \cup I \rangle = \{1, g, g^2, g^3 I, gI, g^2 I, g^3 I : g^6 = 1, I^2 = I\}$ and $\langle G_2 \cup I \rangle = \{1, g, g^2, g^3, g^4, g^5, I, gI, \ldots, g^5 I : g^6 = 1, I^2 = I\}$ . Let $H_1^{'} = \{1, g^2 : g^4 = 1\}$, $H_2^{'} = \{1, g^2, g^4 : g^6 = 1\}$ , $H_1^{''} = \{1, g^2, I, g^2 I : g^4 = 1, I^2 = I\}$ and $H_2^{''} = \{1, g^3 : g^6 = 1\}$. Let $A = \{a_1, a_2\}$ be a set of parameters. Then $(F, A)$ is a soft neutrosophic subbigroup biring over $R_B \langle G \cup I \rangle$, where

$$F(a_1) = 2\mathbb{Z} H_1^{'} \cup \mathbb{R} H_2^{'} ,$$

$$F(a_1) = 8\mathbb{Z} H_1^{''} \cup \mathbb{Q} H_2^{''} .$$

**Theorem 4.2.11.** Let $F, A$ and $(H, A)$ be two soft neutrosophic subbigroup birings over $R_B \langle G \cup I \rangle$. Then their intersection $F, A \cap H, A$ is again a soft neutrosophic subbigroup biring over $R_B \langle G \cup I \rangle$ .

**Proof.** The proof is straightforward.

**Theorem 4.2.12.** Let $F, A$ and $H, B$ be two soft neutrosophic subbigroup birings over $R_B \langle G \cup I \rangle$. If $A \cap B = \phi$, then $F, A \cup H, B$ is a soft neutrosophic subbigroup biring over $R_B \langle G \cup I \rangle$.

**Proof.** This is straightforward.

**Proposition 4.2.6.** Let $(F, A)$ and $(K, B)$ be two soft neutosophic subbigroup birings over $R_B \langle G \cup I \rangle$. Then

1. Their extended intersection $(F, A) \cap_E (K, B)$ is a soft neutrosophic subbigroup biring over $R_B \langle G \cup I \rangle$.
2. Their restricted intersection $(F, A) \cap_R (K, B)$ is a soft neutrosophic subbigroup biring over $R_B \langle G \cup I \rangle$.
3. Their *AND* operation $(F, A) \vee (K, B)$ is a soft neutrosophic subbigroup biring over $R_B \langle G \cup I \rangle$.

**Proof.** This is straightforward.

**Remark 4.2.7.** Let $(F, A)$ and $(K, B)$ be two soft neutosophic subbigroup birings over $R_B \langle G \cup I \rangle$. Then

1. Their restricted union $(F, A) \cup_R (K, B)$ is not a soft neutrosophic subbigroup biring over $R_B \langle G \cup I \rangle$.
2. Their extended union $(F, A) \cup_E (K, B)$ is not a soft neutrosophic subbigroup biring over $R_B \langle G \cup I \rangle$.
3. Their *OR* operation $(F, A) \vee (K, B)$ is a soft neutrosophic subbigroup biring over $R_B \langle G \cup I \rangle$.

One can easily show these remarks by the help of examples.

**Definition 4.2.6.** Let $R_B \langle G \cup I \rangle = R_1 \langle G_1 \cup I \rangle \cup R_2 \langle G_2 \cup I \rangle$ be a neutrosophic bigroup biring. Let $A$ be a set of parameters and $(F, A)$ be a soft set over $R_B \langle G \cup I \rangle$. Then $(F, A)$ is called a soft subbiring over $R_B \langle G \cup I \rangle$ if and only if $F(a)$ is a subbiring of $R_B \langle G \cup I \rangle$ for all $a \in A$.

**Theorem 4.2.13.** Let $F, A$ and $(H, A)$ be two soft subbirings over $R_B \langle G \cup I \rangle$. Then their intersection $F, A \cap H, A$ is again a soft subbiring over $R_B \langle G \cup I \rangle$.

**Proof.** The proof is straightforward.

**Theorem 4.2.14.** Let $F, A$ and $H, B$ be two soft subbirings over $R_B \langle G \cup I \rangle$. If $A \cap B = \phi$, then $F, A \cup H, B$ is a soft subbiring over $R_B \langle G \cup I \rangle$.

**Proof.** This is straightforward.

**Proposition 4.2.7.** Let $(F, A)$ and $(K, B)$ be two soft subbirings over $R_B \langle G \cup I \rangle$. Then

1. Their extended intersection $(F, A) \cap_E (K, B)$ is a soft subbiring over $R_B \langle G \cup I \rangle$.
2. Their restricted intersection $(F, A) \cap_R (K, B)$ is a soft subbiring over $R_B \langle G \cup I \rangle$.
3. Their $AND$ operation $(F, A) \vee (K, B)$ is a soft subbiring over $R_B \langle G \cup I \rangle$.

**Proof.** These are straightforward.

**165**

**Remark 4.2.8.** Let $(F, A)$ and $(K, B)$ be two soft subbirings over $R_B \langle G \cup I \rangle$. Then

1. Their restricted union $(F, A) \cup_R (K, B)$ is not a soft subbiring over $R_B \langle G \cup I \rangle$.
2. Their extended union $(F, A) \cup_E (K, B)$ is not a soft subbiring over $R_B \langle G \cup I \rangle$.
3. Their $OR$ operation $(F, A) \vee (K, B)$ is a soft subbiring over $R_B \langle G \cup I \rangle$.

One can easily show these remarks by the help of examples.

**Definition 4.2.7.** Let $R_B \langle G \cup I \rangle = R_1 \langle G_1 \cup I \rangle \cup R_2 \langle G_2 \cup I \rangle$ be a neutrosophic bigroup biring. Then $(F, A)$ is called an absolute soft neutosophic bigroup biring if $F(a) = R_B \langle G \cup I \rangle$, for all $a \in A$.

**Definition 4.2.8.** Let $(F, A)$ and $(K, B)$ be two soft neutrosophic bigroup birings over $R_B \langle G \cup I \rangle$. Then $(K, B)$ is called soft subneutrosophic bigroup biring of $(F, A)$, if

1. $B \subseteq A$, and
2. $K(a)$ is a subneutrosophic bigroup biring of $F(a)$ for all $a \in B$.

**Definition 4.2.9.** Let $R_B \langle G \cup I \rangle$ be a neutrosophic bigroup biring. Let $A$ be a set of parameters and $(F, A)$ be a soft set over $R_B \langle G \cup I \rangle$. Then $(F, A)$ is called soft neutosophic biideal if and only if $F(a)$ is a neutrosophic biideal of $R_B \langle G \cup I \rangle$ for all $a \in A$.

Similarly one can define soft neutrosophic left biideal and soft neutrosophic right biideal over $R_B \langle G \cup I \rangle$.

**Example 4.2.5.** Let $R_B \langle G \cup I \rangle = \mathbb{Z} \langle G_1 \cup I \rangle \cup \mathbb{Z} \langle G_2 \cup I \rangle$ be a neutrosophic bigroup biring. Let $A = \{a_1, a_2, a_3\}$ be a set of parameters. Then $(F, A)$ is a soft neutrosophic biideal over $R_B \langle G \cup I \rangle$, where

$$F(a_1) = 2\mathbb{Z} \langle G_1 \cup I \rangle \cup 4\mathbb{Z} \langle G_2 \cup I \rangle,$$

$$F(a_2) = 4\mathbb{Z} \langle G_1 \cup I \rangle \cup 6\mathbb{Z} \langle G_2 \cup I \rangle,$$

$$F(a_3) = 6\mathbb{Z} \langle G_1 \cup I \rangle \cup 8\mathbb{Z} \langle G_2 \cup I \rangle.$$

**Theorem 4.2.15.** All soft neutrosophic biideals are trivially soft subneutrosophic bigroup birings but the converse is not true in general.

We can easily establish the converse by the help of example.

**Proposition 4.2.8.** Let $(F, A)$ and $(K, B)$ be two soft neutosophic biideals over $R_B \langle G \cup I \rangle$. Then

1. Their extended intersection $(F, A) \cap_E (K, B)$ is a soft neutrosophic biideal over $R_B \langle G \cup I \rangle$.
2. Their restricted intersection $(F, A) \cap_R (K, B)$ is a soft neutrosophic biideal over $R_B \langle G \cup I \rangle$.
3. Their $AND$ operation $(F, A) \vee (K, B)$ is a soft neutrosophic biideal over $R_B \langle G \cup I \rangle$.

**Proof.** This is straightforward.

**Remark 4.2.9.** Let $(F, A)$ and $(K, B)$ be two soft neutrosophic biideal over $R_B \langle G \cup I \rangle$. Then

1. Their restricted union $(F, A) \cup_R (K, B)$ is not a soft neutrosophic biideal over $R_B \langle G \cup I \rangle$.
2. Their extended union $(F, A) \cup_E (K, B)$ is not a soft neutrosophic biideal over $R_B \langle G \cup I \rangle$.
3. Their $OR$ operation $(F, A) \vee (K, B)$ is a soft neutrosophic biideal over $R_B \langle G \cup I \rangle$.

One can easily show these remarks by the help of examples.

**Definition 4.2.10.** Let $R_B \langle G \cup I \rangle$ be a neutrosophic bigroup biring. Let $A$ be a set of parameters and $(F, A)$ be a soft set over $R \langle G \cup I \rangle$. Then $(F, A)$ is called soft pseudo neutosophic biideal if and only if $F(a)$ is a pseudo neutrosophic biideal of $R_B \langle G \cup I \rangle$ for all $a \in A$.

**Theorem 4.2.16.** All soft pseudo neutrosophic biideals are trivially soft subneutrosophic bigroup birings but the converse is not true in general.

We can easily establish the converse by the help of example.

**Theorem 4.2.17.** All soft pseudo neutrosophic biideals are trivially soft neutrosophic biideals but the converse is not true in general.

We can easily establish the converse by the help of example.

**Proposition 4.2.9.** Let $(F, A)$ and $(K, B)$ be two soft pseudo neutosophic biideals over $R_B \langle G \cup I \rangle$. Then

**168**

1. Their extended intersection $(F, A) \cap_E (K, B)$ is a soft pseudo neutrosophic biideals over $R_b \langle G \cup I \rangle$.

2. Their restricted intersection $(F, A) \cap_R (K, B)$ is a soft pseudo neutrosophic biideals over $R_B \langle G \cup I \rangle$.

3. Their $AND$ operation $(F, A) \vee (K, B)$ is a soft pseudo neutrosophic biideals over $R_B \langle G \cup I \rangle$.

**Proof.** This is straightforward.

**Remark 4.2.10.** Let $(F, A)$ and $(K, B)$ be two soft pseudo neutosophic biideals over $R_B \langle G \cup I \rangle$. Then

1. Their restricted union $(F, A) \cup_R (K, B)$ is not a soft pseudo neutrosophic biideals over $R_B \langle G \cup I \rangle$.

2. Their extended union $(F, A) \cup_E (K, B)$ is not a soft pseudo neutrosophic biideals over $R_B \langle G \cup I \rangle$.

3. Their $OR$ operation $(F, A) \vee (K, B)$ is a soft pseudo neutrosophic biideals over $R_B \langle G \cup I \rangle$.

One can easily show these remarks by the help of examples.

**Definition 4.2.10.** Let $(F, A)$ and $(K, B)$ be two soft neutrosophic bigroup birings over $R_B \langle G \cup I \rangle$. Then $(K, B)$ is called soft neutrosophic biideal of $(F, A)$, if

1. $B \subseteq A$, and
2. $K(a)$ is a neutrosophic biideal of $F(a)$ for all $a \in B$.

**Theorem 4.2.18.** Every soft neutrosophic biideal of the soft neutrosophic bigroup biring over a neutrosophic bigroup biring is trivially a soft subneutrosophic bigroup biring.

In the proceeding section soft neutrosophic group ring over a neutrosophic group ring is to extend soft neutrosophic N-group N-ring over a neutrosophic N-group N-ring.

### 4.3    Soft Neutrosophic N-group N-ring

In this last section of the present chapter, we define soft neutrosophic N-group N-ring over a neutrosophic N-group N-ring. This is the generalization of soft neutrosophic group ring. We also studied some basic properties and many example to illustrate this theory.

We now proceed on to define soft neutrosophic N-group N-ring over a neutrosophic N-group N-ring.

**Definition 4.3.1.** Let $N\left(R\langle G \cup I\rangle\right) = R_1\langle G_1 \cup I\rangle \cup R_2\langle G_2 \cup I\rangle \cup \ldots \cup R_n\langle G_n \cup I\rangle$ be a neutrosophic bigroup biring over a biring $R = R_1 \cup R_2 \cup \ldots \cup R_n$ such that $n \geq 3$. Let $A$ be a set of parameters and let $(F, A)$ be a soft set over $N\left(R\langle G \cup I\rangle\right)$. Then $(F, A)$ is called a soft neutrosophic N-group N-ring over $N\left(R\langle G\cup I\rangle\right)$ if and only if $F(a)$ is a subneutrosophic N-group N-ring of $N\left(R\langle G\cup I\rangle\right)$ for all $a\in A$.

See the following example for illustration.

**Example 4.3.1.** Let $N\left(R\langle G\cup I\rangle\right)=\mathbb{Q}\langle G_1\cup I\rangle\cup\mathbb{R}\langle G_2\cup I\rangle\cup\mathbb{Z}\langle G_3\cup I\rangle$ be a neutrosophic 3-group 3-ring, where $R=\mathbb{Q}\cup\mathbb{R}\cup\mathbb{Z}$ and

$\langle G_1\cup I\rangle=\{1,g,g^2,g^3,g^4,g^5,I,gI,...,g^5I:g^6=1,I^2=I\}$,

$\langle G_2\cup I\rangle=\{1,g,g^2,g^3,I,gI,g^2I,g^3I:g^4=1,I^2=I\}$ and

$\langle G_3\cup I\rangle=\{1,g,g^2,g^3,g^4,...g^7,I,gI,g^2I,...,g^7I:g^8=1,I^2=I\}$.

Let $\langle H_1^{'}\cup I\rangle\cup\langle H_2^{'}\cup I\rangle\cup\langle H_3^{'}\cup I\rangle=\{1,g^3:g^6=1\}\cup\{1,g^2:g^4=1\}\cup\{1,g^4:g^8=1\}$,

$\langle H_1^{''}\cup I\rangle\cup\langle H_2^{''}\cup I\rangle\cup\langle H_3^{''}\cup I\rangle=\{1,g^3,I,g^3I:g^6=1,I^2=I\}\cup\{1,g^2,I,g^2I\}\cup\{1,g^2,I,g^2I:g^8=1\}$.
Let $A=\{a_1,a_2\}$ be a set of parameters. Then $(F,A)$ is a soft neutrosophic 3-group 3-ring over $N\left(R\langle G\cup I\rangle\right)$, where

$$F(a_1)=\mathbb{Q}\langle H_1^{'}\cup I\rangle\cup\mathbb{R}\langle H_2^{'}\cup I\rangle\cup\mathbb{Z}\langle H_3^{'}\cup I\rangle,$$

$$F(a_2)=\mathbb{Q}\langle H_1^{''}\cup I\rangle\cup\mathbb{R}\langle H_2^{''}\cup I\rangle\cup\mathbb{Z}\langle H_3^{''}\cup I\rangle.$$

**Theorem 4.3.1.** Let $F,A$ and $(H,A)$ be two soft neutrosophic N-group N-rings over $N\left(R\langle G\cup I\rangle\right)$. Then their intersection $F,A\cap H,A$ is again a soft neutrosophic N-group N-ring over $N\left(R\langle G\cup I\rangle\right)$.

The proof is straightforward, so left as an exercise for the interested readers.

**Theorem 4.3.2.** Let $(F, A)$ and $(H, B)$ be two soft neutrosophic N-group N-rings over $N\big(R\langle G \cup I\rangle\big)$. If $A \cap B = \phi$, then $(F, A) \cup (H, B)$ is a soft neutrosophic N-group N-ring over $N\big(R\langle G \cup I\rangle\big)$.

**Proof.** This is straightforward.

**Proposition 4.3.1.** Let $(F, A)$ and $(K, B)$ be two soft neutosophic N-group N-rings over $N\big(R\langle G \cup I\rangle\big)$. Then

1. Their extended intersection $(F, A) \cap_E (K, B)$ is again a soft neutrosophic N-group N-ring over $N\big(R\langle G \cup I\rangle\big)$.

2. Their restricted intersection $(F, A) \cap_R (K, B)$ is again a soft neutrosophic N-group N-ring over $R\langle G \cup I\rangle$.

3. Their $AND$ operation $(F, A) \vee (K, B)$ is again a soft neutrosophic N-group N-ring over $R\langle G \cup I\rangle$.

**Proof.** This is straightforward.

**Remark 4.3.1.** Let $(F, A)$ and $(K, B)$ be two soft neutosophic N-group N-rings over $N\big(R\langle G \cup I\rangle\big)$. Then

1. Their restricted union $(F, A) \cup_R (K, B)$ is not a soft neutrosophic N-group N-ring over $N\big(R\langle G \cup I\rangle\big)$.

2. Their extended union $(F, A) \cup_E (K, B)$ is not a soft neutrosophic N-group N-ring over $N\big(R\langle G \cup I\rangle\big)$.

3. Their $OR$ operation $(F, A) \vee (K, B)$ is not a soft neutrosophic N-group N-ring over $N\big(R\langle G \cup I\rangle\big)$.

**172**

One can establish these remarks by the help of following examples.

**Definition 4.3.2.** Let $N(R\langle G \cup I\rangle)$ be a neutrosophic N-group N-ring. Let $A$ be a set of parameters and $(F,A)$ be a soft set over $N(R\langle G \cup I\rangle)$. Then $(F,A)$ is called soft neutosophic N-group sub N-ring over $N(R\langle G \cup I\rangle)$ if and only if $F(a)$ is subneutrosophic N-group sub N-ring of $N(R\langle G \cup I\rangle)$ for all $a \in A$.

**Example 4.3.2.** Let $N(R\langle G \cup I\rangle) = \mathbb{R}\langle G_1 \cup I\rangle \cup \mathbb{C}\langle G_2 \cup I\rangle \cup \mathbb{Z}\langle G_3 \cup I\rangle$ be a neutrosophic 3-group 3-ring, where

$\langle G_1 \cup I\rangle = \{1, g, g^2, g^3, g^4, g^5, I, gI, ..., g^5 I : g^6 = 1, I^2 = I\}$,

$\langle G_2 \cup I\rangle = \{1, g, g^2, g^3, I, gI, g^2 I, g^3 I : g^4 = 1, I^2 = I\}$ and

$\langle G_3 \cup I\rangle = \{1, g, g^2, g^3, g^4, g^5, g^6, g^7, I, gI, ..., g^7 I : g^8 = 1, I^2 = I\}$.

Let $A = \{a_1, a_2\}$ be a set of parameters. Then $(F,A)$ is a soft neutrosophic 3-group sub 3-ring over $N(R\langle G \cup I\rangle)$, where

$$F(a_1) = \mathbb{Q}\langle G_1 \cup I\rangle \cup \mathbb{R}\langle G_2 \cup I\rangle \cup 2\mathbb{Z}\langle G_3 \cup I\rangle,$$

$$F(a_2) = \mathbb{Z}\langle G_1 \cup I\rangle \cup \mathbb{Q}\langle G_2 \cup I\rangle \cup 3\mathbb{Z}\langle G_3 \cup I\rangle.$$

**Theorem 4.3.3.** Let $F,A$ and $(H,A)$ be two soft neutrosophic N-group sub N-rings over $N(R\langle G \cup I\rangle)$. Then their intersection $F,A \cap H,A$ is again a soft neutrosophic N-group sub N-rings over $N(R\langle G \cup I\rangle)$.

**Proof.** The proof is straightforward.

**173**

**Theorem 4.3.4.** Let $F,A$ and $H,B$ be two soft neutrosophic N-group sub N-rings over $N\left(R\langle G\cup I\rangle\right)$. If $A\cap B=\phi$, then $F,A\ \cup\ H,B$ is a soft neutrosophic N-group sub N-ring over $N\left(R\langle G\cup I\rangle\right)$.

**Proof.** This is straightforward.

**Remark 4.3.2.** Let $(F,A)$ and $(K,B)$ be two soft neutosophic N-group sub N-rings over $N\left(R\langle G\cup I\rangle\right)$. Then

1. Their extended union $(F,A)\cup_E(K,B)$ is not a soft neutrosophic N-group sub N-ring over $N\left(R\langle G\cup I\rangle\right)$.
2. Their restricted union $(F,A)\cup_R(K,B)$ is not a soft neutrosophic N-group sub N-ring over $N\left(R\langle G\cup I\rangle\right)$.
3. Their $OR$ operation $(F,A)\vee(K,B)$ is not a soft neutrosophic N-group sub N-ring over $N\left(R\langle G\cup I\rangle\right)$.

One can easily show these remaks by the help of examples.

**Proposition 4.3.2.** Let $(F,A)$ and $(K,B)$ be two soft neutosophic N-group sub N-rings over $N\left(R\langle G\cup I\rangle\right)$. Then

1. Their restricted intersection $(F,A)\cap_R(K,B)$ is a soft neutrosophic N-group sub N-ring over $N\left(R\langle G\cup I\rangle\right)$.
2. Their extended intersection $(F,A)\cap_E(K,B)$ is a soft neutrosophic N-group sub N-ring over $N\left(R\langle G\cup I\rangle\right)$.
3. Their $AND$ operation $(F,A)\vee(K,B)$ is a soft neutrosophic N-group

sub N-ring over $N\left(R\langle G \cup I\rangle\right)$.

**Proof.** These are straightforward.

**Definition 4.3.3.** Let $N\left(R\langle G \cup I\rangle\right)$ be a neutrosophic bigroup biring. Let $A$ be a set of parameters and $(F,A)$ be a soft set over $R\langle G \cup I\rangle$. Then $(F,A)$ is called soft mixed neutosophic N-group N-ring over $N\left(R\langle G \cup I\rangle\right)$ if for some $a \in A$, $F(a)$ is subneutrosophic N-group sub N-ring and for the remaining $a \in A$, $F(a)$ is a subneutrosophic N-group N-ring of $N\left(R\langle G \cup I\rangle\right)$.

**Example 4.3.3.** Let $N\left(R\langle G \cup I\rangle\right)=\mathbb{Z}\langle G_1 \cup I\rangle \cup \mathbb{C}\langle G_2 \cup I\rangle \cup \mathbb{R}\langle G_3 \cup I\rangle$ be a neutrosophic 3-group 3-ring, where

$\langle G_1 \cup I\rangle=\{1,g,g^2,g^3,g^4,g^5,I,gI,...,g^5I:g^6=1,I^2=I\}$ ,

$\langle G_2 \cup I\rangle=\{1,g,g^2,g^3,I,gI,g^2I,g^{53}I:g^4=1,I^2=I\}$ and

$\langle G_3 \cup I\rangle=\{1,g,g^2,g^3,g^4,g^5,g^6,g^7,I,gI,...,g^7I:g^8=1,I^2=I\}$ .

Let $A=\{a_1,a_2,a_3,a_4\}$ be a set of parameters. Then $(F,A)$ is a soft mixed neutrosophic 3-group 3-ring over $N(R\langle G \cup I\rangle)$, where

$$F(a_1)=2\mathbb{Z}\langle G_1 \cup I\rangle \cup \mathbb{R}\langle G_2 \cup I\rangle \cup \mathbb{Q}\langle G_3 \cup I\rangle ,$$

$$F(a_2)=4\mathbb{Z}\langle G_1 \cup I\rangle \cup \mathbb{Q}\langle G_2 \cup I\rangle \cup \mathbb{Z}\langle G_3 \cup I\rangle ,$$

$$F(a_3)=\mathbb{Z}\langle H_1^{'} \cup I\rangle \cup \mathbb{C}\langle H_2^{'} \cup I\rangle \cup \mathbb{R}\langle H_3^{'} \cup I\rangle ,$$

$$F(a_4)=\mathbb{Z}\langle H_1^{''} \cup I\rangle \cup \mathbb{C}\langle H_2^{''} \cup I\rangle \cup \mathbb{R}\langle H_3^{''} \cup I\rangle .$$

Where $\langle H_1' \cup I \rangle = \{1, g^2 : g^4 = 1\}$, $\langle H_2' \cup I \rangle = \{1, g^3, I, g^3I : g^6 = 1, I^2 = I\}$,

$\langle H_3' \cup I \rangle = \{1, g^4, I, g^4I : g^8 = 1, I^2 = I\}$, $\langle H_1'' \cup I \rangle = \{1, g^2, I, g^2I : g^4 = 1, I^2 = I\}$,

$\langle H_2'' \cup I \rangle = \{1, g^2, g^4, I, g^2I, g^4I : g^6 = 1, I^2 = I\}$ and $\langle H_3'' \cup I \rangle = \{1, I, g^4I : g^8 = 1, I^2 = I\}$.

**Theorem 4.3.5.** Let $F, A$ and $(H, A)$ be two soft mixed neutrosophic N-group N-rings over $N\ R\langle G \cup I \rangle$. Then their intersection $F, A \cap H, A$ is again a soft mixed neutrosophic N-group N-ring over $N\ R\langle G \cup I \rangle$.

**Proof.** The proof is straightforward.

**Theorem 4.3.6.** Let $F, A$ and $H, B$ be two soft mixed neutrosophic N-group N-rings over $N\ R\langle G \cup I \rangle$. If $A \cap B = \phi$, then $F, A \cup H, B$ is a soft mixed neutrosophic N-group N-ring over $N\ R\langle G \cup I \rangle$.

**Proof.** This is straightforward.

**Remark 4.3.3.** Let $R\langle G \cup I \rangle$ be a neutrosophic N-group N-ring. Then $N\ R\langle G \cup I \rangle$ can have soft neutrosophic N-group N-ring, soft neutrosophic N-group sub N-ring and soft mixed neutrosophic N-group N-ring over $N\ R\langle G \cup I \rangle$.

**Proof:** It is obvious.

**Proposition 4.3.3.** Let $(F, A)$ and $(K, B)$ be two soft mixed neutosophic N-group N-rings over $N\ R\langle G \cup I \rangle$. Then

**176**

1. Their extended intersection $(F,A) \cap_E (K,B)$ is a soft mixed neutrosophic N-group N-ring over $N \ R \langle G \cup I \rangle$ .

2. Their restricted intersection $(F,A) \cap_R (K,B)$ is a soft mixed neutrosophic N-group N-ring over $N \ R \langle G \cup I \rangle$ .

3. Their $AND$ operation $(F,A) \vee (K,B)$ is a soft mixed neutrosophic N-group N-ring over $N \ R \langle G \cup I \rangle$ .

**Proof.** This is straightforward.

**Remark 4.3.4.** Let $(F,A)$ and $(K,B)$ be two soft mixed neutosophic N-group N-rings over $N \ R \langle G \cup I \rangle$ . Then

1. Their restricted union $(F,A) \cup_R (K,B)$ is not a soft mixed neutrosophic N-group N-ring over $N \ R \langle G \cup I \rangle$ .

2. Their extended union $(F,A) \cup_E (K,B)$ is not a soft mixed neutrosophic N-group N-ring over $N \ R \langle G \cup I \rangle$ .

3. Their $OR$ operation $(F,A) \vee (K,B)$ is a soft mixed neutrosophic N-group N-ring over $N \ R \langle G \cup I \rangle$ .

One can easily show these remarks by the help of examples.

**Definition 4.3.4.** Let $N \ R \langle G \cup I \rangle$ be a neutrosophic N-group N-ring. Let $A$ be a set of parameters and $(F,A)$ be a soft set over $N \ R \langle G \cup I \rangle$ . Then $(F,A)$ is called soft neutosophic sub N-ring if and only if $F(a)$ is neutrosophic sub N-ring of $N \ R \langle G \cup I \rangle$ for all $a \in A$ .

**Example 4.3.4.** Let $N\left(R\langle G \cup I\rangle\right) = \mathbb{Z}\langle G_1 \cup I\rangle \cup \mathbb{C}\langle G_2 \cup I\rangle \cup \mathbb{R}\langle G_3 \cup I\rangle$ be a neutrosophic 3-group 3-ring. Let $A = \{a_1, a_2\}$ be a set of parameters. Then $(F, A)$ is a soft neutrosophic sub 3-ring over $N\ R\langle G \cup I\rangle$ , where

$$F(a_1) = \langle 2\mathbb{Z} \cup I\rangle \cup \langle \mathbb{R} \cup I\rangle \cup \langle \mathbb{Q} \cup I\rangle,$$

$$F(a_2) = \langle 4\mathbb{Z} \cup I\rangle \cup \langle \mathbb{Q} \cup I\rangle \cup \langle \mathbb{Z} \cup I\rangle.$$

**Theorem 4.3.7.** Let $F, A$ and $(H, A)$ be two soft neutrosophic sub N-rings over $N\ R\langle G \cup I\rangle$ . Then their intersection $F, A \cap H, A$ is again a soft neutrosophic sub N-ring over $N\ R\langle G \cup I\rangle$ .

**Proof.** The proof is straightforward.

**Theorem 4.3.8.** Let $F, A$ and $H, B$ be two soft neutrosophic sub N-rings over $N\ R\langle G \cup I\rangle$ . If $A \cap B = \phi$, then $F, A \cup H, B$ is also a soft neutrosophic sub N-ring over $N\ R\langle G \cup I\rangle$ .

**Proof.** This is straightforward.

**Proposition 4.3.4.** Let $(F, A)$ and $(K, B)$ be two soft neutosophic sub N-rings over $N\ R\langle G \cup I\rangle$ . Then

1. Their extended intersection $(F, A) \cap_E (K, B)$ is a soft neutrosophic sub N-ring over $N\ R\langle G \cup I\rangle$ .
2. Their restricted intersection $(F, A) \cap_R (K, B)$ is a soft neutrosophic sub N-ring over $N\ R\langle G \cup I\rangle$ .

3. Their *AND* operation $(F,A) \vee (K,B)$ is a soft neutrosophic sub N-ring over $N\ R\langle G \cup I \rangle$ .

**Proof.** This is straightforward.

**Remark 4.3.5.** Let $(F,A)$ and $(K,B)$ be two soft neutosophic sub N-rings over $N\ R\langle G \cup I \rangle$ . Then

1. Their restricted union $(F,A) \cup_R (K,B)$ is not a soft neutrosophic sub N-ring over $N\ R\langle G \cup I \rangle$ .
2. Their extended union $(F,A) \cup_E (K,B)$ is not a soft neutrosophic sub N-ring over $N\ R\langle G \cup I \rangle$ .
3. Their *OR* operation $(F,A) \vee (K,B)$ is a soft neutrosophic sub N-ring over $N\ R\langle G \cup I \rangle$ .

One can easily show these remarks by the help of examples.

**Definition 4.3.5.** Let $N\ R\langle G \cup I \rangle$ be a neutrosophic N-group N-ring. Let $A$ be a set of parameters and $(F,A)$ be a soft set over $N\ R\langle G \cup I \rangle$ . Then $(F,A)$ is called soft pseudo neutosophic sub N-ring if and only if $F(a)$ is a pseudo neutrosophic sub N-ring of $N\ R\langle G \cup I \rangle$ for all $a \in A$ .

**Theorem 4.3.9.** Let $F,A$ and $(H,A)$ be two soft pseudo neutrosophic sub N-rings over $N\ R\langle G \cup I \rangle$ . Then their intersection $F,A \cap H,A$ is again a soft pseudo neutrosophic sub N-ring over $N\ R\langle G \cup I \rangle$ .

**Proof.** The proof is straightforward.

**Theorem 4.3.10.** Let $F, A$ and $H, B$ be two soft pseudo neutrosophic sub N-rings over $N\ R\langle G \cup I\rangle$ . If $A \cap B = \phi$, then $F, A\ \cup\ H, B$ is a soft pseudo neutrosophic sub N-ring over $N\ R\langle G \cup I\rangle$ .

**Proof.** This is straightforward.

**Proposition 4.3.5.** Let $(F, A)$ and $(K, B)$ be two soft pseudo neutosophic sub N-rings over $N\ R\langle G \cup I\rangle$ . Then

1. Their extended intersection $(F, A) \cap_E (K, B)$ is a soft pseudo neutrosophic sub N-ring over $N\ R\langle G \cup I\rangle$ .
2. Their restricted intersection $(F, A) \cap_R (K, B)$ is a soft pseudo neutrosophic sub N-ring over $N\ R\langle G \cup I\rangle$ .
3. Their $AND$ operation $(F, A) \vee (K, B)$ is a soft pseudo neutrosophic sub N-ring over $N\ R\langle G \cup I\rangle$ .

**Proof.** This is straightforward.

**Remark 4.3.6.** Let $(F, A)$ and $(K, B)$ be two soft pseudo neutosophic sub N-rings over $N\ R\langle G \cup I\rangle$ . Then

1. Their restricted union $(F, A) \cup_R (K, B)$ is not a soft pseudo neutrosophic sub N-ring over $N\ R\langle G \cup I\rangle$ .
2. Their extended union $(F, A) \cup_E (K, B)$ is not a soft pseudo neutrosophic sub N-ring over $N\ R\langle G \cup I\rangle$ .

3. Their $OR$ operation $(F, A) \vee (K, B)$ is a soft pseudo neutrosophic sub N-ring over $_N R \langle G \cup I \rangle$.

One can easily show these remarks by the help of examples.

**Definition 4.3.6.** Let $_N R \langle G \cup I \rangle$ be a neutrosophic N-group N-ring. Let $A$ be a set of parameters and $(F, A)$ be a soft set over $_N R \langle G \cup I \rangle$. Then $(F, A)$ is called soft neutosophic sub N-group N-ring if and only if $F(a)$ is neutrosophic sub N-group N-ring of $_N R \langle G \cup I \rangle$ for all $a \in A$.

**Example 4.3.5.** Let $N(R\langle G \cup I \rangle) = \mathbb{Z}\langle G_1 \cup I \rangle \cup \mathbb{C}\langle G_2 \cup I \rangle \cup \mathbb{R}\langle G_3 \cup I \rangle$ be a neutrosophic 3-group 3-ring, where
$\langle G_1 \cup I \rangle = \{1, g, g^2, g^3 I, gI, g^2 I, g^3 I : g^6 = 1, I^2 = I\}$ and
$\langle G_2 \cup I \rangle = \{1, g, g^2, g^3, g^4, g^5, I, gI, ..., g^5 I : g^6 = 1, I^2 = I\}$ and
$\langle G_3 \cup I \rangle = \{1, g, g^2, g^3, g^4, g^5, g^6, g^6, I, gI, ..., g^7 I : g^8 = 1, I^2 = I\}$. Let $H_1^{'} = \{1, g^2 : g^4 = 1\}$,
$H_2^{'} = \{1, g^2, g^4 : g^6 = 1\}$, $H_3^{'} = \{1, g^2, g^4, g^6 : g^8 = 1\}$, $H_1^{''} = \{1, g^2, I, g^2 I : g^4 = 1, I^2 = I\}$,
$H_2^{''} = \{1, g^3 : g^6 = 1\}$ and $H_3^{''} = \{1, g^4, I, g^4 I : g^8 = 1, I^2 = I\}$. Let $A = \{a_1, a_2\}$ be a set of parameters. Then $(F, A)$ is a soft neutrosophic sub 3-group 3-ring over $N(R\langle G \cup I \rangle)$, where

$$F(a_1) = 2\mathbb{Z}H_1^{'} \cup \mathbb{R}H_2^{'} \cup \mathbb{Q}H_3^{'},$$

$$F(a_1) = 8\mathbb{Z}H_1^{''} \cup \mathbb{Q}H_2^{''} \cup \mathbb{Z}H_3^{''}.$$

**Theorem 4.3.11.** Let $F, A$ and $(H, A)$ be two soft neutrosophic sub N-group N-rings over $N(R\langle G \cup I \rangle)$. Then their intersection $F,A \cap H,A$ is again a soft neutrosophic sub N-group N-ring over $N\left(R\langle G \cup I\rangle\right)$ .

**Proof.** The proof is straightforward.

**Theorem 4.3.12.** Let $F,A$ and $H,B$ be two soft neutrosophic sub N-group N-rings over $N\left(R\langle G \cup I\rangle\right)$. If $A \cap B = \phi$, then $F,A \cup H,B$ is a soft neutrosophic sub N-group N-ring over $R\langle G \cup I\rangle$.

**Proof.** This is straightforward.

**Proposition 4.3.6.** Let $(F,A)$ and $(K,B)$ be two soft neutosophic sub N-group N-rings over $N\left(R\langle G \cup I\rangle\right)$. Then

1. Their extended intersection $(F,A) \cap_E (K,B)$ is a soft neutrosophic sub N-group N-ring over $N\left(R\langle G \cup I\rangle\right)$.
2. Their restricted intersection $(F,A) \cap_R (K,B)$ is a soft neutrosophic sub N-group N-ring over $N\left(R\langle G \cup I\rangle\right)$.
3. Their *AND* operation $(F,A) \vee (K,B)$ is a soft neutrosophic sub N-group N-ring over $N\left(R\langle G \cup I\rangle\right)$.

**Proof.** This is straightforward.

**Remark 4.3.7.** Let $(F,A)$ and $(K,B)$ be two soft neutosophic sub N-group N-rings over $N\left(R\langle G \cup I\rangle\right)$. Then

1. Their restricted union $(F,A) \cup_R (K,B)$ is not a soft neutrosophic sub

N-group N-ring over $N\left(R\langle G\cup I\rangle\right)$.

2. Their extended union $(F,A)\cup_E(K,B)$ is not a soft neutrosophic sub N-group N-ring over $N\left(R\langle G\cup I\rangle\right)$.

3. Their *OR* operation $(F,A)\vee(K,B)$ is a soft neutrosophic sub N-group N-ring over $N\left(R\langle G\cup I\rangle\right)$.

One can easily show these remarks by the help of examples.

**Definition 4.3.7.** Let $N\left(R\langle G\cup I\rangle\right)$ be a neutrosophic N-group N-ring. Let $A$ be a set of parameters and $(F,A)$ be a soft set over $N\left(R\langle G\cup I\rangle\right)$. Then $(F,A)$ is called a soft sub N-ring over $N\left(R\langle G\cup I\rangle\right)$ if and only if $F(a)$ is a sub N-ring of $R\langle G\cup I\rangle$ for all $a\in A$ .

**Theorem 4.3.13.** Let $F,A$ and $(H,A)$ be two soft sub N-rings over $N\left(R\langle G\cup I\rangle\right)$. Then their intersection $F,A\cap H,A$ is again a soft sub N-ring over $N\left(R\langle G\cup I\rangle\right)$ .

**Proof.** The proof is straightforward.

**Theorem 4.3.14.** Let $F,A$ and $H,B$ be two soft sub N-rings over $N\left(R\langle G\cup I\rangle\right)$. If $A\cap B=\phi$, then $F,A\cup H,B$ is a soft sub N-ring over $N\left(R\langle G\cup I\rangle\right)$.

**Proof.** This is straightforward.

**Proposition 4.3.7.** Let $(F, A)$ and $(K, B)$ be two soft sub N-rings over $N\big(R\langle G \cup I\rangle\big)$. Then

1. Their extended intersection $(F, A) \cap_E (K, B)$ is a soft sub N-ring over $N\big(R\langle G \cup I\rangle\big)$.

2. Their restricted intersection $(F, A) \cap_R (K, B)$ is a soft sub N-ring over $N\big(R\langle G \cup I\rangle\big)$.

3. Their $AND$ operation $(F, A) \vee (K, B)$ is a soft sub N-ring over $N\big(R\langle G \cup I\rangle\big)$.

**Proof.** These are straightforward.

**Remark 4.3.8.** Let $(F, A)$ and $(K, B)$ be two soft sub N-rings over $N\big(R\langle G \cup I\rangle\big)$. Then

1. Their restricted union $(F, A) \cup_R (K, B)$ is not a soft sub N-ring over $N\big(R\langle G \cup I\rangle\big)$.

2. Their extended union $(F, A) \cup_E (K, B)$ is not a soft sub N-ring over $N\big(R\langle G \cup I\rangle\big)$.

3. Their $OR$ operation $(F, A) \vee (K, B)$ is a soft sub N-ring over $N\big(R\langle G \cup I\rangle\big)$.

One can easily show these remarks by the help of examples.

**Definition 4.3.8.** Let $N\left(R\langle G\cup I\rangle\right)$ be a neutrosophic N-group N-ring. Then $(F,A)$ is called an absolute soft neutosophic N-group N-ring if $F(a)=N\left(R\langle G\cup I\rangle\right)$, for all $a\in A$.

**Definition 4.3.9.** Let $(F,A)$ and $(K,B)$ be two soft neutrosophic N-group N-rings over $N\left(R\langle G\cup I\rangle\right)$. Then $(K,B)$ is called soft subneutrosophic N-group N-ring of $(F,A)$, if

1. $B\subseteq A$, and
2. $K(a)$ is a subneutrosophic N-group N-ring of $F(a)$ for all $a\in B$.

**Definition 4.3.10.** Let $N\left(R\langle G\cup I\rangle\right)$ be a neutrosophic N-group N-ring. Let $A$ be a set of parameters and $(F,A)$ be a soft set over $N\left(R\langle G\cup I\rangle\right)$. Then $(F,A)$ is called soft neutosophic N-ideal if and only if $F(a)$ is a neutrosophic N-ideal of $N\left(R\langle G\cup I\rangle\right)$ for all $a\in A$.

Similarly one can define soft neutrosophic left N-ideal and soft neutrosophic right N-ideal over $N\left(R\langle G\cup I\rangle\right)$.

**Example 4.3.6.** Let $N\left(R\langle G\cup I\rangle\right)=\mathbb{Z}\langle G_1\cup I\rangle\cup\mathbb{Z}\langle G_2\cup I\rangle\cup\mathbb{Z}\langle G_3\cup I\rangle$ be a neutrosophic 3-group 3-ring. Let $A=\{a_1,a_2,a_3\}$ be a set of parameters. Then $(F,A)$ is a soft neutrosophic 3-ideal over $N\left(R\langle G\cup I\rangle\right)$, where

$$F(a_1)=2\mathbb{Z}\langle G_1\cup I\rangle\cup 4\mathbb{Z}\langle G_2\cup I\rangle\cup 6\mathbb{Z}\langle G_3\cup I\rangle,$$

$$F(a_2)=4\mathbb{Z}\langle G_1\cup I\rangle\cup 6\mathbb{Z}\langle G_2\cup I\rangle\cup 8\mathbb{Z}\langle G_3\cup I\rangle,$$

$$F(a_3)=6\mathbb{Z}\langle G_1\cup I\rangle\cup 8\mathbb{Z}\langle G_2\cup I\rangle\cup 10\mathbb{Z}\langle G_3\cup I\rangle.$$

**Theorem 4.3.15.** All soft neutrosophic N-ideals are trivially soft subneutrosophic N-group N-rings but the converse is not true in general.

We can easily establish the converse by the help of example.

**Proposition 4.3.8.** Let $(F,A)$ and $(K,B)$ be two soft neutosophic N-ideals over $N\left(R\langle G \cup I\rangle\right)$. Then

1. Their extended intersection $(F,A) \cap_E (K,B)$ is a soft neutrosophic N-ideal over $N\left(R\langle G \cup I\rangle\right)$.
2. Their restricted intersection $(F,A) \cap_R (K,B)$ is a soft neutrosophic N-ideal over $N\left(R\langle G \cup I\rangle\right)$.
3. Their *AND* operation $(F,A) \vee (K,B)$ is a soft neutrosophic N-ideal over $N\left(R\langle G \cup I\rangle\right)$.

**Proof.** This is straightforward.

**Remark 4.3.9.** Let $(F,A)$ and $(K,B)$ be two soft neutosophic N-ideal over $N\left(R\langle G \cup I\rangle\right)$. Then

1. Their restricted union $(F,A) \cup_R (K,B)$ is not a soft neutrosophic N-ideal over $N\left(R\langle G \cup I\rangle\right)$.
2. Their extended union $(F,A) \cup_E (K,B)$ is not a soft neutrosophic N-ideal over $N\left(R\langle G \cup I\rangle\right)$.
3. Their *OR* operation $(F,A) \vee (K,B)$ is a soft neutrosophic N-ideal over $N\left(R\langle G \cup I\rangle\right)$.

One can easily show these remarks by the help of examples.

**Definition 4.3.11.** Let $N\left(R\langle G\cup I\rangle\right)$ be a neutrosophic N-group N-ring. Let $A$ be a set of parameters and $(F,A)$ be a soft set over $N\left(R\langle G\cup I\rangle\right)$. Then $(F,A)$ is called soft pseudo neutosophic N-ideal if and only if $F(a)$ is a pseudo neutrosophic N-ideal of $N\left(R\langle G\cup I\rangle\right)$ for all $a\in A$.

**Theorem 4.3.16.** All soft pseudo neutrosophic N-ideals are trivially soft subneutrosophic N-group N-rings but the converse is not true in general.

We can easily establish the converse by the help of example.

**Theorem 4.3.17.** All soft pseudo neutrosophic N-ideals are trivially soft neutrosophic N-ideals but the converse is not true in general.

We can easily establish the converse by the help of example.

**Proposition 4.3.9.** Let $(F,A)$ and $(K,B)$ be two soft pseudo neutosophic N-ideals over $N\left(R\langle G\cup I\rangle\right)$. Then

1. Their extended intersection $(F,A)\cap_E (K,B)$ is a soft pseudo neutrosophic N-ideals over $N\left(R\langle G\cup I\rangle\right)$.
2. Their restricted intersection $(F,A)\cap_R (K,B)$ is a soft pseudo neutrosophic N-ideals over $N\left(R\langle G\cup I\rangle\right)$.
3. Their $AND$ operation $(F,A)\vee (K,B)$ is a soft pseudo neutrosophic N-ideals over $N\left(R\langle G\cup I\rangle\right)$.

**Proof.** This is straightforward.

**Remark 4.3.10.** Let $(F, A)$ and $(K, B)$ be two soft pseudo neutosophic Ni-ideals over $N(R\langle G \cup I\rangle)$. Then

1. Their restricted union $(F, A) \cup_R (K, B)$ is not a soft pseudo neutrosophic N-ideals over $N(R\langle G \cup I\rangle)$.

2. Their extended union $(F, A) \cup_E (K, B)$ is not a soft pseudo neutrosophic N-ideals over $N(R\langle G \cup I\rangle)$.

3. Their $OR$ operation $(F, A) \vee (K, B)$ is a soft pseudo neutrosophic N-ideals over $N(R\langle G \cup I\rangle)$.

One can easily show these remarks by the help of examples.

**Definition 4.3.12.** Let $(F, A)$ and $(K, B)$ be two soft neutrosophic N-group N-rings over $N(R\langle G \cup I\rangle)$. Then $(K, B)$ is called soft neutrosophic N-ideal of $(F, A)$, if

1. $B \subseteq A$, and
2. $K(a)$ is a neutrosophic N-ideal of $F(a)$ for all $a \in B$.

**Theorem 4.3.18.** Every soft neutrosophic N-ideal of the soft neutrosophic N-group N-ring over a neutrosophic N-group N-ring is trivially a soft subneutrosophic N-group N-ring.

# Chapter No. 5

# Soft Neutrosophic Semigroup Ring and Their Generalization

In this chapter we take the important neutrosophic algebraic structures neutrosophic semigroup rings to define soft neutrosophic semigroup rings over neutrosophic semigroup rings. We also extend this theory and give the general concept of soft neutrosophic N-semigroup N-rings over neutrosophic N-semigroup N-rings. Some core properties are also studied here and many illustrative examples are given.

We proceed to define soft neutrosophic semigroup ring over a neutrosophic semigroup ring.

## 5.1    Soft Neutrosophic Semigroup Ring

This section is about to introduces soft neutrosophic semigroup ring over a neutrosophic semigroup ring. This is more generalized from soft neutroosphic group ring over a neutrosophic group ring. This is the parameterized collection of subneutrosophic semigroup rings of a neutroosphic semigroup ring. We also give some basic and fundamental results with illustrative examples.

**Definition 5.1.1.** Let $R\langle S \cup I \rangle$ be a neutrosophic semigroup ring. Let $A$ be a set of parameters and $(F, A)$ be a soft set over $R\langle S \cup I \rangle$. Then $(F, A)$ is called soft neutosophic semigroup ring if and only if $F(a)$ is subneutrosophic semigroup ring of $R\langle S \cup I \rangle$ for all $a \in A$.

This situation can be explained in the following example.

**Example 5.1.1.** Let $\mathbb{Q}\langle Z^+ \cup \{0\} \cup \{I\} \rangle$ be a neutrosophic semigroup ring, where $\mathbb{Q} =$ field of rationals and $\langle S \cup I \rangle = \langle Z^+ \cup \{0\} \cup \{I\} \rangle$ be a neutrosophic semigroup under $_+$. Let $\langle H_1 \cup I \rangle = \langle 2Z^+ \cup I \rangle$ and $\langle H_2 \cup I \rangle = \langle 3Z^+ \cup I \rangle$. Let $A = \{a_1, a_2\}$ be a set of parameters. Then $(F, A)$ is a soft neutrosophic semigroup ring over $R\langle S \cup I \rangle$, where

$$F(a_1) = \mathbb{Q}\langle H_1 \cup I \rangle = \langle 2Z^+ \cup I \rangle,$$

$$F(a_2) = \langle 3Z^+ \cup I \rangle.$$

Now some characterization of soft neutrosophic semigroup ring over a neutrosophic semigroup ring is presented here.

**Theorem 5.1.1.** Let $F, A$ and $(H, A)$ be two soft neutrosophic semigroup rings over $R\langle S \cup I \rangle$. Then their intersection $F, A \cap H, A$ is again a soft neutrosophic semigroup ring over $R\langle S \cup I \rangle$.

The proof is straightforward, so left as an exercise for the interested readers.

**Theorem 5.1.2.** Let $F, A$ and $H, B$ be two soft neutrosophic semigroup rings over $R\langle S \cup I \rangle$. If $A \cap B = \phi$, then $F, A \cup H, B$ is a soft neutrosophic semigroup ring over $R\langle S \cup I \rangle$.

**Theorem 5.1.3.** If $R\langle S \cup I \rangle$ is a commutative neutrosophic semigroup ring. Then the soft neutrosophic semigroup ring $(F, A)$ is also commutative soft neutrosophic semigroup ring.

**Proof:** It is trivial.

**Proposition 5.1.1.** Let $(F, A)$ and $(K, B)$ be two soft neutosophic semigroup rings over $R\langle S \cup I \rangle$. Then

1. Their extended intersection $(F, A) \cap_E (K, B)$ is again a soft neutrosophic semigroup ring over $R\langle S \cup I \rangle$.
2. Their restricted intersection $(F, A) \cap_R (K, B)$ is again a soft neutrosophic semigroup ring over $R\langle S \cup I \rangle$.
3. Their $AND$ operation $(F, A) \vee (K, B)$ is again a soft neutrosophic semigroup ring over $R\langle S \cup I \rangle$.

**Remark 5.1.1.** Let $(F, A)$ and $(K, B)$ be two soft neutosophic semigroup rings over $R\langle S \cup I \rangle$. Then

1. Their restricted union $(F, A) \cup_R (K, B)$ is not a soft neutrosophic semigroup ring over $R\langle S \cup I \rangle$.
2. Their extended union $(F, A) \cup_E (K, B)$ is not a soft neutrosophic semigroup ring over $R\langle S \cup I \rangle$.

**191**

3. Their $OR$ operation $(F,A) \vee (K,B)$ is not a soft neutrosophic semigroup ring over $R\langle S \cup I \rangle$.

One can easily prove these remarks by the help of examples.

**Definition 5.1.2.** Let $R\langle S \cup I \rangle$ be a neutrosophic semigroup ring. Let $A$ be a set of parameters and $(F,A)$ be a soft set over $R\langle S \cup I \rangle$. Then $(F,A)$ is called soft neutosophic semigroup subring if and only if $F(a)$ is subneutrosophic semigroup subring of $R\langle S \cup I \rangle$ for all $a \in A$.

**Example 5.1.2.** Let $\mathbb{C}\langle Z^+ \cup \{0\} \cup \{I\} \rangle$ be a neutrosophic semigroup ring, where $\mathbb{C} =$ field of complex numbers and $\langle Z^+ \cup \{0\} \cup \{I\} \rangle$ is a neutrosophic semigroup ring under $+$. Let $A = \{a_1, a_2, a_3\}$ be a set of parameters. Then $(F,A)$ is a soft neutrosophic semigroup subring over $\mathbb{C}\langle Z^+ \cup \{0\} \cup \{I\} \rangle$, where

$$F(a_1) = \mathbb{R}\langle Z^+ \cup \{0\} \cup \{I\} \rangle,$$

$$F(a_2) = \mathbb{Q}\langle Z^+ \cup \{0\} \cup \{I\} \rangle,$$

$$F(a_3) = \mathbb{Z}\langle Z^+ \cup \{0\} \cup \{I\} \rangle.$$

**Theorem 5.1.4.** Let $F,A$ and $(H,A)$ be two soft neutrosophic semigroup subrings over $R\langle S \cup I \rangle$. Then their intersection $F,A \cap H,A$ is again a soft neutrosophic semigroup subrings over $R\langle S \cup I \rangle$.

**Proof.** The proof is straightforward.

**Theorem 5.1.5.** Let $F, A$ and $H, B$ be two soft neutrosophic semigroup subrings over $R\langle S \cup I \rangle$. If $A \cap B = \phi$, then $F, A \cup H, B$ is a soft neutrosophic semigroup subring over $R\langle S \cup I \rangle$.

**Proof.** This is straightforward.

**Remark 5.1.2.** Let $(F, A)$ and $(K, B)$ be two soft neutosophic semigroup subrings over $R\langle S \cup I \rangle$. Then

1. Their extended union $(F, A) \cup_E (K, B)$ is not a soft neutrosophic semigroup subring over $R\langle S \cup I \rangle$.
2. Their restricted union $(F, A) \cup_R (K, B)$ is not a soft neutrosophic semigroup subring over $R\langle S \cup I \rangle$.
3. Their $OR$ operation $(F, A) \vee (K, B)$ is not a soft neutrosophic semigroup subring over $R\langle S \cup I \rangle$.

One can easily show these remarks by the help of examples.

**Proposition 5.1.2.** Let $(F, A)$ and $(K, B)$ be two soft neutosophic semigroup subrings over $R\langle S \cup I \rangle$. Then

1. Their restricted intersection $(F, A) \cap_R (K, B)$ is a soft neutrosophic semigroup subring over $R\langle S \cup I \rangle$.
2. Their extended intersection $(F, A) \cap_E (K, B)$ is a soft neutrosophic semigroup subring over $R\langle S \cup I \rangle$.
3. Their $AND$ operation $(F, A) \vee (K, B)$ is a soft neutrosophic semigroup subring over $R\langle S \cup I \rangle$.

**Definition 5.1.3.** Let $R\langle S \cup I \rangle$ be a neutrosophic semigroup ring. Let $A$ be a set of parameters and $(F, A)$ be a soft set over $R\langle S \cup I \rangle$. Then $(F, A)$ is called soft mixed neutosophic semigroup ring if for some $a \in A$, $F(a)$ is subneutrosophic semigroup subring and for the remaining $a \in A$, $F(a)$ is a subneutrosophic semigroup ring of $R\langle S \cup I \rangle$.

**Example 5.1.3.** Let $R\langle S \cup I \rangle = \mathbb{C}\langle Z^+ \cup \{0\} \cup \{I\}\rangle$ be a neutrosophic semigroup ring, where $\mathbb{C} =$ field of complex numbers and $\langle S \cup I \rangle = \langle Z^+ \cup \{0\} \cup \{I\}\rangle$ be a neutrosophic semigroup under $+$. Let $\langle H_1 \cup I \rangle = \langle 2Z^+ \cup I \rangle$ and $\langle H_2 \cup I \rangle = \langle 3Z^+ \cup I \rangle$. Let $A = \{a_1, a_2, a_3, a_4, a_5\}$ be a set of parameters. Then $(F, A)$ is a soft mixed neutrosophic semigroup ring over $\mathbb{C}\langle G \cup I \rangle$, where

$$F(a_1) = \mathbb{R}\langle G \cup I \rangle,$$

$$F(a_2) = \mathbb{Q}\langle G \cup I \rangle,$$

$$F(a_3) = \mathbb{Z}\langle G \cup I \rangle,$$

$$F(a_4) = \mathbb{C}\langle H_1 \cup I \rangle,$$

$$F(a_5) = \mathbb{C}\langle H_2 \cup I \rangle.$$

**Theorem 5.1.6.** Let $F, A$ and $(H, A)$ be two soft mixed neutrosophic semigroup rings over $R\langle S \cup I \rangle$. Then their intersection $F, A \cap H, A$ is again a soft mixed neutrosophic semigroup ring over $R\langle S \cup I \rangle$.

The proof is straightforward, so left as an exercise for the readers.

**Theorem 5.1.7.** Let $F, A$ and $H, B$ be two soft mixed neutrosophic semigroup rings over $R\langle S \cup I \rangle$. If $A \cap B = \phi$, then $F, A \cup H, B$ is a soft mixed neutrosophic semigroup ring over $R\langle S \cup I \rangle$.

**Proof.** This is straightforward.

**Remark 5.1.3.** Let $R\langle S \cup I \rangle$ be a neutrosophic semigroup ring. Then $R\langle S \cup I \rangle$ can have soft neutrosophic semigroup ring, soft neutrosophic semigroup subring and soft mixed neutrosophic semigroup ring over $R\langle S \cup I \rangle$.

**Proof:** It is obvious.

**Proposition 5.1.3.** Let $(F, A)$ and $(K, B)$ be two soft mixed neutosophic semigroup rings over $R\langle S \cup I \rangle$. Then

1. Their extended intersection $(F, A) \cap_E (K, B)$ is a soft mixed neutrosophic semigroup ring over $R\langle S \cup I \rangle$.
2. Their restricted intersection $(F, A) \cap_R (K, B)$ is a soft mixed neutrosophic semigroup ring over $R\langle S \cup I \rangle$.
3. Their $AND$ operation $(F, A) \vee (K, B)$ is a soft mixed neutrosophic semigroup ring over $R\langle S \cup I \rangle$.

**Proof.** This is straightforward.

**Remark 5.1.4.** Let $(F, A)$ and $(K, B)$ be two soft mixed neutosophic semigroup rings over $R\langle S \cup I \rangle$. Then

1. Their restricted union $(F, A) \cup_R (K, B)$ is not a soft mixed neutrosophic semigroup ring over $R\langle S \cup I \rangle$.
2. Their extended union $(F, A) \cup_E (K, B)$ is not a soft mixed neutrosophic semigroup ring over $R\langle S \cup I \rangle$.
3. Their $OR$ operation $(F, A) \vee (K, B)$ is a soft mixed neutrosophic semigroup ring over $R\langle S \cup I \rangle$.

One can easily show these remarks by the help of examples.

**Definition 5.1.4.** Let $R\langle S \cup I \rangle$ be a neutrosophic semigroup ring. Let $A$ be a set of parameters and $(F, A)$ be a soft set over $R\langle S \cup I \rangle$. Then $(F, A)$ is called soft neutosophic subring if and only if $F(a)$ is neutrosophic subring of $R\langle S \cup I \rangle$ for all $a \in A$.

**Example 5.1.4.** Let $R\langle S \cup I \rangle = \mathbb{C}\langle Z^+ \cup \{0\} \cup \{I\} \rangle$ be a neutrosophic semigroup ring, where $\mathbb{C} = $ field of complex numbers and $\langle Z^+ \cup \{0\} \cup \{I\} \rangle$ be a neutrosophic semigroup under $+$. Let $A = \{a_1, a_2, a_3\}$ be a set of parameters. Then $(F, A)$ is a soft neutrosophic subring over $\mathbb{C}\langle S \cup I \rangle$, where

$$F(a_1) = \langle \mathbb{R} \cup I \rangle,$$

$$F(a_2) = \langle \mathbb{Q} \cup I \rangle,$$

$$F(a_3) = \langle \mathbb{Z} \cup I \rangle.$$

**Theorem 5.1.8.** Let $F,A$ and $(H,A)$ be two soft neutrosophic subrings over $R\langle S \cup I\rangle$. Then their intersection $F,A \cap H,A$ is again a soft neutrosophic subring over $R\langle S \cup I\rangle$.

**Proof.** The proof is straightforward.

**Theorem 5.1.9.** Let $F,A$ and $H,B$ be two soft neutrosophic subrings over $R\langle S \cup I\rangle$. If $A \cap B = \phi$, then $F,A \cup H,B$ is also a soft neutrosophic subring over $R\langle S \cup I\rangle$.

**Proof.** This is straightforward.

**Proposition 5.1.4.** Let $(F,A)$ and $(K,B)$ be two soft neutosophic subrings over $R\langle S \cup I\rangle$. Then

1. Their extended intersection $(F,A) \cap_E (K,B)$ is a soft neutrosophic subring over $R\langle S \cup I\rangle$.
2. Their restricted intersection $(F,A) \cap_R (K,B)$ is a soft neutrosophic subring over $R\langle S \cup I\rangle$.
3. Their $AND$ operation $(F,A) \vee (K,B)$ is a soft neutrosophic subring over $R\langle S \cup I\rangle$.

**Proof.** This is straightforward.

**Remark 5.1.5.** Let $(F,A)$ and $(K,B)$ be two soft neutosophic subrings over $R\langle S \cup I\rangle$. Then

1. Their restricted union $(F,A) \cup_R (K,B)$ is not a soft neutrosophic subring over $R\langle S \cup I \rangle$.

2. Their extended union $(F,A) \cup_E (K,B)$ is not a soft neutrosophic subring over $R\langle S \cup I \rangle$.

3. Their *OR* operation $(F,A) \vee (K,B)$ is a soft neutrosophic subring over $R\langle S \cup I \rangle$.

One can easily show these remarks by the help of examples.

**Definition 5.1.5.** Let $R\langle S \cup I \rangle$ be a neutrosophic semigroup ring. Let $A$ be a set of parameters and $(F,A)$ be a soft set over $R\langle S \cup I \rangle$. Then $(F,A)$ is called soft pseudo neutrosophic subring if and only if $F(a)$ is a pseudo neutrosophic subring of $R\langle S \cup I \rangle$ for all $a \in A$ .

**Theorem 5.1.10.** Let $F,A$ and $(H,A)$ be two soft pseudo neutrosophic subrings over $R\langle S \cup I \rangle$. Then their intersection $F,A \cap H,A$ is again a soft pseudo neutrosophic subring over $R\langle S \cup I \rangle$ .

**Proof.** The proof is straightforward.

**Theorem 5.1.11.** Let $F,A$ and $H,B$ be two soft pseudo neutrosophic subrings over $R\langle S \cup I \rangle$. If $A \cap B = \phi$, then $F,A \cup H,B$ is a soft pseudo neutrosophic subring over $R\langle S \cup I \rangle$.

**Proof.** This is straightforward.

**Proposition 5.1.5.** Let $(F,A)$ and $(K,B)$ be two soft pseudo neutosophic subrings over $R\langle S \cup I \rangle$. Then

1. Their extended intersection $(F,A) \cap_E (K,B)$ is a soft pseudo neutrosophic subring over $R\langle S \cup I \rangle$.
2. Their restricted intersection $(F,A) \cap_R (K,B)$ is a soft pseudo neutrosophic subring over $R\langle S \cup I \rangle$.
3. Their $AND$ operation $(F,A) \vee (K,B)$ is a soft pseudo neutrosophic subring over $R\langle S \cup I \rangle$.

**Proof.** This is straightforward.

**Remark 5.1.6.** Let $(F,A)$ and $(K,B)$ be two soft pseudo neutosophic subrings over $R\langle S \cup I \rangle$. Then

1. Their restricted union $(F,A) \cup_R (K,B)$ is not a soft pseudo neutrosophic subring over $R\langle S \cup I \rangle$.
2. Their extended union $(F,A) \cup_E (K,B)$ is not a soft pseudo neutrosophic subring over $R\langle S \cup I \rangle$.
3. Their $OR$ operation $(F,A) \vee (K,B)$ is a soft pseudo neutrosophic subring over $R\langle S \cup I \rangle$.

One can easily show these remarks by the help of examples.

**Definition 5.1.6.** Let $R\langle S \cup I \rangle$ be a neutrosophic semigroup ring. Let $A$ be a set of parameters and $(F,A)$ be a soft set over $R\langle S \cup I \rangle$. Then $(F,A)$ is called soft neutosophic subsemigroup ring if and only if $F(a)$ is neutrosophic subsemigroup ring of $R\langle S \cup I\rangle$ for all $a \in A$ .

**Example 5.1.5.** Let $R\langle S \cup I\rangle = \mathbb{C}\langle Z^{+} \cup \{0\} \cup \{I\}\rangle$ be a neutrosophic semigroup ring, where $\mathbb{C} =$ field of complex numbers and $\langle S \cup I\rangle = \langle Z^{+} \cup \{0\} \cup \{I\}\rangle$. Let $H_1 = 2Z^{+} \cup \{0\}$ and $H_2 = 5Z^{+} \cup \{0\}$ be subsemigroups of the neutrosophic semigroup $\langle S \cup I\rangle$.

Let $A = \{a_1, a_2\}$ be a set of parameters. Then $(F,A)$ is a soft neutrosophic subsemigroup ring over $\mathbb{C}\langle S \cup I\rangle$, where

$$F(a_1) = \mathbb{R}H_1,$$

$$F(a_2) = \mathbb{Q}H_2 .$$

**Theorem 5.1.12.** Let $F, A$ and $(H,A)$ be two soft neutrosophic subsemigroup rings over $R\langle S \cup I\rangle$. Then their intersection $F, A \cap H, A$ is again a soft neutrosophic subsemigroup ring over $R\langle S \cup I\rangle$ .

**Proof.** The proof is straightforward.

**Theorem 5.1.13.** Let $F, A$ and $H, B$ be two soft neutrosophic subsemigroup rings over $R\langle S \cup I\rangle$. If $A \cap B = \phi$, then $F, A \cup H, B$ is a soft neutrosophic subsemigroup ring over $R\langle S \cup I\rangle$.

**Proof.** This is straightforward.

**Proposition 5.1.6.** Let $(F, A)$ and $(K, B)$ be two soft neutosophic subsemigroup rings over $R\langle S \cup I \rangle$. Then

1. Their extended intersection $(F, A) \cap_E (K, B)$ is a soft neutrosophic subsemigroup ring over $R\langle S \cup I \rangle$.
2. Their restricted intersection $(F, A) \cap_R (K, B)$ is a soft neutrosophic subsemigroup ring over $R\langle S \cup I \rangle$.
3. Their *AND* operation $(F, A) \vee (K, B)$ is a soft neutrosophic subsemigroup ring over $R\langle S \cup I \rangle$.

**Proof.** This is straightforward.

**Remark 5.1.7.** Let $(F, A)$ and $(K, B)$ be two soft neutosophic subsemigroup rings over $R\langle S \cup I \rangle$. Then

1. Their restricted union $(F, A) \cup_R (K, B)$ is not a soft neutrosophic subsemigroup ring over $R\langle S \cup I \rangle$.
2. Their extended union $(F, A) \cup_E (K, B)$ is not a soft neutrosophic subsemigroup ring over $R\langle S \cup I \rangle$.
3. Their *OR* operation $(F, A) \vee (K, B)$ is a soft neutrosophic subsemigroup ring over $R\langle S \cup I \rangle$.

One can easily show these remarks by the help of examples.

**Definition 5.1.7.** Let $R\langle S \cup I \rangle$ be a neutrosophic semigroup ring. Let $A$ be a set of parameters and $(F, A)$ be a soft set over $R\langle S \cup I \rangle$. Then $(F, A)$ is called a soft subring over $R\langle S \cup I \rangle$ if and only if $F(a)$ is a subring of $R\langle S \cup I \rangle$ for all $a \in A$.

**Example 5.1.6.** Let $\mathbb{Z}_2\langle S \cup I \rangle$ be a neutrosophic semigroup ring, where $\langle S \cup I \rangle = \{1, g, g^2, g^3, I, gI, g^2I, g^3I : g^4 = 1, I^2 = I\}$.

Let $A = \{a_1, a_2\}$ be a set of parameters. Then $(F, A)$ is a soft subring over $\mathbb{Z}_2\langle S \cup I \rangle$, where

$$F(a_1) = \{0, 1 + g^2\},$$

$$F(a_2) = \{0, 1 + g, g + g^3, 1 + g^3\}.$$

**Theorem 5.1.15.** Let $F, A$ and $(H, A)$ be two soft subrings over $R\langle S \cup I \rangle$. Then their intersection $F, A \cap H, A$ is again a soft subring over $R\langle S \cup I \rangle$.

The proof is straightforward, so left as an exercise for the readers.

**Theorem 5.1.16.** Let $F, A$ and $H, B$ be two soft subrings over $R\langle S \cup I \rangle$. If $A \cap B = \phi$, then $F, A \cup H, B$ is a soft subring over $R\langle S \cup I \rangle$.

**Proof.** This is straightforward.

**Proposition 5.1.7.** Let $(F, A)$ and $(K, B)$ be two soft subrings over $R\langle S \cup I \rangle$. Then

1. Their extended intersection $(F, A) \cap_E (K, B)$ is a soft subring over

$R \langle S \cup I \rangle$.

2. Their restricted intersection $(F, A) \cap_R (K, B)$ is a soft subring over $R \langle S \cup I \rangle$.

3. Their *AND* operation $(F, A) \vee (K, B)$ is a soft subring over $R \langle S \cup I \rangle$.

**Proof.** These are straightforward.

**Remark 5.1.8.** Let $(F, A)$ and $(K, B)$ be two soft subrings over $R \langle S \cup I \rangle$. Then

1. Their restricted union $(F, A) \cup_R (K, B)$ is not a soft subring over $R \langle S \cup I \rangle$.

2. Their extended union $(F, A) \cup_E (K, B)$ is not a soft subring over $R \langle S \cup I \rangle$.

3. Their *OR* operation $(F, A) \vee (K, B)$ is a soft subring over $R \langle S \cup I \rangle$.

One can easily show these remarks by the help of examples.

**Definition 5.1.8.** Let $R \langle S \cup I \rangle$ be a neutrosophic semigroup ring. Then $(F, A)$ is called an absolute soft neutosophic semigroup ring if $F(a) = \mathrm{R} \langle S \cup I \rangle$, for all $a \in A$.

**Definition 5.1.9.** Let $(F, A)$ and $(K, B)$ be two soft neutrosophic semigroup rings over $R \langle S \cup I \rangle$. Then $(K, B)$ is called soft subneutrosophic semigroup ring of $(F, A)$, if

1. $B \subseteq A$, and
2. $K(a)$ is a subneutrosophic semigroup ring of $F(a)$ for all $a \in B$.

**Example 5.1.7.** Let $\mathbb{C}\langle \mathbb{Z}^+ / \{0\} \cup I \rangle$ be a neutrosophic semigroup ring, where $\mathbb{C} =$ field of complex numbers and $\langle \mathbb{Z}^+ / \{0\} \cup I \rangle$.

Let $A = \{a_1, a_2, a_3, a_4, a_5, a_6, a_7\}$ be a set of parameters. Then $(F, A)$ is a soft neutrosophic semigroup ring over $\mathbb{C}\langle \mathbb{Z}^+ / \{0\} \cup I \rangle$, where

$$F(a_1) = \mathbb{R}\langle \mathbb{Z}^+ / \{0\} \cup I \rangle,$$

$$F(a_2) = \mathbb{Q}\langle \mathbb{Z}^+ / \{0\} \cup I \rangle,$$

$$F(a_3) = \mathbb{Z}\langle \mathbb{Z}^+ / \{0\} \cup I \rangle,$$

$$F(a_4) = \mathbb{C}\langle 2\mathbb{Z}^+ / \{0\} \cup I \rangle,$$

$$F(a_5) = \mathbb{C}\langle 3\mathbb{Z}^+ / \{0\} \cup I \rangle,$$

$$F(a_6) = \mathbb{C}\langle 4\mathbb{Z}^+ / \{0\} \cup I \rangle,$$

$$F(a_7) = \mathbb{C}\langle 5\mathbb{Z}^+ / \{0\} \cup I \rangle.$$

Let $B = \{a_1, a_2, a_5, a_7\}$ be a set of parameters. Then $(H, B)$ is a soft subneutrosophic semigroup ring of $(F, A)$, where

$$H(a_1) = \mathbb{Q}\langle \mathbb{Z}^+ / \{0\} \cup I \rangle,$$

$$H(a_2) = \mathbb{Z}\langle \mathbb{Z}^+ / \{0\} \cup I \rangle,$$

$$H(a_5) = \mathbb{C}\langle 6\mathbb{Z}^+ / \{0\} \cup I \rangle,$$

$$H(a_7) = \mathbb{C}\langle 10\mathbb{Z}^+ / \{0\} \cup I \rangle.$$

**Definition 5.1.10.** Let $R\langle S \cup I\rangle$ be a neutrosophic semigroup ring. Let $A$ be a set of parameters and $(F,A)$ be a soft set over $R\langle S \cup I\rangle$. Then $(F,A)$ is called soft neutosophic ideal if and only if $F(a)$ is a neutrosophic ideal of $R\langle S \cup I\rangle$ for all $a \in A$.

Similarly one can define soft neutrosophic left ideal and soft neutrosophic right ideal over $R\langle S \cup I\rangle$.

**Example 5.1.8.** Let $R\langle S \cup I\rangle$ be a neutrosophic semigroup ring, where $R = \mathbb{Z}$. Let $A = \{a_1, a_2, a_3\}$ be a set of parameters. Then $(F,A)$ is a soft neutrosophic ideal over $R\langle S \cup I\rangle$, where

$$F(a_1) = 2\mathbb{Z}\langle S \cup I\rangle,$$

$$F(a_2) = 4\mathbb{Z}\langle S \cup I\rangle,$$

$$F(a_3) = 6\mathbb{Z}\langle S \cup I\rangle.$$

**Theorem 5.1.17.** All soft neutrosophic ideals are trivially soft neutrosophic semigroup rings but the converse is not true in general.

We can easily establish the converse by the help of example.

**Proposition 5.1.8.** Let $(F,A)$ and $(K,B)$ be two soft neutosophic ideals over $R\langle S \cup I\rangle$. Then

1. Their extended intersection $(F,A) \cap_E (K,B)$ is a soft neutrosophic ideal over $R\langle S \cup I\rangle$.
2. Their restricted intersection $(F,A) \cap_R (K,B)$ is a soft neutrosophic ideal over $R\langle S \cup I\rangle$.

3. Their *AND* operation $(F,A) \vee (K,B)$ is a soft neutrosophic ideal over $R\langle S \cup I \rangle$.

**Proof.** This is straightforward.

**Remark 5.1.9.** Let $(F,A)$ and $(K,B)$ be two soft neutosophic ideals over $R\langle S \cup I \rangle$. Then

1. Their restricted union $(F,A) \cup_R (K,B)$ is not a soft neutrosophic ideal over $R\langle S \cup I \rangle$.
2. Their extended union $(F,A) \cup_E (K,B)$ is not a soft neutrosophic ideal over $R\langle S \cup I \rangle$.
3. Their *OR* operation $(F,A) \vee (K,B)$ is a soft neutrosophic ideal over $R\langle S \cup I \rangle$.

One can easily show these remarks by the help of examples.

**Definition 5.1.11.** Let $(F,A)$ and $(K,B)$ be two soft neutrosophic semigroup rings over $R\langle S \cup I \rangle$. Then $(K,B)$ is called soft neutrosophic ideal of $(F,A)$, if

1. $B \subseteq A$, and
2. $K(a)$ is a neutrosophic ideal of $F(a)$ for all $a \in B$.

**Example 5.1.9.** Let $R\langle S \cup I \rangle$ be a neutrosophic semigroup ring, where $R = \mathbb{Z}$. Let $A = \{a_1, a_2, a_3\}$ be a set of parameters. Then $(F,A)$ is a soft neutrosophic semigroup ring over $R\langle S \cup I \rangle$, where

$$F(a_1) = 2\mathbb{Z}\langle S \cup I \rangle,$$

$$F(a_2) = 4\mathbb{Z}\langle S \cup I \rangle,$$

$$F(a_3) = 6\mathbb{Z}\langle S \cup I \rangle.$$

Let $B = \{a_1, a_2\}$ be another set of parameters. Then clearly $(H, B)$ is a soft neutrosophic ideal of $(F, A)$, where

$$H(a_1) = 8\mathbb{Z}\langle S \cup I \rangle,$$

$$H(a_2) = 12\mathbb{Z}\langle S \cup I \rangle.$$

**Theorem 5.1.18.** Every soft neutrosophic ideal of the soft neutrosophic semigroup ring over a neutrosophic semigroup ring is trivially a soft subneutrosophic semigroup ring.

In the proceeding section, the authors give the definitions and notions about soft neutrosophic bisemigroup birings over neutrosophic bisemigroup birings.

## 5.2   Soft Neutrosophic Bisemigroup Biring

In this section, the authors introduced soft neutrosophic bisemigroup birings over neutrosophic bisemigroup birings. This means that we are going towards the generalization of soft neutrosophic semigroup rings.

We also give some basic and fundamental properties and characteristics of soft neutrosophic bisemigroup birings over neutrosophic bisemigroup birings with many illustrative examples.

**Definition 5.2.1.** Let $R_B \langle S \cup I \rangle = R_1 \langle S_1 \cup I \rangle \cup R_2 \langle S_2 \cup I \rangle$ be a neutrosophic bisemigroup biring over a biring $R_B = R_1 \cup R_2$. Let $A$ be a set of parameters and let $(F, A)$ be a soft set over $R_B \langle S \cup I \rangle$. Then $(F, A)$ is called a soft neutrosophic bisemigroup biring over $R_B \langle S \cup I \rangle$ if and only if $F(a)$ is a subneutrosophic bisemigroup biring of $R_B \langle S \cup I \rangle$ for all $a \in A$.

See the following example for better understanding.

**Example 5.2.1.** Let $R_B \langle S \cup I \rangle = \mathbb{Q} \langle S_1 \cup I \rangle \cup \mathbb{R} \langle S_2 \cup I \rangle$ be a neutrosophic bisemigroup biring, where $R_B = \mathbb{Q} \cup \mathbb{R}$ and
$\langle S_1 \cup I \rangle = \{1, g, g^2, g^3, g^4, g^5, I, gI, ..., g^5 I : g^6 = 1, I^2 = I\}$ is a neutrosophic semigroup and $\langle S_2 \cup I \rangle = \langle Z^+ \setminus \{0\} \cup I \rangle$ be a neutrosophic semigroup under $+$. Let
$\langle H_1^{'} \cup I \rangle \cup \langle H_2^{'} \cup I \rangle = \{1, g^3 : g^6 = 1\} \cup \langle 2Z^+ \setminus \{0\} \cup I \rangle$,
$\langle H_1^{''} \cup I \rangle \cup \langle H_2^{''} \cup I \rangle = \{1, g^3, I, g^3 I : g^6 = 1, I^2 = I\} \cup \langle 3Z^+ \setminus \{0\} \cup I \rangle$,
$\langle H_1^{'''} \cup I \rangle \cup \langle H_2^{'''} \cup I \rangle = \{1, g^2, g^4 : g^6 = 1, I^2 = I\} \cup \langle 4Z^+ \setminus \{0\} \cup I \rangle$ and
$\langle H_1^{''''} \cup I \rangle \cup \langle H_2^{''''} \cup I \rangle = \{1, g^2, g^4, I, g^2 I, g^4 I : g^6 = 1, I^2 = I\} \cup \langle 5Z^+ \setminus \{0\} \cup I \rangle$. Let
$A = \{a_1, a_2, a_3, a_4\}$ be a set of parameters. Then $(F, A)$ is a soft neutrosophic bisemigroup biring over $R_B \langle S \cup I \rangle$, where

$$F(a_1) = \mathbb{Q} \langle H_1^{'} \cup I \rangle \cup \mathbb{R} \langle H_2^{'} \cup I \rangle,$$

$$F(a_2) = \mathbb{Q} \left\langle H^{'}_{\ 1} \cup I \right\rangle \cup \mathbb{R} \left\langle H^{'}_{\ 2} \cup I \right\rangle,$$

$$F(a_3) = \mathbb{Q} \left\langle H^{''}_{\ 1} \cup I \right\rangle \cup \mathbb{R} \left\langle H^{''}_{\ 2} \cup I \right\rangle,$$

$$F(a_4) = \mathbb{Q} \left\langle H^{'''}_{\ 1} \cup I \right\rangle \cup \mathbb{R} \left\langle H^{'''}_{\ 2} \cup I \right\rangle.$$

**Theorem 5.2.1.** Let $F, A$ and $(H,A)$ be two soft neutrosophic bisemigroup birings over $R_B \left\langle S \cup I \right\rangle$. Then their intersection $F, A \cap H, A$ is again a soft neutrosophic bisemigroup biring over $R_B \left\langle S \cup I \right\rangle$.

The proof is straightforward, so we suggest the readers to prove it by yourself as an exercise.

**Theorem 5.2.2.** Let $F, A$ and $H, B$ be two soft neutrosophic bisemigroup birings over $R_B \left\langle S \cup I \right\rangle$. If $A \cap B = \phi$, then $F, A \cup H, B$ is a soft neutrosophic bisemigroup biring over $R_B \left\langle S \cup I \right\rangle$.

**Proof.** This is straightforward.

**Proposition 5.2.1.** Let $(F, A)$ and $(K, B)$ be two soft neutosophic bisemigroup birings over a neutrosophic bisemigroup biring $R_B \left\langle S \cup I \right\rangle$. Then

1. Their extended intersection $(F, A) \cap_E (K, B)$ is again a soft neutrosophic bisemigroup biring over $R_B \left\langle S \cup I \right\rangle$.
2. Their restricted intersection $(F, A) \cap_R (K, B)$ is again a soft neutrosophic bisemigroup biring over $R_B \left\langle G \cup I \right\rangle$.
3. Their $AND$ operation $(F, A) \vee (K, B)$ is again a soft neutrosophic

bisemigroup biring over $R_B \langle G \cup I \rangle$.

**Proof.** This is straightforward.

**Remark 5.2.1.** Let $(F, A)$ and $(K, B)$ be two soft neutosophic bisemigroup birings over $R_B \langle S \cup I \rangle$. Then

1. Their restricted union $(F, A) \cup_R (K, B)$ is not a soft neutrosophic bisemigroup biring over $R_B \langle S \cup I \rangle$.
2. Their extended union $(F, A) \cup_E (K, B)$ is not a soft neutrosophic bisemigroup biring over $R_B \langle S \cup I \rangle$.
3. Their $OR$ operation $(F, A) \vee (K, B)$ is not a soft neutrosophic bisemigroup biring over $R_B \langle S \cup I \rangle$.

One can establish these remarks by the help of following examples.

**Definition 5.2.2.** Let $R_B \langle S \cup I \rangle = R_1 \langle S_1 \cup I \rangle \cup R_2 \langle S_2 \cup I \rangle$ be a neutrosophic bisemigroup biring. Let $A$ be a set of parameters and $(F, A)$ be a soft set over $R_B \langle S \cup I \rangle$. Then $(F, A)$ is called soft neutosophic bisemigroup subbiring over $R_B \langle S \cup I \rangle$ if and only if $F(a)$ is subneutrosophic bisemigroup subbiring of $R_B \langle S \cup I \rangle$ for all $a \in A$.

**Example 5.2.2.** Let $R_B \langle S \cup I \rangle = \mathbb{R} \langle S_1 \cup I \rangle \cup \mathbb{C} \langle S_2 \cup I \rangle$ be a neutrosophic bisemigroup biring, where $\langle S_1 \cup I \rangle = \{1, g, g^2, g^3, g^4, g^5, I, gI, ..., g^5 I : g^6 = 1, I^2 = I\}$ and $\langle S_2 \cup I \rangle = \{1, g, g^2, g^3, I, gI, g^2 I, g^3 I : g^4 = 1, I^2 = I\}$ are neutrosophic semigroups. Let $A = \{a_1, a_2\}$ be a set of parameters. Then $(F, A)$ is a soft neutrosophic bisemigroup subbiring over $R_B \langle S \cup I \rangle$, where

$$F(a_1) = \mathbb{Q}\langle S_1 \cup I \rangle \cup \mathbb{R}\langle S_2 \cup I \rangle,$$

$$F(a_2) = \mathbb{Z}\langle S_1 \cup I \rangle \cup \mathbb{Q}\langle S_2 \cup I \rangle.$$

**Theorem 5.2.3.** Let $F,A$ and $(H,A)$ be two soft neutrosophic bisemigroup subbirings over $R_B\langle S \cup I \rangle$. Then their intersection $F,A \cap H,A$ is again a soft neutrosophic bisemigroup subbirings over $R_B\langle S \cup I \rangle$.

The proof is straightforward, so we suggest the readers to prove it by yourself as an exercise.

**Theorem 5.2.4.** Let $F,A$ and $H,B$ be two soft neutrosophic bisemigroup subbirings over $R_B\langle S \cup I \rangle$. If $A \cap B = \phi$, then $F,A \cup H,B$ is a soft neutrosophic bisemigroup subbiring over $R_B\langle S \cup I \rangle$.

**Proof.** This is straightforward.

**Remark 5.2.2.** Let $(F,A)$ and $(K,B)$ be two soft neutosophic bisemigroup subbirings over $R_B\langle S \cup I \rangle$. Then

1. Their extended union $(F,A) \cup_E (K,B)$ is not a soft neutrosophic bisemigroup subbiring over $R_B\langle S \cup I \rangle$.
2. Their restricted union $(F,A) \cup_R (K,B)$ is not a soft neutrosophic bisemigroup subbiring over $R_B\langle S \cup I \rangle$.
3. Their $OR$ operation $(F,A) \vee (K,B)$ is not a soft neutrosophic bisemigroup subbiring over $R_B\langle S \cup I \rangle$.

One can easily show these remaks by the help of examples.

**Proposition 5.2.2.** Let $(F, A)$ and $(K, B)$ be two soft neutosophic bisemigroup subbirings over $R_B \langle S \cup I \rangle$. Then

1. Their restricted intersection $(F, A) \cap_R (K, B)$ is a soft neutrosophic bisemigroup subbring over $R_B \langle S \cup I \rangle$.
2. Their extended intersection $(F, A) \cap_E (K, B)$ is a soft neutrosophic bisemigroup subbiring over $R_B \langle S \cup I \rangle$.
3. Their *AND* operation $(F, A) \vee (K, B)$ is a soft neutrosophic bisemigroup subbring over $R_B \langle S \cup I \rangle$.

**Proof.** These are straightforward.

**Definition 5.2.3.** Let $R_B \langle S \cup I \rangle = R_1 \langle S_1 \cup I \rangle \cup R_2 \langle S_2 \cup I \rangle$ be a neutrosophic bisemigroup biring. Let $A$ be a set of parameters and $(F, A)$ be a soft set over $R_B \langle S \cup I \rangle$. Then $(F, A)$ is called soft mixed neutrosophic bisemigroup biring over $R_B \langle S \cup I \rangle$ if for some $a \in A$, $F(a)$ is subneutrosophic bisemigroup subbiring and for the remaining $a \in A$, $F(a)$ is a subneutrosophic bisemigroup biring of $R_B \langle S \cup I \rangle$.

**Example 5.2.3.** Let $R_B \langle S \cup I \rangle = \mathbb{Z} \langle S_1 \cup I \rangle \cup \mathbb{C} \langle S_2 \cup I \rangle$ be a neutrosophic bisemigroup biring, where $\langle S_1 \cup I \rangle = \{1, g, g^2, g^3, g^4, g^5, I, gI, ..., g^5 I : g^6 = 1, I^2 = I\}$ and $\langle S_2 \cup I \rangle = \{1, g, g^2, g^3, I, gI, g^2 I, g^{53} I : g^4 = 1, I^2 = I\}$. Let $A = \{a_1, a_2, a_3, a_4\}$ be a set of parameters. Then $(F, A)$ is a soft mixed neutrosophic bisemigroup biring over $R_B \langle S \cup I \rangle$, where

$$F(a_1) = 2\mathbb{Z}\langle S_1 \cup I \rangle \cup \mathbb{R}\langle S_2 \cup I \rangle,$$

$$F(a_2) = 4\mathbb{Z}\langle S_1 \cup I \rangle \cup \mathbb{Q}\langle S_2 \cup I \rangle,$$

$$F(a_3) = \mathbb{Z}\langle H'_1 \cup I \rangle \mathbb{C}\langle H'_2 \cup I \rangle,$$

$$F(a_4) = \mathbb{Z}\langle H^{''}_1 \cup I \rangle \mathbb{C}\langle H^{''}_2 \cup I \rangle.$$

**Where** $\langle H'_1 \cup I \rangle = \{1, g^2 : g^4 = 1\}$, $\langle H'_2 \cup I \rangle = \{1, g^3, I, g^3I : g^6 = 1, I^2 = I\}$, $\langle H^{''}_1 \cup I \rangle = \{1, g^2, I, g^2I : g^4 = 1, I^2 = I\}$ **and** $\langle H^{''}_2 \cup I \rangle = \{1, g^2, g^4, I, g^2I, g^4I : g^6 = 1, I^2 = I\}$.

**Theorem 5.2.5.** Let $\langle F, A \rangle$ and $(H, A)$ be two soft mixed neutrosophic bisemigroup birings over $R_B \langle S \cup I \rangle$. Then their intersection $\langle F, A \rangle \cap \langle H, A \rangle$ is again a soft mixed neutrosophic bisemigroup biring over $R_B \langle S \cup I \rangle$.

**Proof.** The proof is straightforward.

**Theorem 5.2.6.** Let $\langle F, A \rangle$ and $\langle H, B \rangle$ be two soft mixed neutrosophic bisemigroup birings over $R_B \langle S \cup I \rangle$. If $A \cap B = \phi$, then $\langle F, A \rangle \cup \langle H, B \rangle$ is a soft mixed neutrosophic bisemigroup biring over $R_B \langle S \cup I \rangle$.

**Remark 5.2.3.** Let $R_B \langle S \cup I \rangle$ be a neutrosophic bisemigroup biring. Then $R_B \langle S \cup I \rangle$ can have soft neutrosophic bisemigroup biring, soft neutrosophic bisemigroup subbiring and soft mixed neutrosophic bisemigroup biring over $R_B \langle S \cup I \rangle$.

It is obvious, so left as an exercise for the readers.

**Proposition 5.2.3.** Let $(F, A)$ and $(K, B)$ be two soft mixed neutosophic bisemigroup birings over $R_B \langle S \cup I \rangle$. Then

1. Their extended intersection $(F, A) \cap_E (K, B)$ is a soft mixed neutrosophic bisemigroup biring over $R_B \langle S \cup I \rangle$.

2. Their restricted intersection $(F, A) \cap_R (K, B)$ is a soft mixed neutrosophic bisemigroup biring over $R_B \langle S \cup I \rangle$.

3. Their *AND* operation $(F, A) \vee (K, B)$ is a soft mixed neutrosophic bisemigroup biring over $R_B \langle S \cup I \rangle$.

**Proof.** This is straightforward.

**Remark 5.2.4.** Let $(F, A)$ and $(K, B)$ be two soft mixed neutosophic bisemigroup birings over $R_B \langle S \cup I \rangle$. Then

1. Their restricted union $(F, A) \cup_R (K, B)$ is not a soft mixed neutrosophic bisemigroup biring over $R_B \langle S \cup I \rangle$.

2. Their extended union $(F, A) \cup_E (K, B)$ is not a soft mixed neutrosophic bisemigroup biring over $R_B \langle S \cup I \rangle$.

3. Their *OR* operation $(F, A) \vee (K, B)$ is a soft mixed neutrosophic bisemigroup biring over $R_B \langle S \cup I \rangle$.

One can easily show these remarks by the help of examples.

**Definition 5.2.4.** Let $R_B \langle S \cup I \rangle = R_1 \langle S_1 \cup I \rangle \cup R_2 \langle S_2 \cup I \rangle$ be a neutrosophic bisemigroup biring. Let $A$ be a set of parameters and $(F, A)$ be a soft set over $R_B \langle S \cup I \rangle$. Then $(F, A)$ is called soft neutosophic subbiring if and only if $F(a)$ is neutrosophic subbiring of $R_B \langle S \cup I \rangle$ for all $a \in A$.

**Example 5.2.4.** Let $R_B \langle S \cup I \rangle = \mathbb{Z} \langle S_1 \cup I \rangle \cup \mathbb{C} \langle S_1 \cup I \rangle$ be a neutrosophic bisemigroup biring. Let $A = \{a_1, a_2, a_3\}$ be a set of parameters. Then $(F, A)$ is a soft neutrosophic subbiring over $R_B \langle S \cup I \rangle$, where

$$F(a_1) = \langle 2\mathbb{Z} \cup I \rangle \cup \langle \mathbb{R} \cup I \rangle,$$

$$F(a_2) = \langle 4\mathbb{Z} \cup I \rangle \cup \langle \mathbb{Q} \cup I \rangle,$$

$$F(a_3) = \langle 6\mathbb{Z} \cup I \rangle \cup \langle \mathbb{Z} \cup I \rangle.$$

**Theorem 5.2.7.** Let $F, A$ and $(H, A)$ be two soft neutrosophic subbirings over $R_B \langle S \cup I \rangle$. Then their intersection $F, A \cap H, A$ is again a soft neutrosophic subbiring over $R_B \langle S \cup I \rangle$.

**Proof.** The proof is straightforward.

**Theorem 5.2.8.** Let $F, A$ and $H, B$ be two soft neutrosophic subbirings over $R_B \langle S \cup I \rangle$. If $A \cap B = \phi$, then $F, A \cup H, B$ is also a soft neutrosophic subbiring over $R_B \langle S \cup I \rangle$.

**Proof.** This is straightforward.

**Proposition 5.2.4.** Let $(F, A)$ and $(K, B)$ be two soft neutosophic subbirings over $R_B \langle S \cup I \rangle$. Then

1. Their extended intersection $(F,A) \cap_E (K,B)$ is a soft neutrosophic subbiring over $R_B \langle S \cup I \rangle$.

2. Their restricted intersection $(F,A) \cap_R (K,B)$ is a soft neutrosophic subbiring over $R_B \langle S \cup I \rangle$.

3. Their *AND* operation $(F,A) \vee (K,B)$ is a soft neutrosophic subbiring over $R_B \langle S \cup I \rangle$.

**Proof.** This is straightforward.

**Remark 5.2.5.** Let $(F,A)$ and $(K,B)$ be two soft neutosophic subbirings over $R_B \langle S \cup I \rangle$. Then

1. Their restricted union $(F,A) \cup_R (K,B)$ is not a soft neutrosophic subbiring over $R_B \langle S \cup I \rangle$.

2. Their extended union $(F,A) \cup_E (K,B)$ is not a soft neutrosophic subbiring over $R_B \langle S \cup I \rangle$.

3. Their *OR* operation $(F,A) \vee (K,B)$ is a soft neutrosophic subbiring over $R_B \langle S \cup I \rangle$.

One can easily show these remarks by the help of examples.

**Definition 5.2.5.** Let $R_B \langle S \cup I \rangle = R_1 \langle S_1 \cup I \rangle \cup R_2 \langle S_2 \cup I \rangle$ be a neutrosophic bisemigroup biring. Let $A$ be a set of parameters and $(F,A)$ be a soft set over $R_B \langle S \cup I \rangle$. Then $(F,A)$ is called soft pseudo neutosophic subbiring if and only if $F(a)$ is a pseudo neutrosophic subbiring of $R_B \langle S \cup I \rangle$ for all $a \in A$.

**Theorem 5.2.9.** Let $(F, A)$ and $(H, A)$ be two soft pseudo neutrosophic subbirings over $R_B \langle S \cup I \rangle$. Then their intersection $(F, A) \cap (H, A)$ is again a soft pseudo neutrosophic subbiring over $R_B \langle S \cup I \rangle$.

**Proof.** The proof is straightforward.

**Theorem 5.2.10.** Let $(F, A)$ and $(H, B)$ be two soft pseudo neutrosophic subbirings over $R_B \langle S \cup I \rangle$. If $A \cap B = \phi$, then $(F, A) \cup (H, B)$ is a soft pseudo neutrosophic subbiring over $R_B \langle S \cup I \rangle$.

**Proof.** This is straightforward.

**Proposition 5.2.5.** Let $(F, A)$ and $(K, B)$ be two soft pseudo neutosophic subbirings over $R_B \langle S \cup I \rangle$. Then

1. Their extended intersection $(F, A) \cap_E (K, B)$ is a soft pseudo neutrosophic subbiring over $R_B \langle S \cup I \rangle$.
2. Their restricted intersection $(F, A) \cap_R (K, B)$ is a soft pseudo neutrosophic subbiring over $R_b \langle S \cup I \rangle$.
3. Their $AND$ operation $(F, A) \vee (K, B)$ is a soft pseudo neutrosophic subbiring over $R_B \langle S \cup I \rangle$.

**Proof.** This is straightforward.

**Remark 5.2.6.** Let $(F, A)$ and $(K, B)$ be two soft pseudo neutosophic subbirings over $R_B \langle S \cup I \rangle$. Then

1. Their restricted union $(F, A) \cup_R (K, B)$ is not a soft pseudo

neutrosophic subbiring over $R_B \langle S \cup I \rangle$.

2. Their extended union $(F, A) \cup_E (K, B)$ is not a soft pseudo neutrosophic subbbiring over $R_B \langle S \cup I \rangle$.

3. Their $OR$ operation $(F, A) \vee (K, B)$ is a soft pseudo neutrosophic subbiring over $R_B \langle S \cup I \rangle$.

One can easily show these remarks by the help of examples.

**Definition 5.2.6.** Let $R_B \langle S \cup I \rangle = R_1 \langle S_1 \cup I \rangle \cup R_2 \langle S_2 \cup I \rangle$ be a neutrosophic bisemigroup biring. Let $A$ be a set of parameters and $(F, A)$ be a soft set over $R_B \langle S \cup I \rangle$. Then $(F, A)$ is called soft neutosophic subbisemigroup biring if and only if $F(a)$ is neutrosophic subbisemigroup biring of $R_B \langle S \cup I \rangle$ for all $a \in A$.

**Example 5.2.5.** Let $R_B \langle S \cup I \rangle = \mathbb{Z} \langle S_1 \cup I \rangle \cup \mathbb{C} \langle S_2 \cup I \rangle$ be a neutrosophic bisemigroup biring, where $\langle S_1 \cup I \rangle = \{1, g, g^2, g^3 I, gI, g^2 I, g^3 I : g^6 = 1, I^2 = I\}$ and $\langle S_2 \cup I \rangle = \{1, g, g^2, g^3, g^4, g^5, I, gI, ..., g^5 I : g^6 = 1, I^2 = I\}$. Let $H_1^{'} = \{1, g^2 : g^4 = 1\}$, $H_2^{'} = \{1, g^2, g^4 : g^6 = 1\}$, $H_1^{''} = \{1, g^2, I, g^2 I : g^4 = 1, I^2 = I\}$ and $H_2^{''} = \{1, g^3 : g^6 = 1\}$ are neutrosophic subsemigroups respectively. Let $A = \{a_1, a_2\}$ be a set of parameters. Then $(F, A)$ is a soft neutrosophic subbisemigroup biring over $R_B \langle S \cup I \rangle$, where

$$F(a_1) = 2\mathbb{Z} H_1^{'} \cup \mathbb{R} H_2^{'},$$

$$F(a_1) = 8\mathbb{Z} H_1^{''} \cup \mathbb{Q} H_2^{''}.$$

**Theorem 5.2.11.** Let $F, A$ and $(H, A)$ be two soft neutrosophic subbisemigroup birings over $R_B \langle S \cup I \rangle$. Then their intersection $F, A \cap H, A$ is again a soft neutrosophic subbisemigroup biring over $R_B \langle S \cup I \rangle$.

**Proof.** The proof is straightforward.

**Theorem5.2.12.** Let $F, A$ and $H, B$ be two soft neutrosophic subbisemigroup birings over $R_B \langle S \cup I \rangle$. If $A \cap B = \phi$, then $F, A \cup H, B$ is a soft neutrosophic subbisemigroup biring over $R_B \langle S \cup I \rangle$.

**Proof.** This is straightforward.

**Proposition 5.2.6.** Let $(F, A)$ and $(K, B)$ be two soft neutosophic subbisemigroup birings over $R_B \langle S \cup I \rangle$. Then

1. Their extended intersection $(F, A) \cap_E (K, B)$ is a soft neutrosophic subbisemigroup biring over $R_B \langle S \cup I \rangle$.
2. Their restricted intersection $(F, A) \cap_R (K, B)$ is a soft neutrosophic subbisemigroup biring over $R_B \langle S \cup I \rangle$.
3. Their $AND$ operation $(F, A) \vee (K, B)$ is a soft neutrosophic subbisemigroup biring over $R_B \langle S \cup I \rangle$.

**Proof.** This is straightforward.

**Remark 5.2.7.** Let $(F, A)$ and $(K, B)$ be two soft neutosophic subbisemigroup birings over $R_B \langle S \cup I \rangle$. Then

1. Their restricted union $(F,A) \cup_R (K,B)$ is not a soft neutrosophic subbisemigroup biring over $R_B \langle S \cup I \rangle$.

2. Their extended union $(F,A) \cup_E (K,B)$ is not a soft neutrosophic subbisemigroup biring over $R_B \langle S \cup I \rangle$.

3. Their $OR$ operation $(F,A) \vee (K,B)$ is a soft neutrosophic subbisemigroup biring over $R_B \langle S \cup I \rangle$.

One can easily show these remarks by the help of examples.

**Definition 5.2.7.** Let $R_B \langle S \cup I \rangle = R_1 \langle S_1 \cup I \rangle \cup R_2 \langle S_2 \cup I \rangle$ be a neutrosophic bisemigroup biring. Let $A$ be a set of parameters and $(F,A)$ be a soft set over $R_B \langle S \cup I \rangle$. Then $(F,A)$ is called a soft subbiring over $R_B \langle S \cup I \rangle$ if and only if $F(a)$ is a subbiring of $R_B \langle S \cup I \rangle$ for all $a \in A$.

**Theorem 5.2.13.** Let $F,A$ and $(H,A)$ be two soft subbirings over $R_B \langle S \cup I \rangle$. Then their intersection $F,A \cap H,A$ is again a soft subbiring over $R_B \langle S \cup I \rangle$.

**Proof.** The proof is straightforward.

**Theorem 5.2.14.** Let $F,A$ and $H,B$ be two soft subbirings over $R_B \langle S \cup I \rangle$. If $A \cap B = \phi$, then $F,A \cup H,B$ is a soft subbiring over $R_B \langle S \cup I \rangle$.

**Proof.** This is straightforward.

**Proposition 5.2.7.** Let $(F,A)$ and $(K,B)$ be two soft subbirings over $R_B \langle S \cup I \rangle$. Then

1. Their extended intersection $(F,A) \cap_E (K,B)$ is a soft subbiring over $R_B \langle S \cup I \rangle$.

2. Their restricted intersection $(F,A) \cap_R (K,B)$ is a soft subbiring over $R_B \langle S \cup I \rangle$.

3. Their $AND$ operation $(F,A) \vee (K,B)$ is a soft subbiring over $R_B \langle S \cup I \rangle$.

**Proof.** These are straightforward.

**Remark 5.2.8.** Let $(F,A)$ and $(K,B)$ be two soft subbirings over $R_B \langle S \cup I \rangle$. Then

1. Their restricted union $(F,A) \cup_R (K,B)$ is not a soft subbiring over $R_B \langle S \cup I \rangle$.

2. Their extended union $(F,A) \cup_E (K,B)$ is not a soft subbiring over $R_B \langle S \cup I \rangle$.

3. Their $OR$ operation $(F,A) \vee (K,B)$ is a soft subbiring over $R_B \langle S \cup I \rangle$.

One can easily show these remarks by the help of examples.

**Definition 5.2.8.** Let $R_B \langle S \cup I \rangle = R_1 \langle S_1 \cup I \rangle \cup R_2 \langle S_2 \cup I \rangle$ be a neutrosophic bigroup biring. Then $(F,A)$ is called an absolute soft neutosophic bigroup biring if $F(a) = R_B \langle G \cup I \rangle$, for all $a \in A$.

**Definition 5.2.9.** Let $(F,A)$ and $(K,B)$ be two soft neutrosophic bisemigroup birings over $R_B \langle S \cup I \rangle$. Then $(K,B)$ is called soft subneutrosophic bisemigroup biring of $(F,A)$, if

1. $B \subseteq A$, and
2. $K(a)$ is a subneutrosophic bisemigroup biring of $F(a)$ for all $a \in B$.

**Definition 5.2.10.** Let $R_B \langle S \cup I \rangle$ be a neutrosophic bisemigroup biring. Let $A$ be a set of parameters and $(F, A)$ be a soft set over $R_B \langle S \cup I \rangle$. Then $(F, A)$ is called soft neutrosophic biideal if and only if $F(a)$ is a neutrosophic biideal of $R_B \langle S \cup I \rangle$ for all $a \in A$.

Similarly one can define soft neutrosophic left biideal and soft neutrosophic right biideal over $R_B \langle S \cup I \rangle$.

**Example 5.2.6.** Let $R_B \langle S \cup I \rangle = \mathbb{Z} \langle S_1 \cup I \rangle \cup \mathbb{Z} \langle S_2 \cup I \rangle$ be a neutrosophic bisemigroup biring. Let $A = \{a_1, a_2, a_3\}$ be a set of parameters. Then $(F, A)$ is a soft neutrosophic biideal over $R_B \langle S \cup I \rangle$, where

$$F(a_1) = 2\mathbb{Z} \langle S_1 \cup I \rangle \cup 4\mathbb{Z} \langle S_2 \cup I \rangle,$$

$$F(a_2) = 4\mathbb{Z} \langle S_1 \cup I \rangle \cup 6\mathbb{Z} \langle S_2 \cup I \rangle,$$

$$F(a_3) = 6\mathbb{Z} \langle S_1 \cup I \rangle \cup 8\mathbb{Z} \langle S_2 \cup I \rangle.$$

**Theorem 5.2.15.** All soft neutrosophic biideals are trivially soft subneutrosophic bisemigroup birings but the converse is not true in general.

We can easily establish the converse by the help of example.

**Proposition 5.2.8.** Let $(F, A)$ and $(K, B)$ be two soft neutosophic biideals over $R_B \langle S \cup I \rangle$. Then

1. Their extended intersection $(F, A) \cap_E (K, B)$ is a soft neutrosophic biideal over $R_B \langle S \cup I \rangle$.

2. Their restricted intersection $(F, A) \cap_R (K, B)$ is a soft neutrosophic biideal over $R_B \langle S \cup I \rangle$.

3. Their *AND* operation $(F, A) \vee (K, B)$ is a soft neutrosophic biideal over $R_B \langle S \cup I \rangle$.

**Proof.** This is straightforward.

**Remark 5.2.9.** Let $(F, A)$ and $(K, B)$ be two soft neutrosophic biideal over $R_B \langle S \cup I \rangle$. Then

1. Their restricted union $(F, A) \cup_R (K, B)$ is not a soft neutrosophic biideal over $R_B \langle S \cup I \rangle$.

2. Their extended union $(F, A) \cup_E (K, B)$ is not a soft neutrosophic biideal over $R_B \langle S \cup I \rangle$.

3. Their *OR* operation $(F, A) \vee (K, B)$ is a soft neutrosophic biideal over $R_B \langle S \cup I \rangle$.

One can easily show these remarks by the help of examples.

**Definition 5.2.11.** Let $R_B \langle S \cup I \rangle$ be a neutrosophic bisemigroup biring. Let $A$ be a set of parameters and $(F, A)$ be a soft set over $R \langle S \cup I \rangle$. Then $(F, A)$ is called soft pseudo neutosophic biideal if and only if $F(a)$ is a pseudo neutrosophic biideal of $R_B \langle S \cup I \rangle$ for all $a \in A$.

**Theorem 5.2.16.** All soft pseudo neutrosophic biideals are trivially soft subneutrosophic bisemigroup birings but the converse is not true in general.

We can easily establish the converse by the help of example.

**Theorem 5.2.17.** All soft pseudo neutrosophic biideals are trivially soft neutrosophic biideals but the converse is not true in general.

We can easily establish the converse by the help of example.

**Proposition 5.2.9.** Let $(F, A)$ and $(K, B)$ be two soft pseudo neutosophic biideals over $R_B \langle S \cup I \rangle$. Then

1. Their extended intersection $(F, A) \cap_E (K, B)$ is a soft pseudo neutrosophic biideals over $R_B \langle S \cup I \rangle$.

2. Their restricted intersection $(F, A) \cap_R (K, B)$ is a soft pseudo neutrosophic biideals over $R_B \langle S \cup I \rangle$.

3. Their $AND$ operation $(F, A) \vee (K, B)$ is a soft pseudo neutrosophic biideals over $R_B \langle S \cup I \rangle$.

**Proof.** This is straightforward.

**Remark 5.2.10.** Let $(F, A)$ and $(K, B)$ be two soft pseudo neutosophic biideals over $R_B \langle S \cup I \rangle$. Then

1. Their restricted union $(F, A) \cup_R (K, B)$ is not a soft pseudo neutrosophic biideals over $R_B \langle S \cup I \rangle$.

2. Their extended union $(F, A) \cup_E (K, B)$ is not a soft pseudo neutrosophic biideals over $R_B \langle S \cup I \rangle$.

3. Their $OR$ operation $(F,A) \vee (K,B)$ is not a soft pseudo neutrosophic biideals over $R_B \langle S \cup I \rangle$.

One can easily show these remarks by the help of examples.

**Definition 5.2.12.** Let $(F,A)$ and $(K,B)$ be two soft neutrosophic bisemigroup birings over $R_B \langle S \cup I \rangle$. Then $(K,B)$ is called soft neutrosophic biideal of $(F,A)$, if

1. $B \subseteq A$, and
2. $K(a)$ is a neutrosophic biideal of $F(a)$ for all $a \in B$.

**Theorem 5.2.18.** Every soft neutrosophic biideal of the soft neutrosophic bisemigroup biring over a neutrosophic bisemigroup biring is trivially a soft subneutrosophic bisemigroup biring.

In this final section, we finally extend soft neutrosophic semigroup ring to soft neutrosophic N-semigroup N-ring.

## 5.3    Soft Neutrosophic N-semigroup N-ring

In this section, the authors introduced soft neutrosophic N-semigroup N-rings over neutrosophic N-semigroup N-rings. This is the generalization of soft neutrosophic semigroup ring over neutrosophic semigroup ring. Here we also coined some basic and fundamental properties and other important notions of soft neutrosophic N-semigroup N-rings with sufficient amount of illustrative examples.

We now proceed to define soft neutrosophic N-semigroup N-ring over a neutrosophic N-semigroup N-ring as follows.

**Definition 5.3.1.** Let $N\left(R\langle S \cup I\rangle\right) = R_1\langle S_1 \cup I\rangle \cup R_2\langle S_2 \cup I\rangle \cup ... \cup R_n\langle S_n \cup I\rangle$ be a neutrosophic N-semigroup N-ring over $R = R_1 \cup R_2 \cup ... \cup R_n$ such that $n \geq 3$. Let $A$ be a set of parameters and let $(F, A)$ be a soft set over $N\left(R\langle S \cup I\rangle\right)$. Then $(F, A)$ is called a soft neutrosophic N-semigroup N-ring over $N\left(R\langle S \cup I\rangle\right)$ if and only if $F(a)$ is a subneutrosophic N-semigroup N-ring of $N\left(R\langle S \cup I\rangle\right)$ for all $a \in A$.

This situation can be shown as follows in the example.

**Example 5.3.1.** Let $N\left(R\langle S \cup I\rangle\right) = \mathbb{Q}\langle S_1 \cup I\rangle \cup \mathbb{R}\langle S_2 \cup I\rangle \cup \mathbb{Z}\langle S_3 \cup I\rangle$ be a neutrosophic 3-semigroup 3-ring, where $R = \mathbb{Q} \cup \mathbb{R} \cup \mathbb{Z}$ and $\langle S_1 \cup I\rangle = \{1, g, g^2, g^3, g^4, g^5, I, gI, ..., g^5I : g^6 = 1, I^2 = I\}$, $\langle S_2 \cup I\rangle = \{1, g, g^2, g^3, I, gI, g^2I, g^3I : g^4 = 1, I^2 = I\}$ and $\langle S_3 \cup I\rangle = \langle Z^+ \setminus \{0\} \cup I\rangle$ are neutrosophic semigroups. Let
$\langle H_1^\cdot \cup I\rangle \cup \langle H_2^\cdot \cup I\rangle \cup \langle H_3^\cdot \cup I\rangle = \{1, g^3 : g^6 = 1\} \cup \{1, g^2 : g^4 = 1\} \cup \langle 2Z^+ \setminus \{0\} \cup I\rangle$,
$\langle H_1^{\cdot\cdot} \cup I\rangle \cup \langle H_2^{\cdot\cdot} \cup I\rangle \cup \langle H_3^{\cdot\cdot} \cup I\rangle = \{1, g^3, I, g^3I : g^6 = 1, I^2 = I\} \cup \{1, g^2, I, g^2 I\} \cup \langle 3Z^+ \setminus \{0\} \cup I\rangle$
Let $A = \{a_1, a_2\}$ be a set of parameters. Then $(F, A)$ is a soft neutrosophic 3-semigroup 3-ring over $N\left(R\langle S \cup I\rangle\right)$, where

$$F(a_1) = \mathbb{Q}\left\langle H_1^{'} \cup I \right\rangle \cup \mathbb{R}\left\langle H_2^{'} \cup I \right\rangle \cup \mathbb{Z}\left\langle H_3^{'} \cup I \right\rangle,$$

$$F(a_2) = \mathbb{Q}\left\langle H_1^{''} \cup I \right\rangle \cup \mathbb{R}\left\langle H_2^{''} \cup I \right\rangle \cup \mathbb{Z}\left\langle H_3^{''} \cup I \right\rangle.$$

**Theorem 5.3.1.** Let $F, A$ and $(H, A)$ be two soft neutrosophic N-semigroup N-rings over $N\left(R\langle S \cup I \rangle\right)$. Then their intersection $F, A \cap H, A$ is again a soft neutrosophic N-semigroup N-ring over $N\left(R\langle S \cup I \rangle\right)$.

The proof is straightforward, so we suggest the readers to prove it by yourself as an exercise.

**Theorem 5.3.2.** Let $F, A$ and $H, B$ be two soft neutrosophic N-semigroup N-rings over $N\left(R\langle S \cup I \rangle\right)$. If $A \cap B = \phi$, then $F, A \cup H, B$ is a soft neutrosophic N-semigroup N-ring over $N\left(R\langle S \cup I \rangle\right)$.

This is straightforward, so left as an exercise.

**Proposition 5.3.1.** Let $(F, A)$ and $(K, B)$ be two soft neutosophic N-semigroup N-rings over $N\left(R\langle S \cup I \rangle\right)$. Then

1. Their extended intersection $(F, A) \cap_E (K, B)$ is again a soft neutrosophic N-semigroup N-ring over $N\left(R\langle S \cup I \rangle\right)$.
2. Their restricted intersection $(F, A) \cap_R (K, B)$ is again a soft neutrosophic N-semigroup N-ring over $R\langle S \cup I \rangle$.
3. Their $AND$ operation $(F, A) \vee (K, B)$ is again a soft neutrosophic N-semigroup N-ring over $R\langle S \cup I \rangle$.

This is straightforward, so left as an exercise for the readers.

**Remark 5.3.1.** Let $(F,A)$ and $(K,B)$ be two soft neutosophic N-semigroup N-rings over $N\big(R\langle S\cup I\rangle\big)$. Then

1. Their restricted union $(F,A)\cup_R (K,B)$ is not a soft neutrosophic N-semigroup N-ring over $N\big(R\langle S\cup I\rangle\big)$.
2. Their extended union $(F,A)\cup_E (K,B)$ is not a soft neutrosophic N-semigroup N-ring over $N\big(R\langle S\cup I\rangle\big)$.
3. Their $OR$ operation $(F,A)\vee (K,B)$ is not a soft neutrosophic N-semigroup N-ring over $N\big(R\langle S\cup I\rangle\big)$.

One can establish these remarks by the help of following examples.

**Definition 5.3.2.** Let $N\big(R\langle S\cup I\rangle\big)$ be a neutrosophic N-semigroup N-ring. Let $A$ be a set of parameters and $(F,A)$ be a soft set over $N\big(R\langle S\cup I\rangle\big)$. Then $(F,A)$ is called soft neutosophic N-semigroup sub N-ring over $N\big(R\langle S\cup I\rangle\big)$ if and only if $F(a)$ is subneutrosophic N-semigroup sub N-ring of $N\big(R\langle S\cup I\rangle\big)$ for all $a\in A$.

**Example 5.3.2.** Let $N\big(R\langle S\cup I\rangle\big)=\mathbb{R}\langle S_1\cup I\rangle\cup\mathbb{C}\langle S_2\cup I\rangle\cup\mathbb{Z}\langle S_3\cup I\rangle$ be a neutrosophic 3-semigroup 3-ring, where
$\langle S_1\cup I\rangle=\{1,g,g^2,g^3,g^4,g^5,I,gI,...,g^5I:g^6=1,I^2=I\}$,
$\langle S_2\cup I\rangle=\{1,g,g^2,g^3,I,gI,g^2I,g^3I:g^4=1,I^2=I\}$ and $\langle G_3\cup I\rangle=\langle Z^+\setminus\{0\}\cup I\rangle$ are neutrosophic semigroups.

Let $A = \{a_1, a_2\}$ be a set of parameters. Then $(F, A)$ is a soft neutrosophic 3-semigroup sub 3-ring over $N(R\langle S \cup I \rangle)$, where

$$F(a_1) = \mathbb{Q}\langle S_1 \cup I \rangle \cup \mathbb{R}\langle S_2 \cup I \rangle \cup 2\mathbb{Z}\langle S_3 \cup I \rangle,$$

$$F(a_2) = \mathbb{Z}\langle S_1 \cup I \rangle \cup \mathbb{Q}\langle S_2 \cup I \rangle \cup 3\mathbb{Z}\langle S_3 \cup I \rangle.$$

**Theorem 5.3.3.** Let $F, A$ and $(H, A)$ be two soft neutrosophic N-semigroup sub N-rings over $N(R\langle S \cup I \rangle)$. Then their intersection $F, A \cap H, A$ is again a soft neutrosophic N-semigroup sub N-rings over $N(R\langle S \cup I \rangle)$.

**Proof.** The proof is straightforward.

**Theorem 5.3.4.** Let $F, A$ and $H, B$ be two soft neutrosophic N-semigroup sub N-rings over $N(R\langle S \cup I \rangle)$. If $A \cap B = \phi$, then $F, A \cup H, B$ is a soft neutrosophic N-semigroup sub N-ring over $N(R\langle S \cup I \rangle)$.

**Proof.** This is straightforward.

**Remark 5.3.2.** Let $(F, A)$ and $(K, B)$ be two soft neutosophic N-semigroup sub N-rings over $N(R\langle S \cup I \rangle)$. Then

1. Their extended union $(F, A) \cup_E (K, B)$ is not a soft neutrosophic N-semigroup sub N-ring over $N(R\langle S \cup I \rangle)$.
2. Their restricted union $(F, A) \cup_R (K, B)$ is not a soft neutrosophic N-semigroup sub N-ring over $N(R\langle S \cup I \rangle)$.

3. Their *OR* operation $(F,A) \vee (K,B)$ is not a soft neutrosophic N-semigroup sub N-ring over $N(R\langle S \cup I \rangle)$.

One can easily show these remaks by the help of examples.

**Proposition 5.3.2.** Let $(F,A)$ and $(K,B)$ be two soft neutosophic N-semigroup sub N-rings over $N(R\langle S \cup I \rangle)$. Then

1. Their restricted intersection $(F,A) \cap_R (K,B)$ is a soft neutrosophic N-semigroup sub N-ring over $N(R\langle S \cup I \rangle)$.

2. Their extended intersection $(F,A) \cap_E (K,B)$ is a soft neutrosophic N-semigroup sub N-ring over $N(R\langle S \cup I \rangle)$.

3. Their *AND* operation $(F,A) \vee (K,B)$ is a soft neutrosophic N-semigroup sub N-ring over $N(R\langle S \cup I \rangle)$.

**Proof.** These are straightforward.

**Definition 5.3.3.** Let $N(R\langle S \cup I \rangle)$ be a neutrosophic bisemigroup biring. Let $A$ be a set of parameters and $(F,A)$ be a soft set over $R\langle S \cup I \rangle$. Then $(F,A)$ is called soft mixed neutosophic N-semigroup N-ring over $N(R\langle S \cup I \rangle)$ if for some $a \in A$, $F(a)$ is subneutrosophic N-semigroup sub N-ring and for the remaining $a \in A$, $F(a)$ is a subneutrosophic N-semigroup N-ring of $N(R\langle S \cup I \rangle)$.

**Example 5.3.3.** Let $N(R\langle S \cup I \rangle) = \mathbb{Z}\langle S_1 \cup I \rangle \cup \mathbb{C}\langle S_2 \cup I \rangle \cup \mathbb{R}\langle S_3 \cup I \rangle$ be a neutrosophic 3-semigroup 3-ring, where
$\langle S_1 \cup I \rangle = \{1, g, g^2, g^3, g^4, g^5, I, gI, ..., g^5I : g^6 = 1, I^2 = I\}$,

$\langle S_2 \cup I \rangle = \{1, g, g^2, g^3, I, gI, g^2I, g^{53}I : g^4 = 1, I^2 = I\}$ and $\langle S_3 \cup I \rangle = \langle Z^+ \setminus \{0\} \cup I \rangle$ are neutrosophic semigroups.

Let $A = \{a_1, a_2, a_3, a_4\}$ be a set of parameters. Then $(F, A)$ is a soft mixed neutrosophic 3-semigroup 3-ring over $N(R\langle S \cup I \rangle)$, where

$$F(a_1) = 2\mathbb{Z}\langle S_1 \cup I \rangle \cup \mathbb{R}\langle S_2 \cup I \rangle \cup \mathbb{Q}\langle S_3 \cup I \rangle,$$

$$F(a_2) = 4\mathbb{Z}\langle S_1 \cup I \rangle \cup \mathbb{Q}\langle S_2 \cup I \rangle \cup \mathbb{Z}\langle S_3 \cup I \rangle,$$

$$F(a_3) = \mathbb{Z}\langle H_1^{'} \cup I \rangle \cup \mathbb{C}\langle H_2^{'} \cup I \rangle \cup \mathbb{R}\langle H_3^{'} \cup I \rangle,$$

$$F(a_4) = \mathbb{Z}\langle H_1^{''} \cup I \rangle \cup \mathbb{C}\langle H_2^{''} \cup I \rangle \cup \mathbb{R}\langle H_3^{''} \cup I \rangle.$$

Where $\langle H_1^{'} \cup I \rangle = \{1, g^2 : g^4 = 1\}$, $\langle H_2^{'} \cup I \rangle = \{1, g^3, I, g^3I : g^6 = 1, I^2 = I\}$, $\langle H_3^{'} \cup I \rangle = \langle 2Z^+ \setminus \{0\} \cup I \rangle$, $\langle H_1^{''} \cup I \rangle = \{1, g^2, I, g^2I : g^4 = 1, I^2 = I\}$, $\langle H_2^{''} \cup I \rangle = \{1, g^2, g^4, I, g^2I, g^4I : g^6 = 1, I^2 = I\}$ and $\langle H_3^{''} \cup I \rangle = \langle 3Z^+ \setminus \{0\} \cup I \rangle$.

**Theorem 5.3.5.** Let $(F, A)$ and $(H, A)$ be two soft mixed neutrosophic N-semigroup N-rings over $N(R\langle S \cup I \rangle)$. Then their intersection $(F, A) \cap (H, A)$ is again a soft mixed neutrosophic N-semigroup N-ring over $N(R\langle S \cup I \rangle)$.

**Proof.** The proof is straightforward.

**Theorem 5.3.6.** Let $(F, A)$ and $(H, B)$ be two soft mixed neutrosophic N-semigroup N-rings over $N(R\langle S \cup I \rangle)$. If $A \cap B = \phi$, then $(F, A) \cup (H, B)$ is a soft mixed neutrosophic N-semigroup N-ring over $N(R\langle S \cup I \rangle)$.

**Proof.** This is straightforward.

**Remark 5.3.3.** Let $R\langle S \cup I \rangle$ be a neutrosophic N-semigroup N-ring. Then $N\ R\langle S \cup I \rangle$ can have soft neutrosophic N-semigroup N-ring, soft neutrosophic N-semigroup sub N-ring and soft mixed neutrosophic N-semigroup N-ring over $N\ R\langle S \cup I \rangle$.

**Proof:** It is obvious.

**Proposition 5.3.3.** Let $(F, A)$ and $(K, B)$ be two soft mixed neutosophic N-semigroup N-rings over $N\ R\langle S \cup I \rangle$. Then

1. Their extended intersection $(F, A) \cap_E (K, B)$ is a soft mixed neutrosophic N-semigroup N-ring over $N\ R\langle S \cup I \rangle$.
2. Their restricted intersection $(F, A) \cap_R (K, B)$ is a soft mixed neutrosophic N-semigroup N-ring over $N\ R\langle S \cup I \rangle$.
3. Their $AND$ operation $(F, A) \vee (K, B)$ is a soft mixed neutrosophic N-semigroup N-ring over $N\ R\langle S \cup I \rangle$.

**Proof.** This is straightforward.

**Remark 5.3.4.** Let $(F, A)$ and $(K, B)$ be two soft mixed neutosophic N-semigroup N-rings over $N\ R\langle S \cup I \rangle$. Then

1. Their restricted union $(F, A) \cup_R (K, B)$ is not a soft mixed neutrosophic N-semigroup N-ring over $N\ R\langle S \cup I \rangle$.
2. Their extended union $(F, A) \cup_E (K, B)$ is not a soft mixed

neutrosophic N-semigroup N-ring over $N \ R\langle S \cup I \rangle$ .

3. Their $OR$ operation $(F,A) \vee (K,B)$ is a soft mixed neutrosophic N-semigroup N-ring over $N \ R\langle S \cup I \rangle$ .

One can easily show these remarks by the help of examples.

**Definition 5.3.4.** Let $N \ R\langle S \cup I \rangle$ be a neutrosophic N-semigroup N-ring. Let $A$ be a set of parameters and $(F,A)$ be a soft set over $N \ R\langle S \cup I \rangle$ . Then $(F,A)$ is called soft neutosophic sub N-ring if and only if $F(a)$ is neutrosophic sub N-ring of $N \ R\langle S \cup I \rangle$ for all $a \in A$ .

**Example 5.3.4.** Let $N\left(R\langle S \cup I \rangle\right) = \mathbb{Z}\langle S_1 \cup I \rangle \cup \mathbb{C}\langle S_2 \cup I \rangle \cup \mathbb{R}\langle S_3 \cup I \rangle$ be a neutrosophic 3-semigroup 3-ring. Let $A = \{a_1, a_2\}$ be a set of parameters. Then $(F,A)$ is a soft neutrosophic sub 3-ring over $N \ R\langle S \cup I \rangle$ , where

$$F(a_1) = \langle 2\mathbb{Z} \cup I \rangle \cup \langle \mathbb{R} \cup I \rangle \cup \langle \mathbb{Q} \cup I \rangle,$$

$$F(a_2) = \langle 4\mathbb{Z} \cup I \rangle \cup \langle \mathbb{Q} \cup I \rangle \cup \langle \mathbb{Z} \cup I \rangle.$$

**Theorem 5.3.7.** Let $F,A$ and $(H,A)$ be two soft neutrosophic sub N-rings over $N \ R\langle S \cup I \rangle$ . Then their intersection $F,A \cap H,A$ is again a soft neutrosophic sub N-ring over $N \ R\langle S \cup I \rangle$ .

**Proof.** The proof is straightforward.

**Theorem 5.3.8.** Let $F,A$ and $H,B$ be two soft neutrosophic sub N-rings over $N \ R\langle S \cup I \rangle$ . If $A \cap B = \phi$, then $F,A \cup H,B$ is also a soft neutrosophic sub N-ring over $N \ R\langle S \cup I \rangle$ .

**Proof.** This is straightforward.

**Proposition 5.3.4.** Let $(F, A)$ and $(K, B)$ be two soft neutosophic sub N-rings over $N\ R\langle S \cup I\rangle$ . Then

1. Their extended intersection $(F, A) \cap_E (K, B)$ is a soft neutrosophic sub N-ring over $N\ R\langle S \cup I\rangle$ .

2. Their restricted intersection $(F, A) \cap_R (K, B)$ is a soft neutrosophic sub N-ring over $N\ R\langle S \cup I\rangle$ .

3. Their $AND$ operation $(F, A) \vee (K, B)$ is a soft neutrosophic sub N-ring over $N\ R\langle S \cup I\rangle$ .

**Proof.** This is straightforward.

**Remark 5.3.5.** Let $(F, A)$ and $(K, B)$ be two soft neutosophic sub N-rings over $N\ R\langle S \cup I\rangle$ . Then

1. Their restricted union $(F, A) \cup_R (K, B)$ is not a soft neutrosophic sub N-ring over $N\ R\langle S \cup I\rangle$ .

2. Their extended union $(F, A) \cup_E (K, B)$ is not a soft neutrosophic sub N-ring over $N\ R\langle S \cup I\rangle$ .

3. Their $OR$ operation $(F, A) \vee (K, B)$ is a soft neutrosophic sub N-ring over $N\ R\langle S \cup I\rangle$ .

One can easily show these remarks by the help of examples.

**Definition 5.3.5.** Let $N\ R\langle S \cup I \rangle$ be a neutrosophic N-semigroup N-ring. Let $A$ be a set of parameters and $(F,A)$ be a soft set over $N\ R\langle S \cup I \rangle$. Then $(F,A)$ is called soft pseudo neutosophic sub N-ring if and only if $F(a)$ is a pseudo neutrosophic sub N-ring of $N\ R\langle S \cup I \rangle$ for all $a \in A$.

**Theorem 5.3.9.** Let $F,A$ and $(H,A)$ be two soft pseudo neutrosophic sub N-rings over $N\ R\langle S \cup I \rangle$. Then their intersection $F,A \cap H,A$ is again a soft pseudo neutrosophic sub N-ring over $N\ R\langle S \cup I \rangle$.

**Proof.** The proof is straightforward.

**Theorem 5.3.10.** Let $F,A$ and $H,B$ be two soft pseudo neutrosophic sub N-rings over $N\ R\langle S \cup I \rangle$. If $A \cap B = \phi$, then $F,A \cup H,B$ is a soft pseudo neutrosophic sub N-ring over $N\ R\langle S \cup I \rangle$.

The proof is straightforward, so we suggest the readers to prove it by yourself as an exercise.

**Proposition 5.3.5.** Let $(F,A)$ and $(K,B)$ be two soft pseudo neutosophic sub N-rings over $N\ R\langle S \cup I \rangle$. Then

1. Their extended intersection $(F,A) \cap_E (K,B)$ is a soft pseudo neutrosophic sub N-ring over $N\ R\langle S \cup I \rangle$.
2. Their restricted intersection $(F,A) \cap_R (K,B)$ is a soft pseudo neutrosophic sub N-ring over $N\ R\langle S \cup I \rangle$.

3. Their *AND* operation $(F,A) \vee (K,B)$ is a soft pseudo neutrosophic sub N-ring over $N \ R\langle S \cup I \rangle$ .

**Proof.** This is straightforward.

**Remark 5.3.6.** Let $(F,A)$ and $(K,B)$ be two soft pseudo neutosophic sub N-rings over $N \ R\langle S \cup I \rangle$ . Then

1. Their restricted union $(F,A) \cup_R (K,B)$ is not a soft pseudo neutrosophic sub N-ring over $N \ R\langle S \cup I \rangle$ .
2. Their extended union $(F,A) \cup_E (K,B)$ is not a soft pseudo neutrosophic sub N-ring over $N \ R\langle S \cup I \rangle$ .
3. Their *OR* operation $(F,A) \vee (K,B)$ is a soft pseudo neutrosophic sub N-ring over $N \ R\langle S \cup I \rangle$ .

One can easily show these remarks by the help of examples.

**Definition 5.3.6.** Let $N \ R\langle S \cup I \rangle$ be a neutrosophic N-semigroup N-ring. Let $A$ be a set of parameters and $(F,A)$ be a soft set over $N \ R\langle S \cup I \rangle$ . Then $(F,A)$ is called soft neutosophic sub N-semigroup N-ring if and only if $F(a)$ is neutrosophic sub N-semigroup N-ring of $N \ R\langle S \cup I \rangle$ for all $a \in A$ .

**Example 5.3.5.** Let $N\left(R\langle S \cup I \rangle\right) = \mathbb{Z}\langle S_1 \cup I \rangle \cup \mathbb{C}\langle S_2 \cup I \rangle \cup \mathbb{R}\langle S_3 \cup I \rangle$ be a neutrosophic 3-semigroup 3-ring, where

$\langle S_1 \cup I \rangle = \{1, g, g^2, g^3 I, gI, g^2 I, g^3 I : g^6 = 1, I^2 = I \}$ ,

$\langle S_2 \cup I \rangle = \{1, g, g^2, g^3, g^4, g^5, I, gI, ..., g^5 I : g^6 = 1, I^2 = I \}$ and $\langle S_3 \cup I \rangle = \langle Z^+ \setminus \{0\} \cup I \rangle$ .

Let $H_1^{'} = \{1, g^2 : g^4 = 1\}$, $H_2^{'} = \{1, g^2, g^4 : g^6 = 1\}$, $H_3^{'} = \langle 2Z^+ \setminus \{0\} \cup I \rangle$,

$H_1^{''} = \{1, g^2, I, g^2 I : g^4 = 1, I^2 = I\}$, $H_2^{''} = \{1, g^3 : g^6 = 1\}$ and $H_3^{''} = \langle 3Z^+ \setminus \{0\} \cup I \rangle$. Let $A = \{a_1, a_2\}$ be a set of parameters. Then $(F, A)$ is a soft neutrosophic sub 3-semigroup 3-ring over $N\left(R\langle S \cup I \rangle\right)$, where

$$F(a_1) = 2\mathbb{Z}H_1^{'} \cup \mathbb{R}H_2^{'} \cup \mathbb{Q}H_3^{'},$$

$$F(a_1) = 8\mathbb{Z}H_1^{''} \cup \mathbb{Q}H_2^{''} \cup \mathbb{Z}H_3^{''}.$$

**Theorem 5.3.11.** Let $F, A$ and $(H, A)$ be two soft neutrosophic sub N-semigroup N-rings over $N\left(R\langle S \cup I \rangle\right)$. Then their intersection $F, A \cap H, A$ is again a soft neutrosophic sub N-semigroup N-ring over $N\left(R\langle S \cup I \rangle\right)$.

**Proof.** The proof is straightforward.

**Theorem 5.3.12.** Let $F, A$ and $H, B$ be two soft neutrosophic sub N-semigroup N-rings over $N\left(R\langle S \cup I \rangle\right)$. If $A \cap B = \phi$, then $F, A \cup H, B$ is a soft neutrosophic sub N-semigroup N-ring.

**Proof.** This is straightforward.

**Proposition 5.3.6.** Let $(F, A)$ and $(K, B)$ be two soft neutosophic sub N-semigroup N-rings over $N\left(R\langle S \cup I \rangle\right)$. Then

1. Their extended intersection $(F, A) \cap_E (K, B)$ is a soft neutrosophic sub N-semigroup N-ring over $N\left(R\langle S \cup I \rangle\right)$.
2. Their restricted intersection $(F, A) \cap_R (K, B)$ is a soft neutrosophic

sub N-semigroup N-ring over $N(R\langle S \cup I\rangle)$.

3. Their *AND* operation $(F,A) \vee (K,B)$ is a soft neutrosophic sub N-semigroup N-ring over $N(R\langle S \cup I\rangle)$.

**Proof.** This is straightforward.

**Remark 5.3.7.** Let $(F,A)$ and $(K,B)$ be two soft neutosophic sub N-semigroup N-rings over $N(R\langle S \cup I\rangle)$. Then

1. Their restricted union $(F,A) \cup_R (K,B)$ is not a soft neutrosophic sub N-semigroup N-ring over $N(R\langle S \cup I\rangle)$.

2. Their extended union $(F,A) \cup_E (K,B)$ is not a soft neutrosophic sub N-semigroup N-ring over $N(R\langle S \cup I\rangle)$.

3. Their *OR* operation $(F,A) \vee (K,B)$ is a soft neutrosophic sub N-semigroup N-ring over $N(R\langle S \cup I\rangle)$.

One can easily show these remarks by the help of examples.

**Definition 5.3.7.** Let $N(R\langle S \cup I\rangle)$ be a neutrosophic N-semigroup N-ring. Let $A$ be a set of parameters and $(F,A)$ be a soft set over $N(R\langle S \cup I\rangle)$. Then $(F,A)$ is called a soft sub N-ring over $N(R\langle S \cup I\rangle)$ if and only if $F(a)$ is a sub N-ring of $N(R\langle S \cup I\rangle)$ for all $a \in A$ .

**Theorem 5.3.13.** Let $F,A$ and $(H,A)$ be two soft sub N-rings over $N(R\langle S \cup I\rangle)$. Then their intersection $F,A \cap H,A$ is again a soft sub N-ring over $N(R\langle S \cup I\rangle)$ .

**Proof.** The proof is straightforward.

**Theorem 5.3.14.** Let $(F, A)$ and $(H, B)$ be two soft sub N-rings over $N\big(R\langle S \cup I \rangle\big)$. If $A \cap B = \phi$, then $(F, A) \cup (H, B)$ is a soft sub N-ring over $N\big(R\langle S \cup I \rangle\big)$.

**Proof.** This is straightforward.

**Proposition 5.3.7.** Let $(F, A)$ and $(K, B)$ be two soft sub N-rings over $N\big(R\langle S \cup I \rangle\big)$. Then

1. Their extended intersection $(F, A) \cap_E (K, B)$ is a soft sub N-ring over $N\big(R\langle S \cup I \rangle\big)$.

2. Their restricted intersection $(F, A) \cap_R (K, B)$ is a soft sub N-ring over $N\big(R\langle S \cup I \rangle\big)$.

3. Their $AND$ operation $(F, A) \vee (K, B)$ is a soft sub N-ring over $N\big(R\langle S \cup I \rangle\big)$.

**Proof.** These are straightforward.

**Remark 5.3.8.** Let $(F, A)$ and $(K, B)$ be two soft sub N-rings over $N\big(R\langle S \cup I \rangle\big)$. Then

1. Their restricted union $(F, A) \cup_R (K, B)$ is not a soft sub N-ring over $N\big(R\langle S \cup I \rangle\big)$.

2. Their extended union $(F, A) \cup_E (K, B)$ is not a soft sub N-ring over $N\big(R\langle S \cup I \rangle\big)$.

3. Their *OR* operation $(F, A) \vee (K, B)$ is a soft sub N-ring over $N\big(R\langle S \cup I\rangle\big)$.

One can easily show these remarks by the help of examples.

**Definition 5.3.8.** Let $N\big(R\langle S \cup I\rangle\big)$ be a neutrosophic N-semigroup N-ring. Then $(F, A)$ is called an absolute soft neutosophic N-semigroup N-ring if $F(a) = N\big(R\langle S \cup I\rangle\big)$, for all $a \in A$.

**Definition 5.3.9.** Let $(F, A)$ and $(K, B)$ be two soft neutrosophic N-semigroup N-rings over $N\big(R\langle S \cup I\rangle\big)$. Then $(K, B)$ is called soft subneutrosophic N-semigroup N-ring of $(F, A)$, if

1. $B \subseteq A$, and
2. $K(a)$ is a subneutrosophic N-semigroup N-ring of $F(a)$ for all $a \in B$.

**Definition 5.3.10.** Let $N\big(R\langle S \cup I\rangle\big)$ be a neutrosophic N-semigroup N-ring. Let $A$ be a set of parameters and $(F, A)$ be a soft set over $N\big(R\langle S \cup I\rangle\big)$. Then $(F, A)$ is called soft neutosophic N-ideal if and only if $F(a)$ is a neutrosophic N-ideal of $N\big(R\langle S \cup I\rangle\big)$ for all $a \in A$.

Similarly one can define soft neutrosophic left N-ideal and soft neutrosophic right N-ideal over $N\big(R\langle S \cup I\rangle\big)$.

**Example 5.3.6.** Let $N\big(R\langle S \cup I\rangle\big) = \mathbb{Z}\langle S_1 \cup I\rangle \cup \mathbb{Z}\langle S_2 \cup I\rangle \cup \mathbb{Z}\langle S_3 \cup I\rangle$ be a neutrosophic 3-semigroup 3-ring. Let $A = \{a_1, a_2, a_3\}$ be a set of parameters. Then $(F, A)$ is a soft neutrosophic 3-ideal over $N\big(R\langle S \cup I\rangle\big)$, where

$$F(a_1) = 2\mathbb{Z}\langle S_1 \cup I\rangle \cup 4\mathbb{Z}\langle S_2 \cup I\rangle \cup 6\mathbb{Z}\langle S_3 \cup I\rangle,$$

$$F(a_2) = 4\mathbb{Z}\langle S_1 \cup I\rangle \cup 6\mathbb{Z}\langle S_2 \cup I\rangle \cup 8\mathbb{Z}\langle S_3 \cup I\rangle,$$

$$F(a_3) = 6\mathbb{Z}\langle S_1 \cup I\rangle \cup 8\mathbb{Z}\langle S_2 \cup I\rangle \cup 10\mathbb{Z}\langle S_3 \cup I\rangle.$$

**Theorem 5.3.15.** All soft neutrosophic N-ideals are trivially soft subneutrosophic N-semigroup N-rings but the converse is not true in general.

We can easily establish the converse by the help of example.

**Proposition 5.3.8.** Let $(F, A)$ and $(K, B)$ be two soft neutosophic N-ideals over $N\left(R\langle S \cup I\rangle\right)$. Then

1. Their extended intersection $(F, A) \cap_E (K, B)$ is a soft neutrosophic N-ideal over $N\left(R\langle S \cup I\rangle\right)$.
2. Their restricted intersection $(F, A) \cap_R (K, B)$ is a soft neutrosophic N-ideal over $N\left(R\langle S \cup I\rangle\right)$.
3. Their $AND$ operation $(F, A) \vee (K, B)$ is a soft  neutrosophic N-ideal over $N\left(R\langle S \cup I\rangle\right)$.

**Proof.** This is straightforward.

**Remark 5.3.9.** Let $(F, A)$ and $(K, B)$ be two soft neutosophic N-ideal over $N\left(R\langle S \cup I\rangle\right)$. Then

1. Their restricted union $(F, A) \cup_R (K, B)$ is not a soft neutrosophic N-ideal over $N\left(R\langle S \cup I\rangle\right)$.

2. Their extended union $(F,A) \cup_E (K,B)$ is not a soft neutrosophic N-ideal over $N\big(R\langle S \cup I\rangle\big)$.

3. Their $OR$ operation $(F,A) \vee (K,B)$ is a soft neutrosophic N-ideal over $N\big(R\langle S \cup I\rangle\big)$.

One can easily show these remarks by the help of examples.

**Definition 5.3.11.** Let $N\big(R\langle S \cup I\rangle\big)$ be a neutrosophic N-semigroup N-ring. Let $A$ be a set of parameters and $(F,A)$ be a soft set over $N\big(R\langle S \cup I\rangle\big)$. Then $(F,A)$ is called soft pseudo neutosophic N-ideal if and only if $F(a)$ is a pseudo neutrosophic N-ideal of $N\big(R\langle S \cup I\rangle\big)$ for all $a \in A$.

**Theorem 5.3.16.** All soft pseudo neutrosophic N-ideals are trivially soft subneutrosophic N-semigroup N-rings but the converse is not true in general.

We can easily establish the converse by the help of example.

**Theorem 5.3.17.** All soft pseudo neutrosophic N-ideals are trivially soft neutrosophic N-ideals but the converse is not true in general.

We can easily establish the converse by the help of example.

**Proposition 5.3.9.** Let $(F,A)$ and $(K,B)$ be two soft pseudo neutosophic N-ideals over $N\big(R\langle S \cup I\rangle\big)$. Then

1. Their extended intersection $(F,A) \cap_E (K,B)$ is a soft pseudo neutrosophic N-ideals over $N\big(R\langle S \cup I\rangle\big)$.

2. Their restricted intersection $(F,A) \cap_R (K,B)$ is a soft pseudo

neutrosophic N-ideals over $N\left(R\langle S\cup I\rangle\right)$.

3. Their *AND* operation $(F,A)\vee(K,B)$ is a soft pseudo neutrosophic N-ideals over $N\left(R\langle S\cup I\rangle\right)$.

**Proof.** This is straightforward.

**Remark 5.3.10.** Let $(F,A)$ and $(K,B)$ be two soft pseudo neutosophic N-ideals over $N\left(R\langle S\cup I\rangle\right)$. Then

1. Their restricted union $(F,A)\cup_R(K,B)$ is not a soft pseudo neutrosophic N-ideals over $N\left(R\langle S\cup I\rangle\right)$.

2. Their extended union $(F,A)\cup_E(K,B)$ is not a soft pseudo neutrosophic N-ideals over $N\left(R\langle S\cup I\rangle\right)$.

3. Their *OR* operation $(F,A)\vee(K,B)$ is a soft pseudo neutrosophic N-ideals over $N\left(R\langle S\cup I\rangle\right)$.

One can easily show these remarks by the help of examples.

**Definition 5.2.12.** Let $(F,A)$ and $(K,B)$ be two soft neutrosophic N-semigroup N-rings over $N\left(R\langle S\cup I\rangle\right)$. Then $(K,B)$ is called soft neutrosophic N-ideal of $(F,A)$, if

1. $B\subseteq A$, and
2. $K(a)$ is a neutrosophic N-ideal of $F(a)$ for all $a\in B$.

**Theorem 5.3.18.** Every soft neutrosophic N-ideal of the soft neutrosophic N-semigroup N-ring over a neutrosophic N-semigroup N-ring is trivially a soft subneutrosophic N-semigroup N-ring.

The proof is straightforward, so we suggest the readers to prove it by yourself as an exercise.

# Chapter No. 6

# Soft Mixed Neutrosophic N-Algebraic Structures

In this chapter, the authors for the first time introduced soft mixed neutrosophic N-algebraic structures over mixed neutrosophic N-algebraic structures which are basically the approximated collection of mixed neutrosophic sub N-algebraic structures of a mixed neutrosophic N-algebraic structure. We also define soft dual neutrosophic N-algebraic structures over dual neutrosophic N-algebraic structures, soft weak neutrosophic N-algebraic structures over weak neutrosophic N-algebraic structures etc. We also give some basic and fundamental properties and characterization of soft mixed neutrosophic N-algebraic structures with many illustrative examples.

We now define soft mixed neutrosophic n-algebraic structures over a mixed neutrosophic N-algebraic structures.

**Definition 6.1.1.** Let $\langle M \cup I \rangle$ be a mixed neutrosophic $N$-algebraic structure and let $(F, A)$ soft set over $\langle M \cup I \rangle$. Then $(F, A)$ is called a soft mixed neutrosophic $N$-algebraic structure if and only if $F(a)$ is a mixed neutrosophic sub $N$-algebraic structure of $\langle M \cup I \rangle$ for all $a \in A$.

This situation can be explained in the following example.

**Example 6.1.1.** Let $\left\{ \langle M \cup I \rangle = M_1 \cup M_2 \cup M_3 \cup M_4 \cup M_5, *_1, *_2, *_3, *_4, *_5 \right\}$ be a mixed neutosophic $5$-structure, where

$M_1 = \langle \mathbb{Z}_3 \cup I \rangle$, a neutrosophic group under multiplication modolu $3$,

$M_2 = \langle \mathbb{Z}_6 \cup I \rangle$, a neutrosophic semigroup under multiplication modolu $6$,

$M_3 = \{0, 1, 2, 3, 1I, 2I, 3I,$, a neutrosophic groupoid under multiplication modolu $4\}$,

$M_4 = S_3$, and $M_5 = \{Z_{10}$, a semigroup under multiplication modolu $10\}$.

Let $A = \{a_1, a_2, a_3\}$ be a set of parameters and let $(F, A)$ be a soft set over $\langle M \cup I \rangle$, where

$$F(a_1) = \{1, I\} \cup \{0, 3, 3I\} \cup \{0, 2, 2I\} \cup A_3 \cup \{0, 2, 4, 6, 8\},$$

$$F(a_2) = \{2, I\} \cup \{0, 2, 4, 2I, 4I\} \cup \{0, 2, 2I\} \cup A_3 \cup \{0, 5\},$$

$$F(a_3) = \{1, 2\} \cup \{0, 3\} \cup \{0, 2\} \cup A_3 \cup \{0, 2, 4, 6, 8\}.$$

Clearly $(F, A)$ is a soft mixed neutrosophic $N$-algebraic structure over $\langle M \cup I \rangle$.

We now give some characterization of soft mixed neutrosophic N-algebraic structures.

**Theorem 6.1.1.** Let $(F,A)$ and $(H,A)$ be two soft mixed neutrosophic $N$-algebraic structures over $\langle M \cup I \rangle$. Then their intersection $(F,A) \cap (H,A)$ is again a soft mixed neutrosophic $N$-algebraic structure over $\langle M \cup I \rangle$.

The proof is straightforward, so left as an exercise for the readers.

**Theorem 6.1.2.** Let $(F,A)$ and $(H,B)$ be two soft mixed neutrosophic $N$-algebraic stuctures over $\langle M \cup I \rangle$. If $A \cap B = \phi$, then $(F,A) \cup (H,B)$ is a soft mixed neutrosophic $N$-algebraic structure over $\langle M \cup I \rangle$.

**Proof.** The proof is straightforward.

**Proposition 6.1.1.** Let $(F,A)$ and $(K,C)$ be two soft mixed neutrosophic $N$-algebraic structure over $\langle M \cup I \rangle$. Then

1. Their extended intersection $(F,A) \cap_E (K,C)$ is a soft mixed neutrosophic $N$-algebraic structure over $\langle M \cup I \rangle$.
2. Their restricted intersection $(F,A) \cap_R (K,C)$ is a soft mixed neutrosophic $N$-algebraic structure over $\langle M \cup I \rangle$.
3. Their $AND$ operation $(F,A) \wedge (K,C)$ is a soft mixed neutrosophic $N$-algebraic structure over $\langle M \cup I \rangle$.

**Remark 6.1.1.** Let $(F,A)$ and $(K,C)$ be two soft mixed neutrosophic $N$-algebraic structure over $\langle M \cup I \rangle$. Then

1. Their restricted  union $(F,A) \cup_R (K,C)$ may not be a soft mixed neutrosophic $N$-algebraic structure  over $\langle M \cup I \rangle$.
2. Their extended union  $(F,A) \cup_E (K,C)$ may not be a soft mxed neutrosophic $N$-algebraic structure  over $\langle M \cup I \rangle$.
3. Their *OR* operation $(F,A) \vee (K,C)$ may not be a soft mixed neutrosophic $N$-algebraic structure  over $\langle M \cup I \rangle$.

To establish the above remark, see the following Example.

**Example 6.1.2.**  Let $\left\{ \langle M \cup I \rangle = M_1 \cup M_2 \cup M_3 \cup M_4 \cup M_5, *_1, *_2, *_3, *_4, *_5 \right\}$ be a mixed neutosophic $5$-structure, where

$M_1 = \langle \mathbb{Z}_3 \cup I \rangle$, a neutrosophic group under multiplication modolu $3$,

$M_2 = \langle \mathbb{Z}_6 \cup I \rangle$, a neutrosophic semigroup under multiplication modolu $6$,

$M_3 = \{0,1,2,3,1I,2I,3I,$, a neutrosophic groupoid under multiplication modolu $4\}$,

$M_4 = S_3$, and $M_5 = \{Z_{10}$, a semigroup under multiplication modolu $10\}$.

Let $A = \{a_1, a_2, a_3\}$ be a set of parameters. Then $(F,A)$ is a soft mixed neutrosophic  $N$-algebraic structure  over $\langle M \cup I \rangle$, where

$$F(a_1) = \{1, I\} \cup \{0, 3, 3I\} \cup \{0, 2, 2I\} \cup A_3 \cup \{0, 2, 4, 6, 8\},$$

$$F(a_2) = \{2, I\} \cup \{0, 2, 4, 2I, 4I\} \cup \{0, 2, 2I\} \cup A_3 \cup \{0, 5\},$$

$$F(a_3) = \{1, 2\} \cup \{0, 3\} \cup \{0, 2\} \cup A_3 \cup \{0, 2, 4, 6, 8\}.$$

Let $B = \{a_1, a_4\}$ be a set of parameters. Let $(K, C)$ be another soft mixed neutrosophic $N$-algebraic structure over $\langle M \cup I \rangle$, where

$$K(a_1) = \{1, I\} \cup \{0, 3I\} \cup \{0, 2, 2I\} \cup A_3 \cup \{0, 2, 4, 6, 8\},$$

$$K(a_2) = \{1, 2\} \cup \{0, 3I\} \cup \{0, 2I\} \cup A_3 \cup \{0, 5\}.$$

Let $C = A \cap B = \{a_1\}$. Then the restricted union $(F, A) \cup_R (K, B) = (H, C)$, where

$$H(a_1) = F(a_1) \cup K(a_1) = \{1, I, 2\} \cup \{0, 3I\} \cup \{0, 2, 2I\} \cup A_3 \cup \{0, 2, 4, 5, 6, 8\}.$$

Thus clearly in $H(a_1)$, $\{1, I, 2\}, \{0, 2, 4, 5, 6, 8\}$ are not subgroups. This shows that $(H, C)$ is not a soft mixed neutrosophic $N$-algebraic structure over $\langle M \cup I \rangle$.

Similarly one can easily show 2 and 3 by the help of examples.

**Definition 6.1.2.** Let $\langle D \cup I \rangle$ be a mixed dual neutrosophic $N$-algebraic structure and let $(F, A)$ soft set over $\langle D \cup I \rangle$. Then $(F, A)$ is called a soft mixed dual neutrosophic $N$-algebraic structure if and only if $F(a)$ is a mixed dual neutrosophic sub $N$-algebraic structure of $\langle D \cup I \rangle$ for all $a \in A$.

**Example 6.1.3.** Let $\left\{ \langle D \cup I \rangle = D_1 \cup D_2 \cup D_3 \cup D_4 \cup D_5, *_1, *_2, *_3, *_4, *_5 \right\}$ be a mixed dual neutosophic 5-structure, where

$D_1 = L_7(4)$, $D_2 = S_4$, $D_3 = \{Z_{10}$, a semigroup under multiplication modulo 10$\}$,

$D_4 = \{0,1,2,3$ , a groupoid under multiplication modulo $4\}$ ,

$D_5 = \langle L_7(4) \cup I \rangle$ .

Let $A = \{a_1, a_2\}$ be a set of parameters and let $(F, A)$ be a soft set over $\langle D \cup I \rangle$, where

$$F(a_1) = \{e, 2\} \cup A_4 \cup \{0, 2, 4, 6, 8\} \cup \{0, 2\} \cup \{e, eI, 2, 2I\} ,$$

$$F(a_2) = \{e, 3\} \cup S_3 \cup \{0, 5\} \cup \{0, 2\} \cup \{e, eI, 3, 3I\} .$$

Clearly $(F, A)$ is a soft mixed dual neutrosophic $N$ -structure over $\langle D \cup I \rangle$.

**Theorem 6.1.3.** If $\langle D \cup I \rangle$ is a mixed dual neutrosophic $N$ -algebraic structure. Then $(F, A)$ over $\langle W \cup I \rangle$ is also a soft mixed dual neutrosophic $N$ -algebraic structure.

**Theorem 6.1.4.** Let $(F, A)$ and $(H, A)$ be two soft mixed dual neutrosophic $N$ -algebraic structures over $\langle D \cup I \rangle$. Then their intersection $(F, A) \cap (H, A)$ is again a soft mixed dual neutrosophic $N$ -algebraic structure over $\langle D \cup I \rangle$.

**Proof.** The proof is straightforward.

**Theorem 6.1.5.** Let $(F, A)$ and $(H, B)$ be two soft mixed dual neutrosophic $N$ -algebraic stuctures over $\langle D \cup I \rangle$. If $A \cap B = \phi$, then $(F, A) \cup (H, B)$ is a soft mixed dual neutrosophic $N$ -algebraic structure over $\langle M \cup I \rangle$.

**Proof.** The proof is straightforward.

**Proposition 6.1.2.** Let $(F,A)$ and $(K,C)$ be two soft mixed dual neutrosophic $N$-algebraic structure over $\langle D \cup I \rangle$. Then

1. Their extended intersection $(F,A) \cap_E (K,C)$ is a soft mixed dual neutrosophic $N$-algebraic structure over $\langle D \cup I \rangle$.

2. Their restricted intersection $(F,A) \cap_R (K,C)$ is a soft mixed dual neutrosophic $N$-algebraic structure over $\langle D \cup I \rangle$.

3. Their *AND* operation $(F,A) \wedge (K,C)$ is a soft mixed dual neutrosophic $N$-algebraic structure over $\langle D \cup I \rangle$.

**Remark 6.1.2.** Let $(F,A)$ and $(K,C)$ be two soft mixed Dual neutrosophic $N$-algebraic structure over $\langle D \cup I \rangle$. Then

1. Their restricted union $(F,A) \cup_R (K,C)$ may not be a soft mixed dual neutrosophic $N$-algebraic structure over $\langle D \cup I \rangle$.

2. Their extended union $(F,A) \cup_E (K,C)$ may not be a soft mixed dual neutrosophic $N$-algebraic structure over $\langle D \cup I \rangle$.

3. Their *OR* operation $(F,A) \vee (K,C)$ may not be a soft mixed dual neutrosophic $N$-algebraic structure over $\langle D \cup I \rangle$.

One can easily establish the above remarks by the help of examples.

**Definition 6.1.3.** Let $\langle W \cup I \rangle$ be a weak mixed neutrosophic $N$-algebraic structure and let $(F,A)$ soft set over $\langle W \cup I \rangle$. Then $(F,A)$ is called a soft weak mixed neutrosophic $N$-algebraic structure if and only if $F(a)$ is a weak mixed neutrosophic sub $N$-structure of $\langle W \cup I \rangle$ for all $a \in A$.

**Theorem 6.1.6.** If $\langle W \cup I \rangle$ is a weak mixed neutrosophic $N$-algebraic structure. Then $(F,A)$ over $\langle W \cup I \rangle$ is also a soft weak mixed neutrosophic $N$-algebraic structure.

The restricted intersection, extended intersection and the $AND$ operation of two soft weak mixed neutrosophic $N$-algebraic strucures are again soft weak mixed neutrosophic $N$-algebraic structures.

The restricted union, extended union and the $OR$ operation of two soft weak mixed neutrosophic $N$-algebraic strucures may not be soft weak mixed neutrosophic $N$-algebraic structures.

**Definition 6.1.4.** Let $\langle V \cup I \rangle$ be a weak mixed dual neutrosophic $N$-algebraic structure and let $(F,A)$ soft set over $\langle V \cup I \rangle$. Then $(F,A)$ is called a soft weak mixed dual neutrosophic $N$-algebraic structure if and only if $F(a)$ is a weak mixed dual neutrosophic sub $N$-structure of $\langle V \cup I \rangle$ for all $a \in A$.

**Theorem 6.1.7.** If $\langle V \cup I \rangle$ is a weak mixed dual neutrosophic $N$-algebraic structure. Then $(F,A)$ over $\langle V \cup I \rangle$ is also a soft weak mixed dual neutrosophic $N$-algebraic structure.

The restricted intersection, extended intersection and the $AND$ operation of two soft weak mixed dual neutrosophic $N$-algebraic strucures are again soft weak mixed dual neutrosophic $N$-algebraic structures.

The restricted union, extended union and the *OR* operation of two soft weak mixed dual neutrosophic $N$-algebraic strucures may not be soft weak mixed dual neutrosophic $N$-algebraic structures.

**Definition 6.1.5.** Let $(F,A)$ and $(H,C)$ be two soft mixed neutrosophic $N$-algebraic structures over $\langle M \cup I \rangle$. Then $(H,C)$ is called soft mixed neutrosophic sub $N$-algebraic structure of $(F,A)$, if

1. $C \subseteq A$.
2. $H(a)$ is a mixed neutrosophic sub $N$-structure of $F(a)$ for all $a \in C$.

It is important to note that a soft mixed neutrosophic $N$-algebraic structure can have soft weak mixed neutrosophic sub $N$-algebraic structure. But a soft weak mixed neutrosophic sub $N$-strucure cannot in general have a soft mixed neutrosophic $N$-stucture.

**Definition 6.1.6.** Let $\langle W \cup I \rangle$ be a weak mixed neutrosophic $N$-algebraic structure and let $(F,A)$ soft set over $\langle W \cup I \rangle$. Then $(F,A)$ is called a soft weak mixed deficit neutrosophic $N$-algebraic structure if and only if $F(a)$ is a weak mixed deficit neutrosophic sub $N$-structure of $\langle W \cup I \rangle$ for all $a \in A$.

**Proposition 6.1.3.** Let $(F,A)$ and $(K,C)$ be two soft weak mixed deficit neutrosophic $N$-algebraic structure over $\langle W \cup I \rangle$. Then

1. Their extended intersection $(F,A) \cap_E (K,C)$ is a soft weak mixed deficit neutrosophic $N$-algebraic structure over $\langle M \cup I \rangle$.
2. Their restricted intersection $(F,A) \cap_R (K,C)$ is a soft weak mixed

deficit neutrosophic $N$-algebraic structure over $\langle M \cup I \rangle$.

3. Their *AND* operation $(F,A) \wedge (K,C)$ is a soft weak mixed deficit neutrosophic $N$-algebraic structure over $\langle M \cup I \rangle$.

**Remark 6.1.3.** Let $(F,A)$ and $(K,C)$ be two soft weak mixed deficit neutrosophic $N$-algebraic structure over $\langle W \cup I \rangle$. Then

1. Their restricted union $(F,A) \cup_R (K,C)$ may not be a soft weak mixed deficit neutrosophic $N$-algebraic structure over $\langle W \cup I \rangle$.

2. Their extended union $(F,A) \cup_E (K,C)$ may not be a soft weak mixed deficit neutrosophic $N$-algebraic structure over $\langle W \cup I \rangle$.

3. Their *OR* operation $(F,A) \vee (K,C)$ may not be a soft weak mixed deficit neutrosophic $N$-algebraic structure over $\langle Ws \cup I \rangle$.

One can easily establish the above remarks by the help of examples.

**Definition 6.1.7.** Let $\langle M \cup I \rangle$ be a mixed neutrosophic $N$-algebraic structure and let $(F,A)$ soft set over $\langle M \cup I \rangle$. Then $(F,A)$ is called a soft Lagrange mixed neutrosophic $N$-algebraic structure if and only if $F(a)$ is a Lagrange mixed neutrosophic sub $N$-structure of $\langle M \cup I \rangle$ for all $a \in A$.

**Theorem 6.1.8.** If $\langle M \cup I \rangle$ is a Lagrange mixed neutrosophic $N$-algebraic structure. Then $(F,A)$ over $\langle M \cup I \rangle$ is also a soft Lagrange mixed neutrosophic $N$-algebraic structure.

**Remark 6.1.4.** Let $(F,A)$ and $(K,C)$ be two soft Lagrange mixed neutrosophic $N$-algebraic structures over $\langle M \cup I \rangle$. Then

1. Their restricted union $(F,A) \cup_R (K,C)$ may not be a soft Lagrange mixed neutrosophic $N$-algebraic structure over $\langle M \cup I \rangle$.

2. Their extended union $(F,A) \cup_E (K,C)$ may not be a soft Lagrange mixed neutrosophic $N$-algebraic structure over $\langle M \cup I \rangle$.

3. Their *OR* operation $(F,A) \vee (K,C)$ may not be a soft Lagrange mixed neutrosophic $N$-algebraic structure over $\langle M \cup I \rangle$.

4. Their extended intersection $(F,A) \cap_E (K,C)$ may not be a soft Lagrange mixed neutrosophic $N$-algebraic structure over $\langle M \cup I \rangle$.

5. Their restricted intersection $(F,A) \cap_R (K,C)$ may not be a soft Lagrange mixed neutrosophic $N$-algebraic structure over $\langle M \cup I \rangle$.

6. Their *AND* operation $(F,A) \wedge (K,C)$ may not be a soft Lagrange mixed neutrosophic $N$-algebraic structure over $\langle M \cup I \rangle$.

One can easily establish the above remarks by the help of examples.

Now on similar lines, we can define soft Lagrange weak deficit mixed neutrosophic $N$-algebraic structures.

**Definition 6.1.8.** Let $\langle M \cup I \rangle$ be a mixed neutrosophic $N$-algebraic structure and let $(F,A)$ soft set over $\langle M \cup I \rangle$. Then $(F,A)$ is called a soft weak Lagrange mixed neutrosophic $N$-algebraic structure if and only if $F(a)$ is not a Lagrange mixed neutrosophic sub $N$-structure of $\langle M \cup I \rangle$ for some $a \in A$.

**Remark 6.1.5.** Let $(F,A)$ and $(K,C)$ be two soft weak Lagrange mixed neutrosophic $N$-algebraic structures over $\langle M \cup I \rangle$. Then

1. Their restricted union $(F,A) \cup_R (K,C)$ may not be a soft weak Lagrange mixed neutrosophic $N$-algebraic structure over $\langle M \cup I \rangle$.

2. Their extended union $(F,A) \cup_E (K,C)$ may not be a soft weak Lagrange mixed neutrosophic $N$-algebraic structure over $\langle M \cup I \rangle$.

3. Their $OR$ operation $(F,A) \vee (K,C)$ may not be a soft weak Lagrange mixed neutrosophic $N$-algebraic structure over $\langle M \cup I \rangle$.

4. Their extended intersection $(F,A) \cap_E (K,C)$ may not be a soft weak Lagrange mixed neutrosophic $N$-algebraic structure over $\langle M \cup I \rangle$.

5. Their restricted intersection $(F,A) \cap_R (K,C)$ may not be a soft weak Lagrange mixed neutrosophic $N$-algebraic structure over $\langle M \cup I \rangle$.

6. Their $AND$ operation $(F,A) \wedge (K,C)$ may not be a soft weak Lagrange mixed neutrosophic $N$-algebraic structure over $\langle M \cup I \rangle$.

One can easily establish the above remarks by the help of examples.

Similarly we can define soft weak Lagrange weak deficit mixed neutrosophic $N$-algebraic structures.

**Definition 6.1.9.** Let $\langle M \cup I \rangle$ be a mixed neutrosophic $N$-algebraic structure and let $(F,A)$ soft set over $\langle M \cup I \rangle$. Then $(F,A)$ is called a soft Lagrange free mixed neutrosophic $N$-algebraic structure if and only if $F(a)$ is not a Lagrange mixed neutrosophic sub $N$-structure of $\langle M \cup I \rangle$ for all $a \in A$.

**Theorem 6.1.9.** If $\langle M \cup I \rangle$ is a Lagrange free mixed neutrosophic $N$-algebraic structure. Then $(F,A)$ over $\langle M \cup I \rangle$ is also a soft Lagrange free mixed neutrosophic $N$-algebraic structure.

**Remark 6.1.6.** Let $(F,A)$ and $(K,C)$ be two soft Lagrange free mixed neutrosophic $N$-algebraic structures over $\langle M \cup I \rangle$. Then

1. Their restricted union $(F,A) \cup_R (K,C)$ may not be a soft Lagrange free mixed neutrosophic $N$-algebraic structure over $\langle M \cup I \rangle$.
2. Their extended union $(F,A) \cup_E (K,C)$ may not be a soft Lagrange free mixed neutrosophic $N$-algebraic structure over $\langle M \cup I \rangle$.
3. Their $OR$ operation $(F,A) \vee (K,C)$ may not be a soft Lagrange free mixed neutrosophic $N$-algebraic structure over $\langle M \cup I \rangle$.
4. Their extended intersection $(F,A) \cap_E (K,C)$ may not be a soft Lagrange free mixed neutrosophic $N$-algebraic structure over $\langle M \cup I \rangle$.
5. Their restricted intersection $(F,A) \cap_R (K,C)$ may not be a soft Lagrange free mixed neutrosophic $N$-algebraic structure over $\langle M \cup I \rangle$.
6. Their $AND$ operation $(F,A) \wedge (K,C)$ may not be a soft Lagrange free mixed neutrosophic $N$-algebraic structure over $\langle M \cup I \rangle$.

One can easily establish the above remarks by the help of examples.

Similarly we can define soft Lagrange free weak deficit mixed neutrosophic $N$-algebraic structures.

**Chapter No. 7**

# SUGGESTED PROBLEMS

In this chapter a number of problems about the soft neutrosophic algebraic structures and their generalization are presented.

These problems are as follows:

1. Can one define soft P-Sylow neutrosophic groupoid over a neutrosophic groupoid?

2. How one can define soft weak Sylow neutrosophic groupoid and soft Sylow free neutrosophic groupoid over a neutrosophic groupoid?

3. Define soft Sylow neutrosophic subgroupoid of a soft Sylow neutrosophic groupoid over a neutrosophic groupoid?

4. What are soft Sylow pseudo neutrosophic groupoids, weak Sylow

pseudo neutrosophic groupoids and free Sylow pseudo neutrosophic group over a neutrosophic groupoids?

5. How one can define soft Lagrange Sylow neutrosophic groupoid over a neutrosophic groupoid?

6. What is the difference between soft Sylow neutrosophic groupoid and soft Lagrange Sylow neutrosophic groupoid?

7. Can one define soft Cauchy neutrosophic groupoids and soft semi Cauchy neutrosophic groupoids over a neutrosophic groupoids?

8. Define soft strong Cauchy neutrosophic groupoid over a neutrosophic groupoid?

9. Define all the properties and results in above problems in case of soft neutrosophic bigroupoid over a neutrosophic bigroupoid?

10. Also define these notions for soft neutrosophic N-groups over neutrosophic N-groups?

11. Give some examples of soft neutrosophic groupoids over neutrosophic groupoids?

12. Write a soft neutrosophic groupoid over a neutrosophic groupoid which has soft neutrosophic subgroupoids, soft neutrosophic strong subgroupoids and soft subgroupoids?

13. Also do these calculations for a soft neutrosophic bigroupoid over a neutrosophic bigroupoid and for a soft neutrosophic N-groupoid over a neutrosophic N-groupoid respectively?

14. Can one define a soft neutrosophic cyclic groupoid over a neutrosophic groupoid?

15. What are the soft neutrosophic idempotent groupoids over a neutrosophic groupoids?  Give some of their examples?

16. What is a soft weakly neutrosophic idempotent groupoid over a neutrosophic groupoid?

17. Give 3 examples of soft neutrosophic ring over neutrosophic rings?

18. Also give some examples of soft neutrosophic N-rings over neutrosophic N-rings?

19. Can one define a soft neutrosophic maximal and a soft neutrosophic minimal ideal over a neutrosophic ring?

20. How can one define a soft neutrosophic biideal over a neutrosophic biring? Also give some examples of it?

21. Give some examples of soft neutrosophic birings and soft neutrosophic N-rings?

22. Also consider the above exercise 22, find soft neutrosophic 5-ideal of a soft neutrosophic 5-ring?

23. Give different examples of soft neutrosophic 3-ring, soft neutrosophic 4-ring and soft neutrosophic 6-ring over a neutrosophic 3-ring, a neutrosophic 4-ring and a neutrosophic 6-ring which are finite respectively?

24. Can we define soft Lagrange neutrosophic ring over a neutrosophic ring?

25. When a soft neutrosophic ring will be a soft Lagrange free neutrosophic ring over a neutrosophic ring?

26. What are soft Lagrange neutrosophic birings over neutrosophic birings?

27. Give some example of soft neutrosophic birings over neutrosophic birings to illustrate this situation?

28. Also define soft neutrosophic biideal over a neutrosophic biring with examples?

29. Give some examples of soft neutrosophic biideal of a soft neutrosophic ring?

30. What are soft Lagrange neutrosophic N-rings over neutrosophic N-rings?

31. What will be the soft neutrosophic ring over a neutrosophic ring of integers?

32. Give some example of soft neutrosophic N-rings over neutrosophic N-rings to illustrate this situation?

33. Define soft neutrosophic field over a neutrosophic complex fields?

34. Give an example of soft neutrosophic field over the neutrosophic real number field?

35. Give some different examples of soft neutrosophic bifields over neutrosophic bifields?

36. Also find soft neutrosophic subbifields in above exercise 35.

37. What is the difference between soft neutrosophic field and soft neutrosophic ring? Also give some examples?

38. Give some examples of soft neutrosophic rings, soft neutrosophic birings and soft neutrosophic 4-rings over neutrosophic rings, neutrosophic birings and neutrosophic 4-rings respectively?

39. Can we define soft Lagrange neutrosophic field over a neutrosophic field?

40. If yes, give some examples. If no, give the reason.

41. Give some examples of soft neutrosophic group rings over neutrosophic group rings?

42. Give some examples of commutative soft neutrosophic group rings over neutrosophic group rings?

43. Can one define soft strong neutrosophic group ring over a neutrosophic group ring?

44. Give some examples of it?

45. Give two examples of soft neutrosophic group ring, three examples of soft neutrosophic bigroup biring and also 2 examples of soft neutrosophic N-group N-ring?

46. Define soft strong neutrosophic right ideal of the soft neutrosophic group ring ?

47. Also give 2 examples of above exercise 46?

48. Under what condition a soft neutrosophic group ring will be a soft strong neutrosophic group ring?

49. Give some examples of soft neutrosophic bigroup biring over neutrosophic bigroup biring?

50. Can one define soft Lagrange neutrosophic group ring over a neutrosophic group ring?

51. Give examples of soft neutrosophic 3-group 3-ring?

52. Find all the soft neutrosophic 3-ideals over a soft neutrosophic 3-group 3-ring ?

53. Give some examples of soft neutrosophic semigroup ring over a neutrosophic semigroup ring?

54. Give some examples of soft mixed neutrosophic semigroup ring over a neutrosophic semigroup ring?

55. Give some examples of soft neutrosophic subring over a neutrosophic semigroup ring?

56. Give some examples of soft neutrosophic subsemigroup ring over a neutrosophic semigroup ring?

57. Give 3 examples of soft neutrosophic bisemigroup biring over a neutrosophic bisemigroup biring?

58. Can one define soft strong neutrosophic semigroup ring over a neutrosophic semigroup ring? Also give some examples?

59. Can one define soft strong neutrosophic ideal over a neutrosophic semigroup ring?

60. Give examples soft neutrosophic 3-semigroup 3-rings?

61. Give some examples of Soft neutrosophic mixed 4-algebraic structure over a neutrosophic mixed 4-structure?

62. Give some examples of Soft neutrosophic mixed 5-algebraic structure over a neutrosophic mixed 4-structure?

63. Give some examples of Soft neutrosophic deficit 4-algebraic structure over a neutrosophic mixed 4-structure?

64. Give some examples of Soft strong neutrosophic mixed 3-algebraic structure over a neutrosophic mixed 3-structure?

65. Give soft neutrosophic mixed sub 3-algebraic structures of a soft neutrosophic mixed 3-algebraic structure over a neutrosophic mixed 3-algebraic structure?

# References

An elaborative list of references is given for the readers to study the notions of neutrosophic set theory and soft set theory in innovative directions that several researchers carried out.

# ABOUT THE AUTHORS

❖ **Mumtaz Ali** has recently completed the degree of M. Phil in the Department of Mathematics, Quaid-e-Azam University Isalmabad. He recently started the research on soft set theory, neutrosophic set theory, algebraic structures including Smarandache algebraic structures, fuzzy theory, coding/ communication theory.

He can be contacted at the following E- mails:

mumtazali770@yahoo.com,

bloomy_boy2006@yahoo.com

❖ **Dr. Florentin Smarandache** is a Professor of Mathematics at the University of New Mexico in USA. He published over 75 books and 100 articles and notes in mathematics, physics, philosophy, psychology, literature, rebus. In mathematics his research is in number theory, non-Euclidean geometry, synthetic geometry, algebraic structures, statistics, neutrosophic logic and set (generalizations of fuzzy logic and set respectively), neutrosophic probability (generalization of classical and imprecise probability). Also, small contributions to nuclear and particle physics, information fusion, neutrosophy (a generalization of dialectics), law of sensations and stimuli, etc.

He can be contacted at the following E-mails:

smarand@unm.edu

fsmarandache@gmail.com

❖ **Dr. Muhammad Shabir** is a Professor in the Department of Mathematics, Quaid-e-Azam University Isalmabad. His interests are in fuzzy and neutrosophic theories, soft set theory, algebraic structures. In the past decade he has guided 6 Ph.D scholars in different fields of Algebras and logics. Currently, 8 Ph. D scholars are working under his guidance. He has written more than 120 research papers and he is also an author of the 2 book.

He can be contacted at the following E-mail:

**mshbirbhatti@yahoo.co.uk**